\documentclass[preprint,12pt]{elsarticle}
\usepackage{amsmath,amssymb,amsfonts,amsthm,array,enumerate,extarrows}
\usepackage[all]{xy}
\usepackage{graphicx}
\usepackage{hyperref}
\usepackage{mathrsfs}
\usepackage{tikz}
\hypersetup{hypertex=true,
colorlinks=true,
linkcolor=blue,
anchorcolor=blue,
citecolor=blue}
\usepackage[nameinlink]{cleveref}
\usepackage{wrapfig}
\usepackage{enumitem}
\usepackage[a4paper, left = 2.5cm, right = 2.5cm, top = 2.54cm, bottom = 2.54cm]{geometry}
\usepackage{mdframed}
\newtheorem{theorem}{\indent Theorem}[section]
\newtheorem{lemma}{\indent Lemma}[section]

\newtheorem{definition}{\indent Def}

\newtheorem{remark}{\indent Remark}[section]

\usepackage{amssymb}
\usepackage{amsmath}

\journal{}

\begin{document}

\begin{frontmatter}



\title{Addition Theorems for Real Vector Spherical Harmonics and Explicit Matrix Representations of the Quasi-Periodic Elastic Single Layer Potential} 


\author[1]{Xin Feng\corref{cor1}} 
\ead{fengxin24@mails.jlu.edu.cn}

\affiliation[1]{organization={School of Mathematics, Jilin University},
            addressline={Qianjin Street 2699},
            city={Changchun},
            postcode={130012},
            state={Jilin},
            country={China}}
\cortext[cor1]{Corresponding author}
\begin{abstract}

This paper develops a multipole expansion method for the quasi-periodic elastic single layer potential $\mathcal{S}_D^{\alpha,0}$ associated with the Kelvin tensor in one-dimensional periodic arrays. A key step in this approach is the derivation of translation addition theorems for the real vector spherical harmonics $V_{lm}$, $W_{lm}$, and $X_{lm}$. These addition theorems enable the exact calculation of all matrix entries of $\mathcal{S}_D^{\alpha,0}$ in closed form. By working entirely within the spherical harmonic basis, the proposed analytical method overcomes the poor convergence and mesh-dependent issues commonly caused by the direct surface discretization of weakly singular kernels. Additionally, the involved infinite sums are evaluated exactly using polylogarithm functions, which eliminates the need for series truncation. As an application, the integral equation $\mathcal{S}_D^{\alpha,0}[f]=\varphi$ is reduced to a linear system. This framework is further extended to dimer geometries consisting of two disjoint balls in each cell, where the off-diagonal matrices are explicitly formulated via the Lerch transcendent.\end{abstract}


\end{frontmatter}




\section{Introduction}
In the study of wave propagation in periodic resonator structures, layer potential techniques combined with Floquet-Bloch theory provide a powerful analytical framework. In the acoustic setting, a multipole expansion method based on the scalar spherical harmonic basis has been developed in \citet{ammari2020topologically} to analyze subwavelength band gaps and topological edge modes in one-dimensional resonator chains. A key ingredient of that approach is the addition theorem for scalar spherical waves (see  \cite[Appendix~A]{ammari2020topologically}), which allows the quasi-periodic single layer potential to be represented explicitly as a matrix in the spherical harmonic basis, analytically computable lattice sums.

The present paper develops the elastic analogue of this framework. The governing operator is now the Lamé system, whose fundamental solution is the Kelvin tensor, and the natural basis on the boundary of a spherical resonator consists of the three families of real vector spherical harmonics $V_{lm}$, $W_{lm}$, $X_{lm}$ rather than scalar spherical harmonics \cite{StammXiang2022}. The central difficulty is that addition theorems for these vector harmonics under translation are not available in the literature, yet they are indispensable for computing the matrix representation of the quasi-periodic elastic single layer potential $\mathcal{S}_D^{\alpha,0}$.

A further motivation for this approach is numerical. The kernel of the elastic single layer potential $\mathcal{S}_D^{\alpha,0}$ is weakly singular, and direct discretization of the operator on the sphere via surface meshes leads to bad convergence and numerical results that are sensitive to the choice of mesh. By working entirely in the vector spherical harmonic basis and deriving analytical expressions for all matrix entries, our method is completely mesh-free and free of numerical discretization error. Moreover, the infinite lattice sums that appear in the entries of $\mathbf{M}(\alpha)$ are evaluated in closed form via polylogarithm functions\cite{lewin1981polylogarithms} and, in the dimer case, the Lerch transcendent\cite{Lerch1887}, thereby avoiding the need for any truncation of infinite series. 

This gap is filled by deriving addition theorems for all three families of real vector spherical harmonics, and applying them to compute in closed form all entries of the matrix $\mathbf{M}(\alpha)$ of $\mathcal{S}_D^{\alpha,0}$ in the vector spherical harmonic basis. The framework is further extended to the dimer geometry, with coupling matrices $\mathbf{M}^{12}$ and $\mathbf{M}^{21}$ expressed through the Lerch transcendent. These results provide the analytic foundation for a multipole expansion method for elastic subwavelength resonator chains \cite{Ammari2024, RenChenGaoLi2025SubwavelengthBandgaps}, directly paralleling and generalizing the scalar approach of \cite{ammari2020topologically}.

\section{Preliminaries}
In this section, the single layer potential for the elastic Kelvin Green tensor, real vector spherical harmonics and the addition theorems for solid harmonics are introduced. The underlying space is $\mathbb{R}^3$ throughout this discussion.
Denote by $B_\rho(x_0)$ the ball centered at $x_0$ with radius $\rho$ and $\partial B_\rho(x_0)$ its boundary. In particular, let $\mathbb{S}^2$ denote the unit sphere in $\mathbb{R}^3$.
\subsection{Kelvin tensor and single layer potential}
 Define the Lamé operator corresponding to the Lamé constants $\lambda, \mu$  by
\begin{equation}
\mathcal{L}^{\lambda, \mu}\boldsymbol{u} := -\mu\Delta \boldsymbol{u} -(\lambda+ \mu)\nabla\nabla\cdot \boldsymbol{u}.
\end{equation}
This definition is consistent with  \citet[P1, Eq~1.1]{StammXiang2022}, \citet[P96, section~5.2]{steinbach2008numerical}  and \citet{ammari2015mathematical,ammari2018mathematical,ammari2007polarization} up to a sign. It is worth  mentioning that one should reverse the sign when using the results in this paper if the operator being considered is $-\mathcal{L}^{\lambda, \mu}$.
Meanwhile, the fundamental solution $\boldsymbol{G}(x)$ is given by \cite{ammari2015mathematical,ammari2018mathematical,steinbach2008numerical,ammari2007polarization}
\begin{equation}\label{eq-Kelvin-tensor}\boldsymbol{G}_{ij}(x):=\dfrac{1}{8\pi\vert x\vert }\cdot \left(\dfrac{\lambda+3\mu}{\lambda+2\mu}\delta_{ij}+\dfrac{\lambda+\mu}{\lambda+2\mu}\dfrac{x_ix_j}{\vert x\vert^2}\right),
\end{equation}
called the Kelvin solution tensor, where $\delta_{ij}$ is the Kronecker symbol.\\
Using the Kelvin tensor \eqref{eq-Kelvin-tensor}, the single layer potential on $L^2(\partial B_\rho(x_0))^3$ is defined by
\begin{equation}
\mathcal{S}[\psi](x):=\int_{\partial B_\rho (x_0)} \boldsymbol{G}(x-y)\psi(y)d\sigma(y),\quad x\in\mathbb{R}^3.
\end{equation}
This operator is well-defined because of the weak singularity of $\boldsymbol{G}$\cite{ammari2007polarization,ammari2015mathematical}.

Define the quasi-periodic single layer potential by
\begin{equation}\label{eq-def-S_D}
\mathcal{S}_D^{\alpha,0}[\psi](x) := \int_{\partial D} \boldsymbol{G}^{\alpha,0}(x-y) \psi(y) \mathrm{d}\sigma(y), \quad x \in \partial D,
\end{equation}
where $D=B_\rho(0)$ with $\rho<1/2$ and $\boldsymbol{G}^{\alpha,0}(x-y)=\sum\limits_{n\in \mathbb{Z}}\boldsymbol{G}(x-y-\mathbf{a})e^{in\alpha}$ with $\mathbf{a}=(n,0,0)$.\\
The operator $\mathcal{S}_D^{\alpha,0}[\psi](x)$ can be decomposed as 
\begin{equation}\label{eq-operator-S}\mathcal{S}_D^{\alpha,0}[\psi](x)=\mathcal{S}_D[\psi](x)+\sum\limits_{n\neq 0}\mathcal{S}_{D+n}[\psi](x)e^{in\alpha},\end{equation}
where the operator $(L^2(\partial D))^3\xlongrightarrow{\mathcal{S}_{D+n}}(L^2(\partial D))^3$ is defined by
 $$\mathcal{S}_{D+n}[\psi](x)=\int_{\partial D}\boldsymbol{G}(x-y-\mathbf{a})\psi(y)d\sigma(y).$$
By making the change of variable $y'=y+\mathbf{a}$, one obtains:
$$\mathcal{S}_{D+n}[\psi](x)=\int_{\partial D+\mathbf{a}}\boldsymbol{G}(x-y)\psi(y-\mathbf{a})d\sigma(y).$$ 
It can be regarded as an operator from $(L^2(\partial D+\mathbf{a}))^3$ to $(L^2(\partial D))^3$. 

\subsection{Real vector  spherical harmonics}
The real  spherical harmonics on $\mathbb{S}^2$ are well known. To distinguish the real and complex spherical harmonics, their indices are placed in different positions. For instance, $Y_l^m$ is a complex spherical harmonic and $Y_{lm}$ is the real one. Furthermore, the relation between them is given by \cite{romanowski2008transformation}:
\begin{equation}\label{eq-relation-complex-real-Y}Y_{lm}=\left\{
\begin{aligned}
\dfrac{1}{\sqrt{2}}(Y_l^{-m}+(-1)^m Y_l^m)&, \quad m>0,\\
Y_l^0\qquad\qquad &,\quad m=0,\\
\dfrac{i}{\sqrt{2}}(Y_l^m-(-1)^mY_l^{-m})&,\quad m<0.
\end{aligned}\right.\end{equation}
The real vector spherical harmonics are defined by 
$$\begin{aligned}
V_{lm}(\hat{\mathbf{r}})&:=\nabla_SY_{lm}(\hat{\mathbf{r}})-(l+1)Y_{lm}(\hat{\mathbf{r}})\hat{\mathbf{r}},\\
W_{lm}(\hat{\mathbf{r}})&:=\nabla_S Y_{lm}(\hat{\mathbf{r}})+lY_{lm}(\hat{\mathbf{r}})\hat{\mathbf{r}},\\
X_{lm}(\hat{\mathbf{r}})&:=\hat{\mathbf{r}}\times\nabla_SY_{lm}(\hat{\mathbf{r}})
\end{aligned}$$
where $\nabla_S=\hat{\theta}\dfrac{\partial}{\partial \theta}+\hat{\phi}\dfrac{1}{\sin \theta}\dfrac{\partial}{\partial \phi}$ is the surface gradient operator and $r$, $\theta$, and $\phi$ are the spherical polar coordinates of the 3-dimensional vector $\mathbf{r}$\cite{StammXiang2022}. The notations $Y_{lm}^1,Y_{lm}^2,Y_{lm}^3$ are used to represent $V_{lm}, W_{lm},X_{lm}$ respectively in this paper. Meanwhile, the complex vector spherical harmonics can be derived from the relation \eqref{eq-relation-complex-real-Y}:
$$\begin{aligned}
V_l^m(\hat{\mathbf{r}})&=\nabla_SY_l^m(\hat{\mathbf{r}})-(l+1)Y_l^m(\hat{\mathbf{r}})\hat{\mathbf{r}},\\
W_l^m(\hat{\mathbf{r}})&=\nabla_S Y_l^m(\hat{\mathbf{r}})+lY_l^m(\hat{\mathbf{r}})\hat{\mathbf{r}},\\
X_l^m(\hat{\mathbf{r}})&=\hat{\mathbf{r}}\times\nabla_SY_l^m(\hat{\mathbf{r}}).
\end{aligned}$$
Precisely,
$$    V_\lambda^\mu=\left\{
\begin{aligned}
\dfrac{(-1)^\mu}{\sqrt{2}}(V_{\lambda\mu}+i V_{\lambda,-\mu})&, \quad \mu>0,\\
V_{\lambda 0}\qquad\qquad &,\quad \mu=0,\\
\dfrac{1}{\sqrt{2}}(V_{\lambda,-\mu}-iV_{\lambda\mu})&,\quad \mu<0,
\end{aligned}\right.    $$
$$    W_\lambda^\mu=\left\{
\begin{aligned}
\dfrac{(-1)^\mu}{\sqrt{2}}(W_{\lambda\mu}+i W_{\lambda,-\mu})&, \quad \mu>0,\\
W_{\lambda 0}\qquad\qquad &,\quad \mu=0,\\
\dfrac{1}{\sqrt{2}}(W_{\lambda,-\mu}-iW_{\lambda\mu})&,\quad \mu<0,
\end{aligned}\right.    $$
and
$$    X_\lambda^\mu=\left\{
\begin{aligned}
\dfrac{(-1)^\mu}{\sqrt{2}}(X_{\lambda\mu}+i X_{\lambda,-\mu})&, \quad \mu>0,\\
X_{\lambda 0}\qquad\qquad &,\quad \mu=0,\\
\dfrac{1}{\sqrt{2}}(X_{\lambda,-\mu}-iX_{\lambda\mu})&,\quad \mu<0.
\end{aligned}\right.    $$

\subsection{Solid harmonics and their addition theorems}

Define two complex harmonics:
\begin{equation}\label{solid-harmonics}R_l^m(\mathbf{r})=\sqrt{\dfrac{4\pi}{2l+1}}\vert \mathbf{r}\vert^lY_l^m(\hat{\mathbf{r}}),\quad I_l^m(\mathbf{r})=\sqrt{\dfrac{4\pi}{2l+1}}\dfrac{Y_l^m(\hat{\mathbf{r}})}{\vert \mathbf{r}\vert^{l+1}}.\end{equation}
The translation addition theorem of the regular solid harmonic gives a finite expansion\cite{caola1978solid,tough1977properties}:
$$
R_{l}^{m}(\mathbf{r} + \mathbf{a}) = \sum_{\lambda=0}^{l} {2l \choose 2\lambda}^{1/2} \sum_{\mu=-\lambda}^{\lambda} R_{\lambda}^{\mu}(\mathbf{r}) R_{l-\lambda}^{m-\mu}(\mathbf{a}) \langle \lambda, \mu; l - \lambda, m - \mu | l m \rangle,
$$
where the Clebsch-Gordan coefficient is given by
$$
\langle \lambda, \mu; l - \lambda, m - \mu | l m \rangle =  {l + m\choose \lambda + \mu}^{1/2} {l - m\choose\lambda - \mu}^{1/2} {2l \choose 2\lambda}^{-1/2}.
$$

The similar expansion for irregular solid harmonics gives an infinite series\cite{caola1978solid,tough1977properties}:
$$
I_{l}^{m}(\mathbf{r} + \mathbf{a}) = \sum_{\lambda=0}^{\infty} {2l + 2\lambda + 1 \choose 2\lambda}^{1/2} \sum_{\mu=-\lambda}^{\lambda} R_{\lambda}^{\mu}(\mathbf{r}) I_{l+\lambda}^{m-\mu}(\mathbf{a}) \langle \lambda, \mu; l + \lambda, m - \mu | l m \rangle
$$
with \(|\mathbf{r}| \leq |\mathbf{a}|\). The quantity between pointed brackets is again a Clebsch-Gordan coefficient,
$$
\langle \lambda, \mu; l + \lambda, m - \mu | l m \rangle = (-1)^{\lambda+\mu} {l + \lambda - m + \mu\choose \lambda + \mu}^{1/2}{l + \lambda + m - \mu\choose\lambda - \mu} ^{1/2} {2l + 2\lambda + 1 \choose 2\lambda}^{-1/2}.
$$

\section{Addition theorem for vector spherical harmonics}

For the complex valued vector spherical harmonics $V_l^m$, $W_l^m$ and $X_l^m$, their addition theorems are derived with the aid of the following lemmas:

\begin{lemma}The solid harmonics $R_l^m(\mathbf{r})$ and $I_l^m(\mathbf{r})$ are defined by \eqref{solid-harmonics} as before. Then
\begin{enumerate}
\item $\nabla I_l^m(\mathbf{r})=\sqrt{\dfrac{4\pi}{2l+1}}r^{-l-2}V_l^m(\hat{\mathbf{r}})$.
\item $\nabla R_\lambda^\mu(\mathbf{r})=\sqrt{\dfrac{4\pi}{2\lambda+1}}r^{\lambda-1}W_\lambda^\mu(\hat{\mathbf{r}})$.
\item $\mathbf{r}\times \nabla I_l^m(\mathbf{r})=\sqrt{\dfrac{4\pi}{2l+1}} r^{-l-1}X_l^m(\hat{\mathbf{r}})$
\item $\mathbf{r}\times \nabla R_\lambda^\mu(\mathbf{r})=\sqrt{\dfrac{4\pi}{2\lambda+1}}r^\lambda X_\lambda^\mu(\hat{\mathbf{r}})$.
\item $\mathbf{a}\times \nabla R_\lambda^\mu(\mathbf{r})=\sqrt{\dfrac{4\pi}{2\lambda+1}}r^{\lambda-1}(\mathbf{a}\times W_\lambda^\mu(\hat{\mathbf{r}}))$
\end{enumerate}
\end{lemma}
\begin{theorem}Let $\mathbf{r}'=\mathbf{r}+\mathbf{a}$ with $\vert \mathbf{a}\vert>r$. It can be derived smoothly that
\begin{equation}\label{eq-addcomplex-V}{r'}^{-l-2}V_l^m(\hat{\mathbf{r}'})
=\sum\limits_{\lambda=1}^\infty\sum\limits_{\mu=-\lambda}^\lambda
(-1)^{\lambda+\mu}
\sqrt{\dfrac{2l+1}{2\lambda+1}}
r^{\lambda-1}
I_{l+\lambda}^{m-\mu}(\mathbf{a})\sqrt{{l+\lambda+\mu-m\choose \lambda+\mu}{l+\lambda+m-\mu\choose \lambda-\mu}}W_\lambda^\mu(\hat{\mathbf{r}}),\end{equation}
$${r'}^{l-1}W_l^m(\hat{\mathbf{r}'})=\sum\limits_{\lambda=0}^l\sum\limits_{\mu=-\lambda}^\lambda
\sqrt{\dfrac{2l+1}{2\lambda+1}}r^{\lambda-1}R_{l-\lambda}^{m-\mu}(\mathbf{a})\sqrt{{l+m\choose \lambda+\mu}{l-m\choose \lambda-\mu}}W_\lambda^\mu(\hat{\mathbf{r}}),$$
\begin{equation}\label{eq-addcomplex-X}\begin{aligned}{r'}^{-l-1}X_l^m(\hat{\mathbf{r}'})=\sum\limits_{\lambda=1}^\infty\sum\limits_{\mu=-\lambda}^\lambda(-1)^{\lambda+\mu}\sqrt{\dfrac{2l+1}{2\lambda+1}} r^{\lambda-1}I_{l+\lambda}^{m-\mu}(\mathbf{a})&
\sqrt{{l+\lambda+\mu-m\choose \lambda+\mu}{l+\lambda+m-\mu\choose \lambda-\mu}}\\
&\left[r X_\lambda^\mu(\hat{\mathbf{r}})+\mathbf{a}\times W_\lambda^\mu(\hat{\mathbf{r}})\right],
\end{aligned}
\end{equation}
where $r=\vert \mathbf{r}\vert $ and $r'=\vert \mathbf{r}+\mathbf{a}\vert$.
\end{theorem}
These results can be moved to real valued spherical harmonics case. 
\subsection{Addition theorem for ${r'}^{l-1}W_{lm}$ }
To obtain the addition theorem for ${r'}^{l-1}W_{lm}$, for instance, one can apply the addition theorem for the complex-valued harmonics $W_l^m$ and, conversely, transform the terms $W_\lambda^\mu$ back into the real-valued $W_{\lambda\mu}$. Then the addition theorem is obtained, while keeping complex-valued terms $R_{l-\lambda}^{m-\mu}$ and discarding $W_0^0$ which vanishes.
\begin{theorem}
Denote the coefficient $\sqrt{\dfrac{2l+1}{2\lambda+1}}\sqrt{{l+m\choose \lambda+\mu}{l-m\choose \lambda-\mu}}R_{l-\lambda}^{m-\mu}(\mathbf{a})$ by $A_{l,\lambda}^{m,\mu}(\mathbf{a})$. Then, when $m>0$,
$$\begin{aligned}
{r'}^{l-1}W_{lm}(\mathbf{r}')
=&\sum\limits_{\lambda=1}^l\sum\limits_{\mu=-\lambda}^\lambda r^{\lambda-1}\dfrac{1}{\sqrt{2}}\left(A_{l,\lambda}^{-m,\mu}(\mathbf{a})+(-1)^mA_{l,\lambda}^{m,\mu}(\mathbf{a})\right)W_\lambda^\mu(\hat{\mathbf{r}})\\
=&\sum\limits_{\lambda=1}^l\sum\limits_{\mu=-\lambda}^{-1} r^{\lambda-1}\dfrac{1}{2}\left(A_{l,\lambda}^{-m,\mu}(\mathbf{a})+(-1)^mA_{l,\lambda}^{m,\mu}(\mathbf{a})\right)[W_{\lambda,-\mu}(\hat{\mathbf{r}})-iW_{\lambda,\mu}(\hat{\mathbf{r}})]\\
&+\sum\limits_{\lambda=1}^l\sum\limits_{\mu=1}^\lambda r^{\lambda-1}\dfrac{(-1)^{\mu}}{2}\left(A_{l,\lambda}^{-m,\mu}(\mathbf{a})+(-1)^mA_{l,\lambda}^{m,\mu}(\mathbf{a})\right)[W_{\lambda,\mu}(\hat{\mathbf{r}})+iW_{\lambda,-\mu}(\hat{\mathbf{r}})]\\
&+\sum\limits_{\lambda=1}^l r^{\lambda-1}\dfrac{1}{\sqrt{2}}\left(A_{l,\lambda}^{-m,0}(\mathbf{a})+(-1)^mA_{l,\lambda}^{m,0}(\mathbf{a})\right)W_{\lambda,0}(\hat{\mathbf{r}}).
\end{aligned}$$
Similarly, when $m<0$,
$$\begin{aligned}
{r'}^{l-1}W_{lm}(\mathbf{r}')
=&\sum\limits_{\lambda=1}^l\sum\limits_{\mu=-\lambda}^\lambda r^{\lambda-1}\dfrac{i}{\sqrt{2}}\left(A_{l,\lambda}^{m,\mu}(\mathbf{a})-(-1)^mA_{l,\lambda}^{-m,\mu}(\mathbf{a})\right)W_\lambda^\mu(\hat{\mathbf{r}})\\
=&\sum\limits_{\lambda=1}^l\sum\limits_{\mu=-\lambda}^{-1} r^{\lambda-1}\dfrac{i}{2}\left(A_{l,\lambda}^{m,\mu}(\mathbf{a})-(-1)^mA_{l,\lambda}^{-m,\mu}(\mathbf{a})\right)[W_{\lambda,-\mu}(\hat{\mathbf{r}})-iW_{\lambda,\mu}(\hat{\mathbf{r}})]\\
&+\sum\limits_{\lambda=1}^l\sum\limits_{\mu=1}^\lambda r^{\lambda-1}\dfrac{(-1)^{\mu}}{2}\cdot i\left(A_{l,\lambda}^{m,\mu}(\mathbf{a})-(-1)^mA_{l,\lambda}^{-m,\mu}(\mathbf{a})\right)[W_{\lambda,\mu}(\hat{\mathbf{r}})+iW_{\lambda,-\mu}(\hat{\mathbf{r}})]\\
&+\sum\limits_{\lambda=1}^l r^{\lambda-1}\dfrac{i}{\sqrt{2}}\left(A_{l,\lambda}^{m,0}(\mathbf{a})-(-1)^mA_{l,\lambda}^{-m,0}(\mathbf{a})\right)W_{\lambda,0}(\hat{\mathbf{r}}).
\end{aligned}$$
When $m=0$, 
$$\begin{aligned}
{r'}^{l-1}W_{l0}(\mathbf{r}')=&{r'}^{l-1}W_l^0(\mathbf{r}')=\sum\limits_{\lambda=1}^l\sum\limits_{\mu=-\lambda}^\lambda
r^{\lambda-1}A_{l,\lambda}^{0,\mu}(\mathbf{a})W_\lambda^\mu(\hat{\mathbf{r}})\\
=&\sum\limits_{\lambda=1}^l\sum\limits_{\mu=-\lambda}^{-1}
r^{\lambda-1}\dfrac{1}{\sqrt{2}}A_{l,\lambda}^{0,\mu}(\mathbf{a})[W_{\lambda,-\mu}(\hat{\mathbf{r}})-iW_{\lambda,\mu}(\hat{\mathbf{r}})]\\
&+\sum\limits_{\lambda=1}^l\sum\limits_{\mu=1}^{\lambda}
r^{\lambda-1}\dfrac{(-1)^\mu}{\sqrt{2}}A_{l,\lambda}^{0,\mu}(\mathbf{a})[W_{\lambda,\mu}(\hat{\mathbf{r}})+iW_{\lambda,-\mu}(\hat{\mathbf{r}})]\\
&+\sum\limits_{\lambda=1}^lA_{l,\lambda}^{0,0}(\mathbf{a})W_{\lambda,0}(\hat{\mathbf{r}}).
\end{aligned}$$
\end{theorem}

In this way, the addition theorems for real-valued vector spherical harmonics are obtained as follows.
\subsection{Addition theorem for ${r'}^{-l-2}V_l^m$}
\begin{theorem}
Denote the coefficient in the addition theorem for ${r'}^{-l-2}V_{lm}$ by $H_{l,\lambda}^{m,\mu}$, i.e.,
$$\boxed{H_{l,\lambda}^{m,\mu}:=(-1)^{\lambda+\mu}
\sqrt{\dfrac{2l+1}{2\lambda+1}}
I_{l+\lambda}^{m-\mu}(\mathbf{a})
\sqrt{{l+\lambda+\mu-m\choose \lambda+\mu}{l+\lambda+m-\mu\choose \lambda-\mu}}.}$$
When $m>0$,
$$\begin{aligned}{r'}^{-l-2}V_{lm}(\mathbf{r}')
=&\sum\limits_{\lambda=1}^\infty\sum\limits_{\mu=-\lambda}^\lambda r^{\lambda-1} \dfrac{1}{\sqrt{2}}\left(H_{l,\lambda}^{-m,\mu}(\mathbf{a})+(-1)^mH_{l,\lambda}^{m,\mu}(\mathbf{a})\right)W_\lambda^\mu(\hat{\mathbf{r}})\\
=&\sum\limits_{\lambda=1}^\infty\sum\limits_{\mu=-\lambda}^{-1}r^{\lambda-1} \dfrac{1}{2}\left(H_{l,\lambda}^{-m,\mu}(\mathbf{a})+(-1)^mH_{l,\lambda}^{m,\mu}(\mathbf{a})\right)[W_{\lambda,-\mu}(\hat{\mathbf{r}})-iW_{\lambda,\mu}(\hat{\mathbf{r}})]\\
&+\sum\limits_{\lambda=1}^\infty\sum\limits_{\mu=1}^{\lambda}r^{\lambda-1} \dfrac{(-1)^\mu}{2}\left(H_{l,\lambda}^{-m,\mu}(\mathbf{a})+(-1)^mH_{l,\lambda}^{m,\mu}(\mathbf{a})\right)[W_{\lambda,\mu}(\hat{\mathbf{r}})+iW_{\lambda,-\mu}(\hat{\mathbf{r}})]\\
&+\sum\limits_{\lambda=1}^\infty r^{\lambda-1} \dfrac{1}{\sqrt{2}}\left(H_{l,\lambda}^{-m,0}(\mathbf{a})+(-1)^mH_{l,\lambda}^{m,0}(\mathbf{a})\right)W_{\lambda,0}(\hat{\mathbf{r}})
\end{aligned}$$

When $m<0$,
$$\begin{aligned}{r'}^{-l-2}V_{lm}(\mathbf{r}')& =\sum\limits_{\lambda=1}^\infty\sum\limits_{\mu=-\lambda}^\lambda r^{\lambda-1} \dfrac{i}{\sqrt{2}}\left(H_{l,\lambda}^{m,\mu}(\mathbf{a})-(-1)^mH_{l,\lambda}^{-m,\mu}(\mathbf{a})\right)W_\lambda^\mu(\hat{\mathbf{r}})\\
&=\sum\limits_{\lambda=1}^\infty\sum\limits_{\mu=-\lambda}^{-1} r^{\lambda-1}\dfrac{i}{2}\bigl(H_{l,\lambda}^{m,\mu}(\mathbf{a})-(-1)^mH_{l,\lambda}^{-m,\mu}(\mathbf{a})\bigr)[W_{\lambda,-\mu}(\hat{\mathbf{r}})-iW_{\lambda,\mu}(\hat{\mathbf{r}})]\\
&\quad+\sum\limits_{\lambda=1}^\infty\sum\limits_{\mu=1}^\lambda r^{\lambda-1}\dfrac{(-1)^{\mu}}{2}\cdot i\bigl(H_{l,\lambda}^{m,\mu}(\mathbf{a})-(-1)^mH_{l,\lambda}^{-m,\mu}(\mathbf{a})\bigr)[W_{\lambda,\mu}(\hat{\mathbf{r}})+iW_{\lambda,-\mu}(\hat{\mathbf{r}})]\\
&\quad+\sum\limits_{\lambda=1}^\infty r^{\lambda-1}\dfrac{i}{\sqrt{2}}\bigl(H_{l,\lambda}^{m,0}(\mathbf{a})-(-1)^mH_{l,\lambda}^{-m,0}(\mathbf{a})\bigr)W_{\lambda,0}(\hat{\mathbf{r}}).
\end{aligned}$$
When $m=0$,
$$\begin{aligned}{r'}^{-l-2}V_{lm}(\mathbf{r}')&=\sum\limits_{\lambda=1}^\infty\sum\limits_{\mu=-\lambda}^\lambda r^{\lambda-1} H_{l,\lambda}^{0,\mu}(\mathbf{a})W_\lambda^\mu(\hat{\mathbf{r}})\\
&=\sum\limits_{\lambda=1}^\infty\sum\limits_{\mu=-\lambda}^{-1}
r^{\lambda-1}\dfrac{1}{\sqrt{2}}H_{l,\lambda}^{0,\mu}(\mathbf{a})[W_{\lambda,-\mu}(\hat{\mathbf{r}})-iW_{\lambda,\mu}(\hat{\mathbf{r}})]\\
&\quad+\sum\limits_{\lambda=1}^\infty\sum\limits_{\mu=1}^{\lambda}
r^{\lambda-1}\dfrac{(-1)^\mu}{\sqrt{2}}H_{l,\lambda}^{0,\mu}(\mathbf{a})[W_{\lambda,\mu}(\hat{\mathbf{r}})+iW_{\lambda,-\mu}(\hat{\mathbf{r}})]\\
&\quad+\sum\limits_{\lambda=1}^\infty r^{\lambda-1}H_{l,\lambda}^{0,0}(\mathbf{a})\,W_{\lambda,0}(\hat{\mathbf{r}}).
\end{aligned}$$
\end{theorem}

\subsection{Addition theorem for ${r'}^{-l}V_{lm}$ and ${r'}^{-l}W_{lm}$ }
By Eq.\eqref{eq-addcomplex-V}, it follows that
$$\begin{aligned}
{r'}^{-l}V_l^m(\hat{\mathbf{r}'})
&={r'}^{2}\cdot{r'}^{-l-2}V_l^m(\hat{\mathbf{r}'})={r'}^{2}\cdot  \sum\limits_{\lambda=1}^\infty\sum\limits_{\mu=-\lambda}^{\lambda}H_{l,\lambda}^{m,\mu}(\mathbf{a}) r^{\lambda-1}W_\lambda^\mu(\hat{\mathbf{r}})\\
&=\sum\limits_{\lambda=1}^\infty\sum\limits_{\mu=-\lambda}^{\lambda}H_{l,\lambda}^{m,\mu}(\mathbf{a})((r^{\lambda+1}+a^2r^{\lambda-1})W_\lambda^\mu(\hat{\mathbf{r}})+2r^{\lambda}\hat{\mathbf{r}}\cdot \mathbf{a}W_\lambda^\mu(\hat{\mathbf{r}})).\end{aligned}$$
From  \citet[Eq.~26]{haber2017tensor}, it follows that
\begin{lemma}
$\hat{\mathbf{r}}\cdot \mathbf{a}=\sqrt{\dfrac{4\pi}{3}}\sum\limits_{q=-1}^{1}(-1)^q a_{-q}Y_1^q(\hat{\mathbf{r}}),$
where $a_{-q}:=\mathbf{a}\cdot \hat{\mathbf{e}}_{-q}$, more precisely,
\begin{equation}\label{eq-a-q}\mathbf{a}\cdot \hat{\mathbf{e}}_{-q}=\left\{
\begin{aligned}
\dfrac{1}{\sqrt{2}}(a_x-ia_y)&,\quad q=1,\\
a_z\quad&,\quad q=0\\
-\dfrac{1}{\sqrt{2}}(a_x+ia_y)&,\quad q=-1
\end{aligned}\right.\end{equation}
\end{lemma}
\begin{remark}When  $  \mathbf{b} = (-n, 0, 0)  $,
$$b_{-q} =
\begin{cases}
\dfrac{n}{\sqrt{2}} & q = -1 \\
0 & q = 0 \\
-\dfrac{n}{\sqrt{2}} & q = 1
\end{cases}$$
\end{remark}
From \citet[Eq.~4, Eq.~46, Eq.~88]{haber2017tensor}, they imply that
\begin{equation}\label{eq-W-vecY}W_\lambda^\mu(\hat{\mathbf{r}})=\sqrt{\lambda(2\lambda+1)}\vec{Y}_{\lambda,\lambda-1,\mu}(\hat{\mathbf{r}}),\end{equation}
\begin{equation}\label{eq-vecY-expand}\vec{Y}_{\lambda,\lambda-1,\mu}(\hat{\mathbf{r}})=\sum\limits_{m_1=-1}^{1}\langle \lambda-1,\mu-m_1; 1, m_1\mid \lambda\mu\rangle Y_{\lambda-1}^{\mu-m_1}(\hat{\mathbf{r}})\chi_{1,m_1},\end{equation}
$$Y_1^q(\hat{\mathbf{r}})Y_{\lambda-1}^{\mu-m_1}(\hat{\mathbf{r}})=\sum\limits_{\vert\lambda-2\vert \leqslant l'\leqslant\lambda,\atop l'\neq\lambda-1} \sqrt{\dfrac{3(2\lambda-1)}{4\pi(2l'+1)}}\langle 1,q;\lambda-1,\mu-m_1\mid l',q+\mu-m_1\rangle\langle 1,0;\lambda-1,0\mid l',0\rangle Y_{l'}^{q+\mu-m_1}(\hat{\mathbf{r}}).$$
Therefore, it follows that
\begin{equation}\label{eq-raWlambdamu}\begin{aligned}
(\hat{\mathbf{r}}\cdot\mathbf{a})W_\lambda^\mu(\hat{\mathbf{r}})
&=\sum\limits_{q=-1}^{1}\sum\limits_{m_1=-1}^{1}\sum\limits_{\vert\lambda-2\vert \leqslant l'\leqslant\lambda,\atop l'\neq\lambda-1}\sqrt{\dfrac{\lambda(2\lambda+1)(2\lambda-1)}{2l'+1}}(-1)^qa_{-q}\\
\langle\lambda-1,\mu-m_1;1,m_1\mid \lambda,\mu\rangle&\langle 1,q;\lambda-1,\mu-m_1\mid l',q+\mu-m_1\rangle\langle 1,0;\lambda-1,0\mid l',0\rangle Y_{l'}^{q+\mu-m_1}(\hat{\mathbf{r}})\chi_{1,m_1}.
\end{aligned}\end{equation}
\begin{remark}If $\lambda=1$, then $l'$ only takes value on $1$.
\end{remark}
Now the last term $ Y_{l'}^{q+\mu-m_1}(\hat{\mathbf{r}})\chi_{1,m_1}$ should be recoupled. From \cite[Eq.~4]{haber2017tensor} and the orthonormal relation:
$$\sum\limits_{j,m}\langle l,k;1,m_s\mid j,m\rangle\langle l,k';1,m_s'\mid j,m\rangle=\delta_{kk'}\delta_{m_sm_s'},$$ 
it is obtained that
\begin{equation}\label{eq-Ykl'}Y_{l'}^k(\hat{\mathbf{r}})\chi_{1,m_1}= \sum\limits_{\vert l'-1\vert\leqslant j\leqslant \vert l'+1\vert}\langle l',k;1,m_1\mid j,k+m_1\rangle \vec{Y}_{j,l',k+m_1}(\hat{\mathbf{r}}).\end{equation}
\begin{remark}If $l'=0$, then $j$ only takes value on $1$. And if this is the case,  $\vec{Y}_{1,0,k+m_1}(\hat{\mathbf{r}})=\dfrac{1}{\sqrt{3}}W_1^{k+m_1}(\hat{\mathbf{r}})$ by \cite[Eq.~46]{haber2017tensor}.
\end{remark}
By \cite[Eq~45, Eq~46, Eq~47]{haber2017tensor}, $(\hat{\mathbf{r}}\cdot\mathbf{a})W_\lambda^\mu(\hat{\mathbf{r}})$ is the linear combination of vector harmonics since
$$\vec{Y}_{l'-1,l',k+m_1}=\dfrac{1}{\sqrt{l'(2l'-1)}}V_{l'-1}^{k+m_1},$$
$$\vec{Y}_{l',l',k+m_1}=\dfrac{-i}{\sqrt{l'(l'+1)}}X_{l'}^{k+m_1},$$
$$\vec{Y}_{l'+1,l',k+m_1}=\dfrac{1}{\sqrt{(l'+1)(2l'+3)}}W_{l'+1}^{k+m_1}.$$
Combining Eq~\eqref{eq-raWlambdamu} and Eq~\eqref{eq-Ykl'} where $k=q+\mu-m_1$, it follows that
\begin{theorem}
$$\begin{aligned}{r'}^{-l}V_l^m(\hat{\mathbf{r}'})
&=\sum\limits_{\lambda=1}^\infty \sum\limits_{\mu=-\lambda}^{\lambda}H_{l,\lambda}^{m,\mu}(\mathbf{a})(r^{\lambda+1}+a^2r^{\lambda-1})W_\lambda^\mu(\hat{\mathbf{r}})\\
&+\sum\limits_{\lambda=1}^\infty \sum\limits_{\mu=-\lambda}^{\lambda}\sum\limits_{q=-1}^{1}\sum\limits_{m_1=-1}^{1}\sum\limits_{\vert\lambda-2\vert \leqslant l'\leqslant\lambda,\atop l'\neq\lambda-1}\sum\limits_{\vert l'-1\vert\leqslant j\leqslant \vert l'+1\vert}
H_{l,\lambda}^{m,\mu}(\mathbf{a})2r^{\lambda}K_{l',j,\lambda}^{m_1,\mu,q}a_{-q}\vec{Y}_{j,l',q+\mu}(\hat{\mathbf{r}})
\end{aligned}$$
where $\boxed{K_{l',j,\lambda}^{m_1,\mu,q}}=\sqrt{\dfrac{\lambda(2\lambda+1)(2\lambda-1)}{2l'+1}}(-1)^q\langle\lambda-1,\mu-m_1;1,m_1\mid \lambda,\mu\rangle\langle 1,q;\lambda-1,\mu-m_1\mid l',q+\mu-m_1\rangle\langle 1,0;\lambda-1,0\mid l',0\rangle\langle l',q+\mu-m_1;1,m_1\mid j,q+\mu\rangle  $ and 
$$\vec{Y}_{l'-1,l',q+\mu}=\dfrac{1}{\sqrt{l'(2l'-1)}}V_{l'-1}^{q+\mu},$$
$$\vec{Y}_{l',l',q+\mu}=\dfrac{-i}{\sqrt{l'(l'+1)}}X_{l'}^{q+\mu},$$
$$\vec{Y}_{l'+1,l',q+\mu}=\dfrac{1}{\sqrt{(l'+1)(2l'+3)}}W_{l'+1}^{q+\mu}.$$
\end{theorem}
Given the addition theorem of ${r'}^{-l}V_l^m(\hat{\mathbf{r}'})$, it is not difficult to derive the addition theorem of ${r'}^{-l}W_l^m(\hat{\mathbf{r}'})$.
Note that
\begin{equation}\label{eq-W-V=Y}W_l^m(\hat{\mathbf{r}})-V_l^m(\hat{\mathbf{r}})=(2l+1)Y_l^m(\hat{\mathbf{r}})\hat{\mathbf{r}}.\end{equation}
It follows that
$${r'}^{-l}W_l^m(\hat{\mathbf{r}'})={r'}^{-l}V_l^m(\hat{\mathbf{r}'})+{r'}^{-l}(2l+1)Y_l^m(\hat{\mathbf{r}'})\hat{\mathbf{r}'}.$$
The term ${r'}^{-l}Y_l^m(\hat{\mathbf{r}'})\hat{\mathbf{r}'}$ can be decomposed as
\begin{equation}\label{eq-tornapart}
\begin{aligned}{r'}^{-l}Y_l^m(\hat{\mathbf{r}'})\hat{\mathbf{r}'}&=\sqrt{\dfrac{2l+1}{4\pi}}\mathbf{r}I_l^m(\hat{\mathbf{r}'})+\sqrt{\dfrac{2l+1}{4\pi}}\mathbf{a}I_l^m(\hat{\mathbf{r}'}).
\end{aligned}\end{equation}
Using the addition theorem of $I_l^m$:
$$I_l^m(\mathbf{r}+\mathbf{a})=\sum\limits_{\lambda=0}^\infty\sum\limits_{\mu=-\lambda}^\lambda\sqrt{\dfrac{2\lambda+1}{2l+1}}H_{l,\lambda}^{m,\mu}(\mathbf{a})R_\lambda^\mu(\mathbf{r}),$$
the first term becomes
$$\sqrt{\dfrac{2l+1}{4\pi}}\mathbf{r}I_l^m(\hat{\mathbf{r}'})=\sum\limits_{\lambda=0}^\infty\sum\limits_{\mu=-\lambda}^{\lambda}\dfrac{r^{\lambda+1}}{2\lambda+1}H_{l,\lambda}^{m,\mu}(\mathbf{a})(W_\lambda^\mu(\hat{\mathbf{r}})-V_\lambda^\mu(\hat{\mathbf{r}}))$$
by the Eq~\eqref{eq-W-V=Y} and the definition of $R_\lambda^\mu$.
\begin{lemma}Under the definition of $a_{-q}$ in Eq~\eqref{eq-a-q}, the constant vector $\mathbf{a}$ can be expanded as 
\begin{equation}\label{eq-a-expand}\mathbf{a}=\sum\limits_{q=-1}^1(-1)^qa_{-q}\chi_{1,q}.\end{equation}
\end{lemma}
By the addition theorem of $I_l^m$ and Eq~\eqref{eq-Ykl'}, the second term in Eq~\eqref{eq-tornapart} can be expressed as
$$\begin{aligned}
\sqrt{\dfrac{2l+1}{4\pi}}\mathbf{a}I_l^m(\hat{\mathbf{r}'})&=\sum\limits_{\lambda=0}^\infty\sum\limits_{\mu=-\lambda}^\lambda r^{\lambda} H_{l,\lambda}^{m,\mu}(\mathbf{a}) Y_\lambda^\mu(\hat{\mathbf{r}})\mathbf{a}\\
&=\sum\limits_{\lambda=0}^\infty\sum\limits_{\mu=-\lambda}^\lambda\sum\limits_{q=-1}^1\sum\limits_{\vert\lambda-1\vert\leqslant j\leqslant \vert\lambda+1\vert}(-1)^qa_{-q} r^{\lambda} H_{l,\lambda}^{m,\mu}(\mathbf{a}) \langle \lambda,\mu;1,q\mid j,\mu+q\rangle\vec{Y}_{j,\lambda,q+\mu}(\hat{\mathbf{r}}).
\end{aligned}$$
\begin{theorem}
$$\begin{aligned}
{r'}^{-l}W_l^m(\hat{\mathbf{r}'})
=&\sum\limits_{\lambda=1}^\infty \sum\limits_{\mu=-\lambda}^{\lambda}H_{l,\lambda}^{m,\mu}(\mathbf{a})(r^{\lambda+1}+a^2r^{\lambda-1})W_\lambda^\mu(\hat{\mathbf{r}})\\
&+\sum\limits_{\lambda=1}^\infty \sum\limits_{\mu=-\lambda}^{\lambda}\sum\limits_{q=-1}^{1}\sum\limits_{m_1=-1}^{1}\sum\limits_{\vert\lambda-2\vert \leqslant l'\leqslant\lambda,\atop l'\neq\lambda-1}\sum\limits_{\vert l'-1\vert\leqslant j\leqslant \vert l'+1\vert}
H_{l,\lambda}^{m,\mu}(\mathbf{a})2r^{\lambda}K_{l',j,\lambda}^{m_1,\mu,q}a_{-q}\vec{Y}_{j,l',q+\mu}(\hat{\mathbf{r}})\\
&+(2l+1)\sum\limits_{\lambda=0}^\infty\sum\limits_{\mu=-\lambda}^{\lambda}\dfrac{r^{\lambda+1}}{2\lambda+1}H_{l,\lambda}^{m,\mu}(\mathbf{a})(W_\lambda^\mu(\hat{\mathbf{r}})-V_\lambda^\mu(\hat{\mathbf{r}}))\\
&+(2l+1)\sum\limits_{\lambda=0}^\infty\sum\limits_{\mu=-\lambda}^\lambda\sum\limits_{q=-1}^1\sum\limits_{\vert \lambda-1\vert\leqslant j\leqslant \vert\lambda+1\vert}(-1)^{q}a_{-q} r^\lambda H_{l,\lambda}^{m,\mu}(\mathbf{a})\langle \lambda,\mu;1,q\mid j,\mu+q\rangle\vec{Y}_{j,\lambda,q+\mu}(\hat{\mathbf{r}}).
\end{aligned}$$
\end{theorem}

Thus, the real forms of addition theorems are given by
\begin{theorem}
When $m>0$,
$$\begin{aligned}
&{r'}^{-l}V_{lm}(\hat{\mathbf{r}'})=\sum_{\lambda=1}^\infty\sum_{\mu=-\lambda}^\lambda\dfrac{1}{\sqrt{2}}\bigl(H_{l,\lambda}^{-m,\mu}(\mathbf{a})+(-1)^m H_{l,\lambda}^{m,\mu}(\mathbf{a})\bigr)(r^{\lambda+1}+a^2 r^{\lambda-1})W_\lambda^\mu(\hat{\mathbf{r}})\\
&+\sum_{\lambda=1}^\infty\sum_{\mu=-\lambda}^\lambda\sum_{q=-1}^1\sum_{m_1=-1}^1\sum_{\substack{|\lambda-2|\leqslant k\leqslant\lambda\\k\neq\lambda-1}}\sum_{|k-1|\leqslant j\leqslant|k+1|}
2r^\lambda a_{-q}\dfrac{1}{\sqrt{2}}\bigl(H_{l,\lambda}^{-m,\mu}(\mathbf{a})+(-1)^m H_{l,\lambda}^{m,\mu}(\mathbf{a})\bigr) K_{k,j,\lambda}^{m_1,\mu,q}\vec{Y}_{j,k,q+\mu}(\hat{\mathbf{r}})
\end{aligned}$$
When $m<0$,
$$\begin{aligned}
&{r'}^{-l}V_{lm}(\hat{\mathbf{r}'})=\sum_{\lambda=1}^\infty\sum_{\mu=-\lambda}^\lambda\dfrac{i}{\sqrt{2}}\bigl(H_{l,\lambda}^{m,\mu}(\mathbf{a})-(-1)^m H_{l,\lambda}^{-m,\mu}(\mathbf{a})\bigr)(r^{\lambda+1}+a^2 r^{\lambda-1})W_\lambda^\mu(\hat{\mathbf{r}})\\
&+\sum_{\lambda=1}^\infty\sum_{\mu=-\lambda}^\lambda\sum_{q=-1}^1\sum_{m_1=-1}^1\sum_{\substack{|\lambda-2|\leqslant k\leqslant\lambda\\k\neq\lambda-1}}\sum_{|k-1|\leqslant j\leqslant|k+1|}
\dfrac{i}{\sqrt{2}}\bigl(H_{l,\lambda}^{m,\mu}(\mathbf{a})-(-1)^m H_{l,\lambda}^{-m,\mu}(\mathbf{a})\bigr)2r^\lambda K_{k,j,\lambda}^{m_1,\mu,q}a_{-q}\vec{Y}_{j,k,q+\mu}(\hat{\mathbf{r}}).
\end{aligned}$$
When $m=0$,
$$\begin{aligned}
&{r'}^{-l}V_{l0}(\hat{\mathbf{r}'})=\sum_{\lambda=1}^\infty\sum_{\mu=-\lambda}^\lambda H_{l,\lambda}^{0,\mu}(\mathbf{a})(r^{\lambda+1}+a^2 r^{\lambda-1})W_\lambda^\mu(\hat{\mathbf{r}})\\
&+\sum_{\lambda=1}^\infty\sum_{\mu=-\lambda}^\lambda\sum_{q=-1}^1\sum_{m_1=-1}^1\sum_{\substack{|\lambda-2|\leqslant k\leqslant\lambda\\k\neq\lambda-1}}\sum_{|k-1|\leqslant j\leqslant|k+1|}
H_{l,\lambda}^{0,\mu}(\mathbf{a})\,2r^\lambda K_{k,j,\lambda}^{m_1,\mu,q}a_{-q}\vec{Y}_{j,k,q+\mu}(\hat{\mathbf{r}}).
\end{aligned}$$
where
$$W_\lambda^\mu=\left\{
\begin{aligned}
\dfrac{(-1)^\mu}{\sqrt{2}}(W_{\lambda\mu}+i W_{\lambda,-\mu})&, \quad \mu>0,\\
W_{\lambda 0}\qquad\qquad &,\quad \mu=0,\\
\dfrac{1}{\sqrt{2}}(W_{\lambda,-\mu}-iW_{\lambda\mu})&,\quad \mu<0.
\end{aligned}\right.
$$
and 
$$\vec{Y}_{k-1,k,q+\mu}=
\left\{
\begin{aligned}
\dfrac{(-1)^{q+\mu}}{ \sqrt{2k(2k-1)} }\bigl(V_{k-1,q+\mu}+i\,V_{k-1,-(q+\mu)}\bigr)&,q+\mu>0\\
\dfrac{1}{\sqrt{k(2k-1)}}\,V_{k-1,0}\quad &,q+\mu=0\\
\dfrac{1}{ \sqrt{2k(2k-1)} }\bigl(V_{k-1,-(q+\mu)}-i\,V_{k-1,q+\mu}\bigr)&,q+\mu<0
\end{aligned}
\right.$$
$$\vec{Y}_{k,k,q+\mu}=
\left\{\begin{aligned}\dfrac{-i(-1)^{q+\mu}}{ \sqrt{2k(k+1)} }\bigl(X_{k,q+\mu}+i\,X_{k,-(q+\mu)}\bigr)&,q+\mu>0\\
\dfrac{-i}{\sqrt{k(k+1)}}\,X_{k0}\quad &,q+\mu=0\\
\dfrac{-i}{ \sqrt{2k(k+1)} }\bigl(X_{k,-(q+\mu)}-i\,X_{k,q+\mu}\bigr)&,q+\mu<0
\end{aligned}\right.
$$
$$\vec{Y}_{k+1,k,q+\mu}=
\left\{\begin{aligned}\dfrac{(-1)^{q+\mu}}{ \sqrt{2(k+1)(2k+3)} }\bigl(W_{k+1,q+\mu}+i\,W_{k+1,-(q+\mu)}\bigr)&,q+\mu>0\\
\dfrac{1}{\sqrt{(k+1)(2k+3)}}\,W_{k+1,0}\quad &,q+\mu=0\\
\dfrac{1}{ \sqrt{2(k+1)(2k+3)} }\bigl(W_{k+1,-(q+\mu)}-i\,W_{k+1,q+\mu}\bigr)&,q+\mu<0
\end{aligned}\right.$$
\end{theorem}

\begin{theorem}
When $m>0$,
$$\begin{aligned}&{r'}^{-l}W_{lm}(\hat{\mathbf{r}'})
=\sum_{\lambda=1}^\infty\sum_{\mu=-\lambda}^\lambda\dfrac{1}{\sqrt{2}}\bigl(H_{l,\lambda}^{-m,\mu}(\mathbf{a})+(-1)^m H_{l,\lambda}^{m,\mu}(\mathbf{a})\bigr)(r^{\lambda+1}+a^2 r^{\lambda-1})W_\lambda^\mu(\hat{\mathbf{r}})\\
+&\sum_{\lambda=1}^\infty\sum_{\mu=-\lambda}^\lambda\sum_{q=-1}^1\sum_{m_1=-1}^1\sum_{\substack{|\lambda-2|\leqslant k\leqslant\lambda\\k\neq\lambda-1}}\sum_{|k-1|\leqslant j\leqslant|k+1|}
\dfrac{1}{\sqrt{2}}2r^\lambda a_{-q}\bigl(H_{l,\lambda}^{-m,\mu}(\mathbf{a})+(-1)^m H_{l,\lambda}^{m,\mu}(\mathbf{a})\bigr) K_{k,j,\lambda}^{m_1,\mu,q}\vec{Y}_{j,k,q+\mu}(\hat{\mathbf{r}})\\
+&(2l+1)\sum_{\lambda=0}^\infty\sum_{\mu=-\lambda}^\lambda\dfrac{r^{\lambda+1}}{2\lambda+1}\dfrac{1}{\sqrt{2}}\bigl(H_{l,\lambda}^{-m,\mu}(\mathbf{a})+(-1)^m H_{l,\lambda}^{m,\mu}(\mathbf{a})\bigr)(W_\lambda^\mu(\hat{\mathbf{r}})-V_\lambda^\mu(\hat{\mathbf{r}}))\\
+&(2l+1)\sum_{\lambda=0}^\infty\sum_{\mu=-\lambda}^\lambda\sum_{q=-1}^1\sum_{|\lambda-1|\leqslant j\leqslant|\lambda+1|}
\dfrac{(-1)^qa_{-q}r^\lambda}{\sqrt{2}}\bigl(H_{l,\lambda}^{-m,\mu}(\mathbf{a})+(-1)^m H_{l,\lambda}^{m,\mu}(\mathbf{a})\bigr)  \\
&\langle\lambda,\mu;1,q\mid j,\mu+q\rangle\vec{Y}_{j,\lambda,q+\mu}(\hat{\mathbf{r}}).
\end{aligned}$$
When $m<0$,
$$ \begin{aligned}
&{r'}^{-l}W_{lm}(\hat{\mathbf{r}'})
=\sum_{\lambda=1}^\infty\sum_{\mu=-\lambda}^\lambda\dfrac{i}{\sqrt{2}}\bigl(H_{l,\lambda}^{m,\mu}(\mathbf{a})-(-1)^m H_{l,\lambda}^{-m,\mu}(\mathbf{a})\bigr)(r^{\lambda+1}+a^2 r^{\lambda-1})W_\lambda^\mu(\hat{\mathbf{r}})\\
+&\sum_{\lambda=1}^\infty\sum_{\mu=-\lambda}^\lambda\sum_{q=-1}^1\sum_{m_1=-1}^1\sum_{\substack{|\lambda-2|\leqslant k\leqslant\lambda\\k\neq\lambda-1}}\sum_{|k-1|\leqslant j\leqslant|k+1|}
\dfrac{i}{\sqrt{2}}\bigl(H_{l,\lambda}^{m,\mu}(\mathbf{a})-(-1)^m H_{l,\lambda}^{-m,\mu}(\mathbf{a})\bigr)2r^\lambda K_{k,j,\lambda}^{m_1,\mu,q}a_{-q}\vec{Y}_{j,k,q+\mu}(\hat{\mathbf{r}})\\
+&(2l+1)\sum_{\lambda=0}^\infty\sum_{\mu=-\lambda}^\lambda\dfrac{r^{\lambda+1}}{2\lambda+1}\dfrac{i}{\sqrt{2}}\bigl(H_{l,\lambda}^{m,\mu}(\mathbf{a})-(-1)^m H_{l,\lambda}^{-m,\mu}(\mathbf{a})\bigr)(W_\lambda^\mu(\hat{\mathbf{r}})-V_\lambda^\mu(\hat{\mathbf{r}}))\\
+&(2l+1)\sum_{\lambda=0}^\infty\sum_{\mu=-\lambda}^\lambda\sum_{q=-1}^1\sum_{|\lambda-1|\leqslant j\leqslant|\lambda+1|}
\dfrac{i(-1)^q a_{-q} r^\lambda}{\sqrt{2}}\bigl(H_{l,\lambda}^{m,\mu}(\mathbf{a})-(-1)^m H_{l,\lambda}^{-m,\mu}(\mathbf{a})\bigr)\\&\langle\lambda,\mu;1,q\mid j,\mu+q\rangle\vec{Y}_{j,\lambda,q+\mu}(\hat{\mathbf{r}}).
\end{aligned}    $$
When $   m=0   $,
$$    \begin{aligned}
{r'}^{-l}W_{l0}(\hat{\mathbf{r}'})
=&\sum_{\lambda=1}^\infty\sum_{\mu=-\lambda}^\lambda H_{l,\lambda}^{0,\mu}(\mathbf{a})(r^{\lambda+1}+a^2 r^{\lambda-1})W_\lambda^\mu(\hat{\mathbf{r}})\\
&+\sum_{\lambda=1}^\infty\sum_{\mu=-\lambda}^\lambda\sum_{q=-1}^1\sum_{m_1=-1}^1\sum_{\substack{|\lambda-2|\leqslant k\leqslant\lambda\\k\neq\lambda-1}}\sum_{|k-1|\leqslant j\leqslant|k+1|}
H_{l,\lambda}^{0,\mu}(\mathbf{a})\,2r^\lambda K_{k,j,\lambda}^{m_1,\mu,q}a_{-q}\vec{Y}_{j,k,q+\mu}(\hat{\mathbf{r}})\\
&+(2l+1)\sum_{\lambda=0}^\infty\sum_{\mu=-\lambda}^\lambda\dfrac{r^{\lambda+1}}{2\lambda+1}H_{l,\lambda}^{0,\mu}(\mathbf{a})(W_\lambda^\mu(\hat{\mathbf{r}})-V_\lambda^\mu(\hat{\mathbf{r}}))\\
&+(2l+1)\sum_{\lambda=0}^\infty\sum_{\mu=-\lambda}^\lambda\sum_{q=-1}^1\sum_{|\lambda-1|\leqslant j\leqslant|\lambda+1|}
H_{l,\lambda}^{0,\mu}(\mathbf{a})(-1)^q a_{-q} r^\lambda\langle\lambda,\mu;1,q\mid j,\mu+q\rangle\vec{Y}_{j,\lambda,q+\mu}(\hat{\mathbf{r}}).
\end{aligned}    $$
where
$$    \vec{Y}_{\lambda-1,\lambda,q+\mu}=
\left\{
\begin{aligned}
\dfrac{(-1)^{q+\mu}}{ \sqrt{2\lambda(2\lambda-1)} }\bigl(V_{\lambda-1,q+\mu}+i\,V_{\lambda-1,-(q+\mu)}\bigr)&,q+\mu>0\\
\dfrac{1}{\sqrt{\lambda(2\lambda-1)}}\,V_{\lambda-1,0}\quad &,q+\mu=0\\
\dfrac{1}{ \sqrt{2\lambda(2\lambda-1)} }\bigl(V_{\lambda-1,-(q+\mu)}-i\,V_{\lambda-1,q+\mu}\bigr)&,q+\mu<0
\end{aligned}
\right.    $$
$$    \vec{Y}_{\lambda,\lambda,q+\mu}=
\left\{\begin{aligned}\dfrac{-i(-1)^{q+\mu}}{ \sqrt{2\lambda(\lambda+1)} }\bigl(X_{\lambda,q+\mu}+i\,X_{\lambda,-(q+\mu)}\bigr)&,q+\mu>0\\
\dfrac{-i}{\sqrt{\lambda(\lambda+1)}}\,X_{\lambda0}\quad &,q+\mu=0\\
\dfrac{-i}{ \sqrt{2\lambda(\lambda+1)} }\bigl(X_{\lambda,-(q+\mu)}-i\,X_{\lambda,q+\mu}\bigr)&,q+\mu<0
\end{aligned}\right.    $$
$$    \vec{Y}_{\lambda+1,\lambda,q+\mu}=
\left\{\begin{aligned}\dfrac{(-1)^{q+\mu}}{ \sqrt{2(\lambda+1)(2\lambda+3)} }\bigl(W_{\lambda+1,q+\mu}+i\,W_{\lambda+1,-(q+\mu)}\bigr)&,q+\mu>0\\
\dfrac{1}{\sqrt{(\lambda+1)(2\lambda+3)}}\,W_{\lambda+1,0}\quad &,q+\mu=0\\
\dfrac{1}{ \sqrt{2(\lambda+1)(2\lambda+3)} }\bigl(W_{\lambda+1,-(q+\mu)}-i\,W_{\lambda+1,q+\mu}\bigr)&,q+\mu<0
\end{aligned}\right.    $$
\end{theorem}

\subsection{Addition theorem for ${r'}^{-l-1}X_{lm}$}
Note that   the coefficient of the addition theorem for $X_l^m$ that multiplies the bracket $  [r\,X_\lambda^\mu+\mathbf{a}\times W_\lambda^\mu]  $ is identical to the coefficient that appeared in the $  V_{lm}  $ case. Therefore, the following addition theorem is obtained.\begin{theorem}
When $ m>0 $,
$$\begin{aligned}
&{r'}^{-l-1}X_{lm}(\mathbf{r}')\\
=&\sum\limits_{\lambda=1}^\infty\sum\limits_{\mu=-\lambda}^{-1} r^{\lambda-1}\dfrac{1}{2}\bigl(H_{l,\lambda}^{-m,\mu}(\mathbf{a})+(-1)^mH_{l,\lambda}^{m,\mu}(\mathbf{a})\bigr)\Bigl\{r[X_{\lambda,-\mu}(\hat{\mathbf{r}})-iX_{\lambda,\mu}(\hat{\mathbf{r}})]+\mathbf{a}\times[W_{\lambda,-\mu}(\hat{\mathbf{r}})-iW_{\lambda,\mu}(\hat{\mathbf{r}})]\Bigr\}\\
+&\sum\limits_{\lambda=1}^\infty\sum\limits_{\mu=1}^\lambda r^{\lambda-1}\dfrac{(-1)^{\mu}}{2}\bigl(H_{l,\lambda}^{-m,\mu}(\mathbf{a})+(-1)^mH_{l,\lambda}^{m,\mu}(\mathbf{a})\bigr)\Bigl\{r[X_{\lambda,\mu}(\hat{\mathbf{r}})+iX_{\lambda,-\mu}(\hat{\mathbf{r}})]+\mathbf{a}\times[W_{\lambda,\mu}(\hat{\mathbf{r}})+iW_{\lambda,-\mu}(\hat{\mathbf{r}})]\Bigr\}\\
&\quad+\sum\limits_{\lambda=1}^\infty r^{\lambda-1}\dfrac{1}{\sqrt{2}}\bigl(H_{l,\lambda}^{-m,0}(\mathbf{a})+(-1)^mH_{l,\lambda}^{m,0}(\mathbf{a})\bigr)\Bigl\{r\,X_{\lambda,0}(\hat{\mathbf{r}})+\mathbf{a}\times W_{\lambda,0}(\hat{\mathbf{r}})\Bigr\}.
\end{aligned}$$
When $  m<0  $,
$$\begin{aligned}
&{r'}^{-l-1}X_{lm}(\mathbf{r}')\\
=&\sum\limits_{\lambda=1}^\infty\sum\limits_{\mu=-\lambda}^{-1} r^{\lambda-1}\dfrac{i}{2}\bigl(H_{l,\lambda}^{m,\mu}(\mathbf{a})-(-1)^mH_{l,\lambda}^{-m,\mu}(\mathbf{a})\bigr)\Bigl\{r[X_{\lambda,-\mu}(\hat{\mathbf{r}})-iX_{\lambda,\mu}(\hat{\mathbf{r}})]+\mathbf{a}\times[W_{\lambda,-\mu}(\hat{\mathbf{r}})-iW_{\lambda,\mu}(\hat{\mathbf{r}})]\Bigr\}\\
+&\sum\limits_{\lambda=1}^\infty\sum\limits_{\mu=1}^\lambda r^{\lambda-1}\dfrac{(-1)^{\mu}}{2}\cdot i\bigl(H_{l,\lambda}^{m,\mu}(\mathbf{a})-(-1)^mH_{l,\lambda}^{-m,\mu}(\mathbf{a})\bigr)\Bigl\{r[X_{\lambda,\mu}(\hat{\mathbf{r}})+iX_{\lambda,-\mu}(\hat{\mathbf{r}})]\\&+\mathbf{a}\times[W_{\lambda,\mu}(\hat{\mathbf{r}})+
iW_{\lambda,-\mu}(\hat{\mathbf{r}})]\Bigr\}\\
&\quad+\sum\limits_{\lambda=1}^\infty r^{\lambda-1}\dfrac{i}{\sqrt{2}}\bigl(H_{l,\lambda}^{m,0}(\mathbf{a})-(-1)^mH_{l,\lambda}^{-m,0}(\mathbf{a})\bigr)\Bigl\{r\,X_{\lambda,0}(\hat{\mathbf{r}})+\mathbf{a}\times W_{\lambda,0}(\hat{\mathbf{r}})\Bigr\}.
\end{aligned}$$
When $  m=0  $,
$$\begin{aligned}
{r'}^{-l-1}X_{l0}(\mathbf{r}')
&=\sum\limits_{\lambda=1}^\infty\sum\limits_{\mu=-\lambda}^{-1}
r^{\lambda-1}\dfrac{1}{\sqrt{2}}H_{l,\lambda}^{0,\mu}(\mathbf{a})\Bigl\{r[X_{\lambda,-\mu}(\hat{\mathbf{r}})-iX_{\lambda,\mu}(\hat{\mathbf{r}})]+\mathbf{a}\times[W_{\lambda,-\mu}(\hat{\mathbf{r}})-iW_{\lambda,\mu}(\hat{\mathbf{r}})]\Bigr\}\\
&\quad+\sum\limits_{\lambda=1}^\infty\sum\limits_{\mu=1}^{\lambda}
r^{\lambda-1}\dfrac{(-1)^\mu}{\sqrt{2}}H_{l,\lambda}^{0,\mu}(\mathbf{a})\Bigl\{r[X_{\lambda,\mu}(\hat{\mathbf{r}})+iX_{\lambda,-\mu}(\hat{\mathbf{r}})]+\mathbf{a}\times[W_{\lambda,\mu}(\hat{\mathbf{r}})+iW_{\lambda,-\mu}(\hat{\mathbf{r}})]\Bigr\}\\
&\quad+\sum\limits_{\lambda=1}^\infty H_{l,\lambda}^{0,0}(\mathbf{a})\Bigl\{r\,X_{\lambda,0}(\hat{\mathbf{r}})+\mathbf{a}\times W_{\lambda,0}(\hat{\mathbf{r}})\Bigr\}.
\end{aligned}$$
\end{theorem}

Note that it is difficult to compute the inner product involving the cross product $\mathbf{a}\times W_{\lambda}^\mu$ in Eq.~\eqref{eq-addcomplex-X}.
Therefore, the cross product must be converted to the linear combination of vector spherical harmonics. 
\begin{lemma}For the spherical basis 
$$\chi_{1,\pm 1}=\dfrac{1}{\sqrt{2}}\begin{pmatrix}\mp 1\\-i\\0\end{pmatrix},\quad \chi_{1,0}=\begin{pmatrix}\mp 1\\0\\0\end{pmatrix},$$
the formula for their  cross product is given by
$$\chi_{1, m} \times \chi_{1, n} = i \cdot \text{sgn}(m - n) \chi_{1, m+n}.$$
\end{lemma}
Using Eq~\eqref{eq-W-vecY}, Eq~\eqref{eq-vecY-expand}, Eq~\eqref{eq-a-expand}, and Eq~\eqref{eq-Ykl'}, it yields that
$$\begin{aligned}\mathbf{a}\times W_\lambda^{\mu}(\hat{\mathbf{r}})&=\sum\limits_{q=-1}^1\sum\limits_{m_1=-1}^1\sum\limits_{\vert\lambda-2\vert\leqslant j\leqslant\vert \lambda\vert}i (-1)^q\text{sgn}(q-m_1)a_{-q}\sqrt{\lambda(2\lambda+1)}\\
&\langle\lambda-1,\mu-m_1;1,m_1\mid \lambda,\mu\rangle
\langle \lambda-1,\mu-m_1;1,q+m_1\mid j,\mu+q \rangle \vec{Y}_{j,\lambda-1,q+\mu}(\hat{\mathbf{r}})
\end{aligned}$$
Therefore, the addition theorem for $X_l^m$ \eqref{eq-addcomplex-X} reads

\begin{equation}\label{eq-addthm-X}
\begin{aligned}&{r'}^{-l-1}X_l^m(\hat{\mathbf{r}'})
=\sum\limits_{\lambda=1}^\infty\sum\limits_{\mu=-\lambda}^\lambda(-1)^{\lambda+\mu}\sqrt{\dfrac{2l+1}{2\lambda+1}} r^{\lambda}I_{l+\lambda}^{m-\mu}(\mathbf{a})
\sqrt{{l+\lambda+\mu-m\choose \lambda+\mu}{l+\lambda+m-\mu\choose \lambda-\mu}}X_\lambda^\mu(\hat{\mathbf{r}})\\
&+\sum\limits_{\lambda=1}^\infty\sum\limits_{\mu=-\lambda}^\lambda\sum\limits_{q=-1}^1\sum\limits_{m_1=-1}^1\sum\limits_{\vert\lambda-2\vert\leqslant j\leqslant\vert \lambda\vert}
i (-1)^{\lambda+\mu+q}\text{sgn}(q-m_1)a_{-q}\sqrt{\lambda(2l+1)} r^{\lambda-1}I_{l+\lambda}^{m-\mu}(\mathbf{a})\\
&\sqrt{{l+\lambda+\mu-m\choose \lambda+\mu}{l+\lambda+m-\mu\choose \lambda-\mu}}
\langle\lambda-1,\mu-m_1;1,m_1\mid \lambda,\mu\rangle\\&
\langle \lambda-1,\mu-m_1;1,q+m_1\mid j,\mu+q \rangle \vec{Y}_{j,\lambda-1,q+\mu}(\hat{\mathbf{r}})
\end{aligned}
\end{equation}

Denote the coefficient 
$$\begin{aligned}\boxed{L_{l,j,\lambda}^{m,\mu,q,m_1}(\mathbf{a})}=
&i (-1)^{\lambda+\mu+q}\text{sgn}(q-m_1)a_{-q}\sqrt{\lambda(2l+1)} I_{l+\lambda}^{m-\mu}(\mathbf{a}) 
\sqrt{{l+\lambda+\mu-m\choose \lambda+\mu}{l+\lambda+m-\mu\choose \lambda-\mu}}\\
&\langle\lambda-1,\mu-m_1;1,m_1\mid \lambda,\mu\rangle
\langle \lambda-1,\mu-m_1;1,q+m_1\mid j,\mu+q \rangle
\end{aligned}$$

Then  Eq~\eqref{eq-addthm-X} becomes

\begin{equation}
\begin{aligned}{r'}^{-l-1}X_l^m(\hat{\mathbf{r}'})
=&\sum\limits_{\lambda=1}^\infty\sum\limits_{\mu=-\lambda}^\lambda r^{\lambda} H_{l,\lambda}^{m,\mu}(\mathbf{a})X_\lambda^\mu(\hat{\mathbf{r}})\\
&+\sum\limits_{\lambda=1}^\infty\sum\limits_{\mu=-\lambda}^\lambda\sum\limits_{q=-1}^1\sum\limits_{m_1=-1}^1\sum\limits_{\vert\lambda-2\vert\leqslant j\leqslant\vert \lambda\vert}
r^{\lambda-1}L_{l,j,\lambda}^{m,\mu,q,m_1}(\mathbf{a})\vec{Y}_{j,\lambda-1,q+\mu}(\hat{\mathbf{r}})
\end{aligned}
\end{equation}
Hence, the real valued addition theorem for ${r'}^{-l-1}X_l^m(\hat{\mathbf{r}'})$ is given by

\begin{theorem}
When $m>0$,

\begin{equation}
\begin{aligned}{r'}^{-l-1}X_{lm}(\hat{\mathbf{r}'})
=&\sum\limits_{\lambda=1}^\infty\sum\limits_{\mu=-\lambda}^\lambda r^{\lambda} \dfrac{1}{\sqrt{2}}(H_{l,\lambda}^{-m,\mu}(\mathbf{a})+(-1)^m H_{l,\lambda}^{m,\mu}(\mathbf{a}))X_\lambda^\mu(\hat{\mathbf{r}})\\
&+\sum\limits_{\lambda=1}^\infty\sum\limits_{\mu=-\lambda}^\lambda\sum\limits_{q=-1}^1\sum\limits_{m_1=-1}^1\sum\limits_{\vert\lambda-2\vert\leqslant j\leqslant\vert \lambda\vert}
r^{\lambda-1}\dfrac{1}{\sqrt{2}}(L_{l,j,\lambda}^{-m,\mu,q,m_1}(\mathbf{a})\\&+(-1)^mL_{l,j,\lambda}^{m,\mu,q,m_1}(\mathbf{a}))\vec{Y}_{j,\lambda-1,q+\mu}(\hat{\mathbf{r}})
\end{aligned}
\end{equation}
When $m<0$,

\begin{equation}
\begin{aligned}{r'}^{-l-1}X_{lm}(\hat{\mathbf{r}'})
=&\sum\limits_{\lambda=1}^\infty\sum\limits_{\mu=-\lambda}^\lambda r^{\lambda} \dfrac{i}{\sqrt{2}}(H_{l,\lambda}^{m,\mu}(\mathbf{a})-(-1)^m H_{l,\lambda}^{-m,\mu}(\mathbf{a}))X_\lambda^\mu(\hat{\mathbf{r}})\\
&+\sum\limits_{\lambda=1}^\infty\sum\limits_{\mu=-\lambda}^\lambda\sum\limits_{q=-1}^1\sum\limits_{m_1=-1}^1\sum\limits_{\vert\lambda-2\vert\leqslant j\leqslant\vert \lambda\vert}
r^{\lambda-1}\dfrac{i}{\sqrt{2}}(L_{l,j,\lambda}^{m,\mu,q,m_1}(\mathbf{a})\\
&-(-1)^mL_{l,j,\lambda}^{-m,\mu,q,m_1}(\mathbf{a}))\vec{Y}_{j,\lambda-1,q+\mu}(\hat{\mathbf{r}})
\end{aligned}
\end{equation}

When $m=0$,

\begin{equation}
\begin{aligned}{r'}^{-l-1}X_{l0}(\hat{\mathbf{r}'})
=&\sum\limits_{\lambda=1}^\infty\sum\limits_{\mu=-\lambda}^\lambda r^{\lambda} H_{l,\lambda}^{0,\mu}(\mathbf{a})X_\lambda^\mu(\hat{\mathbf{r}})\\
&+\sum\limits_{\lambda=1}^\infty\sum\limits_{\mu=-\lambda}^\lambda\sum\limits_{q=-1}^1\sum\limits_{m_1=-1}^1\sum\limits_{\vert\lambda-2\vert\leqslant j\leqslant\vert \lambda\vert}
r^{\lambda-1}L_{l,j,\lambda}^{0,\mu,q,m_1}(\mathbf{a})\vec{Y}_{j,\lambda-1,q+\mu}(\hat{\mathbf{r}})
\end{aligned}
\end{equation}

where
$$    \vec{Y}_{\lambda-2,\lambda-1,q+\mu}=
\left\{
\begin{aligned}
\dfrac{(-1)^{q+\mu}}{ \sqrt{2(\lambda-1)(2\lambda-3)} }\bigl(V_{\lambda-2,q+\mu}+i\,V_{\lambda-2,-(q+\mu)}\bigr)&,q+\mu>0\\
\dfrac{1}{\sqrt{(\lambda-1)(2\lambda-3)}}\,V_{\lambda-2,0}\quad &,q+\mu=0\\
\dfrac{1}{ \sqrt{2(\lambda-1)(2\lambda-3)} }\bigl(V_{\lambda-2,-(q+\mu)}-i\,V_{\lambda-2,q+\mu}\bigr)&,q+\mu<0
\end{aligned}
\right.    $$
$$    \vec{Y}_{\lambda-1,\lambda-1,q+\mu}=
\left\{\begin{aligned}\dfrac{-i(-1)^{q+\mu}}{ \sqrt{2(\lambda-1)(\lambda)} }\bigl(X_{\lambda-1,q+\mu}+i\,X_{\lambda-1,-(q+\mu)}\bigr)&,q+\mu>0\\
\dfrac{-i}{\sqrt{(\lambda-1)(\lambda)}}\,X_{\lambda-1,0}\quad &,q+\mu=0\\
\dfrac{-i}{ \sqrt{2(\lambda-1)(\lambda)} }\bigl(X_{\lambda-1,-(q+\mu)}-i\,X_{\lambda-1,q+\mu}\bigr)&,q+\mu<0
\end{aligned}\right.    $$
$$    \vec{Y}_{\lambda,\lambda-1,q+\mu}=
\left\{\begin{aligned}\dfrac{(-1)^{q+\mu}}{ \sqrt{2(\lambda)(2\lambda+1)} }\bigl(W_{\lambda,q+\mu}+i\,W_{\lambda,-(q+\mu)}\bigr)&,q+\mu>0\\
\dfrac{1}{\sqrt{(\lambda)(2\lambda+1)}}\,W_{\lambda,0}\quad &,q+\mu=0\\
\dfrac{1}{ \sqrt{2(\lambda)(2\lambda+1)} }\bigl(W_{\lambda,-(q+\mu)}-i\,W_{\lambda,q+\mu}\bigr)&,q+\mu<0
\end{aligned}\right.    $$
\end{theorem}

\section{Inner product results}
Define the inner product over $\mathbb{S}^2$ by
$$(f,g)=\int_{\mathbb{S}^2}f\cdot g.$$

By \citet{StammXiang2022}, the orthogonal properties are given by
$$\begin{aligned}
&\int_{\mathbb{S}^2} V_{l m} \cdot W_{l' m'} = 0, 
\qquad  
\int_{\mathbb{S}^2} W_{l m} \cdot X_{l' m'} = 0,\qquad  \int_{\mathbb{S}^2} X_{l m} \cdot V_{l' m'} = 0, \\
&\int_{\mathbb{S}^2} V_{l m} \cdot V_{l' m'} 
= \delta_{ll'} \delta_{m m'} (l+1)(2l+1), \quad 
 \int_{\mathbb{S}^2} W_{l m} \cdot W_{l' m'} 
= \delta_{ll'} \delta_{m m'} l(2l+1), \\
&\int_{\mathbb{S}^2} X_{l m} \cdot X_{l' m'} 
= \delta_{ll'} \delta_{m m'} l(l+1).
\end{aligned}$$

\subsection{Main lemma}
The following two lemmas are proved in \citet{StammXiang2022}:
\begin{lemma}The family of vector spherical harmonics gives a complete basis of $L^2(\partial B_\rho(x_0))^3$ and any real function $f\in L^2(\partial B_\rho(x_0))^3$ can be represented as
$$f(x)=\sum\limits_{l=0}^\infty\sum\limits_{m=-l}^m \sum\limits_{k=1}^3F_{lm}^kY_{lm}^k\left(\dfrac{x-x_0}{\rho}\right)$$

\end{lemma}

\begin{lemma}\label{lem-main}Let $\underline{Y_{lm}}=(V_{lm},W_{lm},X_{lm}).$ Then\\
when $ |x| > 1$, it holds that 
$$
(\mathcal{S}\underline{Y_{l m}})(x)
=
\underline{Y_{l m}}(\hat{x})
A^{out}_{\mathcal{S},l}(x),
$$
where the matrix $ A^{out}_{\mathcal{S},l}(x)$ is given by
$$
A^{out}_{\mathcal{S},l}(x)
=
\begin{bmatrix}
\dfrac{(3l+1)\mu+l\lambda}
{(2l+3)(2l+1)\mu(2\mu+\lambda)}
|x|^{-l-2}
&
\dfrac{l(\mu+\lambda)}
{2(2l+1)\mu(2\mu+\lambda)}
\left(|x|^{-l-2}-|x|^{-l}\right)
&
0
\\[12pt]
0
&
\dfrac{(3l+2)\mu+(l+1)\lambda}
{(2l-1)(2l+1)\mu(2\mu+\lambda)}
|x|^{-l}
&
0
\\[12pt]
0
&
0
&
\dfrac{1}{(2l+1)\mu}
|x|^{-l-1}
\end{bmatrix}.
$$
\end{lemma}
Therefore, applying the Lemma for  $\partial D+(n,0,0)=\partial B_\rho(n,0,0)$ with $\rho<1/2$, $n\in\mathbb{Z}$ and $n\neq 0$:
\begin{equation}\label{eq-S-D+nY}\mathcal{S}_{D+n}[Y^k_{lm}]=\rho \sum\limits_{j=1}^3Y^j_{lm}(\widehat{x-\mathbf{a}})\left[A_{\mathcal{S},l}^{out}\left(\dfrac{x-\mathbf{a}}{\rho}\right)\right]_j^{k},\end{equation}
where the  superscript $k$ of $ A_{\mathcal{S},l}^{out}\left(\dfrac{x-\mathbf{a}}{\rho}\right) $ stands for the $k$th column of the matrix and the subscript $j$ the $j$th row. Let $a_{ij}$ stand for the coefficient of the entries of $A^{out}_{\mathcal{S},l}(x)$.\\
\subsection{Main results}
It is clear that $$\boxed{(V_{l'm'},\mathcal{S}_{D+n}[V_{lm}])=0},\quad \boxed{(X_{l'm'},\mathcal{S}_{D+n}[V_{lm}])=0},\quad \boxed{(W_{00},\mathcal{S}_{D+n}[V_{lm}])=0}.$$
Let $\mathbf{b}=-\mathbf{a}$.
Now compute $(W_{l'm'},\mathcal{S}_{D+n}[V_{lm}])$ using Eq~\eqref{eq-S-D+nY}.
\begin{theorem}The expression of $\boxed{(W_{l'm'},\mathcal{S}_{D+n}[V_{lm}])}$($l'\geqslant 1$) is given by:
\\
If $m>0$ and $m'>0$,
$$\begin{aligned}(W_{l'm'},\mathcal{S}_{D+n}[V_{lm}])= \rho^{l+3}a_{11}[r^{l'-1}\dfrac{1}{2}(H_{l,l'}^{-m,-m'}(\mathbf{b})+(-1)^mH_{l,l'}^{m,-m'}(\mathbf{b}))l'(2l'+1)\\
+r^{l'-1}\dfrac{(-1)^{m'}}{2}(H_{l,l'}^{-m,m'}(\mathbf{b})+(-1)^mH_{l,l'}^{m,m'}(\mathbf{b}))l'(2l'+1)]
\end{aligned}$$
If $m>0$ and $m'<0$,
$$\begin{aligned}(W_{l'm'},\mathcal{S}_{D+n}[V_{lm}])= \rho^{l+3}a_{11}[r^{l'-1}\dfrac{1}{2}(H_{l,l'}^{-m,m'}(\mathbf{b})+(-1)^mH_{l,l'}^{m,m'}(\mathbf{b}))(-i \cdot l'(2l'+1))\\
+r^{l'-1}\dfrac{(-1)^{m'}}{2}(H_{l,l'}^{-m,-m'}(\mathbf{b})+(-1)^mH_{l,l'}^{m,-m'}(\mathbf{b}))(i\cdot l'(2l'+1))]
\end{aligned}$$
If $m>0$ and $m'=0$,
$$\begin{aligned}(W_{l'm'},\mathcal{S}_{D+n}[V_{lm}])= \rho^{l+3}a_{11}r^{l'-1}\dfrac{1}{\sqrt{2}}(H_{l,l'}^{-m,0}(\mathbf{b})+(-1)^mH_{l,l'}^{m,0}(\mathbf{b}))l'(2l'+1).
\end{aligned}$$

If $m<0$ and $m'>0$,
$$\begin{aligned}(W_{l'm'},\mathcal{S}_{D+n}[V_{lm}])= \rho^{l+3}a_{11}[r^{l'-1}\dfrac{i}{2}(H_{l,l'}^{m,-m'}(\mathbf{b})-(-1)^m H_{l,l'}^{-m,-m'}(\mathbf{b}))l'(2l'+1)\\
+r^{l'-1}\dfrac{(-1)^{m'}}{2}\cdot i (H_{l,l'}^{m,m'}(\mathbf{b})-(-1)^mH_{l,l'}^{-m,m'}(\mathbf{b}))l'(2l'+1)]
\end{aligned}$$

If $m<0$ and $m'<0$,
$$\begin{aligned}(W_{l'm'},\mathcal{S}_{D+n}[V_{lm}])= \rho^{l+3}a_{11}[r^{l'-1}\dfrac{i}{2}(H_{l,l'}^{m,m'}(\mathbf{b})-(-1)^m H_{l,l'}^{-m,m'}(\mathbf{b}))(-i\cdot l'(2l'+1))\\
+r^{l'-1}\dfrac{(-1)^{m'}}{2}\cdot i (H_{l,l'}^{m,-m'}(\mathbf{b})-(-1)^mH_{l,l'}^{-m,-m'}(\mathbf{b}))(i\cdot l'(2l'+1))]
\end{aligned}$$

If $m<0$ and $m'=0$,
$$\begin{aligned}(W_{l'm'},\mathcal{S}_{D+n}[V_{lm}])= \rho^{l+3}a_{11}r^{l'-1}\dfrac{i}{\sqrt{2}}(H_{l,l'}^{m,0}(\mathbf{b})-(-1)^mH_{l,l'}^{-m,0}(\mathbf{b}))l'(2l'+1).\end{aligned}$$

If $m=0$ and $m'>0$,
$$\begin{aligned}(W_{l'm'},\mathcal{S}_{D+n}[V_{lm}])= \rho^{l+3}a_{11}r^{l'-1}
[\dfrac{1}{\sqrt{2}} H_{l,l'}^{0,-m'}(\mathbf{b})l'(2l'+1)+\dfrac{(-1)^{m'}}{\sqrt{2}}H_{l,l'}^{0,m'}(\mathbf{b})l'(2l'+1)]
\end{aligned}$$

If $m=0$ and $m'<0$,
$$\begin{aligned}(W_{l'm'},\mathcal{S}_{D+n}[V_{lm}])= \rho^{l+3}a_{11}r^{l'-1}
[\dfrac{1}{\sqrt{2}} H_{l,l'}^{0,m'}(\mathbf{b})(-i\cdot l'(2l'+1))+\dfrac{(-1)^{m'}}{\sqrt{2}}H_{l,l'}^{0,-m'}(\mathbf{b})(i\cdot l'(2l'+1))]
\end{aligned}$$

If $m=0$ and $m'=0$,
$$\begin{aligned}(W_{l'm'},\mathcal{S}_{D+n}[V_{lm}])= \rho^{l+3}a_{11}r^{l'-1}
[H_{l,l'}^{0,0}(\mathbf{b})l'(2l'+1)]
\end{aligned}$$
where $a_{ij}$ stands for the coefficient of the entries of $A^{out}_{\mathcal{S},l}(x)$.
\end{theorem}

It is clear that $ (W_{0,0},\mathcal{S}_{D+n}[X_{lm}])=0, (X_{00},\mathcal{S}_{D+n}[X_{lm}])=0$.

\begin{theorem}
Entries for {$\boxed{(V_{l'm'},\mathcal{S}_{D+n}[X_{lm}])}$}:\\
When $m>0$, $m'>0$,
$$\begin{aligned}&(V_{l'm'},\mathcal{S}_{D+n}[X_{lm}])\\
=&\rho^{l+2}a_{33}r^{l'+1}\dfrac{\sqrt{(l'+1)(2l'+1)}}{2}\sum\limits_{q=-1}^{1}\sum\limits_{m_1=-1}^{1}[(-1)^{m'}(L_{l,l',l'+2}^{-m,m'-q,q,m_1}(\mathbf{b})+(-1)^mL_{l,l',l'+2}^{m,m'-q,q,m_1}(\mathbf{b}))\\
&+(L_{l,l',l'+2}^{-m,-m'-q,q,m_1}(\mathbf{b})+(-1)^mL_{l,l',l'+2}^{m,-m'-q,q,m_1}(\mathbf{b}))]
\end{aligned}$$

When $m>0$, $m'<0$,
$$\begin{aligned}&(V_{l'm'},\mathcal{S}_{D+n}[X_{lm}])\\
=&\rho^{l+2}a_{33}r^{l'+1}\dfrac{\sqrt{(l'+1)(2l'+1)}}{2}\sum\limits_{q=-1}^{1}\sum\limits_{m_1=-1}^{1}
[i\cdot (-1)^{m'}(L_{l,l',l'+2}^{-m,-m'-q,q,m_1}(\mathbf{b})+(-1)^mL_{l,l',l'+2}^{m,-m'-q,q,m_1}(\mathbf{b}))\\
&-i\cdot (L_{l,l',l'+2}^{-m,m'-q,q,m_1}(\mathbf{b})+(-1)^mL_{l,l',l'+2}^{m,m'-q,q,m_1}(\mathbf{b}))]
\end{aligned}$$
When $m>0$, $m'=0$,
$$\begin{aligned}&(V_{l'm'},\mathcal{S}_{D+n}[X_{lm}])\\
=&\rho^{l+2}a_{33}r^{l'+1}\dfrac{\sqrt{(l'+1)(2l'+1)}}{\sqrt{2}}\sum\limits_{q=-1}^{1}\sum\limits_{m_1=-1}^{1}
[ (L_{l,l',l'+2}^{-m,-q,q,m_1}(\mathbf{b})+(-1)^mL_{l,l',l'+2}^{m,-q,q,m_1}(\mathbf{b}))]
\end{aligned}$$

When $m<0$, $m'>0$,
$$\begin{aligned}&(V_{l'm'},\mathcal{S}_{D+n}[X_{lm}])\\
=&\rho^{l+2}a_{33}r^{l'+1}i\cdot \dfrac{\sqrt{(l'+1)(2l'+1)}}{2}\sum\limits_{q=-1}^{1}\sum\limits_{m_1=-1}^{1}[(-1)^{m'}(L_{l, l',l'+2}^{\textcolor{magenta}{m},m'-q,q,m_1}(\mathbf{b})-(-1)^mL_{l, l',l'+2}^{-m,m'-q,q,m_1}(\mathbf{b}))\\
&+(L_{l, l',l'+2}^{m,-m'-q,q,m_1}(\mathbf{b})-(-1)^mL_{l,l',l'+2}^{-m,-m'-q,q,m_1}(\mathbf{b}))]
\end{aligned}$$

When $m<0$, $m'<0$,
$$\begin{aligned}&(V_{l'm'},\mathcal{S}_{D+n}[X_{lm}])\\
=&\rho^{l+2}a_{33}r^{l'+1} \dfrac{\sqrt{(l'+1)(2l'+1)}}{2}\sum\limits_{q=-1}^{1}\sum\limits_{m_1=-1}^{1}[(L_{l,l',l'+2}^{m,m'-q,q,m_1}(\mathbf{b})-(-1)^mL_{l,l',l'+2}^{-m,m'-q,q,m_1}(\mathbf{b}))\\
&- (-1)^{m'}(L_{l,l',l'+2}^{m,-m'-q,q,m_1}(\mathbf{b})-(-1)^mL_{l,l',l'+2}^{-m,-m'-q,q,m_1}(\mathbf{b}))]
\end{aligned}$$
When $m<0$, $m'=0$,
$$\begin{aligned}&(V_{l'm'},\mathcal{S}_{D+n}[X_{lm}])\\
=&\rho^{l+2}a_{33}r^{l'+1}i\cdot \dfrac{\sqrt{(l'+1)(2l'+1)}}{\sqrt{2}}\sum\limits_{q=-1}^{1}\sum\limits_{m_1=-1}^{1}
[ L_{l,l',l'+2}^{m,-q,q,m_1}(\mathbf{b})-(-1)^mL_{l,l',l'+2}^{-m,-q,q,m_1}(\mathbf{b})]
\end{aligned}$$

When $m=0$, $m'>0$,
$$\begin{aligned}&(V_{l'm'},\mathcal{S}_{D+n}[X_{lm}])\\
=&\rho^{l+2}a_{33}r^{l'+1} \dfrac{\sqrt{(l'+1)(2l'+1)}}{\sqrt{2}}\sum\limits_{q=-1}^{1}\sum\limits_{m_1=-1}^{1}
[(-1)^{m'} L_{l,l',l'+2}^{0,m'-q,q,m_1}(\mathbf{b})+L_{l,l',l'+2}^{0,-m'-q,q,m_1}(\mathbf{b})]
\end{aligned}$$

When $m=0$, $m'<0$,
$$\begin{aligned}&(V_{l'm'},\mathcal{S}_{D+n}[X_{lm}])\\
=&\rho^{l+2}a_{33}r^{l'+1} \dfrac{\sqrt{(l'+1)(2l'+1)}}{\sqrt{2}}\sum\limits_{q=-1}^{1}\sum\limits_{m_1=-1}^{1}
[i\cdot (-1)^{m'}  L_{l,l',l'+2}^{0,-m'-q,q,m_1}(\mathbf{b})-i\cdot L_{l,l',l'+2}^{0,m'-q,q,m_1}(\mathbf{b})]
\end{aligned}$$

When $m=0$, $m'=0$,
$$\begin{aligned}&(V_{l'm'},\mathcal{S}_{D+n}[X_{lm}])\\
=&\rho^{l+2}a_{33}r^{l'+1} \sqrt{(l'+1)(2l'+1)}\sum\limits_{q=-1}^{1}\sum\limits_{m_1=-1}^{1}
[   L_{l,l',l'+2}^{0,-q,q,m_1}(\mathbf{b})]
\end{aligned}$$

\end{theorem}

\begin{remark}Numerical results show that all terms above  vanish.
\end{remark}

\begin{theorem}
Entries for {$\boxed{(W_{l'm'},\mathcal{S}_{D+n}[X_{lm}])}$}:\\
When $m>0, m'>0$,
$$
\begin{aligned}
&(W_{l'm'},\mathcal{S}_{D+n}[X_{lm}])\\
=&\ \rho^{l+2}a_{33}\,r^{l'-1}\,\frac{\sqrt{l'(2l'+1)}}{2}\sum\limits_{q=-1}^{1}\sum\limits_{m_1=-1}^{1}\bigl[(-1)^{m'}\bigl(L_{l,l',l'}^{-m,\,m'-q,\,q,\,m_1}(\mathbf{b})+(-1)^mL_{l,l',l'}^{m,\,m'-q,\,q,\,m_1}(\mathbf{b})\bigr)\\
&\qquad\qquad\qquad\qquad+\bigl(L_{l,l',l'}^{-m,\,-m'-q,\,q,\,m_1}(\mathbf{b})+(-1)^mL_{l,l',l'}^{m,\,-m'-q,\,q,\,m_1}(\mathbf{b})\bigr)\bigr]
\end{aligned}
$$
When  $m>0, m'<0$,
$$
\begin{aligned}
&(W_{l'm'},\mathcal{S}_{D+n}[X_{lm}])\\
=&\ \rho^{l+2}a_{33}\,r^{l'-1}\,\frac{\sqrt{l'(2l'+1)}}{2}\sum\limits_{q=-1}^{1}\sum\limits_{m_1=-1}^{1}\bigl[i(-1)^{m'}\bigl(L_{l,l',l'}^{-m,\,-m'-q,\,q,\,m_1}(\mathbf{b})+(-1)^mL_{l,l',l'}^{m,\,-m'-q,\,q,\,m_1}(\mathbf{b})\bigr)\\
&\qquad\qquad\qquad\qquad-i\bigl(L_{l,l',l'}^{-m,\,m'-q,\,q,\,m_1}(\mathbf{b})+(-1)^mL_{l,l',l'}^{m,\,m'-q,\,q,\,m_1}(\mathbf{b})\bigr)\bigr]
\end{aligned}
$$

When $m>0, m'=0$,
$$
\begin{aligned}
&(W_{l'0},\mathcal{S}_{D+n}[X_{lm}])\\
=&\ \rho^{l+2}a_{33}\,r^{l'-1}\,\frac{\sqrt{l'(2l'+1)}}{\sqrt{2}}\sum\limits_{q=-1}^{1}\sum\limits_{m_1=-1}^{1}\bigl(L_{l,l',l'}^{-m,\,-q,\,q,\,m_1}(\mathbf{b})+(-1)^mL_{l,l',l'}^{m,\,-q,\,q,\,m_1}(\mathbf{b})\bigr)
\end{aligned}
$$

When  $m<0, m'>0$,
$$
\begin{aligned}
&(W_{l'm'},\mathcal{S}_{D+n}[X_{lm}])\\
=&\ \rho^{l+2}a_{33}\,r^{l'-1}\,i\cdot\frac{\sqrt{l'(2l'+1)}}{2}\sum\limits_{q=-1}^{1}\sum\limits_{m_1=-1}^{1}\bigl[(-1)^{m'}\bigl(L_{l,l',l'}^{m,\,m'-q,\,q,\,m_1}(\mathbf{b})-(-1)^mL_{l,l',l'}^{-m,\,m'-q,\,q,\,m_1}(\mathbf{b})\bigr)\\
&\qquad\qquad\qquad\qquad+\bigl(L_{l,l',l'}^{m,\,-m'-q,\,q,\,m_1}(\mathbf{b})-(-1)^mL_{l,l',l'}^{-m,\,-m'-q,\,q,\,m_1}(\mathbf{b})\bigr)\bigr]
\end{aligned}
$$
When $m<0, m'<0$,
$$
\begin{aligned}
&(W_{l'm'},\mathcal{S}_{D+n}[X_{lm}])\\
=&\ \rho^{l+2}a_{33}\,r^{l'-1}\,\frac{\sqrt{l'(2l'+1)}}{2}\sum\limits_{q=-1}^{1}\sum\limits_{m_1=-1}^{1}\bigl[\bigl(L_{l,l',l'}^{m,\,m'-q,\,q,\,m_1}(\mathbf{b})-(-1)^mL_{l,l',l'}^{-m,\,m'-q,\,q,\,m_1}(\mathbf{b})\bigr)\\
&\qquad\qquad\qquad\qquad-(-1)^{m'}\bigl(L_{l,l',l'}^{m,\,-m'-q,\,q,\,m_1}(\mathbf{b})-(-1)^mL_{l,l',l'}^{-m,\,-m'-q,\,q,\,m_1}(\mathbf{b})\bigr)\bigr]
\end{aligned}
$$
When $m<0, m'=0$,
$$
\begin{aligned}
&(W_{l'0},\mathcal{S}_{D+n}[X_{lm}])\\
=&\ \rho^{l+2}a_{33}\,r^{l'-1}\,\frac{i\sqrt{l'(2l'+1)}}{\sqrt{2}}\sum\limits_{q=-1}^{1}\sum\limits_{m_1=-1}^{1}\bigl(L_{l,l',l'}^{m,\,-q,\,q,\,m_1}(\mathbf{b})-(-1)^mL_{l,l',l'}^{-m,\,-q,\,q,\,m_1}(\mathbf{b})\bigr)
\end{aligned}
$$
When $m=0, m'>0$,
$$
\begin{aligned}
&(W_{l'm'},\mathcal{S}_{D+n}[X_{l0}])\\
=&\ \rho^{l+2}a_{33}\,r^{l'-1}\,\frac{\sqrt{l'(2l'+1)}}{\sqrt{2}}\sum\limits_{q=-1}^{1}\sum\limits_{m_1=-1}^{1}\bigl[(-1)^{m'}L_{l,l',l'}^{0,\,m'-q,\,q,\,m_1}(\mathbf{b})+L_{l,l',l'}^{0,\,-m'-q,\,q,\,m_1}(\mathbf{b})\bigr]
\end{aligned}
$$

When  $m=0, m'<0$,
$$
\begin{aligned}
&(W_{l'm'},\mathcal{S}_{D+n}[X_{l0}])\\
=&\ \rho^{l+2}a_{33}\,r^{l'-1}\,\frac{\sqrt{l'(2l'+1)}}{\sqrt{2}}\sum\limits_{q=-1}^{1}\sum\limits_{m_1=-1}^{1}\bigl[i(-1)^{m'}L_{l,l',l'}^{0,\,-m'-q,\,q,\,m_1}(\mathbf{b})-i\,L_{l,l',l'}^{0,\,m'-q,\,q,\,m_1}(\mathbf{b})\bigr]
\end{aligned}
$$

When $m=0, m'=0$,
$$
\begin{aligned}
&(W_{l'0},\mathcal{S}_{D+n}[X_{l0}])\\
=&\ \rho^{l+2}a_{33}\,r^{l'-1}\,\sqrt{l'(2l'+1)}\sum\limits_{q=-1}^{1}\sum\limits_{m_1=-1}^{1}L_{l,l',l'}^{0,\,-q,\,q,\,m_1}(\mathbf{b})
\end{aligned}
$$
\end{theorem}

\begin{theorem}
Entries for {$\boxed{(X_{l'm'},\mathcal{S}_{D+n}[X_{lm}])}$}:\\

When $m>0,\ m'>0$:
$$
\begin{aligned}
&(X_{l'm'},\mathcal{S}_{D+n}[X_{lm}])\\
=&\ \rho^{l+2}a_{33}r^{l'}\Biggl\{
\frac{l'(l'+1)}{2}\Bigl[(-1)^{m'}\bigl(H_{l,l'}^{-m,m'}(\mathbf{b})+(-1)^m H_{l,l'}^{m,m'}(\mathbf{b})\bigr)
+\bigl(H_{l,l'}^{-m,-m'}(\mathbf{b})+(-1)^m H_{l,l'}^{m,-m'}(\mathbf{b})\bigr)\Bigr]\\
&-\frac{i\sqrt{l'(l'+1)}}{2}\sum_{q=-1}^{1}\sum_{m_1=-1}^{1}\Bigl[(-1)^{m'}\bigl(L_{l,l',l'+1}^{-m,m'-q,q,m_1}(\mathbf{b})+(-1)^m L_{l,l',l'+1}^{m,m'-q,q,m_1}(\mathbf{b})\bigr)\\
&\qquad\qquad\qquad\qquad\quad+\bigl(L_{l,l',l'+1}^{-m,-m'-q,q,m_1}(\mathbf{b})+(-1)^m L_{l,l',l'+1}^{m,-m'-q,q,m_1}(\mathbf{b})\bigr)\Bigr]\Biggr\}
\end{aligned}
$$

When $m>0,\ m'<0$:
$$
\begin{aligned}
&(X_{l'm'},\mathcal{S}_{D+n}[X_{lm}])\\
=&\ \rho^{l+2}a_{33}r^{l'}\Biggl\{
\frac{l'(l'+1)}{2}\Bigl[i(-1)^{m'}\bigl(H_{l,l'}^{-m,-m'}(\mathbf{b})+(-1)^m H_{l,l'}^{m,-m'}(\mathbf{b})\bigr)
-i\bigl(H_{l,l'}^{-m,m'}(\mathbf{b})+(-1)^m H_{l,l'}^{m,m'}(\mathbf{b})\bigr)\Bigr]\\
&+\frac{\sqrt{l'(l'+1)}}{2}\sum_{q=-1}^{1}\sum_{m_1=-1}^{1}\Bigl[(-1)^{m'}\bigl(L_{l,l',l'+1}^{-m,-m'-q,q,m_1}(\mathbf{b})+(-1)^m L_{l,l',l'+1}^{m,-m'-q,q,m_1}(\mathbf{b})\bigr)\\
&\qquad\qquad\qquad\qquad\quad-\bigl(L_{l,l',l'+1}^{-m,m'-q,q,m_1}(\mathbf{b})+(-1)^m L_{l,l',l'+1}^{m,m'-q,q,m_1}(\mathbf{b})\bigr)\Bigr]\Biggr\}
\end{aligned}
$$

When $m>0,\ m'=0$:
$$
\begin{aligned}
&(X_{l'0},\mathcal{S}_{D+n}[X_{lm}])\\
=&\ \rho^{l+2}a_{33}r^{l'}\Biggl\{
\frac{l'(l'+1)}{\sqrt{2}}\bigl(H_{l,l'}^{-m,0}(\mathbf{b})+(-1)^m H_{l,l'}^{m,0}(\mathbf{b})\bigr)\\
&-\frac{i\sqrt{l'(l'+1)}}{\sqrt{2}}\sum_{q=-1}^{1}\sum_{m_1=-1}^{1}\bigl(L_{l,l',l'+1}^{-m,-q,q,m_1}(\mathbf{b})+(-1)^m L_{l,l',l'+1}^{m,-q,q,m_1}(\mathbf{b})\bigr)\Biggr\}
\end{aligned}
$$

When $m<0,\ m'>0$:
$$
\begin{aligned}
&(X_{l'm'},\mathcal{S}_{D+n}[X_{lm}])\\
=&\ \rho^{l+2}a_{33}r^{l'}\Biggl\{
\frac{l'(l'+1)}{2}\Bigl[i(-1)^{m'}\bigl(H_{l,l'}^{m,m'}(\mathbf{b})-(-1)^m H_{l,l'}^{-m,m'}(\mathbf{b})\bigr)
+i\bigl(H_{l,l'}^{m,-m'}(\mathbf{b})-(-1)^m H_{l,l'}^{-m,-m'}(\mathbf{b})\bigr)\Bigr]\\
&+\frac{\sqrt{l'(l'+1)}}{2}\sum_{q=-1}^{1}\sum_{m_1=-1}^{1}\Bigl[(-1)^{m'}\bigl(L_{l,l',l'+1}^{m,m'-q,q,m_1}(\mathbf{b})-(-1)^m L_{l,l',l'+1}^{-m,m'-q,q,m_1}(\mathbf{b})\bigr)\\
&\qquad\qquad\qquad\qquad\quad+\bigl(L_{l,l',l'+1}^{m,-m'-q,q,m_1}(\mathbf{b})-(-1)^m L_{l,l',l'+1}^{-m,-m'-q,q,m_1}(\mathbf{b})\bigr)\Bigr]\Biggr\}
\end{aligned}
$$

When $m<0,\ m'<0$:
$$
\begin{aligned}
&(X_{l'm'},\mathcal{S}_{D+n}[X_{lm}])\\
=&\ \rho^{l+2}a_{33}r^{l'}\Biggl\{
\frac{l'(l'+1)}{2}\Bigl[\bigl(H_{l,l'}^{m,m'}(\mathbf{b})-(-1)^m H_{l,l'}^{-m,m'}(\mathbf{b})\bigr)
-(-1)^{m'}\bigl(H_{l,l'}^{m,-m'}(\mathbf{b})-(-1)^m H_{l,l'}^{-m,-m'}(\mathbf{b})\bigr)\Bigr]\\
&+\frac{i\sqrt{l'(l'+1)}}{2}\sum_{q=-1}^{1}\sum_{m_1=-1}^{1}\Bigl[(-1)^{m'}\bigl(L_{l,l',l'+1}^{m,-m'-q,q,m_1}(\mathbf{b})-(-1)^m L_{l,l',l'+1}^{-m,-m'-q,q,m_1}(\mathbf{b})\bigr)\\
&\qquad\qquad\qquad\qquad\quad-\bigl(L_{l,l',l'+1}^{m,m'-q,q,m_1}(\mathbf{b})-(-1)^m L_{l,l',l'+1}^{-m,m'-q,q,m_1}(\mathbf{b})\bigr)\Bigr]\Biggr\}
\end{aligned}
$$

When $m<0,\ m'=0$:
$$
\begin{aligned}
&(X_{l'0},\mathcal{S}_{D+n}[X_{lm}])\\
=&\ \rho^{l+2}a_{33}r^{l'}\Biggl\{
\frac{il'(l'+1)}{\sqrt{2}}\bigl(H_{l,l'}^{m,0}(\mathbf{b})-(-1)^m H_{l,l'}^{-m,0}(\mathbf{b})\bigr)\\
&+\frac{\sqrt{l'(l'+1)}}{\sqrt{2}}\sum_{q=-1}^{1}\sum_{m_1=-1}^{1}\bigl(L_{l,l',l'+1}^{m,-q,q,m_1}(\mathbf{b})-(-1)^m L_{l,l',l'+1}^{-m,-q,q,m_1}(\mathbf{b})\bigr)\Biggr\}
\end{aligned}
$$

When $m=0,\ m'>0$:
$$
\begin{aligned}
&(X_{l'm'},\mathcal{S}_{D+n}[X_{l0}])\\
=&\ \rho^{l+2}a_{33}r^{l'}\Biggl\{
\frac{l'(l'+1)}{\sqrt{2}}\Bigl[(-1)^{m'}H_{l,l'}^{0,m'}(\mathbf{b})+H_{l,l'}^{0,-m'}(\mathbf{b})\Bigr]\\
&-\frac{i\sqrt{l'(l'+1)}}{\sqrt{2}}\sum_{q=-1}^{1}\sum_{m_1=-1}^{1}\Bigl[(-1)^{m'}L_{l,l',l'+1}^{0,m'-q,q,m_1}(\mathbf{b})+L_{l,l',l'+1}^{0,-m'-q,q,m_1}(\mathbf{b})\Bigr]\Biggr\}
\end{aligned}
$$

When $m=0,\ m'<0$:
$$
\begin{aligned}
&(X_{l'm'},\mathcal{S}_{D+n}[X_{l0}])\\
=&\ \rho^{l+2}a_{33}r^{l'}\Biggl\{
\frac{l'(l'+1)}{\sqrt{2}}\Bigl[i(-1)^{m'}H_{l,l'}^{0,-m'}(\mathbf{b})-iH_{l,l'}^{0,m'}(\mathbf{b})\Bigr]\\
&+\frac{\sqrt{l'(l'+1)}}{\sqrt{2}}\sum_{q=-1}^{1}\sum_{m_1=-1}^{1}\Bigl[(-1)^{m'}L_{l,l',l'+1}^{0,-m'-q,q,m_1}(\mathbf{b})-L_{l,l',l'+1}^{0,m'-q,q,m_1}(\mathbf{b})\Bigr]\Biggr\}
\end{aligned}
$$

When $m=0,\ m'=0$:
$$
\begin{aligned}
&(X_{l'0},\mathcal{S}_{D+n}[X_{l0}])\\
=&\ \rho^{l+2}a_{33}r^{l'}\Biggl\{
l'(l'+1)H_{l,l'}^{0,0}(\mathbf{b})
-i\sqrt{l'(l'+1)}\sum_{q=-1}^{1}\sum_{m_1=-1}^{1}L_{l,l',l'+1}^{0,-q,q,m_1}(\mathbf{b})\Biggr\}
\end{aligned}
$$
\end{theorem}

\begin{theorem}
Entries for {$\boxed{(V_{l'm'},\mathcal{S}_{D+n}[W_{lm}])}$}:\\
When $m>0$, $m'>0$,
$$
\begin{aligned}
&(V_{l'm'},\mathcal{S}_{D+n}[W_{lm}])\\
=&\sum\limits_{\lambda\in\{l'+1,l'+3\}}\sum\limits_{q=-1}^{1}\sum\limits_{m_1=-1}^{1}\rho^{l+1} (a_{22}-a_{12})  r^\lambda b_{-q}\sqrt{(l'+1)(2l'+1)}
[(-1)^{m'}(H_{l,\lambda}^{-m,m'-q}(\mathbf{b})\\
&+(-1)^mH_{l,\lambda}^{m,m'-q}(\mathbf{b}))K_{l'+1,l',\lambda}^{m_1,m'-q,q}
+(H_{l,\lambda}^{-m,-m'-q}(\mathbf{b})+(-1)^mH_{l,\lambda}^{m,-m'-q}(\mathbf{b}))K_{l'+1,l',\lambda}^{m_1,-m'-q,q}]\\
&-\rho^{l+1} a_{22}r^{l'+1}\dfrac{(2l+1)(l'+1)}{2}[(-1)^{m'}(H_{l,l'}^{-m,m'}(\mathbf{b})+(-1)^mH_{l,l'}^{m,m'}(\mathbf{b}))+(H_{l,l'}^{-m,-m'}(\mathbf{b})+(-1)^mH_{l,l'}^{m,-m'}(\mathbf{b}))]\\
&+\rho^{l+1} a_{22}r^{l'+1}(2l+1)\dfrac{\sqrt{(l'+1)(2l'+1)}}{2} \sum\limits_{q=-1}^1(-1)^q b_{-q}[(-1)^{m'}(H_{l,l'+1}^{-m,m'-q}(\mathbf{b})+(-1)^mH_{l,l'+1}^{m,m'-q}(\mathbf{b}))\\
&\langle l'+1,m'-q;1,q\mid l',m'\rangle+(H_{l,l'+1}^{-m,-m'-q}(\mathbf{b})+(-1)^mH_{l,l'+1}^{m,-m'-q}(\mathbf{b}))\langle l'+1,-m'-q;1,q\mid l',-m'\rangle]
\end{aligned}$$

When $m>0$, $m'<0$,

$$
\begin{aligned}
&(V_{l'm'},\mathcal{S}_{D+n}[W_{lm}])\\
=&\sum\limits_{\lambda\in\{l'+1,l'+3\}}\sum\limits_{q=-1}^{1}\sum\limits_{m_1=-1}^{1}\rho^{l+1} (a_{22}-a_{12})  r^\lambda b_{-q}\sqrt{(l'+1)(2l'+1)}
[i\cdot (-1)^{m'}(H_{l,\lambda}^{-m,-m'-q}(\mathbf{b})\\
&+(-1)^mH_{l,\lambda}^{m,-m'-q}(\mathbf{b}))K_{l'+1,l',\lambda}^{m_1,-m'-q,q}
-i\cdot (H_{l,\lambda}^{-m,m'-q}(\mathbf{b})+(-1)^mH_{l,\lambda}^{m,m'-q}(\mathbf{b}))K_{l'+1,l',\lambda}^{m_1,m'-q,q}]\\
&-\rho^{l+1} a_{22}r^{l'+1}\dfrac{(2l+1)(l'+1)}{2}[i\cdot (-1)^{m'}(H_{l,l'}^{-m,-m'}(\mathbf{b})+(-1)^mH_{l,l'}^{m,-m'}(\mathbf{b}))\\
&-i\cdot (H_{l,l'}^{-m,m'}(\mathbf{b})+(-1)^mH_{l,l'}^{m,m'}(\mathbf{b}))]\\
&+\rho^{l+1} a_{22}r^{l'+1}(2l+1)\dfrac{\sqrt{(l'+1)(2l'+1)}}{2} \sum\limits_{q=-1}^1(-1)^q b_{-q}[i\cdot (-1)^{m'}(H_{l,l'+1}^{-m,-m'-q}(\mathbf{b})+(-1)^mH_{l,l'+1}^{m,-m'-q}(\mathbf{b}))\\
&\langle l'+1,-m'-q;1,q\mid l',-m'\rangle
-i\cdot (H_{l,l'+1}^{-m,m'-q}(\mathbf{b})+(-1)^mH_{l,l'+1}^{m,m'-q}(\mathbf{b}))\langle l'+1,m'-q;1,q\mid l',m'\rangle]
\end{aligned}$$

When $m>0$, $m'=0$,
$$
\begin{aligned}
&(V_{l'm'},\mathcal{S}_{D+n}[W_{lm}])\\
=&\sum\limits_{\lambda\in\{l'+1,l'+3\}}\sum\limits_{q=-1}^{1}\sum\limits_{m_1=-1}^{1}\rho^{l+1} (a_{22}-a_{12})  r^\lambda b_{-q}\sqrt{2(l'+1)(2l'+1)}
[(H_{l,\lambda}^{-m,-q}(\mathbf{b})+(-1)^mH_{l,\lambda}^{m,-q}(\mathbf{b}))K_{l'+1,l',\lambda}^{m_1,-q,q}]\\
&-\rho^{l+1} a_{22}r^{l'+1}\dfrac{(2l+1)(l'+1)}{\sqrt{2}}[H_{l,l'}^{-m,0}(\mathbf{b})+(-1)^mH_{l,l'}^{m,0}(\mathbf{b})]\\
&+\rho^{l+1} a_{22}r^{l'+1}(2l+1)\dfrac{\sqrt{(l'+1)(2l'+1)}}{\sqrt{2}} \sum\limits_{q=-1}^1(-1)^q b_{-q}(H_{l,l'+1}^{-m,-q}(\mathbf{b})+(-1)^mH_{l,l'+1}^{m,-q}(\mathbf{b}))\\
&\langle l'+1,-q;1,q\mid l',0\rangle
\end{aligned}$$

When $m<0$, $m'>0$,
$$
\begin{aligned}
&(V_{l'm'},\mathcal{S}_{D+n}[W_{lm}])\\
=&\sum\limits_{\lambda\in\{l'+1,l'+3\}}\sum\limits_{q=-1}^{1}\sum\limits_{m_1=-1}^{1}\rho^{l+1} (a_{22}-a_{12})  r^\lambda b_{-q}i\cdot\sqrt{(l'+1)(2l'+1)}
[(-1)^{m'}(H_{l,\lambda}^{m,m'-q}(\mathbf{b})\\
&-(-1)^mH_{l,\lambda}^{-m,m'-q}(\mathbf{b}))K_{l'+1,l',\lambda}^{m_1,m'-q,q}
+(H_{l,\lambda}^{m,-m'-q}(\mathbf{b})-(-1)^mH_{l,\lambda}^{-m,-m'-q}(\mathbf{b}))K_{l'+1,l',\lambda}^{m_1,-m'-q,q}]\\
&-\rho^{l+1} a_{22}r^{l'+1}i\cdot \dfrac{(2l+1)(l'+1)}{2}[(-1)^{m'}(H_{l,l'}^{m,m'}(\mathbf{b})-(-1)^mH_{l,l'}^{-m,m'}(\mathbf{b}))\\
&+(H_{l,l'}^{m,-m'}(\mathbf{b})-(-1)^mH_{l,l'}^{-m,-m'}(\mathbf{b}))]\\
&+\rho^{l+1} a_{22}r^{l'+1}i\cdot (2l+1)\dfrac{\sqrt{(l'+1)(2l'+1)}}{2} \sum\limits_{q=-1}^1(-1)^q b_{-q}[(-1)^{m'}(H_{l,l'+1}^{m,m'-q}(\mathbf{b})-(-1)^mH_{l,l'+1}^{-m,m'-q}(\mathbf{b}))\\
&\langle l'+1,m'-q;1,q\mid l',m'\rangle+(H_{l,l'+1}^{m,-m'-q}(\mathbf{b})-(-1)^mH_{l,l'+1}^{-m,-m'-q}(\mathbf{b}))\langle l'+1,-m'-q;1,q\mid l',-m'\rangle]
\end{aligned}$$

When $m<0$, $m'<0$,

$$
\begin{aligned}
&(V_{l'm'},\mathcal{S}_{D+n}[W_{lm}])\\
=&\sum\limits_{\lambda\in\{l'+1,l'+3\}}\sum\limits_{q=-1}^{1}\sum\limits_{m_1=-1}^{1}\rho^{l+1} (a_{22}-a_{12})  r^\lambda b_{-q}\sqrt{(l'+1)(2l'+1)}
[-(-1)^{m'}(H_{l,\lambda}^{m,-m'-q}(\mathbf{b})\\
&-(-1)^mH_{l,\lambda}^{-m,-m'-q}(\mathbf{b}))K_{l'+1,l',\lambda}^{m_1,-m'-q,q}
+ (H_{l,\lambda}^{m,m'-q}(\mathbf{b})-(-1)^mH_{l,\lambda}^{-m,m'-q}(\mathbf{b}))K_{l'+1,l',\lambda}^{m_1,m'-q,q}]\\
&-\rho^{l+1} a_{22}r^{l'+1}\dfrac{(2l+1)(l'+1)}{2}
[- (-1)^{m'}(H_{l,l'}^{m,-m'}(\mathbf{b})-(-1)^mH_{l,l'}^{-m,-m'}(\mathbf{b}))\\
&+
 (H_{l,l'}^{m,m'}(\mathbf{b})-(-1)^mH_{l,l'}^{-m,m'}(\mathbf{b}))]\\
&+\rho^{l+1} a_{22}r^{l'+1}(2l+1)\dfrac{\sqrt{(l'+1)(2l'+1)}}{2} \sum\limits_{q=-1}^1(-1)^q b_{-q}[-(-1)^{m'}(H_{l,l'+1}^{m,-m'-q}(\mathbf{b})-(-1)^mH_{l,l'+1}^{-m,-m'-q}(\mathbf{b}))\\
&\langle l'+1,-m'-q;1,q\mid l',-m'\rangle
+ (H_{l,l'+1}^{m,m'-q}(\mathbf{b})-(-1)^mH_{l,l'+1}^{-m,m'-q}(\mathbf{b}))\langle l'+1,m'-q;1,q\mid l',m'\rangle]
\end{aligned}$$

When $m<0$, $m'=0$,
$$
\begin{aligned}
&(V_{l'm'},\mathcal{S}_{D+n}[W_{lm}])\\
=&\sum\limits_{\lambda\in\{l'+1,l'+3\}}\sum\limits_{q=-1}^{1}\sum\limits_{m_1=-1}^{1}\rho^{l+1} (a_{22}-a_{12})  r^\lambda b_{-q}i\cdot \sqrt{2(l'+1)(2l'+1)}
[(H_{l,\lambda}^{m,-q}(\mathbf{b})\\
&-(-1)^mH_{l,\lambda}^{-m,-q}(\mathbf{b}))K_{l'+1,l',\lambda}^{m_1,-q,q}]\\
&-\rho^{l+1} a_{22}r^{l'+1}i\cdot \dfrac{(2l+1)(l'+1)}{\sqrt{2}}[(H_{l,l'}^{m,0}(\mathbf{b})-(-1)^mH_{l,l'}^{-m,0}(\mathbf{b})]\\
&+\rho^{l+1} a_{22}r^{l'+1}i\cdot (2l+1)\dfrac{\sqrt{(l'+1)(2l'+1)}}{\sqrt{2}} \sum\limits_{q=-1}^1(-1)^q b_{-q}[ (H_{l,l'+1}^{m,-q}(\mathbf{b})-(-1)^mH_{l,l'+1}^{-m,-q}(\mathbf{b}))\\
&\langle l'+1,-q;1,q\mid l',0\rangle]
\end{aligned}$$

When $m=0$, $m'>0$,
$$
\begin{aligned}
&(V_{l'm'},\mathcal{S}_{D+n}[W_{lm}])\\
=&\sum\limits_{\lambda\in\{l'+1,l'+3\}}\sum\limits_{q=-1}^{1}\sum\limits_{m_1=-1}^{1}\rho^{l+1} (a_{22}-a_{12})  r^\lambda b_{-q}\sqrt{2(l'+1)(2l'+1)}
[(-1)^{m'}H_{l,\lambda}^{0,m'-q}(\mathbf{b}) K_{l'+1,l',\lambda}^{m_1,m'-q,q}\\
&+H_{l,\lambda}^{0,-m'-q}(\mathbf{b})K_{l'+1,l',\lambda}^{m_1,-m'-q,q}]\\
&-\rho^{l+1} a_{22}r^{l'+1}\dfrac{(2l+1)(l'+1)}{\sqrt{2}}[(-1)^{m'}H_{l,l'}^{0,m'}(\mathbf{b})+H_{l,l'}^{0,-m'}(\mathbf{b})]\\
&+\rho^{l+1} a_{22}r^{l'+1}(2l+1)\dfrac{\sqrt{(l'+1)(2l'+1)}}{\sqrt{2}} \sum\limits_{q=-1}^1(-1)^q b_{-q}[(-1)^{m'}H_{l,l'+1}^{0,m'-q}(\mathbf{b})\\
&\langle l'+1,m'-q;1,q\mid l',m'\rangle+(H_{l,l'+1}^{0,-m'-q}(\mathbf{b})\langle l'+1,-m'-q;1,q\mid l',-m'\rangle]
\end{aligned}$$

When $m=0$, $m'<0$,

$$
\begin{aligned}
&(V_{l'm'},\mathcal{S}_{D+n}[W_{lm}])\\
=&\sum\limits_{\lambda\in\{l'+1,l'+3\}}\sum\limits_{q=-1}^{1}\sum\limits_{m_1=-1}^{1}\rho^{l+1} (a_{22}-a_{12})  r^\lambda b_{-q}\sqrt{2(l'+1)(2l'+1)}
[i\cdot (-1)^{m'}H_{l,\lambda}^{0,-m'-q}(\mathbf{b})K_{l'+1,l',\lambda}^{m_1,-m'-q,q}\\
&
-i\cdot H_{l,\lambda}^{0,m'-q}(\mathbf{b})K_{l'+1,l',\lambda}^{m_1,m'-q,q}]\\
&-\rho^{l+1} a_{22}r^{l'+1}\dfrac{(2l+1)(l'+1)}{\sqrt{2}}[i\cdot (-1)^{m'}H_{l,l'}^{0,-m'}(\mathbf{b})-i\cdot H_{l,l'}^{0,m'}(\mathbf{b})]\\
&+\rho^{l+1} a_{22}r^{l'+1}(2l+1)\dfrac{\sqrt{(l'+1)(2l'+1)}}{\sqrt{2}} \sum\limits_{q=-1}^1(-1)^q b_{-q}[i\cdot (-1)^{m'} H_{l,l'+1}^{0,-m'-q}(\mathbf{b})\\
&\langle l'+1,-m'-q;1,q\mid l',-m'\rangle-i\cdot H_{l,l'+1}^{0,m'-q}(\mathbf{b})\langle l'+1,m'-q;1,q\mid l',m'\rangle]
\end{aligned}$$

When $m=0$, $m'=0$,
$$
\begin{aligned}
&(V_{l'm'},\mathcal{S}_{D+n}[W_{lm}])\\
=&\sum\limits_{\lambda\in\{l'+1,l'+3\}}\sum\limits_{q=-1}^{1}\sum\limits_{m_1=-1}^{1}\rho^{l+1} (a_{22}-a_{12})  r^\lambda b_{-q}\sqrt{4(l'+1)(2l'+1)}
[H_{l,\lambda}^{0,-q}(\mathbf{b})K_{l'+1,l',\lambda}^{m_1,-q,q}]\\
&-\rho^{l+1} a_{22}r^{l'+1}(2l+1)(l'+1)[H_{l,l'}^{0,0}(\mathbf{b})]\\
&+\rho^{l+1} a_{22}r^{l'+1}(2l+1)\sqrt{(l'+1)(2l'+1)} \sum\limits_{q=-1}^1(-1)^q b_{-q}[ H_{l,l'+1}^{0,-q}(\mathbf{b})\langle l'+1,-q;1,q\mid l',0\rangle]
\end{aligned}$$

\end{theorem}

\begin{theorem}
Entries for {$\boxed{(X_{l'm'},\mathcal{S}_{D+n}[W_{lm}])}$}:\\
When $m>0$, $m'>0$,
$$
\begin{aligned}
&(X_{l'm'},\mathcal{S}_{D+n}[W_{lm}])\\
=&\sum\limits_{\lambda\in\{l',l'+2\mid l'\geqslant 1\}}\sum\limits_{q=-1}^{1}\sum\limits_{m_1=-1}^{1}-i\cdot \rho^{l+1} (a_{22}-a_{12})  r^\lambda b_{-q}\sqrt{l'(l'+1)}
[(-1)^{m'}(H_{l,\lambda}^{-m,m'-q}(\mathbf{b})\\
&+(-1)^mH_{l,\lambda}^{m,m'-q}(\mathbf{b}))K_{l',l',\lambda}^{m_1,m'-q,q}
+(H_{l,\lambda}^{-m,-m'-q}(\mathbf{b})+(-1)^mH_{l,\lambda}^{m,-m'-q}(\mathbf{b}))K_{l',l',\lambda}^{m_1,-m'-q,q}]\\
&- i\cdot \rho^{l+1} a_{22}r^{l'}(2l+1)\dfrac{\sqrt{l'(l'+1)}}{2} \sum\limits_{q=-1}^1(-1)^q b_{-q}[(-1)^{m'}(H_{l,l'}^{-m,m'-q}(\mathbf{b})+(-1)^mH_{l,l'}^{m,m'-q}(\mathbf{b}))\\
&\langle l',m'-q;1,q\mid l',m'\rangle+(H_{l,l'}^{-m,-m'-q}(\mathbf{b})+(-1)^mH_{l,l'}^{m,-m'-q}(\mathbf{b}))\langle l',-m'-q;1,q\mid l',-m'\rangle]
\end{aligned}$$

When $m>0$, $m'<0$,
$$
\begin{aligned}
&(X_{l'm'},\mathcal{S}_{D+n}[W_{lm}])\\
=&\sum\limits_{\lambda\in\{l',l'+2\mid l'\geqslant 1\}}\sum\limits_{q=-1}^{1}\sum\limits_{m_1=-1}^{1}\rho^{l+1} (a_{22}-a_{12})  r^\lambda b_{-q}\sqrt{l'(l'+1)}
[ (-1)^{m'}(H_{l,\lambda}^{-m,-m'-q}(\mathbf{b})\\
&+(-1)^mH_{l,\lambda}^{m,-m'-q}(\mathbf{b}))K_{l',l',\lambda}^{m_1,-m'-q,q}
- (H_{l,\lambda}^{-m,m'-q}(\mathbf{b})+(-1)^mH_{l,\lambda}^{m,m'-q}(\mathbf{b}))K_{l',l',\lambda}^{m_1,m'-q,q}]\\
&+\rho^{l+1} a_{22}r^{l'}(2l+1)\dfrac{\sqrt{l'(l'+1)}}{2} \sum\limits_{q=-1}^1(-1)^q b_{-q}[ (-1)^{m'}(H_{l,l'}^{-m,-m'-q}(\mathbf{b})+(-1)^mH_{l,l'}^{m,-m'-q}(\mathbf{b}))\\
&\langle l',-m'-q;1,q\mid l',-m'\rangle
- (H_{l,l'}^{-m,m'-q}(\mathbf{b})+(-1)^mH_{l,l'}^{m,m'-q}(\mathbf{b}))\langle l',m'-q;1,q\mid l',m'\rangle]\end{aligned}$$

When $m>0$, $m'=0$,
$$
\begin{aligned}
&(X_{l'm'},\mathcal{S}_{D+n}[W_{lm}])\\
=&\sum\limits_{\lambda\in\{l',l'+2\mid l'\geqslant 1\}}\sum\limits_{q=-1}^{1}\sum\limits_{m_1=-1}^{1} -i\cdot \rho^{l+1} (a_{22}-a_{12})  r^\lambda b_{-q}\sqrt{2l'(l'+1)}
[(H_{l,\lambda}^{-m,-q}(\mathbf{b})+(-1)^mH_{l,\lambda}^{m,-q}(\mathbf{b}))K_{l',l',\lambda}^{m_1,-q,q}]\\
&-i\cdot \rho^{l+1} a_{22}r^{l'}(2l+1)\dfrac{\sqrt{l'(l'+1)}}{\sqrt{2}} \sum\limits_{q=-1}^1(-1)^q b_{-q}[ (H_{l,l'}^{-m,-q}(\mathbf{b})+(-1)^mH_{l,l'}^{m,-q}(\mathbf{b}))\\
&\langle l',-q;1,q\mid l',0\rangle]
\end{aligned}$$

When $m<0$, $m'>0$,
$$
\begin{aligned}
&(X_{l'm'},\mathcal{S}_{D+n}[W_{lm}])\\
=&\sum\limits_{\lambda\in\{l',l'+2\mid l'\geqslant 1\}}\sum\limits_{q=-1}^{1}\sum\limits_{m_1=-1}^{1}\rho^{l+1} (a_{22}-a_{12})  r^\lambda b_{-q}\sqrt{l'(l'+1)}
[(-1)^{m'}(H_{l,\lambda}^{m,m'-q}(\mathbf{b})\\
&-(-1)^mH_{l,\lambda}^{-m,m'-q}(\mathbf{b}))K_{l',l',\lambda}^{m_1,m'-q,q}
+(H_{l,\lambda}^{m,-m'-q}(\mathbf{b})-(-1)^mH_{l,\lambda}^{-m,-m'-q}(\mathbf{b}))K_{l',l',\lambda}^{m_1,-m'-q,q}]\\
& +\rho^{l+1} a_{22}r^{l'} (2l+1)\dfrac{\sqrt{l'(l'+1)}}{2} \sum\limits_{q=-1}^1(-1)^q b_{-q}[(-1)^{m'}(H_{l,l'}^{m,m'-q}(\mathbf{b})-(-1)^mH_{l,l'}^{-m,m'-q}(\mathbf{b}))\\
&\langle l',m'-q;1,q\mid l',m'\rangle+(H_{l,l'}^{m,-m'-q}(\mathbf{b})-(-1)^mH_{l,l'}^{-m,-m'-q}(\mathbf{b}))\langle l',-m'-q;1,q\mid l',-m'\rangle]
\end{aligned}$$

When $m<0$, $m'<0$,
$$
\begin{aligned}
&(X_{l'm'},\mathcal{S}_{D+n}[W_{lm}])\\
=&\sum\limits_{\lambda\in\{l',l'+2\mid l'\geqslant 1\}}\sum\limits_{q=-1}^{1}\sum\limits_{m_1=-1}^{1}\rho^{l+1} (a_{22}-a_{12})  r^\lambda b_{-q}i\cdot \sqrt{l'(l'+1)}
[(-1)^{m'}(H_{l,\lambda}^{m,-m'-q}(\mathbf{b})\\
&-(-1)^mH_{l,\lambda}^{-m,-m'-q}(\mathbf{b}))K_{l',l',\lambda}^{m_1,-m'-q,q}
- (H_{l,\lambda}^{m,m'-q}(\mathbf{b})-(-1)^mH_{l,\lambda}^{-m,m'-q}(\mathbf{b}))K_{l',l',\lambda}^{m_1,m'-q,q}]\\
&-i\cdot \rho^{l+1} a_{22}r^{l'}(2l+1)\dfrac{\sqrt{l'(l'+1)}}{2} \sum\limits_{q=-1}^1(-1)^q b_{-q}[-(-1)^{m'}(H_{l,l'}^{m,-m'-q}(\mathbf{b})-(-1)^mH_{l,l'}^{-m,-m'-q}(\mathbf{b}))\\
&\langle l',-m'-q;1,q\mid l',-m'\rangle
+ (H_{l,l'}^{m,m'-q}(\mathbf{b})-(-1)^mH_{l,l'}^{-m,m'-q}(\mathbf{b}))\langle l',m'-q;1,q\mid l',m'\rangle]
\end{aligned}$$

When $m<0$, $m'=0$,
$$
\begin{aligned}
&(X_{l'm'},\mathcal{S}_{D+n}[W_{lm}])\\
=&\sum\limits_{\lambda\in\{l',l'+2\mid l'\geqslant 1\}}\sum\limits_{q=-1}^{1}\sum\limits_{m_1=-1}^{1}  \rho^{l+1} (a_{22}-a_{12})  r^\lambda b_{-q} \sqrt{2l'(l'+1)}
[(H_{l,\lambda}^{m,-q}(\mathbf{b})-(-1)^mH_{l,\lambda}^{-m,-q}(\mathbf{b}))K_{l',l',\lambda}^{m_1,-q,q}]\\
&+\rho^{l+1} a_{22}r^{l'} (2l+1)\dfrac{\sqrt{l'(l'+1)}}{\sqrt{2}} \sum\limits_{q=-1}^1(-1)^q b_{-q}[ (H_{l,l'}^{m,-q}(\mathbf{b})-(-1)^mH_{l,l'}^{-m,-q}(\mathbf{b}))\\
&\langle l',-q;1,q\mid l',0\rangle]
\end{aligned}$$

When $m=0$, $m'>0$,
$$
\begin{aligned}
&(X_{l'm'},\mathcal{S}_{D+n}[W_{lm}])\\
=&\sum\limits_{\lambda\in\{l',l'+2\mid l'\geqslant 1\}}\sum\limits_{q=-1}^{1}\sum\limits_{m_1=-1}^{1} -i\cdot \rho^{l+1} (a_{22}-a_{12})  r^\lambda b_{-q}\sqrt{2l'(l'+1)}
[(-1)^{m'}H_{l,\lambda}^{0,m'-q}(\mathbf{b}) K_{l',l',\lambda}^{m_1,m'-q,q}\\
&+H_{l,\lambda}^{0,-m'-q}(\mathbf{b})K_{l',l',\lambda}^{m_1,-m'-q,q}]\\
&-i\cdot \rho^{l+1} a_{22}r^{l'}(2l+1)\dfrac{\sqrt{l'(l'+1)}}{\sqrt{2}} \sum\limits_{q=-1}^1(-1)^q b_{-q}[(-1)^{m'}H_{l,l'}^{0,m'-q}(\mathbf{b})\\
&\langle l',m'-q;1,q\mid l',m'\rangle+(H_{l,l'}^{0,-m'-q}(\mathbf{b})\langle l',-m'-q;1,q\mid l',-m'\rangle]
\end{aligned}$$

When $m=0$, $m'<0$,

$$
\begin{aligned}
&(X_{l'm'},\mathcal{S}_{D+n}[W_{lm}])\\
=&\sum\limits_{\lambda\in\{l',l'+2\mid l'\geqslant 1\}}\sum\limits_{q=-1}^{1}\sum\limits_{m_1=-1}^{1}\rho^{l+1} (a_{22}-a_{12})  r^\lambda b_{-q}\sqrt{2l'(l'+1)}
[(-1)^{m'}H_{l,\lambda}^{0,-m'-q}(\mathbf{b})K_{l',l',\lambda}^{m_1,-m'-q,q}\\
&
- H_{l,\lambda}^{0,m'-q}(\mathbf{b})K_{l',l',\lambda}^{m_1,m'-q,q}]\\
&+\rho^{l+1} a_{22}r^{l'}(2l+1)\dfrac{\sqrt{l'(l'+1)}}{\sqrt{2}} \sum\limits_{q=-1}^1(-1)^q b_{-q}[ (-1)^{m'} H_{l,l'}^{0,-m'-q}(\mathbf{b})\\
&\langle l',-m'-q;1,q\mid l',-m'\rangle- H_{l,l'}^{0,m'-q}(\mathbf{b})\langle l',m'-q;1,q\mid l',m'\rangle]
\end{aligned}$$

When $m=0$, $m'=0$,
$$
\begin{aligned}
&(X_{l'm'},\mathcal{S}_{D+n}[W_{lm}])\\
=&\sum\limits_{\lambda\in\{l',l'+2\mid l'\geqslant 1\}}\sum\limits_{q=-1}^{1}\sum\limits_{m_1=-1}^{1} -i\cdot \rho^{l+1} (a_{22}-a_{12})  r^\lambda b_{-q}\sqrt{4l'(l'+1)}
[H_{l,\lambda}^{0,-q}(\mathbf{b})K_{l',l',\lambda}^{m_1,-q,q}]\\
&-i\cdot \rho^{l+1} a_{22}r^{l'}(2l+1)\sqrt{l'(l'+1)} \sum\limits_{q=-1}^1(-1)^q b_{-q}[ H_{l,l'}^{0,-q}(\mathbf{b})\langle l',-q;1,q\mid l',0\rangle]
\end{aligned}$$

\end{theorem}

\begin{theorem}The entries  $\boxed{(W_{l'm'},\mathcal{S}_{D+n}[W_{lm}])}$($l'\geqslant 1$) are given by:
\\
If $m>0$ and $m'>0$,
$$\begin{aligned}&(W_{l'm'},\mathcal{S}_{D+n}[W_{lm}])\\
&=\dfrac{l'(2l'+1)}{2}[ \rho^{l+3}a_{12}r^{l'-1}+\rho^{l+1}(a_{22}-a_{12})(r^{l'+1}+a^2r^{l'-1})+\rho^{l+1}a_{22}r^{l'+1}\dfrac{2l+1}{2l'+1}]\\
&[(-1)^{m'}(H_{l,l'}^{-m,m'}(\mathbf{b})+(-1)^mH_{l,l'}^{m,m'}(\mathbf{b}))+(H_{l,l'}^{-m,-m'}(\mathbf{b})+(-1)^mH_{l,l'}^{m,-m'}(\mathbf{b}))]\\
&+\rho^{l+1}(a_{22}-a_{12})\sqrt{l'(2l'+1)}\sum\limits_{\lambda\in\{l'+1,l'-1\}}\sum\limits_{q=-1}^{q=1}\sum\limits_{m_1=-1}^{1}r^\lambda b_{-q}[(-1)^{m'}(H_{l,\lambda}^{-m,m'-q}(\mathbf{b})\\
&+(-1)^mH_{l,\lambda}^{m,m'-q}(\mathbf{b}))K_{l'-1,l',\lambda}^{m_1,m'-q,q}+(H_{l,\lambda}^{-m,-m'-q}(\mathbf(b))+(-1)^mH_{l,\lambda}^{m,-m'-q}(\mathbf{b}))K_{l'-1,l',\lambda}^{m_1,-m'-q,q}]\\
&+\rho^{l+1}a_{22}(2l+1)\dfrac{\sqrt{l'(2l'+1)}}{2}\sum\limits_{q=-1}^1(-1)^q b_{-q}r^{l'-1}[(-1)^{m'}(H_{l,l'-1}^{-m,m'-q}(\mathbf{b})+(-1)^mH_{l,l'-1}^{m,m'-q}(\mathbf{b}))\\
&\langle l'-1,m'-q;1,q\mid l',m' \rangle+(H_{l,l'-1}^{-m,-m'-q}(\mathbf{b})+(-1)^mH_{l,l'-1}^{m,-m'-q}(\mathbf{b}))\langle l'-1,-m'-q;1,q\mid l',-m' \rangle]
\end{aligned}$$

If $m>0$ and $m'<0$,
$$
\begin{aligned}
&(W_{l'm'},\mathcal{S}_{D+n}[W_{lm}])\\
&=\dfrac{l'(2l'+1)}{2}\Bigl[ \rho^{l+3}a_{12}r^{l'-1}+\rho^{l+1}(a_{22}-a_{12})(r^{l'+1}+a^2r^{l'-1})+\rho^{l+1}a_{22}r^{l'+1}\dfrac{2l+1}{2l'+1}\Bigr]\\
&\quad\Bigl[i\cdot(-1)^{m'}(H_{l,l'}^{-m,-m'}(\mathbf{b})+(-1)^mH_{l,l'}^{m,-m'}(\mathbf{b}))-i\cdot(H_{l,l'}^{-m,m'}(\mathbf{b})+(-1)^mH_{l,l'}^{m,m'}(\mathbf{b}))\Bigr]\\
&+\rho^{l+1}(a_{22}-a_{12})\sqrt{l'(2l'+1)}\sum_{\lambda\in\{l'+1,l'-1\}}\sum_{q=-1}^{1}\sum_{m_1=-1}^{1}r^\lambda b_{-q}\Bigl[i\cdot(-1)^{m'}(H_{l,\lambda}^{-m,-m'-q}(\mathbf{b})\\
&\quad+(-1)^mH_{l,\lambda}^{m,-m'-q}(\mathbf{b}))K_{l'-1,l',\lambda}^{m_1,-m'-q,q}-i\cdot(H_{l,\lambda}^{-m,m'-q}(\mathbf{b})+(-1)^mH_{l,\lambda}^{m,m'-q}(\mathbf{b}))K_{l'-1,l',\lambda}^{m_1,m'-q,q}\Bigr]\\
&+\rho^{l+1}a_{22}(2l+1)\dfrac{\sqrt{l'(2l'+1)}}{2}\sum_{q=-1}^{1}(-1)^q b_{-q}r^{l'-1}\Bigl[i\cdot(-1)^{m'}(H_{l,l'-1}^{-m,-m'-q}(\mathbf{b})+(-1)^mH_{l,l'-1}^{m,-m'-q}(\mathbf{b}))\\
&\quad\langle l'-1,-m'-q;1,q\mid l',-m' \rangle-i\cdot(H_{l,l'-1}^{-m,m'-q}(\mathbf{b})+(-1)^mH_{l,l'-1}^{m,m'-q}(\mathbf{b}))\langle l'-1,m'-q;1,q\mid l',m' \rangle\Bigr]
\end{aligned}
$$

If $m>0$ and $m'=0$,
$$
\begin{aligned}
&(W_{l'm'},\mathcal{S}_{D+n}[W_{lm}])\\
&=\dfrac{l'(2l'+1)}{\sqrt{2}}\Bigl[ \rho^{l+3}a_{12}r^{l'-1}+\rho^{l+1}(a_{22}-a_{12})(r^{l'+1}+a^2r^{l'-1})+\rho^{l+1}a_{22}r^{l'+1}\dfrac{2l+1}{2l'+1}\Bigr]\\
&\quad\Bigl[(H_{l,l'}^{-m,0}(\mathbf{b})+(-1)^mH_{l,l'}^{m,0}(\mathbf{b}))\Bigr]\\
&+\rho^{l+1}(a_{22}-a_{12})\sqrt{2l'(2l'+1)}\sum_{\lambda\in\{l'+1,l'-1\}}\sum_{q=-1}^{1}\sum_{m_1=-1}^{1}r^\lambda b_{-q}\Bigl[(H_{l,\lambda}^{-m,-q}(\mathbf{b})+(-1)^mH_{l,\lambda}^{m,-q}(\mathbf{b}))K_{l'-1,l',\lambda}^{m_1,-q,q}\Bigr]\\
&+\rho^{l+1}a_{22}(2l+1)\dfrac{\sqrt{l'(2l'+1)}}{\sqrt{2}}\sum_{q=-1}^{1}(-1)^q b_{-q}r^{l'-1}\Bigl[(H_{l,l'-1}^{-m,-q}(\mathbf{b})+(-1)^mH_{l,l'-1}^{m,-q}(\mathbf{b}))\langle l'-1,-q;1,q\mid l',0 \rangle\Bigr]
\end{aligned}
$$

If $m<0$ and $m'>0$,
$$
\begin{aligned}
&(W_{l'm'},\mathcal{S}_{D+n}[W_{lm}])\\
&=\dfrac{l'(2l'+1)}{2}\Bigl[ \rho^{l+3}a_{12}r^{l'-1}+\rho^{l+1}(a_{22}-a_{12})(r^{l'+1}+a^2r^{l'-1})+\rho^{l+1}a_{22}r^{l'+1}\dfrac{2l+1}{2l'+1}\Bigr]\\
&\quad i\Bigl[(-1)^{m'}(H_{l,l'}^{m,m'}(\mathbf{b})-(-1)^mH_{l,l'}^{-m,m'}(\mathbf{b}))+(H_{l,l'}^{m,-m'}(\mathbf{b})-(-1)^mH_{l,l'}^{-m,-m'}(\mathbf{b}))\Bigr]\\
&+\rho^{l+1}(a_{22}-a_{12})\sqrt{l'(2l'+1)}\sum_{\lambda\in\{l'+1,l'-1\}}\sum_{q=-1}^{1}\sum_{m_1=-1}^{1}r^\lambda b_{-q}\;i\Bigl[(-1)^{m'}(H_{l,\lambda}^{m,m'-q}(\mathbf{b})\\
&\quad-(-1)^mH_{l,\lambda}^{-m,m'-q}(\mathbf{b}))K_{l'-1,l',\lambda}^{m_1,m'-q,q}+(H_{l,\lambda}^{m,-m'-q}(\mathbf{b})-(-1)^mH_{l,\lambda}^{-m,-m'-q}(\mathbf{b}))K_{l'-1,l',\lambda}^{m_1,-m'-q,q}\Bigr]\\
&+\rho^{l+1}a_{22}(2l+1)\dfrac{\sqrt{l'(2l'+1)}}{2}\sum_{q=-1}^{1}(-1)^q b_{-q}r^{l'-1}\;i\Bigl[(-1)^{m'}(H_{l,l'-1}^{m,m'-q}(\mathbf{b})-(-1)^mH_{l,l'-1}^{-m,m'-q}(\mathbf{b}))\\
&\quad\langle l'-1,m'-q;1,q\mid l',m' \rangle+(H_{l,l'-1}^{m,-m'-q}(\mathbf{b})-(-1)^mH_{l,l'-1}^{-m,-m'-q}(\mathbf{b}))\langle l'-1,-m'-q;1,q\mid l',-m' \rangle\Bigr]
\end{aligned}
$$

If $m<0$ and $m'<0$,
$$
\begin{aligned}
&(W_{l'm'},\mathcal{S}_{D+n}[W_{lm}])\\
&=\dfrac{l'(2l'+1)}{2}\Bigl[ \rho^{l+3}a_{12}r^{l'-1}+\rho^{l+1}(a_{22}-a_{12})(r^{l'+1}+a^2r^{l'-1})+\rho^{l+1}a_{22}r^{l'+1}\dfrac{2l+1}{2l'+1}\Bigr]\\
&\quad\Bigl[-(-1)^{m'}(H_{l,l'}^{m,-m'}(\mathbf{b})-(-1)^mH_{l,l'}^{-m,-m'}(\mathbf{b}))+(H_{l,l'}^{m,m'}(\mathbf{b})-(-1)^mH_{l,l'}^{-m,m'}(\mathbf{b}))\Bigr]\\
&+\rho^{l+1}(a_{22}-a_{12})\sqrt{l'(2l'+1)}\sum_{\lambda\in\{l'+1,l'-1\}}\sum_{q=-1}^{1}\sum_{m_1=-1}^{1}r^\lambda b_{-q}\Bigl[-(-1)^{m'}(H_{l,\lambda}^{m,-m'-q}(\mathbf{b})\\
&\quad-(-1)^mH_{l,\lambda}^{-m,-m'-q}(\mathbf{b}))K_{l'-1,l',\lambda}^{m_1,-m'-q,q}+(H_{l,\lambda}^{m,m'-q}(\mathbf{b})-(-1)^mH_{l,\lambda}^{-m,m'-q}(\mathbf{b}))K_{l'-1,l',\lambda}^{m_1,m'-q,q}\Bigr]\\
&+\rho^{l+1}a_{22}(2l+1)\dfrac{\sqrt{l'(2l'+1)}}{2}\sum_{q=-1}^{1}(-1)^q b_{-q}r^{l'-1}\Bigl[-(-1)^{m'}(H_{l,l'-1}^{m,-m'-q}(\mathbf{b})-(-1)^mH_{l,l'-1}^{-m,-m'-q}(\mathbf{b}))\\
&\quad\langle l'-1,-m'-q;1,q\mid l',-m' \rangle+(H_{l,l'-1}^{m,m'-q}(\mathbf{b})-(-1)^mH_{l,l'-1}^{-m,m'-q}(\mathbf{b}))\langle l'-1,m'-q;1,q\mid l',m' \rangle\Bigr]
\end{aligned}
$$

If $m<0$ and $m'=0$,
$$
\begin{aligned}
&(W_{l'm'},\mathcal{S}_{D+n}[W_{lm}])\\
&=\dfrac{l'(2l'+1)}{\sqrt{2}}\Bigl[ \rho^{l+3}a_{12}r^{l'-1}+\rho^{l+1}(a_{22}-a_{12})(r^{l'+1}+a^2r^{l'-1})+\rho^{l+1}a_{22}r^{l'+1}\dfrac{2l+1}{2l'+1}\Bigr]\\
&\quad i\Bigl[(H_{l,l'}^{m,0}(\mathbf{b})-(-1)^mH_{l,l'}^{-m,0}(\mathbf{b}))\Bigr]\\
&+\rho^{l+1}(a_{22}-a_{12})\sqrt{2l'(2l'+1)}\sum_{\lambda\in\{l'+1,l'-1\}}\sum_{q=-1}^{1}\sum_{m_1=-1}^{1}r^\lambda b_{-q}\;i\Bigl[(H_{l,\lambda}^{m,-q}(\mathbf{b})-(-1)^mH_{l,\lambda}^{-m,-q}(\mathbf{b}))K_{l'-1,l',\lambda}^{m_1,-q,q}\Bigr]\\
&+\rho^{l+1}a_{22}(2l+1)\dfrac{\sqrt{l'(2l'+1)}}{\sqrt{2}}\sum_{q=-1}^{1}(-1)^q b_{-q}r^{l'-1}\;i\Bigl[(H_{l,l'-1}^{m,-q}(\mathbf{b})-(-1)^mH_{l,l'-1}^{-m,-q}(\mathbf{b}))\langle l'-1,-q;1,q\mid l',0 \rangle\Bigr]
\end{aligned}
$$

If $m=0$ and $m'>0$,
$$
\begin{aligned}
&(W_{l'm'},\mathcal{S}_{D+n}[W_{lm}])\\
&=\dfrac{l'(2l'+1)}{\sqrt{2}}\Bigl[ \rho^{l+3}a_{12}r^{l'-1}+\rho^{l+1}(a_{22}-a_{12})(r^{l'+1}+a^2r^{l'-1})+\rho^{l+1}a_{22}r^{l'+1}\dfrac{2l+1}{2l'+1}\Bigr]\\
&\quad\Bigl[(-1)^{m'}H_{l,l'}^{0,m'}(\mathbf{b})+H_{l,l'}^{0,-m'}(\mathbf{b})\Bigr]\\
&+\rho^{l+1}(a_{22}-a_{12})\sqrt{2l'(2l'+1)}\sum_{\lambda\in\{l'+1,l'-1\}}\sum_{q=-1}^{1}\sum_{m_1=-1}^{1}r^\lambda b_{-q}\Bigl[(-1)^{m'}H_{l,\lambda}^{0,m'-q}(\mathbf{b})K_{l'-1,l',\lambda}^{m_1,m'-q,q}\\
&+H_{l,\lambda}^{0,-m'-q}(\mathbf{b})K_{l'-1,l',\lambda}^{m_1,-m'-q,q}\Bigr]\\
&+\rho^{l+1}a_{22}(2l+1)\dfrac{\sqrt{l'(2l'+1)}}{\sqrt{2}}\sum_{q=-1}^{1}(-1)^q b_{-q}r^{l'-1}\Bigl[(-1)^{m'}H_{l,l'-1}^{0,m'-q}(\mathbf{b})\langle l'-1,m'-q;1,q\mid l',m' \rangle\\
&+H_{l,l'-1}^{0,-m'-q}(\mathbf{b})\langle l'-1,-m'-q;1,q\mid l',-m' \rangle\Bigr]
\end{aligned}
$$

If $m=0$ and $m'<0$,
$$
\begin{aligned}
&(W_{l'm'},\mathcal{S}_{D+n}[W_{lm}])\\
&=\dfrac{l'(2l'+1)}{\sqrt{2}}\Bigl[ \rho^{l+3}a_{12}r^{l'-1}+\rho^{l+1}(a_{22}-a_{12})(r^{l'+1}+a^2r^{l'-1})+\rho^{l+1}a_{22}r^{l'+1}\dfrac{2l+1}{2l'+1}\Bigr]\\
&\quad i\Bigl[(-1)^{m'}H_{l,l'}^{0,-m'}(\mathbf{b})-H_{l,l'}^{0,m'}(\mathbf{b})\Bigr]\\
&+\rho^{l+1}(a_{22}-a_{12})\sqrt{2l'(2l'+1)}\sum_{\lambda\in\{l'+1,l'-1\}}\sum_{q=-1}^{1}\sum_{m_1=-1}^{1}r^\lambda b_{-q}\;i\Bigl[(-1)^{m'}H_{l,\lambda}^{0,-m'-q}(\mathbf{b})K_{l'-1,l',\lambda}^{m_1,-m'-q,q}\\
&-H_{l,\lambda}^{0,m'-q}(\mathbf{b})K_{l'-1,l',\lambda}^{m_1,m'-q,q}\Bigr]\\
&+\rho^{l+1}a_{22}(2l+1)\dfrac{\sqrt{l'(2l'+1)}}{\sqrt{2}}\sum_{q=-1}^{1}(-1)^q b_{-q}r^{l'-1}\;i\Bigl[(-1)^{m'}H_{l,l'-1}^{0,-m'-q}(\mathbf{b})\langle l'-1,-m'-q;1,q\mid l',-m' \rangle\\
&-H_{l,l'-1}^{0,m'-q}(\mathbf{b})\langle l'-1,m'-q;1,q\mid l',m' \rangle\Bigr]
\end{aligned}
$$

If $m=0$ and $m'=0$,
$$
\begin{aligned}
&(W_{l'm'},\mathcal{S}_{D+n}[W_{lm}])\\
&=l'(2l'+1)\Bigl[ \rho^{l+3}a_{12}r^{l'-1}+\rho^{l+1}(a_{22}-a_{12})(r^{l'+1}+a^2r^{l'-1})+\rho^{l+1}a_{22}r^{l'+1}\dfrac{2l+1}{2l'+1}\Bigr]H_{l,l'}^{0,0}(\mathbf{b})\\
&+\rho^{l+1}(a_{22}-a_{12})\cdot 2\sqrt{l'(2l'+1)}\sum_{\lambda\in\{l'+1,l'-1\}}\sum_{q=-1}^{1}\sum_{m_1=-1}^{1}r^\lambda b_{-q}\;H_{l,\lambda}^{0,-q}(\mathbf{b})K_{l'-1,l',\lambda}^{m_1,-q,q}\\
&+\rho^{l+1}a_{22}(2l+1)\sqrt{l'(2l'+1)}\sum_{q=-1}^{1}(-1)^q b_{-q}r^{l'-1}\;H_{l,l'-1}^{0,-q}(\mathbf{b})\langle l'-1,-q;1,q\mid l',0 \rangle
\end{aligned}
$$
\end{theorem}

\section{Analytic Summation of $M_{(l'm'p),(lmq)}$}

\begin{lemma}Let $\underline{Y_{lm}}=(Y^1_{lm},Y^2_{lm},Y^3_{lm})=(V_{lm},W_{lm},X_{lm}).$ Then, we have
$$
(\mathcal{S}_D\underline{Y_{l m}})(x)
=
\underline{Y_{l m}}(\hat{x})
A_{\mathcal{S}_D,l},
$$
where the matrix $ A_{\mathcal{S}_D,l}$ is given by
$$
A_{\mathcal{S}_D,l}
=
\begin{bmatrix}
\dfrac{(3l+1)\mu+l\lambda}
{(2l+3)(2l+1)\mu(2\mu+\lambda)}
&
0
&
0
\\[12pt]
0
&
\dfrac{(3l+2)\mu+(l+1)\lambda}
{(2l-1)(2l+1)\mu(2\mu+\lambda)}
&
0
\\[12pt]
0
&
0
&
\dfrac{1}{(2l+1)\mu}
\end{bmatrix}=\text{diag}(\tau_{\mathcal{S}_D,l}^1,\tau_{\mathcal{S}_D,l}^2,\tau_{\mathcal{S}_D,l}^3).
$$
\end{lemma}

Therefore, 
$$( Y_{l'm'}^p,\mathcal{S}_D Y_{lm}^q)= \delta_{ll'}\delta_{mm'}\delta_{pq}\cdot \rho\cdot \tau_{\mathcal{S}_D}^p(l)\cdot \text{norm}_p(l)$$
where  $ \text{norm}_1(l)=  (l+1)(2l+1)  $, $ \text{norm}_2(l)=  l(2l+1)  $, $ \text{norm}_3(l)=l(l+1)  $.

Define the entries 
\begin{equation}\label{eq-M-entries}M_{(l'm'p),(lmq)}(\alpha) = (Y_{l'm'}^p,\mathcal{S}^{\alpha,0}_D Y_{lm}^q)= ( Y_{l'm'}^p,\mathcal{S}_D Y_{lm}^q) + \sum_{n\neq 0}e^{-in\alpha}(Y_{l'm'}^p,\mathcal{S}_{D+n} Y_{lm}^q).\end{equation}
\begin{remark}Since the operator $\mathcal{S}_D^{\alpha,0}$ is complex, there should be a conjugation on $e^{in\alpha}$ when it is extracted from the second slot of inner product, hence $e^{-in\alpha}$.
\end{remark}

\begin{definition}The {\bf polylogarithm function} is defined by a power series in $z$ generalizing the Mercator series, which is also a Dirichlet series in $s$\cite{lewin1981polylogarithms}:
$$\mathrm{Li}_s(z)=\sum\limits_{k=1}^\infty\dfrac{z^k}{k^s}.$$
\end{definition}

\begin{lemma}Denote the series $\sum\limits_{n\neq 0}H_{l,\lambda}^{m,\mu}(\mathbf{b})e^{-in\alpha}$ by $\mathscr{H}_{l,\lambda}^{m,\mu}(\alpha)$. Then
$$\begin{aligned}
\sum\limits_{n\neq 0}H_{l,\lambda}^{m,\mu}(\mathbf{b})e^{-in\alpha}=(-1)^{\lambda+\mu}
\sqrt{\dfrac{2l+1}{2\lambda+1}}&
\sqrt{{l+\lambda+\mu-m\choose \lambda+\mu}{l+\lambda+m-\mu\choose \lambda-\mu}} \sqrt{\dfrac{4\pi}{2(l+\lambda)+1}}\\
&\big(\mathrm{Li}_{l+\lambda+1}(e^{-i\alpha})Y_{l+\lambda}^{m-\mu}(\dfrac{\pi}{2},\pi)+\mathrm{Li}_{l+\lambda+1}(e^{i\alpha})Y_{l+\lambda}^{m-\mu}(\dfrac{\pi}{2},0)\big)
\end{aligned}.$$
\end{lemma}

\begin{theorem}The expression of $\boxed{M_{(l'm'2),(lm1)}(\alpha)}$($l'\geqslant 1$, otherwise this term vanishes) is given by:
\\
If $m>0$ and $m'>0$,
$$\begin{aligned}M_{(l'm'2),(lm1)}(\alpha)= \rho^{l+3}a_{11}[r^{l'-1}\dfrac{1}{2}(\mathscr{H}_{l,l'}^{-m,-m'}(\alpha)+(-1)^m\mathscr{H}_{l,l'}^{m,-m'}(\alpha))l'(2l'+1)\\
+r^{l'-1}\dfrac{(-1)^{m'}}{2}(\mathscr{H}_{l,l'}^{-m,m'}(\alpha)+(-1)^m\mathscr{H}_{l,l'}^{m,m'}(\alpha))l'(2l'+1)]
\end{aligned}$$
If $m>0$ and $m'<0$,
$$\begin{aligned}M_{(l'm'2),(lm1)}(\alpha)= \rho^{l+3}a_{11}[r^{l'-1}\dfrac{1}{2}(\mathscr{H}_{l,l'}^{-m,m'}(\alpha)+(-1)^m\mathscr{H}_{l,l'}^{m,m'}(\alpha))(-i \cdot l'(2l'+1))\\
+r^{l'-1}\dfrac{(-1)^{m'}}{2}(\mathscr{H}_{l,l'}^{-m,-m'}(\alpha)+(-1)^m\mathscr{H}_{l,l'}^{m,-m'}(\alpha))(i\cdot l'(2l'+1))]
\end{aligned}$$
If $m>0$ and $m'=0$,
$$\begin{aligned}M_{(l'm'2),(lm1)}(\alpha)= \rho^{l+3}a_{11}r^{l'-1}\dfrac{1}{\sqrt{2}}(\mathscr{H}_{l,l'}^{-m,0}(\alpha)+(-1)^m\mathscr{H}_{l,l'}^{m,0}(\alpha))l'(2l'+1).
\end{aligned}$$

If $m<0$ and $m'>0$,
$$\begin{aligned}M_{(l'm'2),(lm1)}(\alpha)= \rho^{l+3}a_{11}[r^{l'-1}\dfrac{i}{2}(\mathscr{H}_{l,l'}^{m,-m'}(\alpha)-(-1)^m \mathscr{H}_{l,l'}^{-m,-m'}(\alpha))l'(2l'+1)\\
+r^{l'-1}\dfrac{(-1)^{m'}}{2}\cdot i (\mathscr{H}_{l,l'}^{m,m'}(\alpha)-(-1)^m\mathscr{H}_{l,l'}^{-m,m'}(\alpha))l'(2l'+1)]
\end{aligned}$$

If $m<0$ and $m'<0$,
$$\begin{aligned}M_{(l'm'2),(lm1)}(\alpha)= \rho^{l+3}a_{11}[r^{l'-1}\dfrac{i}{2}(\mathscr{H}_{l,l'}^{m,m'}(\alpha)-(-1)^m \mathscr{H}_{l,l'}^{-m,m'}(\alpha))(-i\cdot l'(2l'+1))\\
+r^{l'-1}\dfrac{(-1)^{m'}}{2}\cdot i (\mathscr{H}_{l,l'}^{m,-m'}(\alpha)-(-1)^m\mathscr{H}_{l,l'}^{-m,-m'}(\alpha))(i\cdot l'(2l'+1))]
\end{aligned}$$

If $m<0$ and $m'=0$,
$$\begin{aligned}M_{(l'm'2),(lm1)}(\alpha)= \rho^{l+3}a_{11}r^{l'-1}\dfrac{i}{\sqrt{2}}(\mathscr{H}_{l,l'}^{m,0}(\alpha)-(-1)^m\mathscr{H}_{l,l'}^{-m,0}(\alpha))l'(2l'+1).\end{aligned}$$

If $m=0$ and $m'>0$,
$$\begin{aligned}M_{(l'm'2),(lm1)}(\alpha)= \rho^{l+3}a_{11}r^{l'-1}
[\dfrac{1}{\sqrt{2}} \mathscr{H}_{l,l'}^{0,-m'}(\alpha)l'(2l'+1)+\dfrac{(-1)^{m'}}{\sqrt{2}}\mathscr{H}_{l,l'}^{0,m'}(\alpha)l'(2l'+1)]
\end{aligned}$$

If $m=0$ and $m'<0$,
$$\begin{aligned}M_{(l'm'2),(lm1)}(\alpha)= \rho^{l+3}a_{11}r^{l'-1}
[\dfrac{1}{\sqrt{2}} \mathscr{H}_{l,l'}^{0,m'}(\alpha)(-i\cdot l'(2l'+1))+\dfrac{(-1)^{m'}}{\sqrt{2}}\mathscr{H}_{l,l'}^{0,-m'}(\alpha)(i\cdot l'(2l'+1))]
\end{aligned}$$

If $m=0$ and $m'=0$,
$$\begin{aligned}M_{(l'm'2),(lm1)}(\alpha)= \rho^{l+3}a_{11}r^{l'-1}
[\mathscr{H}_{l,l'}^{0,0}(\alpha)l'(2l'+1)]
\end{aligned}$$
where $a_{ij}$ stands for the coefficient of the entries of $A^{out}_{\mathcal{S},l}(x)$.
\end{theorem}

\begin{lemma}Denote $\sum\limits_{n\neq 0}L_{l,j,\lambda}^{m,\mu,q,m_1}(\mathbf{b})e^{-in\alpha}$ by $\mathscr{L}_{l,j,\lambda}^{m,\mu,q,m_1}(\alpha)$. Then
$$\begin{aligned}
\mathscr{L}_{l,j,\lambda}^{m,\mu,q,m_1}(\alpha)&=
i (-1)^{\lambda+\mu+q}\text{sgn}(q-m_1)\epsilon_{q}\sqrt{\lambda(2l+1)} 
\sqrt{{l+\lambda+\mu-m\choose \lambda+\mu}{l+\lambda+m-\mu\choose \lambda-\mu}}\\
&\langle\lambda-1,\mu-m_1;1,m_1\mid \lambda,\mu\rangle
\langle \lambda-1,\mu-m_1;1,q+m_1\mid j,\mu+q \rangle \sqrt{\dfrac{4\pi}{2(l+\lambda)+1}}\\
&\big(\mathrm{Li}_{l+\lambda}(e^{-i\alpha})Y_{l+\lambda}^{m-\mu}(\dfrac{\pi}{2},\pi)-\mathrm{Li}_{l+\lambda}(e^{i\alpha})Y_{l+\lambda}^{m-\mu}(\dfrac{\pi}{2},0)\big)
\end{aligned}$$
where $$\epsilon_q=\left\{
\begin{aligned}
\dfrac{1}{\sqrt{2}}&,\quad q=-1\\
0&,\quad q=0\\
-\dfrac{1}{\sqrt{2}}&,\quad,q=1 \end{aligned}\right.$$
\end{lemma}

\begin{theorem}
The expression of {$\boxed{M_{(l'm',2),(lm,3)}(\alpha)}$}($l'\geqslant 1,l\geqslant 1$, otherwise this term vanishes) is given by
When $m>0, m'>0$,
$$
\begin{aligned}
&M_{(l'm',2),(lm,3)}(\alpha)\\
=&\ \rho^{l+2}a_{33}\,r^{l'-1}\,\frac{\sqrt{l'(2l'+1)}}{2}\sum\limits_{q=-1}^{1}\sum\limits_{m_1=-1}^{1}\bigl[(-1)^{m'}\bigl(\mathscr{L}_{l,l',l'}^{-m,\,m'-q,\,q,\,m_1}(\alpha)+(-1)^m\mathscr{L}_{l,l',l'}^{m,\,m'-q,\,q,\,m_1}(\alpha)\bigr)\\
&\qquad\qquad\qquad\qquad+\bigl(\mathscr{L}_{l,l',l'}^{-m,\,-m'-q,\,q,\,m_1}(\alpha)+(-1)^m\mathscr{L}_{l,l',l'}^{m,\,-m'-q,\,q,\,m_1}(\alpha)\bigr)\bigr]
\end{aligned}
$$
When  $m>0, m'<0$,
$$
\begin{aligned}
&M_{(l'm',2),(lm,3)}(\alpha)\\
=&\ \rho^{l+2}a_{33}\,r^{l'-1}\,\frac{\sqrt{l'(2l'+1)}}{2}\sum\limits_{q=-1}^{1}\sum\limits_{m_1=-1}^{1}\bigl[i(-1)^{m'}\bigl(\mathscr{L}_{l,l',l'}^{-m,\,-m'-q,\,q,\,m_1}(\alpha)+(-1)^m\mathscr{L}_{l,l',l'}^{m,\,-m'-q,\,q,\,m_1}(\alpha)\bigr)\\
&\qquad\qquad\qquad\qquad-i\bigl(\mathscr{L}_{l,l',l'}^{-m,\,m'-q,\,q,\,m_1}(\alpha)+(-1)^m\mathscr{L}_{l,l',l'}^{m,\,m'-q,\,q,\,m_1}(\alpha)\bigr)\bigr]
\end{aligned}
$$

When $m>0, m'=0$,
$$
\begin{aligned}
&M_{(l',0,2),(lm,3)}(\alpha)\\
=&\ \rho^{l+2}a_{33}\,r^{l'-1}\,\frac{\sqrt{l'(2l'+1)}}{\sqrt{2}}\sum\limits_{q=-1}^{1}\sum\limits_{m_1=-1}^{1}\bigl(\mathscr{L}_{l,l',l'}^{-m,\,-q,\,q,\,m_1}(\alpha)+(-1)^m\mathscr{L}_{l,l',l'}^{m,\,-q,\,q,\,m_1}(\alpha)\bigr)
\end{aligned}
$$

When  $m<0, m'>0$,
$$
\begin{aligned}
&M_{(l'm',2),(lm,3)}(\alpha)\\
=&\ \rho^{l+2}a_{33}\,r^{l'-1}\,i\cdot\frac{\sqrt{l'(2l'+1)}}{2}\sum\limits_{q=-1}^{1}\sum\limits_{m_1=-1}^{1}\bigl[(-1)^{m'}\bigl(\mathscr{L}_{l,l',l'}^{m,\,m'-q,\,q,\,m_1}(\alpha)-(-1)^m\mathscr{L}_{l,l',l'}^{-m,\,m'-q,\,q,\,m_1}(\alpha)\bigr)\\
&\qquad\qquad\qquad\qquad+\bigl(\mathscr{L}_{l,l',l'}^{m,\,-m'-q,\,q,\,m_1}(\alpha)-(-1)^m\mathscr{L}_{l,l',l'}^{-m,\,-m'-q,\,q,\,m_1}(\alpha)\bigr)\bigr]
\end{aligned}
$$
When $m<0, m'<0$,
$$
\begin{aligned}
&M_{(l'm',2),(lm,3)}(\alpha)\\
=&\ \rho^{l+2}a_{33}\,r^{l'-1}\,\frac{\sqrt{l'(2l'+1)}}{2}\sum\limits_{q=-1}^{1}\sum\limits_{m_1=-1}^{1}\bigl[\bigl(\mathscr{L}_{l,l',l'}^{m,\,m'-q,\,q,\,m_1}(\alpha)-(-1)^m\mathscr{L}_{l,l',l'}^{-m,\,m'-q,\,q,\,m_1}(\alpha)\bigr)\\
&\qquad\qquad\qquad\qquad-(-1)^{m'}\bigl(\mathscr{L}_{l,l',l'}^{m,\,-m'-q,\,q,\,m_1}(\alpha)-(-1)^m\mathscr{L}_{l,l',l'}^{-m,\,-m'-q,\,q,\,m_1}(\alpha)\bigr)\bigr]
\end{aligned}
$$
When $m<0, m'=0$,
$$
\begin{aligned}
&M_{(l',0,2),(lm,3)}(\alpha)\\
=&\ \rho^{l+2}a_{33}\,r^{l'-1}\,\frac{i\sqrt{l'(2l'+1)}}{\sqrt{2}}\sum\limits_{q=-1}^{1}\sum\limits_{m_1=-1}^{1}\bigl(\mathscr{L}_{l,l',l'}^{m,\,-q,\,q,\,m_1}(\alpha)-(-1)^m\mathscr{L}_{l,l',l'}^{-m,\,-q,\,q,\,m_1}(\alpha)\bigr)
\end{aligned}
$$
When $m=0, m'>0$,
$$
\begin{aligned}
&M_{(l'm',2),(l,0,3)}(\alpha)\\
=&\ \rho^{l+2}a_{33}\,r^{l'-1}\,\frac{\sqrt{l'(2l'+1)}}{\sqrt{2}}\sum\limits_{q=-1}^{1}\sum\limits_{m_1=-1}^{1}\bigl[(-1)^{m'}\mathscr{L}_{l,l',l'}^{0,\,m'-q,\,q,\,m_1}(\alpha)+\mathscr{L}_{l,l',l'}^{0,\,-m'-q,\,q,\,m_1}(\alpha)\bigr]
\end{aligned}
$$

When  $m=0, m'<0$,
$$
\begin{aligned}
&M_{(l'm',2),(l,0,3)}(\alpha)\\
=&\ \rho^{l+2}a_{33}\,r^{l'-1}\,\frac{\sqrt{l'(2l'+1)}}{\sqrt{2}}\sum\limits_{q=-1}^{1}\sum\limits_{m_1=-1}^{1}\bigl[i(-1)^{m'}\mathscr{L}_{l,l',l'}^{0,\,-m'-q,\,q,\,m_1}(\alpha)-i\,\mathscr{L}_{l,l',l'}^{0,\,m'-q,\,q,\,m_1}(\alpha)\bigr]
\end{aligned}
$$

When $m=0, m'=0$,
$$
\begin{aligned}
&M_{(l',0,2),(l,0,3)}(\alpha)\\
=&\ \rho^{l+2}a_{33}\,r^{l'-1}\,\sqrt{l'(2l'+1)}\sum\limits_{q=-1}^{1}\sum\limits_{m_1=-1}^{1}\mathscr{L}_{l,l',l'}^{0,\,-q,\,q,\,m_1}(\alpha)
\end{aligned}
$$
\end{theorem}

\begin{lemma}Denote $\sum\limits_{n\neq 0}b_{-q}H_{l,\lambda}^{m,\mu}(\mathbf{b})e^{-in\alpha}$ by $\mathscr{A}_{l,\lambda}^{m,\mu}(\alpha,q)$, then 
$$\begin{aligned}
\sum\limits_{n\neq 0}b_{-q}H_{l,\lambda}^{m,\mu}(\mathbf{b})e^{-in\alpha}
=(-1)^{\lambda+\mu}\epsilon_q
\sqrt{\dfrac{2l+1}{2\lambda+1}}&
\sqrt{{l+\lambda+\mu-m\choose \lambda+\mu}{l+\lambda+m-\mu\choose \lambda-\mu}} \sqrt{\dfrac{4\pi}{2(l+\lambda)+1}}\\
&
\big(\mathrm{Li}_{l+\lambda}(e^{-i\alpha})Y_{l+\lambda}^{m-\mu}(\dfrac{\pi}{2},\pi)-\mathrm{Li}_{l+\lambda}(e^{i\alpha})Y_{l+\lambda}^{m-\mu}(\dfrac{\pi}{2},0)\big)
\end{aligned}.$$

\end{lemma}

\begin{theorem}
The expression of {$\boxed{M_{(l'm',1),(lm,2)}(\alpha)}$}($l\geqslant 1$, otherwise this term vanishes) is given by\\
When $m>0$, $m'>0$,
$$
\begin{aligned}
&M_{(l'm',1),(lm,2)}(\alpha)\\
=&\sum\limits_{\lambda\in\{l'+1,l'+3\}}\sum\limits_{q=-1}^{1}\sum\limits_{m_1=-1}^{1}\rho^{l+1} (a_{22}-a_{12})  r^\lambda \sqrt{(l'+1)(2l'+1)}
[(-1)^{m'}(\mathscr{A}_{l,\lambda}^{-m,m'-q}(\alpha,q)\\
&+(-1)^m\mathscr{A}_{l,\lambda}^{m,m'-q}(\alpha,q))K_{l'+1,l',\lambda}^{m_1,m'-q,q}
+(\mathscr{A}_{l,\lambda}^{-m,-m'-q}(\alpha,q)+(-1)^m\mathscr{A}_{l,\lambda}^{m,-m'-q}(\alpha,q))K_{l'+1,l',\lambda}^{m_1,-m'-q,q}]\\
&-\rho^{l+1} a_{22}r^{l'+1}\dfrac{(2l+1)(l'+1)}{2}[(-1)^{m'}(\mathscr{H}_{l,l'}^{-m,m'}(\alpha)+(-1)^m\mathscr{H}_{l,l'}^{m,m'}(\alpha))\\
&+(\mathscr{H}_{l,l'}^{-m,-m'}(\alpha)+(-1)^m\mathscr{H}_{l,l'}^{m,-m'}(\alpha))]\\
&+\rho^{l+1} a_{22}r^{l'+1}(2l+1)\dfrac{\sqrt{(l'+1)(2l'+1)}}{2} \sum\limits_{q=-1}^1(-1)^q [(-1)^{m'}(\mathscr{A}_{l,l'+1}^{-m,m'-q}(\alpha,q)+(-1)^m\mathscr{A}_{l,l'+1}^{m,m'-q}(\alpha,q))\\
&\langle l'+1,m'-q;1,q\mid l',m'\rangle+(\mathscr{A}_{l,l'+1}^{-m,-m'-q}(\alpha,q)+(-1)^m\mathscr{A}_{l,l'+1}^{m,-m'-q}(\alpha,q))\langle l'+1,-m'-q;1,q\mid l',-m'\rangle]
\end{aligned}$$

When $m>0$, $m'<0$,

$$
\begin{aligned}
&M_{(l'm',1),(lm,2)}(\alpha)\\
=&\sum\limits_{\lambda\in\{l'+1,l'+3\}}\sum\limits_{q=-1}^{1}\sum\limits_{m_1=-1}^{1}\rho^{l+1} (a_{22}-a_{12})  r^\lambda \sqrt{(l'+1)(2l'+1)}
[i\cdot (-1)^{m'}(\mathscr{A}_{l,\lambda}^{-m,-m'-q}(\alpha,q)\\
&+(-1)^m\mathscr{A}_{l,\lambda}^{m,-m'-q}(\alpha,q))K_{l'+1,l',\lambda}^{m_1,-m'-q,q}
-i\cdot (\mathscr{A}_{l,\lambda}^{-m,m'-q}(\alpha,q)+(-1)^m\mathscr{A}_{l,\lambda}^{m,m'-q}(\alpha,q))K_{l'+1,l',\lambda}^{m_1,m'-q,q}]\\
&-\rho^{l+1} a_{22}r^{l'+1}\dfrac{(2l+1)(l'+1)}{2}[i\cdot (-1)^{m'}(\mathscr{H}_{l,l'}^{-m,-m'}(\alpha)+(-1)^m\mathscr{H}_{l,l'}^{m,-m'}(\alpha))\\
&-i\cdot (\mathscr{H}_{l,l'}^{-m,m'}(\alpha)+(-1)^m\mathscr{H}_{l,l'}^{m,m'}(\alpha))]\\
&+\rho^{l+1} a_{22}r^{l'+1}(2l+1)\dfrac{\sqrt{(l'+1)(2l'+1)}}{2} \sum\limits_{q=-1}^1(-1)^q [i\cdot (-1)^{m'}(\mathscr{A}_{l,l'+1}^{-m,-m'-q}(\alpha,q)\\
&+(-1)^m\mathscr{A}_{l,l'+1}^{m,-m'-q}(\alpha,q))\langle l'+1,-m'-q;1,q\mid l',-m'\rangle
-i\cdot (\mathscr{A}_{l,l'+1}^{-m,m'-q}(\alpha,q)\\
&+(-1)^m\mathscr{A}_{l,l'+1}^{m,m'-q}(\alpha,q))\langle l'+1,m'-q;1,q\mid l',m'\rangle]
\end{aligned}$$

When $m>0$, $m'=0$,
$$
\begin{aligned}
&M_{(l',0,1),(lm,2)}(\alpha)\\
=&\sum\limits_{\lambda\in\{l'+1,l'+3\}}\sum\limits_{q=-1}^{1}\sum\limits_{m_1=-1}^{1}\rho^{l+1} (a_{22}-a_{12})  r^\lambda \sqrt{2(l'+1)(2l'+1)}
[(\mathscr{A}_{l,\lambda}^{-m,-q}(\alpha,q)+\\
&(-1)^m\mathscr{A}_{l,\lambda}^{m,-q}(\alpha,q))K_{l'+1,l',\lambda}^{m_1,-q,q}]\\
&-\rho^{l+1} a_{22}r^{l'+1}\dfrac{(2l+1)(l'+1)}{\sqrt{2}}[\mathscr{H}_{l,l'}^{-m,0}(\alpha)+(-1)^m\mathscr{H}_{l,l'}^{m,0}(\alpha)]\\
&+\rho^{l+1} a_{22}r^{l'+1}(2l+1)\dfrac{\sqrt{(l'+1)(2l'+1)}}{\sqrt{2}} \sum\limits_{q=-1}^1(-1)^q (\mathscr{A}_{l,l'+1}^{-m,-q}(\alpha,q)+(-1)^m\mathscr{A}_{l,l'+1}^{m,-q}(\alpha,q))\\
&\langle l'+1,-q;1,q\mid l',0\rangle
\end{aligned}$$

When $m<0$, $m'>0$,
$$
\begin{aligned}
&M_{(l'm',1),(lm,2)}(\alpha)\\
=&\sum\limits_{\lambda\in\{l'+1,l'+3\}}\sum\limits_{q=-1}^{1}\sum\limits_{m_1=-1}^{1}\rho^{l+1} (a_{22}-a_{12})  r^\lambda i\cdot\sqrt{(l'+1)(2l'+1)}
[(-1)^{m'}(\mathscr{A}_{l,\lambda}^{m,m'-q}(\alpha,q)\\
&-(-1)^m\mathscr{A}_{l,\lambda}^{-m,m'-q}(\alpha,q))K_{l'+1,l',\lambda}^{m_1,m'-q,q}
+(\mathscr{A}_{l,\lambda}^{m,-m'-q}(\alpha,q)-(-1)^m\mathscr{A}_{l,\lambda}^{-m,-m'-q}(\alpha,q))K_{l'+1,l',\lambda}^{m_1,-m'-q,q}]\\
&-\rho^{l+1} a_{22}r^{l'+1}i\cdot \dfrac{(2l+1)(l'+1)}{2}[(-1)^{m'}(\mathscr{H}_{l,l'}^{m,m'}(\alpha)-(-1)^m\mathscr{H}_{l,l'}^{-m,m'}(\alpha))\\
&+(\mathscr{H}_{l,l'}^{m,-m'}(\alpha)-(-1)^m\mathscr{H}_{l,l'}^{-m,-m'}(\alpha))]\\
&+\rho^{l+1} a_{22}r^{l'+1}i\cdot (2l+1)\dfrac{\sqrt{(l'+1)(2l'+1)}}{2} \sum\limits_{q=-1}^1(-1)^q [(-1)^{m'}(\mathscr{A}_{l,l'+1}^{m,m'-q}(\alpha,q)-(-1)^m\mathscr{A}_{l,l'+1}^{-m,m'-q}(\alpha,q))\\
&\langle l'+1,m'-q;1,q\mid l',m'\rangle+(\mathscr{A}_{l,l'+1}^{m,-m'-q}(\alpha,q)-(-1)^m\mathscr{A}_{l,l'+1}^{-m,-m'-q}(\alpha,q))\langle l'+1,-m'-q;1,q\mid l',-m'\rangle]
\end{aligned}$$

When $m<0$, $m'<0$,

$$
\begin{aligned}
&M_{(l'm',1),(lm,2)}(\alpha)\\
=&\sum\limits_{\lambda\in\{l'+1,l'+3\}}\sum\limits_{q=-1}^{1}\sum\limits_{m_1=-1}^{1}\rho^{l+1} (a_{22}-a_{12})  r^\lambda \sqrt{(l'+1)(2l'+1)}
[-(-1)^{m'}(\mathscr{A}_{l,\lambda}^{m,-m'-q}(\alpha,q)\\
&-(-1)^m\mathscr{A}_{l,\lambda}^{-m,-m'-q}(\alpha,q))K_{l'+1,l',\lambda}^{m_1,-m'-q,q}
+ (\mathscr{A}_{l,\lambda}^{m,m'-q}(\alpha,q)-(-1)^m\mathscr{A}_{l,\lambda}^{-m,m'-q}(\alpha,q))K_{l'+1,l',\lambda}^{m_1,m'-q,q}]\\
&-\rho^{l+1} a_{22}r^{l'+1}\dfrac{(2l+1)(l'+1)}{2}
[- (-1)^{m'}(\mathscr{H}_{l,l'}^{m,-m'}(\alpha)-(-1)^m\mathscr{H}_{l,l'}^{-m,-m'}(\alpha))\\
&+
 (\mathscr{H}_{l,l'}^{m,m'}(\alpha)-(-1)^m\mathscr{H}_{l,l'}^{-m,m'}(\alpha))]\\
&+\rho^{l+1} a_{22}r^{l'+1}(2l+1)\dfrac{\sqrt{(l'+1)(2l'+1)}}{2} \sum\limits_{q=-1}^1(-1)^q [-(-1)^{m'}(\mathscr{A}_{l,l'+1}^{m,-m'-q}(\alpha,q)\\
&-(-1)^m\mathscr{A}_{l,l'+1}^{-m,-m'-q}(\alpha,q))
\langle l'+1,-m'-q;1,q\mid l',-m'\rangle
+ (\mathscr{A}_{l,l'+1}^{m,m'-q}(\alpha,q)\\
&-(-1)^m\mathscr{A}_{l,l'+1}^{-m,m'-q}(\alpha,q))\langle l'+1,m'-q;1,q\mid l',m'\rangle]
\end{aligned}$$

When $m<0$, $m'=0$,
$$
\begin{aligned}
&M_{(l',0,1),(lm,2)}(\alpha)\\
=&\sum\limits_{\lambda\in\{l'+1,l'+3\}}\sum\limits_{q=-1}^{1}\sum\limits_{m_1=-1}^{1}\rho^{l+1} (a_{22}-a_{12})  r^\lambda i\cdot \sqrt{2(l'+1)(2l'+1)}
[(\mathscr{A}_{l,\lambda}^{m,-q}(\alpha,q)\\
&-(-1)^m\mathscr{A}_{l,\lambda}^{-m,-q}(\alpha,q))K_{l'+1,l',\lambda}^{m_1,-q,q}]\\
&-\rho^{l+1} a_{22}r^{l'+1}i\cdot \dfrac{(2l+1)(l'+1)}{\sqrt{2}}[(\mathscr{H}_{l,l'}^{m,0}(\alpha)-(-1)^m\mathscr{H}_{l,l'}^{-m,0}(\alpha)]\\
&+\rho^{l+1} a_{22}r^{l'+1}i\cdot (2l+1)\dfrac{\sqrt{(l'+1)(2l'+1)}}{\sqrt{2}} \sum\limits_{q=-1}^1(-1)^q [ (\mathscr{A}_{l,l'+1}^{m,-q}(\alpha,q)-(-1)^m\mathscr{A}_{l,l'+1}^{-m,-q}(\alpha,q))\\
&\langle l'+1,-q;1,q\mid l',0\rangle]
\end{aligned}$$

When $m=0$, $m'>0$,
$$
\begin{aligned}
&M_{(l'm',1),(l,0,2)}(\alpha)\\
=&\sum\limits_{\lambda\in\{l'+1,l'+3\}}\sum\limits_{q=-1}^{1}\sum\limits_{m_1=-1}^{1}\rho^{l+1} (a_{22}-a_{12})  r^\lambda \sqrt{2(l'+1)(2l'+1)}
[(-1)^{m'}\mathscr{A}_{l,\lambda}^{0,m'-q}(\alpha,q) K_{l'+1,l',\lambda}^{m_1,m'-q,q}\\
&+\mathscr{A}_{l,\lambda}^{0,-m'-q}(\alpha,q)K_{l'+1,l',\lambda}^{m_1,-m'-q,q}]\\
&-\rho^{l+1} a_{22}r^{l'+1}\dfrac{(2l+1)(l'+1)}{\sqrt{2}}[(-1)^{m'}\mathscr{H}_{l,l'}^{0,m'}(\alpha)+\mathscr{H}_{l,l'}^{0,-m'}(\alpha)]\\
&+\rho^{l+1} a_{22}r^{l'+1}(2l+1)\dfrac{\sqrt{(l'+1)(2l'+1)}}{\sqrt{2}} \sum\limits_{q=-1}^1(-1)^q [(-1)^{m'}\mathscr{A}_{l,l'+1}^{0,m'-q}(\alpha,q)\\
&\langle l'+1,m'-q;1,q\mid l',m'\rangle+(\mathscr{A}_{l,l'+1}^{0,-m'-q}(\alpha,q)\langle l'+1,-m'-q;1,q\mid l',-m'\rangle]
\end{aligned}$$

When $m=0$, $m'<0$,

$$
\begin{aligned}
&M_{(l'm',1),(l,0,2)}(\alpha)\\
=&\sum\limits_{\lambda\in\{l'+1,l'+3\}}\sum\limits_{q=-1}^{1}\sum\limits_{m_1=-1}^{1}\rho^{l+1} (a_{22}-a_{12})  r^\lambda \sqrt{2(l'+1)(2l'+1)}
[i\cdot (-1)^{m'}\mathscr{A}_{l,\lambda}^{0,-m'-q}(\alpha,q)K_{l'+1,l',\lambda}^{m_1,-m'-q,q}\\
&
-i\cdot \mathscr{A}_{l,\lambda}^{0,m'-q}(\alpha,q)K_{l'+1,l',\lambda}^{m_1,m'-q,q}]\\
&-\rho^{l+1} a_{22}r^{l'+1}\dfrac{(2l+1)(l'+1)}{\sqrt{2}}[i\cdot (-1)^{m'}\mathscr{H}_{l,l'}^{0,-m'}(\alpha)-i\cdot \mathscr{H}_{l,l'}^{0,m'}(\alpha)]\\
&+\rho^{l+1} a_{22}r^{l'+1}(2l+1)\dfrac{\sqrt{(l'+1)(2l'+1)}}{\sqrt{2}} \sum\limits_{q=-1}^1(-1)^q [i\cdot (-1)^{m'} \mathscr{A}_{l,l'+1}^{0,-m'-q}(\alpha,q)\\
&\langle l'+1,-m'-q;1,q\mid l',-m'\rangle-i\cdot \mathscr{A}_{l,l'+1}^{0,m'-q}(\alpha,q)\langle l'+1,m'-q;1,q\mid l',m'\rangle]
\end{aligned}$$

When $m=0$, $m'=0$,
$$
\begin{aligned}
&M_{(l',0,1),(l,0,2)}(\alpha)\\
=&\sum\limits_{\lambda\in\{l'+1,l'+3\}}\sum\limits_{q=-1}^{1}\sum\limits_{m_1=-1}^{1}\rho^{l+1} (a_{22}-a_{12})  r^\lambda \sqrt{4(l'+1)(2l'+1)}
[\mathscr{A}_{l,\lambda}^{0,-q}(\alpha,q)K_{l'+1,l',\lambda}^{m_1,-q,q}]\\
&-\rho^{l+1} a_{22}r^{l'+1}(2l+1)(l'+1)[\mathscr{H}_{l,l'}^{0,0}(\alpha)]\\
&+\rho^{l+1} a_{22}r^{l'+1}(2l+1)\sqrt{(l'+1)(2l'+1)} \sum\limits_{q=-1}^1(-1)^q [ \mathscr{A}_{l,l'+1}^{0,-q}(\alpha,q)\langle l'+1,-q;1,q\mid l',0\rangle]
\end{aligned}$$

\end{theorem}

\begin{theorem}
Entries for {$\boxed{M_{(l'm',3),(lm,2)}(\alpha)}$}($l'\geqslant 1,l\geqslant 1$, otherwise this term vanishes):\\
When $m>0$, $m'>0$,
$$
\begin{aligned}
&M_{(l'm',3),(lm,2)}(\alpha)\\
=&\sum\limits_{\lambda\in\{l',l'+2\mid l'\geqslant 1\}}\sum\limits_{q=-1}^{1}\sum\limits_{m_1=-1}^{1}-i\cdot \rho^{l+1} (a_{22}-a_{12})  r^\lambda \sqrt{l'(l'+1)}
[(-1)^{m'}(\mathscr{A}_{l,\lambda}^{-m,m'-q}(\alpha,q)\\
&+(-1)^m\mathscr{A}_{l,\lambda}^{m,m'-q}(\alpha,q))K_{l',l',\lambda}^{m_1,m'-q,q}
+(\mathscr{A}_{l,\lambda}^{-m,-m'-q}(\alpha,q)+(-1)^m\mathscr{A}_{l,\lambda}^{m,-m'-q}(\alpha,q))K_{l',l',\lambda}^{m_1,-m'-q,q}]\\
&- i\cdot \rho^{l+1} a_{22}r^{l'}(2l+1)\dfrac{\sqrt{l'(l'+1)}}{2} \sum\limits_{q=-1}^1(-1)^q [(-1)^{m'}(\mathscr{A}_{l,l'}^{-m,m'-q}(\alpha,q)+(-1)^m\mathscr{A}_{l,l'}^{m,m'-q}(\alpha,q))\\
&\langle l',m'-q;1,q\mid l',m'\rangle+(\mathscr{A}_{l,l'}^{-m,-m'-q}(\alpha,q)+(-1)^m\mathscr{A}_{l,l'}^{m,-m'-q}(\alpha,q))\langle l',-m'-q;1,q\mid l',-m'\rangle]
\end{aligned}$$

When $m>0$, $m'<0$,
$$
\begin{aligned}
&M_{(l'm',3),(lm,2)}(\alpha)\\
=&\sum\limits_{\lambda\in\{l',l'+2\mid l'\geqslant 1\}}\sum\limits_{q=-1}^{1}\sum\limits_{m_1=-1}^{1}\rho^{l+1} (a_{22}-a_{12})  r^\lambda \sqrt{l'(l'+1)}
[ (-1)^{m'}(\mathscr{A}_{l,\lambda}^{-m,-m'-q}(\alpha,q)\\
&+(-1)^m\mathscr{A}_{l,\lambda}^{m,-m'-q}(\alpha,q))K_{l',l',\lambda}^{m_1,-m'-q,q}
- (\mathscr{A}_{l,\lambda}^{-m,m'-q}(\alpha,q)+(-1)^m\mathscr{A}_{l,\lambda}^{m,m'-q}(\alpha,q))K_{l',l',\lambda}^{m_1,m'-q,q}]\\
&+\rho^{l+1} a_{22}r^{l'}(2l+1)\dfrac{\sqrt{l'(l'+1)}}{2} \sum\limits_{q=-1}^1(-1)^q [ (-1)^{m'}(\mathscr{A}_{l,l'}^{-m,-m'-q}(\alpha,q)+(-1)^m\mathscr{A}_{l,l'}^{m,-m'-q}(\alpha,q))\\
&\langle l',-m'-q;1,q\mid l',-m'\rangle
- (\mathscr{A}_{l,l'}^{-m,m'-q}(\alpha,q)+(-1)^m\mathscr{A}_{l,l'}^{m,m'-q}(\alpha,q))\langle l',m'-q;1,q\mid l',m'\rangle]\end{aligned}$$

When $m>0$, $m'=0$,
$$
\begin{aligned}
&M_{(l',0,3),(lm,2)}(\alpha)\\
=&\sum\limits_{\lambda\in\{l',l'+2\mid l'\geqslant 1\}}\sum\limits_{q=-1}^{1}\sum\limits_{m_1=-1}^{1} -i\cdot \rho^{l+1} (a_{22}-a_{12})  r^\lambda \sqrt{2l'(l'+1)}
[(\mathscr{A}_{l,\lambda}^{-m,-q}(\alpha,q)+\\&(-1)^m\mathscr{A}_{l,\lambda}^{m,-q}(\alpha,q))K_{l',l',\lambda}^{m_1,-q,q}]\\
&-i\cdot \rho^{l+1} a_{22}r^{l'}(2l+1)\dfrac{\sqrt{l'(l'+1)}}{\sqrt{2}} \sum\limits_{q=-1}^1(-1)^q [ (\mathscr{A}_{l,l'}^{-m,-q}(\alpha,q)+(-1)^m\mathscr{A}_{l,l'}^{m,-q}(\alpha,q))\\
&\langle l',-q;1,q\mid l',0\rangle]
\end{aligned}$$

When $m<0$, $m'>0$,
$$
\begin{aligned}
&M_{(l'm',3),(lm,2)}(\alpha)\\
=&\sum\limits_{\lambda\in\{l',l'+2\mid l'\geqslant 1\}}\sum\limits_{q=-1}^{1}\sum\limits_{m_1=-1}^{1}\rho^{l+1} (a_{22}-a_{12})  r^\lambda \sqrt{l'(l'+1)}
[(-1)^{m'}(\mathscr{A}_{l,\lambda}^{m,m'-q}(\alpha,q)\\
&-(-1)^m\mathscr{A}_{l,\lambda}^{-m,m'-q}(\alpha,q))K_{l',l',\lambda}^{m_1,m'-q,q}
+(\mathscr{A}_{l,\lambda}^{m,-m'-q}(\alpha,q)-(-1)^m\mathscr{A}_{l,\lambda}^{-m,-m'-q}(\alpha,q))K_{l',l',\lambda}^{m_1,-m'-q,q}]\\
&+ \rho^{l+1} a_{22}r^{l'} (2l+1)\dfrac{\sqrt{l'(l'+1)}}{2} \sum\limits_{q=-1}^1(-1)^q [(-1)^{m'}(\mathscr{A}_{l,l'}^{m,m'-q}(\alpha,q)-(-1)^m\mathscr{A}_{l,l'}^{-m,m'-q}(\alpha,q))\\
&\langle l',m'-q;1,q\mid l',m'\rangle+(\mathscr{A}_{l,l'}^{m,-m'-q}(\alpha,q)-(-1)^m\mathscr{A}_{l,l'}^{-m,-m'-q}(\alpha,q))\langle l',-m'-q;1,q\mid l',-m'\rangle]
\end{aligned}$$

When $m<0$, $m'<0$,
$$
\begin{aligned}
&M_{(l'm',3),(lm,2)}(\alpha)\\
=&\sum\limits_{\lambda\in\{l',l'+2\mid l'\geqslant 1\}}\sum\limits_{q=-1}^{1}\sum\limits_{m_1=-1}^{1}\rho^{l+1} (a_{22}-a_{12})  r^\lambda i\cdot \sqrt{l'(l'+1)}
[(-1)^{m'}(\mathscr{A}_{l,\lambda}^{m,-m'-q}(\alpha,q)\\
&-(-1)^m\mathscr{A}_{l,\lambda}^{-m,-m'-q}(\alpha,q))K_{l',l',\lambda}^{m_1,-m'-q,q}
- (\mathscr{A}_{l,\lambda}^{m,m'-q}(\alpha,q)-(-1)^m\mathscr{A}_{l,\lambda}^{-m,m'-q}(\alpha,q))K_{l',l',\lambda}^{m_1,m'-q,q}]\\
&-i\cdot \rho^{l+1} a_{22}r^{l'}(2l+1)\dfrac{\sqrt{l'(l'+1)}}{2} \sum\limits_{q=-1}^1(-1)^q [-(-1)^{m'}(\mathscr{A}_{l,l'}^{m,-m'-q}(\alpha,q)-(-1)^m\mathscr{A}_{l,l'}^{-m,-m'-q}(\alpha,q))\\
&\langle l',-m'-q;1,q\mid l',-m'\rangle
+ (\mathscr{A}_{l,l'}^{m,m'-q}(\alpha,q)-(-1)^m\mathscr{A}_{l,l'}^{-m,m'-q}(\alpha,q))\langle l',m'-q;1,q\mid l',m'\rangle]
\end{aligned}$$

When $m<0$, $m'=0$,
$$
\begin{aligned}
&M_{(l',0,3),(lm,2)}(\alpha)\\
=&\sum\limits_{\lambda\in\{l',l'+2\mid l'\geqslant 1\}}\sum\limits_{q=-1}^{1}\sum\limits_{m_1=-1}^{1}  \rho^{l+1} (a_{22}-a_{12})  r^\lambda  \sqrt{2l'(l'+1)}
[(\mathscr{A}_{l,\lambda}^{m,-q}(\alpha,q)-(-1)^m\mathscr{A}_{l,\lambda}^{-m,-q}(\alpha,q))K_{l',l',\lambda}^{m_1,-q,q}]\\
&+\rho^{l+1} a_{22}r^{l'} (2l+1)\dfrac{\sqrt{l'(l'+1)}}{\sqrt{2}} \sum\limits_{q=-1}^1(-1)^q [ (\mathscr{A}_{l,l'}^{m,-q}(\alpha,q)-(-1)^m\mathscr{A}_{l,l'}^{-m,-q}(\alpha,q))\\
&\langle l',-q;1,q\mid l',0\rangle]
\end{aligned}$$

When $m=0$, $m'>0$,
$$
\begin{aligned}
&M_{(l'm',3),(l,0,2)}(\alpha)\\
=&\sum\limits_{\lambda\in\{l',l'+2\mid l'\geqslant 1\}}\sum\limits_{q=-1}^{1}\sum\limits_{m_1=-1}^{1} -i\cdot \rho^{l+1} (a_{22}-a_{12})  r^\lambda \sqrt{2l'(l'+1)}
[(-1)^{m'}\mathscr{A}_{l,\lambda}^{0,m'-q}(\alpha,q) K_{l',l',\lambda}^{m_1,m'-q,q}\\
&+\mathscr{A}_{l,\lambda}^{0,-m'-q}(\alpha,q)K_{l',l',\lambda}^{m_1,-m'-q,q}]\\
&-i\cdot \rho^{l+1} a_{22}r^{l'}(2l+1)\dfrac{\sqrt{l'(l'+1)}}{\sqrt{2}} \sum\limits_{q=-1}^1(-1)^q [(-1)^{m'}\mathscr{A}_{l,l'}^{0,m'-q}(\alpha,q)\\
&\langle l',m'-q;1,q\mid l',m'\rangle+(\mathscr{A}_{l,l'}^{0,-m'-q}(\alpha,q)\langle l',-m'-q;1,q\mid l',-m'\rangle]
\end{aligned}$$

When $m=0$, $m'<0$,

$$
\begin{aligned}
&M_{(l'm',3),(l,0,2)}(\alpha)\\
=&\sum\limits_{\lambda\in\{l',l'+2\mid l'\geqslant 1\}}\sum\limits_{q=-1}^{1}\sum\limits_{m_1=-1}^{1}\rho^{l+1} (a_{22}-a_{12})  r^\lambda \sqrt{2l'(l'+1)}
[(-1)^{m'}\mathscr{A}_{l,\lambda}^{0,-m'-q}(\alpha,q)K_{l',l',\lambda}^{m_1,-m'-q,q}\\
&
- \mathscr{A}_{l,\lambda}^{0,m'-q}(\alpha,q)K_{l',l',\lambda}^{m_1,m'-q,q}]\\
&+\rho^{l+1} a_{22}r^{l'}(2l+1)\dfrac{\sqrt{l'(l'+1)}}{\sqrt{2}} \sum\limits_{q=-1}^1(-1)^q [ (-1)^{m'} \mathscr{A}_{l,l'}^{0,-m'-q}(\alpha,q)\\
&\langle l',-m'-q;1,q\mid l',-m'\rangle- \mathscr{A}_{l,l'}^{0,m'-q}(\alpha,q)\langle l',m'-q;1,q\mid l',m'\rangle]
\end{aligned}$$

When $m=0$, $m'=0$,
$$
\begin{aligned}
&M_{(l',0,3),(l,0,2)}(\alpha)\\
=&\sum\limits_{\lambda\in\{l',l'+2\mid l'\geqslant 1\}}\sum\limits_{q=-1}^{1}\sum\limits_{m_1=-1}^{1} -i\cdot \rho^{l+1} (a_{22}-a_{12})  r^\lambda \sqrt{4l'(l'+1)}
[\mathscr{A}_{l,\lambda}^{0,-q}(\alpha,q)K_{l',l',\lambda}^{m_1,-q,q}]\\
&-i\cdot \rho^{l+1} a_{22}r^{l'}(2l+1)\sqrt{l'(l'+1)} \sum\limits_{q=-1}^1(-1)^q [ \mathscr{A}_{l,l'}^{0,-q}(\alpha,q)\langle l',-q;1,q\mid l',0\rangle]
\end{aligned}$$

\end{theorem}

\begin{theorem}
Entries for {$\boxed{M_{(l'm',3),(lm,3)}(\alpha)}$}($l'\geqslant 1,l\geqslant 1$, otherwise this term vanishes):\\

When $m>0,\ m'>0$:
$$
\begin{aligned}
&M_{(l'm',3),(lm,3)}(\alpha)\\
=&\ \rho^{l+2}a_{33}r^{l'}\Biggl\{
\frac{l'(l'+1)}{2}\Bigl[(-1)^{m'}\bigl(\mathscr{H}_{l,l'}^{-m,m'}(\alpha)+(-1)^m \mathscr{H}_{l,l'}^{m,m'}(\alpha)\bigr)
+\bigl(\mathscr{H}_{l,l'}^{-m,-m'}(\alpha)+(-1)^m \mathscr{H}_{l,l'}^{m,-m'}(\alpha)\bigr)\Bigr]\\
&-\frac{i\sqrt{l'(l'+1)}}{2}\sum_{q=-1}^{1}\sum_{m_1=-1}^{1}\Bigl[(-1)^{m'}\bigl(\mathscr{L}_{l,l',l'+1}^{-m,m'-q,q,m_1}(\alpha)+(-1)^m \mathscr{L}_{l,l',l'+1}^{m,m'-q,q,m_1}(\alpha)\bigr)\\
&\qquad\qquad\qquad\qquad\quad+\bigl(\mathscr{L}_{l,l',l'+1}^{-m,-m'-q,q,m_1}(\alpha)+(-1)^m \mathscr{L}_{l,l',l'+1}^{m,-m'-q,q,m_1}(\alpha)\bigr)\Bigr]\Biggr\}\\
&+\delta_{ll'}\delta_{mm'}\rho\dfrac{l(l+1)}{(2l+1)\mu} \end{aligned}
$$

When $m>0,\ m'<0$:
$$
\begin{aligned}
&M_{(l'm',3),(lm,3)}(\alpha)\\
=&\ \rho^{l+2}a_{33}r^{l'}\Biggl\{
\frac{l'(l'+1)}{2}\Bigl[i(-1)^{m'}\bigl(\mathscr{H}_{l,l'}^{-m,-m'}(\alpha)+(-1)^m \mathscr{H}_{l,l'}^{m,-m'}(\alpha)\bigr)
-i\bigl(\mathscr{H}_{l,l'}^{-m,m'}(\alpha)\\
&+(-1)^m \mathscr{H}_{l,l'}^{m,m'}(\alpha)\bigr)\Bigr]\\
&+\frac{\sqrt{l'(l'+1)}}{2}\sum_{q=-1}^{1}\sum_{m_1=-1}^{1}\Bigl[(-1)^{m'}\bigl(\mathscr{L}_{l,l',l'+1}^{-m,-m'-q,q,m_1}(\alpha)+(-1)^m \mathscr{L}_{l,l',l'+1}^{m,-m'-q,q,m_1}(\alpha)\bigr)\\
&\qquad\qquad\qquad\qquad\quad-\bigl(\mathscr{L}_{l,l',l'+1}^{-m,m'-q,q,m_1}(\alpha)+(-1)^m \mathscr{L}_{l,l',l'+1}^{m,m'-q,q,m_1}(\alpha)\bigr)\Bigr]\Biggr\}\\
&+\delta_{ll'}\delta_{mm'}\rho\dfrac{l(l+1)}{(2l+1)\mu}
\end{aligned}
$$

When $m>0,\ m'=0$:
$$
\begin{aligned}
&M_{(l',0,3),(lm,3)}(\alpha)\\
=&\ \rho^{l+2}a_{33}r^{l'}\Biggl\{
\frac{l'(l'+1)}{\sqrt{2}}\bigl(\mathscr{H}_{l,l'}^{-m,0}(\alpha)+(-1)^m \mathscr{H}_{l,l'}^{m,0}(\alpha)\bigr)\\
&-\frac{i\sqrt{l'(l'+1)}}{\sqrt{2}}\sum_{q=-1}^{1}\sum_{m_1=-1}^{1}\bigl(\mathscr{L}_{l,l',l'+1}^{-m,-q,q,m_1}(\alpha)+(-1)^m \mathscr{L}_{l,l',l'+1}^{m,-q,q,m_1}(\alpha)\bigr)\Biggr\}\\
&+\delta_{ll'}\delta_{mm'}\rho\dfrac{l(l+1)}{(2l+1)\mu}
\end{aligned}
$$

When $m<0,\ m'>0$:
$$
\begin{aligned}
&M_{(l'm',3),(lm,3)}(\alpha)\\
=&\ \rho^{l+2}a_{33}r^{l'}\Biggl\{
\frac{l'(l'+1)}{2}\Bigl[i(-1)^{m'}\bigl(\mathscr{H}_{l,l'}^{m,m'}(\alpha)-(-1)^m \mathscr{H}_{l,l'}^{-m,m'}(\alpha)\bigr)\\
&
+i\bigl(\mathscr{H}_{l,l'}^{m,-m'}(\alpha)-(-1)^m \mathscr{H}_{l,l'}^{-m,-m'}(\alpha)\bigr)\Bigr]\\
&+\frac{\sqrt{l'(l'+1)}}{2}\sum_{q=-1}^{1}\sum_{m_1=-1}^{1}\Bigl[(-1)^{m'}\bigl(\mathscr{L}_{l,l',l'+1}^{m,m'-q,q,m_1}(\alpha)-(-1)^m \mathscr{L}_{l,l',l'+1}^{-m,m'-q,q,m_1}(\alpha)\bigr)\\
&\qquad\qquad\qquad\qquad\quad+\bigl(\mathscr{L}_{l,l',l'+1}^{m,-m'-q,q,m_1}(\alpha)-(-1)^m \mathscr{L}_{l,l',l'+1}^{-m,-m'-q,q,m_1}(\alpha)\bigr)\Bigr]\Biggr\}\\
&+\delta_{ll'}\delta_{mm'}\rho\dfrac{l(l+1)}{(2l+1)\mu}\end{aligned}
$$

When $m<0,\ m'<0$:
$$
\begin{aligned}
&M_{(l'm',3),(lm,3)}(\alpha)\\
=&\ \rho^{l+2}a_{33}r^{l'}\Biggl\{
\frac{l'(l'+1)}{2}\Bigl[\bigl(\mathscr{H}_{l,l'}^{m,m'}(\alpha)-(-1)^m \mathscr{H}_{l,l'}^{-m,m'}(\alpha)\bigr)
-(-1)^{m'}\bigl(\mathscr{H}_{l,l'}^{m,-m'}(\alpha)-(-1)^m \mathscr{H}_{l,l'}^{-m,-m'}(\alpha)\bigr)\Bigr]\\
&+\frac{i\sqrt{l'(l'+1)}}{2}\sum_{q=-1}^{1}\sum_{m_1=-1}^{1}\Bigl[(-1)^{m'}\bigl(\mathscr{L}_{l,l',l'+1}^{m,-m'-q,q,m_1}(\alpha)-(-1)^m \mathscr{L}_{l,l',l'+1}^{-m,-m'-q,q,m_1}(\alpha)\bigr)\\
&\qquad\qquad\qquad\qquad\quad-\bigl(\mathscr{L}_{l,l',l'+1}^{m,m'-q,q,m_1}(\alpha)-(-1)^m \mathscr{L}_{l,l',l'+1}^{-m,m'-q,q,m_1}(\alpha)\bigr)\Bigr]\Biggr\}\\
&+\delta_{ll'}\delta_{mm'}\rho\dfrac{l(l+1)}{(2l+1)\mu}\end{aligned}
$$

When $m<0,\ m'=0$:
$$
\begin{aligned}
&M_{(l',0,3),(lm,3)}(\alpha)\\
=&\ \rho^{l+2}a_{33}r^{l'}\Biggl\{
\frac{il'(l'+1)}{\sqrt{2}}\bigl(\mathscr{H}_{l,l'}^{m,0}(\alpha)-(-1)^m \mathscr{H}_{l,l'}^{-m,0}(\alpha)\bigr)\\
&+\frac{\sqrt{l'(l'+1)}}{\sqrt{2}}\sum_{q=-1}^{1}\sum_{m_1=-1}^{1}\bigl(\mathscr{L}_{l,l',l'+1}^{m,-q,q,m_1}(\alpha)-(-1)^m \mathscr{L}_{l,l',l'+1}^{-m,-q,q,m_1}(\alpha)\bigr)\Biggr\}\\
&+\delta_{ll'}\delta_{mm'}\rho\dfrac{l(l+1)}{(2l+1)\mu}
\end{aligned}
$$

When $m=0,\ m'>0$:
$$
\begin{aligned}
&M_{(l'm',3),(l,0,3)}(\alpha)\\
=&\ \rho^{l+2}a_{33}r^{l'}\Biggl\{
\frac{l'(l'+1)}{\sqrt{2}}\Bigl[(-1)^{m'}\mathscr{H}_{l,l'}^{0,m'}(\alpha)+\mathscr{H}_{l,l'}^{0,-m'}(\alpha)\Bigr]\\
&-\frac{i\sqrt{l'(l'+1)}}{\sqrt{2}}\sum_{q=-1}^{1}\sum_{m_1=-1}^{1}\Bigl[(-1)^{m'}\mathscr{L}_{l,l',l'+1}^{0,m'-q,q,m_1}(\alpha)+\mathscr{L}_{l,l',l'+1}^{0,-m'-q,q,m_1}(\alpha)\Bigr]\Biggr\}\\
&+\delta_{ll'}\delta_{mm'}\rho\dfrac{l(l+1)}{(2l+1)\mu}\end{aligned}
$$

When $m=0,\ m'<0$:
$$
\begin{aligned}
&M_{(l'm',3),(l,0,3)}(\alpha)\\
=&\ \rho^{l+2}a_{33}r^{l'}\Biggl\{
\frac{l'(l'+1)}{\sqrt{2}}\Bigl[i(-1)^{m'}\mathscr{H}_{l,l'}^{0,-m'}(\alpha)-i\mathscr{H}_{l,l'}^{0,m'}(\alpha)\Bigr]\\
&+\frac{\sqrt{l'(l'+1)}}{\sqrt{2}}\sum_{q=-1}^{1}\sum_{m_1=-1}^{1}\Bigl[(-1)^{m'}\mathscr{L}_{l,l',l'+1}^{0,-m'-q,q,m_1}(\alpha)-\mathscr{L}_{l,l',l'+1}^{0,m'-q,q,m_1}(\alpha)\Bigr]\Biggr\}\\
&+\delta_{ll'}\delta_{mm'}\rho\dfrac{l(l+1)}{(2l+1)\mu}\end{aligned}
$$

When $m=0,\ m'=0$:
$$
\begin{aligned}
&M_{(l',0,3),(l,0,3)}(\alpha)\\
=&\ \rho^{l+2}a_{33}r^{l'}\Biggl\{
l'(l'+1)\mathscr{H}_{l,l'}^{0,0}(\alpha)
-i\sqrt{l'(l'+1)}\sum_{q=-1}^{1}\sum_{m_1=-1}^{1}\mathscr{L}_{l,l',l'+1}^{0,-q,q,m_1}(\alpha)\Biggr\}\\
&+\delta_{ll'}\delta_{mm'}\rho\dfrac{l(l+1)}{(2l+1)\mu}\end{aligned}
$$
\end{theorem}

\begin{lemma}Define $\mathscr{D}_{l,\lambda}^{m,\mu}(\alpha)=\sum\limits_{n\neq 0} n^2H_{l,\lambda}^{m,\mu}(\mathbf{b})e^{-in\alpha}$. Then
$$\begin{aligned}
&\mathscr{D}_{l,\lambda}^{m,\mu}(\alpha)\\
=&(-1)^{\lambda+\mu}
\sqrt{\dfrac{2l+1}{2\lambda+1}}
\sqrt{{l+\lambda+\mu-m\choose \lambda+\mu}{l+\lambda+m-\mu\choose \lambda-\mu}} \sqrt{\dfrac{4\pi}{2(l+\lambda)+1}}\\
&
\big(\mathrm{Li}_{l+\lambda-1}(e^{-i\alpha})Y_{l+\lambda}^{m-\mu}(\dfrac{\pi}{2},\pi)+\mathrm{Li}_{l+\lambda-1}(e^{i\alpha})Y_{l+\lambda}^{m-\mu}(\dfrac{\pi}{2},0)\big)
\end{aligned}$$

\end{lemma}

\begin{theorem}The entries for  $\boxed{M_{(l'm',2),(lm,2)}(\alpha)}$($l'\geqslant 1,l\geqslant 1$,otherwise this term vanishes) are given by:
\\
If $m>0$ and $m'>0$,
$$\begin{aligned}&M_{(l'm',2),(lm,2)}(\alpha)\\
&=\dfrac{l'(2l'+1)}{2}[ \rho^{l+3}a_{12}r^{l'-1}+\rho^{l+1}(a_{22}-a_{12})(r^{l'+1} )+\rho^{l+1}a_{22}r^{l'+1}\dfrac{2l+1}{2l'+1}]\\
&[(-1)^{m'}(\mathscr{H}_{l,l'}^{-m,m'}(\alpha)+(-1)^m\mathscr{H}_{l,l'}^{m,m'}(\alpha))+(\mathscr{H}_{l,l'}^{-m,-m'}(\alpha)+(-1)^m\mathscr{H}_{l,l'}^{m,-m'}(\alpha))]\\
&+\dfrac{l'(2l'+1)}{2}[\rho^{l+1}(a_{22}-a_{12})r^{l'-1}][(-1)^{m'}(\mathscr{D}_{l,l'}^{-m,m'}(\alpha)+\\
&(-1)^m\mathscr{D}_{l,l'}^{m,m'}(\alpha))+(\mathscr{D}_{l,l'}^{-m,-m'}(\alpha)+(-1)^m\mathscr{D}_{l,l'}^{m,-m'}(\alpha))]\\
&+\rho^{l+1}(a_{22}-a_{12})\sqrt{l'(2l'+1)}\sum\limits_{\lambda\in\{l'+1,l'-1\}}\sum\limits_{q=-1}^{q=1}\sum\limits_{m_1=-1}^{1}r^\lambda [(-1)^{m'}(\mathscr{A}_{l,\lambda}^{-m,m'-q}(\alpha,q)\\
&+(-1)^m\mathscr{A}_{l,\lambda}^{m,m'-q}(\alpha,q))K_{l'-1,l',\lambda}^{m_1,m'-q,q}+(\mathscr{A}_{l,\lambda}^{-m,-m'-q}(\alpha,q)+(-1)^m\mathscr{A}_{l,\lambda}^{m,-m'-q}(\alpha,q))K_{l'-1,l',\lambda}^{m_1,-m'-q,q}]\\
&+\rho^{l+1}a_{22}(2l+1)\dfrac{\sqrt{l'(2l'+1)}}{2}\sum\limits_{q=-1}^1(-1)^q r^{l'-1}[(-1)^{m'}(\mathscr{A}_{l,l'-1}^{-m,m'-q}(\alpha,q)+(-1)^m\mathscr{A}_{l,l'-1}^{m,m'-q}(\alpha,q))\\
&\langle l'-1,m'-q;1,q\mid l',m' \rangle+(\mathscr{A}_{l,l'-1}^{-m,-m'-q}(\alpha,q)+(-1)^m\mathscr{A}_{l,l'-1}^{m,-m'-q}(\alpha,q))\langle l'-1,-m'-q;1,q\mid l',-m' \rangle]\\
&+\delta_{ll'}\delta_{mm'}l(2l+1)\rho a_{22}.
\end{aligned}$$

If $m>0$ and $m'<0$,
$$
\begin{aligned}
&M_{(l'm',2),(lm,2)}(\alpha)\\
&=\dfrac{l'(2l'+1)}{2}\Bigl[ \rho^{l+3}a_{12}r^{l'-1}+\rho^{l+1}(a_{22}-a_{12})(r^{l'+1} )+\rho^{l+1}a_{22}r^{l'+1}\dfrac{2l+1}{2l'+1}\Bigr]\\
&\quad\Bigl[i\cdot(-1)^{m'}(\mathscr{H}_{l,l'}^{-m,-m'}(\alpha)+(-1)^m\mathscr{H}_{l,l'}^{m,-m'}(\alpha))-i\cdot(\mathscr{H}_{l,l'}^{-m,m'}(\alpha)+(-1)^m\mathscr{H}_{l,l'}^{m,m'}(\alpha))\Bigr]\\
&+\dfrac{l'(2l'+1)}{2}[\rho^{l+1}(a_{22}-a_{12})r^{l'-1}]\Bigl[i\cdot(-1)^{m'}(\mathscr{D}_{l,l'}^{-m,-m'}(\alpha)+(-1)^m\mathscr{D}_{l,l'}^{m,-m'}(\alpha))\\
&-i\cdot(\mathscr{D}_{l,l'}^{-m,m'}(\alpha)+(-1)^m\mathscr{D}_{l,l'}^{m,m'}(\alpha))\Bigr]\\
&+\rho^{l+1}(a_{22}-a_{12})\sqrt{l'(2l'+1)}\sum_{\lambda\in\{l'+1,l'-1\}}\sum_{q=-1}^{1}\sum_{m_1=-1}^{1}r^\lambda \Bigl[i\cdot(-1)^{m'}(\mathscr{A}_{l,\lambda}^{-m,-m'-q}(\alpha,q)\\
&\quad+(-1)^m\mathscr{A}_{l,\lambda}^{m,-m'-q}(\alpha,q))K_{l'-1,l',\lambda}^{m_1,-m'-q,q}-i\cdot(\mathscr{A}_{l,\lambda}^{-m,m'-q}(\alpha,q)+(-1)^m\mathscr{A}_{l,\lambda}^{m,m'-q}(\alpha,q))K_{l'-1,l',\lambda}^{m_1,m'-q,q}\Bigr]\\
&+\rho^{l+1}a_{22}(2l+1)\dfrac{\sqrt{l'(2l'+1)}}{2}\sum_{q=-1}^{1}(-1)^q r^{l'-1}\Bigl[i\cdot(-1)^{m'}(\mathscr{A}_{l,l'-1}^{-m,-m'-q}(\alpha,q)+(-1)^m\mathscr{A}_{l,l'-1}^{m,-m'-q}(\alpha,q))\\
&\quad\langle l'-1,-m'-q;1,q\mid l',-m' \rangle-i\cdot(\mathscr{A}_{l,l'-1}^{-m,m'-q}(\alpha,q)+\\
&(-1)^m\mathscr{A}_{l,l'-1}^{m,m'-q}(\alpha,q))\langle l'-1,m'-q;1,q\mid l',m' \rangle\Bigr]\\
&+\delta_{ll'}\delta_{mm'}l(2l+1)\rho a_{22}.
\end{aligned}
$$

If $m>0$ and $m'=0$,
$$
\begin{aligned}
&M_{(l',0,2),(lm,2)}(\alpha)\\
&=\dfrac{l'(2l'+1)}{\sqrt{2}}\Bigl[ \rho^{l+3}a_{12}r^{l'-1}+\rho^{l+1}(a_{22}-a_{12})(r^{l'+1} )+\rho^{l+1}a_{22}r^{l'+1}\dfrac{2l+1}{2l'+1}\Bigr]\\
&\quad\Bigl[(\mathscr{H}_{l,l'}^{-m,0}(\alpha)+(-1)^m\mathscr{H}_{l,l'}^{m,0}(\alpha))\Bigr]\\
&+\dfrac{l'(2l'+1)}{\sqrt{2}}[\rho^{l+1}(a_{22}-a_{12})r^{l'-1}]\Bigl[(\mathscr{D}_{l,l'}^{-m,0}(\alpha)+(-1)^m\mathscr{D}_{l,l'}^{m,0}(\alpha))\Bigr]\\
&+\rho^{l+1}(a_{22}-a_{12})\sqrt{2l'(2l'+1)}\sum_{\lambda\in\{l'+1,l'-1\}}\sum_{q=-1}^{1}\sum_{m_1=-1}^{1}r^\lambda \Bigl[(\mathscr{A}_{l,\lambda}^{-m,-q}(\alpha,q)+(-1)^m\mathscr{A}_{l,\lambda}^{m,-q}(\alpha,q))K_{l'-1,l',\lambda}^{m_1,-q,q}\Bigr]\\
&+\rho^{l+1}a_{22}(2l+1)\dfrac{\sqrt{l'(2l'+1)}}{\sqrt{2}}\sum_{q=-1}^{1}(-1)^q r^{l'-1}\Bigl[(\mathscr{A}_{l,l'-1}^{-m,-q}(\alpha,q)+\\
&(-1)^m\mathscr{A}_{l,l'-1}^{m,-q}(\alpha,q))\langle l'-1,-q;1,q\mid l',0 \rangle\Bigr]\\
&+\delta_{ll'}\delta_{mm'}l(2l+1)\rho a_{22}.
\end{aligned}
$$

If $m<0$ and $m'>0$,
$$
\begin{aligned}
&M_{(l'm',2),(lm,2)}(\alpha)\\
&=\dfrac{l'(2l'+1)}{2}\Bigl[ \rho^{l+3}a_{12}r^{l'-1}+\rho^{l+1}(a_{22}-a_{12})(r^{l'+1} )+\rho^{l+1}a_{22}r^{l'+1}\dfrac{2l+1}{2l'+1}\Bigr]\\
&\quad\cdot i\Bigl[(-1)^{m'}(\mathscr{H}_{l,l'}^{m,m'}(\alpha)-(-1)^m\mathscr{H}_{l,l'}^{-m,m'}(\alpha))+(\mathscr{H}_{l,l'}^{m,-m'}(\alpha)-(-1)^m\mathscr{H}_{l,l'}^{-m,-m'}(\alpha))\Bigr]\\
&+\dfrac{l'(2l'+1)}{2}[\rho^{l+1}(a_{22}-a_{12})r^{l'-1}]\cdot i\Bigl[(-1)^{m'}(\mathscr{D}_{l,l'}^{m,m'}(\alpha)-(-1)^m\mathscr{D}_{l,l'}^{-m,m'}(\alpha))\\
&+(\mathscr{D}_{l,l'}^{m,-m'}(\alpha)-(-1)^m\mathscr{D}_{l,l'}^{-m,-m'}(\alpha))\Bigr]\\
&+\rho^{l+1}(a_{22}-a_{12})\sqrt{l'(2l'+1)}\sum_{\lambda\in\{l'+1,l'-1\}}\sum_{q=-1}^{1}\sum_{m_1=-1}^{1}r^\lambda \;i\Bigl[(-1)^{m'}(\mathscr{A}_{l,\lambda}^{m,m'-q}(\alpha,q)\\
&\quad-(-1)^m\mathscr{A}_{l,\lambda}^{-m,m'-q}(\alpha,q))K_{l'-1,l',\lambda}^{m_1,m'-q,q}+(\mathscr{A}_{l,\lambda}^{m,-m'-q}(\alpha,q)-(-1)^m\mathscr{A}_{l,\lambda}^{-m,-m'-q}(\alpha,q))K_{l'-1,l',\lambda}^{m_1,-m'-q,q}\Bigr]\\
&+\rho^{l+1}a_{22}(2l+1)\dfrac{\sqrt{l'(2l'+1)}}{2}\sum_{q=-1}^{1}(-1)^q r^{l'-1}\;i\Bigl[(-1)^{m'}(\mathscr{A}_{l,l'-1}^{m,m'-q}(\alpha,q)-(-1)^m\mathscr{A}_{l,l'-1}^{-m,m'-q}(\alpha,q))\\
&\quad\langle l'-1,m'-q;1,q\mid l',m' \rangle+(\mathscr{A}_{l,l'-1}^{m,-m'-q}(\alpha,q)\\
&-(-1)^m\mathscr{A}_{l,l'-1}^{-m,-m'-q}(\alpha,q))\langle l'-1,-m'-q;1,q\mid l',-m' \rangle\Bigr]\\
&+\delta_{ll'}\delta_{mm'}l(2l+1)\rho a_{22}.
\end{aligned}
$$

If $m<0$ and $m'<0$,
$$
\begin{aligned}
&M_{(l'm',2),(lm,2)}(\alpha)\\
&=\dfrac{l'(2l'+1)}{2}\Bigl[ \rho^{l+3}a_{12}r^{l'-1}+\rho^{l+1}(a_{22}-a_{12})(r^{l'+1} )+\rho^{l+1}a_{22}r^{l'+1}\dfrac{2l+1}{2l'+1}\Bigr]\\
&\quad\Bigl[-(-1)^{m'}(\mathscr{H}_{l,l'}^{m,-m'}(\alpha)-(-1)^m\mathscr{H}_{l,l'}^{-m,-m'}(\alpha))+(\mathscr{H}_{l,l'}^{m,m'}(\alpha)-(-1)^m\mathscr{H}_{l,l'}^{-m,m'}(\alpha))\Bigr]\\
&+\dfrac{l'(2l'+1)}{2}[\rho^{l+1}(a_{22}-a_{12})r^{l'-1}]\Bigl[-(-1)^{m'}(\mathscr{D}_{l,l'}^{m,-m'}(\alpha)-(-1)^m\mathscr{D}_{l,l'}^{-m,-m'}(\alpha))\\
&+(\mathscr{D}_{l,l'}^{m,m'}(\alpha)-(-1)^m\mathscr{D}_{l,l'}^{-m,m'}(\alpha))\Bigr]\\
&+\rho^{l+1}(a_{22}-a_{12})\sqrt{l'(2l'+1)}\sum_{\lambda\in\{l'+1,l'-1\}}\sum_{q=-1}^{1}\sum_{m_1=-1}^{1}r^\lambda \Bigl[-(-1)^{m'}(\mathscr{A}_{l,\lambda}^{m,-m'-q}(\alpha,q)\\
&\quad-(-1)^m\mathscr{A}_{l,\lambda}^{-m,-m'-q}(\alpha,q))K_{l'-1,l',\lambda}^{m_1,-m'-q,q}+(\mathscr{A}_{l,\lambda}^{m,m'-q}(\alpha,q)-(-1)^m\mathscr{A}_{l,\lambda}^{-m,m'-q}(\alpha,q))K_{l'-1,l',\lambda}^{m_1,m'-q,q}\Bigr]\\
&+\rho^{l+1}a_{22}(2l+1)\dfrac{\sqrt{l'(2l'+1)}}{2}\sum_{q=-1}^{1}(-1)^q r^{l'-1}\Bigl[-(-1)^{m'}(\mathscr{A}_{l,l'-1}^{m,-m'-q}(\alpha,q)-(-1)^m\mathscr{A}_{l,l'-1}^{-m,-m'-q}(\alpha,q))\\
&\quad\langle l'-1,-m'-q;1,q\mid l',-m' \rangle+(\mathscr{A}_{l,l'-1}^{m,m'-q}(\alpha,q)-(-1)^m\mathscr{A}_{l,l'-1}^{-m,m'-q}(\alpha,q))\langle l'-1,m'-q;1,q\mid l',m' \rangle\Bigr]\\
&+\delta_{ll'}\delta_{mm'}l(2l+1)\rho a_{22}.
\end{aligned}
$$

If $m<0$ and $m'=0$,
$$
\begin{aligned}
&M_{(l',0,2),(lm,2)}(\alpha)\\
&=\dfrac{l'(2l'+1)}{\sqrt{2}}\Bigl[ \rho^{l+3}a_{12}r^{l'-1}+\rho^{l+1}(a_{22}-a_{12})(r^{l'+1} )+\rho^{l+1}a_{22}r^{l'+1}\dfrac{2l+1}{2l'+1}\Bigr]\\
&
\cdot i\Bigl[(\mathscr{H}_{l,l'}^{m,0}(\alpha)-(-1)^m\mathscr{H}_{l,l'}^{-m,0}(\alpha))\Bigr]\\
&+\dfrac{l'(2l'+1)}{\sqrt{2}}[\rho^{l+1}(a_{22}-a_{12})r^{l'-1}]\cdot i\Bigl[(\mathscr{D}_{l,l'}^{m,0}(\alpha)-(-1)^m\mathscr{D}_{l,l'}^{-m,0}(\alpha))\Bigr]\\
&+\rho^{l+1}(a_{22}-a_{12})\sqrt{2l'(2l'+1)}\sum_{\lambda\in\{l'+1,l'-1\}}\sum_{q=-1}^{1}\sum_{m_1=-1}^{1}r^\lambda \;i\Bigl[(\mathscr{A}_{l,\lambda}^{m,-q}(\alpha,q)-(-1)^m\mathscr{A}_{l,\lambda}^{-m,-q}(\alpha,q))K_{l'-1,l',\lambda}^{m_1,-q,q}\Bigr]\\
&+\rho^{l+1}a_{22}(2l+1)\dfrac{\sqrt{l'(2l'+1)}}{\sqrt{2}}\sum_{q=-1}^{1}(-1)^q r^{l'-1}\;i\Bigl[(\mathscr{A}_{l,l'-1}^{m,-q}(\alpha,q)\\
&-(-1)^m\mathscr{A}_{l,l'-1}^{-m,-q}(\alpha,q))\langle l'-1,-q;1,q\mid l',0 \rangle\Bigr]\\
&+\delta_{ll'}\delta_{mm'}l(2l+1)\rho a_{22}.
\end{aligned}
$$

If $m=0$ and $m'>0$,
$$
\begin{aligned}
&M_{(l'm',2),(l,0,2)}(\alpha)\\
&=\dfrac{l'(2l'+1)}{\sqrt{2}}\Bigl[ \rho^{l+3}a_{12}r^{l'-1}+\rho^{l+1}(a_{22}-a_{12})(r^{l'+1} )+\rho^{l+1}a_{22}r^{l'+1}\dfrac{2l+1}{2l'+1}\Bigr]\\
&\quad\Bigl[(-1)^{m'}\mathscr{H}_{l,l'}^{0,m'}(\alpha)+\mathscr{H}_{l,l'}^{0,-m'}(\alpha)\Bigr]\\
&+\dfrac{l'(2l'+1)}{\sqrt{2}}[\rho^{l+1}(a_{22}-a_{12})r^{l'-1}]\Bigl[(-1)^{m'}\mathscr{D}_{l,l'}^{0,m'}(\alpha)+\mathscr{D}_{l,l'}^{0,-m'}(\alpha)\Bigr]\\
&+\rho^{l+1}(a_{22}-a_{12})\sqrt{2l'(2l'+1)}\sum_{\lambda\in\{l'+1,l'-1\}}\sum_{q=-1}^{1}\sum_{m_1=-1}^{1}r^\lambda \Bigl[(-1)^{m'}\mathscr{A}_{l,\lambda}^{0,m'-q}(\alpha,q)K_{l'-1,l',\lambda}^{m_1,m'-q,q}\\
&+\mathscr{A}_{l,\lambda}^{0,-m'-q}(\alpha,q)K_{l'-1,l',\lambda}^{m_1,-m'-q,q}\Bigr]\\
&+\rho^{l+1}a_{22}(2l+1)\dfrac{\sqrt{l'(2l'+1)}}{\sqrt{2}}\sum_{q=-1}^{1}(-1)^q r^{l'-1}\Bigl[(-1)^{m'}\mathscr{A}_{l,l'-1}^{0,m'-q}(\alpha,q)\langle l'-1,m'-q;1,q\mid l',m' \rangle\\
&+\mathscr{A}_{l,l'-1}^{0,-m'-q}(\alpha,q)\langle l'-1,-m'-q;1,q\mid l',-m' \rangle\Bigr]\\
&+\delta_{ll'}\delta_{mm'}l(2l+1)\rho a_{22}.
\end{aligned}
$$

If $m=0$ and $m'<0$,
$$
\begin{aligned}
&M_{(l'm',2),(l,0,2)}(\alpha)\\
&=\dfrac{l'(2l'+1)}{\sqrt{2}}\Bigl[ \rho^{l+3}a_{12}r^{l'-1}+\rho^{l+1}(a_{22}-a_{12})(r^{l'+1} )+\rho^{l+1}a_{22}r^{l'+1}\dfrac{2l+1}{2l'+1}\Bigr]\\
&\quad i\Bigl[(-1)^{m'}\mathscr{H}_{l,l'}^{0,-m'}(\alpha)-\mathscr{H}_{l,l'}^{0,m'}(\alpha)\Bigr]\\
&+\dfrac{l'(2l'+1)}{\sqrt{2}}[\rho^{l+1}(a_{22}-a_{12})r^{l'-1}]\cdot i\Bigl[(-1)^{m'}\mathscr{D}_{l,l'}^{0,-m'}(\alpha)-\mathscr{D}_{l,l'}^{0,m'}(\alpha)\Bigr]\\
&+\rho^{l+1}(a_{22}-a_{12})\sqrt{2l'(2l'+1)}\sum_{\lambda\in\{l'+1,l'-1\}}\sum_{q=-1}^{1}\sum_{m_1=-1}^{1}r^\lambda \;i\Bigl[(-1)^{m'}\mathscr{A}_{l,\lambda}^{0,-m'-q}(\alpha,q)K_{l'-1,l',\lambda}^{m_1,-m'-q,q}\\
&-\mathscr{A}_{l,\lambda}^{0,m'-q}(\alpha,q)K_{l'-1,l',\lambda}^{m_1,m'-q,q}\Bigr]\\
&+\rho^{l+1}a_{22}(2l+1)\dfrac{\sqrt{l'(2l'+1)}}{\sqrt{2}}\sum_{q=-1}^{1}(-1)^q r^{l'-1}\;i\Bigl[(-1)^{m'}\mathscr{A}_{l,l'-1}^{0,-m'-q}(\alpha,q)\langle l'-1,-m'-q;1,q\mid l',-m' \rangle\\
&-\mathscr{A}_{l,l'-1}^{0,m'-q}(\alpha,q)\langle l'-1,m'-q;1,q\mid l',m' \rangle\Bigr]\\
&+\delta_{ll'}\delta_{mm'}l(2l+1)\rho a_{22}.
\end{aligned}
$$

If $m=0$ and $m'=0$,
$$
\begin{aligned}
&M_{(l',0,2),(l,0,2)}(\alpha)\\
&=l'(2l'+1)\Bigl[ \rho^{l+3}a_{12}r^{l'-1}+\rho^{l+1}(a_{22}-a_{12})(r^{l'+1} )+\rho^{l+1}a_{22}r^{l'+1}\dfrac{2l+1}{2l'+1}\Bigr]\mathscr{H}_{l,l'}^{0,0}(\alpha)\\
&+{l'(2l'+1)}[\rho^{l+1}(a_{22}-a_{12})r^{l'-1}]\mathscr{D}_{l,l'}^{0,0}(\alpha)\\
&+\rho^{l+1}(a_{22}-a_{12})\cdot 2\sqrt{l'(2l'+1)}\sum_{\lambda\in\{l'+1,l'-1\}}\sum_{q=-1}^{1}\sum_{m_1=-1}^{1}r^\lambda \;\mathscr{A}_{l,\lambda}^{0,-q}(\alpha,q)K_{l'-1,l',\lambda}^{m_1,-q,q}\\
&+\rho^{l+1}a_{22}(2l+1)\sqrt{l'(2l'+1)}\sum_{q=-1}^{1}(-1)^q r^{l'-1}\;\mathscr{A}_{l,l'-1}^{0,-q}(\alpha,q)\langle l'-1,-q;1,q\mid l',0 \rangle\\
&+\delta_{ll'}\delta_{mm'}l(2l+1)\rho a_{22}.
\end{aligned}
$$
\end{theorem}

\begin{theorem}The remaining terms are simple:
$$M_{(l'm',3),(lm,1)}(\alpha)=M_{(l'm',1),(lm,3)}=0,$$
$$M_{(l'm',1),(lm1)}=\delta_{ll'}\delta_{mm'}\rho (l+1)(2l+1)  a_{11}.$$
\end{theorem}

\section{Application}
\subsection{Problem setting}
Given a function $\varphi\in L^2(\partial B_\rho(0))^3$ for $0<\rho<1/2$, the operator equation is set by
\begin{equation}\label{eq-Sf=phi}\mathcal{S}^{\alpha,0}_D[f]=\varphi,\end{equation}
where $D=B_\rho(0)$.
The goal is to determine the $L^2$ function $f$.
\subsection{Solve $f$}
It can be assumed that the coordinates in the expansion of $f$ under the vector spherical harmonics are $\{F_{lm}^k(\alpha)\}$. That means
\begin{equation}\label{eq-f-expression}f(x)=\sum\limits_{l,m,k}F_{lm}^{k}(\alpha)Y_{lm}^k(\hat{x}).\end{equation}
Starting from the equation \eqref{eq-Sf=phi},  taking inner product over the unit sphere $\mathbb{S}^2$ yields:
$$(Y_{l'm'}^p,\mathcal{S}_D^{\alpha,0}[f])=(Y_{l'm'}^p,\varphi),$$
and therefore,
$$\sum\limits_{l,m,k}\overline{F_{lm}^{k}(\alpha)}(Y_{l'm'}^p,\mathcal{S}_D^{\alpha,0}[Y_{lm}^k])=(Y_{l'm'}^p,\varphi).$$
\begin{remark}Since $W_{00}=X_{00}=0$, these terms should be excluded from the basis functions.
\end{remark}
 Now a system of linear equations is obtain:
\begin{equation}\mathbf{M}(\alpha) \overline{\mathbf{F}(\alpha)}=\mathbf{b},
\end{equation}
where $\overline{\boxed{}}$ represents complex conjugation.
The entries of the matrix $\mathbf{M}(\alpha)$ are given by  $M_{(l'm'p),(lmk)}(\alpha)$ in \eqref{eq-M-entries}. The vector $ \mathbf{F}(\alpha)$ is designed as
$$\mathbf{F}(\alpha) = 
\begin{pmatrix}
F_{00}^1 \\
F_{00}^2 \\
F_{00}^3 \\
F_{1,-1}^1 \\ F_{1,-1}^2 \\ F_{1,-1}^3 \\
F_{10}^1 \\ F_{10}^2 \\ F_{10}^3 \\
F_{11}^1 \\ F_{11}^2 \\ F_{11}^3 \\
F_{2,-2}^1 \\ \vdots \\
\vdots \\
F_{L_{\max},L_{\max}}^3
\end{pmatrix}$$
The components of  vector $\mathbf{b}\in\mathbb{R}^N( N=3\sum_{l=0}^{L_{\max}}(2l+1) )$ are expressed as 
$b_{l'm'}^p := \bigl(Y_{l'm'}^p,\varphi\bigr),$ which can be calculated as long as the exact expression of $\varphi$ is given.

 \subsection{Generalize to the dimer case}
 Let $D=D_1\cup D_2$ with $D_1=B_\rho(-d,0,0)$ and $D_2=B_\rho(d,0,0)$.  Here the parameter is designed 
such that these balls are disjoint, say, $d=0.2$ and $\rho=0.1$. Given an $L^2$ function $\varphi$ on $\partial D$, it  is expected that the function $f$ defined on $\partial D$ satisfying
$$\mathcal{S}_D^{\alpha,0}[f]=\varphi,$$
can be determined.
The function $f$ can be decomposed  as 
$$f=f|_{\partial D_1}\cdot \mathbf{1}_{\partial D_1}+f|_{\partial D_2}\cdot \mathbf{1}_{\partial D_2}.$$
Denote these two parts by $g_1$ and $g_2$.
Therefore, 
$$\mathcal{S}_D^{\alpha,0}[f](x)=\int_{\partial D_1}G^{\alpha,0}(x-y)g_1(y)d\sigma+\int_{\partial D_2}G^{\alpha,0}(x-y)g_2 (y)d\sigma=\varphi(x),\quad x\in\partial D.$$
Denote $\mathcal{S}^{\alpha,0}_{D_s\to D_t}[g_s](x):=\mathcal{S}^{\alpha,0}_{D_s}[g_s](x)$ with $x\in \partial D_t$ for $s,t=1,2$.
That is, the notation $\mathcal{S}_{D_s\to D_t}^{\alpha,0}$ means that the operator maps the function on $\partial D_s$ to the function on $\partial D_t$. 

For $x\in\partial D_1$, it holds that
$$\mathcal{S}_{D_1\to D_1}^{\alpha,0}[g_1](x)+\mathcal{S}_{D_2\to D_1}^{\alpha,0}[g_2](x)=\varphi |_{\partial D_1}(x).$$
Similarly, for $x\in\partial D_2$,
 $$\mathcal{S}_{D_1\to D_2}^{\alpha,0}[g_1](x)+\mathcal{S}_{D_2\to D_2}^{\alpha,0}[g_2](x)=\varphi|_{\partial D_2}(x).$$
 For convenience, it is written as
$$\begin{pmatrix}
\mathcal{S}_{D_1\to D_1}^{\alpha,0}&\mathcal{S}_{D_2\to D_1}^{\alpha,0}\\
\mathcal{S}_{D_1\to D_2}^{\alpha,0}&
\mathcal{S}_{D_2\to D_2}^{\alpha,0}\end{pmatrix}\begin{pmatrix}g_1\\g_2\end{pmatrix}=\begin{pmatrix}\varphi |_{\partial D_1}\\\varphi |_{\partial D_2}\end{pmatrix}.$$

Using the substitution $x=x_{i'}+\mathbf{c}_{i'}$, $i'=1,2$, the coordinate $x\in \partial D_{i'}$ can be pulled back into $x_{i'}\in \partial B_\rho(0)$,
where $\mathbf{c}_{1}=(-d,0,0)$ and $\mathbf{c}_2=(d,0,0)$. Also $c_1=-d$ and $c_2=d$. Applying this to the function $\varphi$ and defining the inner product over $\mathbb{S}^2$, one obtains:
$$b_{l'm',(i')}^{p}:=(Y_{l'm'}^p,\varphi|_{\partial D_{i'}}).$$
These form a vector called $\mathbf{b}_{i'}\in\mathbb{R}^{N}$ $( N=3\sum_{l=0}^{L_{\max}}(2l+1) )$ for $i'=1,2$.\\

 The function $g_{i'}(x)=g_{i'}(\mathbf{c}_{i'}+\rho \widehat{x_{i'}})$ can be regarded as a function with respect to $\widehat{x_{i'}}\in\mathbb{S}^2$. 
 Thus,  it is assumed that 
 $$\boxed{g_{i'}(\mathbf{c}_{i'}+\rho \widehat{x_{i'}})=\sum\limits_{l,m,k}F_{l,m,(i')}^{k}(\alpha) Y_{lm}^k(\widehat{x_{i'}}).}$$
For $x\in \partial D_t$, 
$$\begin{aligned}
\mathcal{S}_{D_s\to D_t}^{\alpha,0}[g_{s}](x)&
=\int_{\partial D_s}\boldsymbol{G}^{\alpha,0}(x-y)g_{s}(y)d\sigma(y)\\
&=\int_{\partial D_s}\boldsymbol{G}^{\alpha,0}(x-y)g_{s}(\mathbf{c}_{s}+\rho \widehat{y_s})d\sigma(y)\\
&=\int_{\partial B_\rho(0)}\boldsymbol{G}^{\alpha,0}(\mathbf{c}_t+\rho \widehat{x_t}-\mathbf{c}_s-\rho\widehat{y_s})g_{s}(\mathbf{c}_{s}+\rho \widehat{y_s})d\sigma(y_s)\\
&=\int_{\partial B_\rho(0)}\boldsymbol{G}^{\alpha,0}(\mathbf{c}_t+\rho \widehat{x_t}-\mathbf{c}_s-\rho\widehat{y_s})\sum\limits_{l,m,k}F_{l,m,(s)}^{k}(\alpha)Y_{lm}^{k}(\widehat{y_s})d\sigma(y_s).
\end{aligned}$$
It is observed that when $t=s$, $$\mathcal{S}_{D_s\to D_t}^{\alpha,0}[g_{s}](x)=\sum\limits_{l,m,k}F_{l,m,(s)}^{k}(\alpha)\mathcal{S}_D^{\alpha,0}[Y_{lm}^{k}](x_t),$$
where the operator $\mathcal{S}_D^{\alpha,0}$ is defined in Eq~\eqref{eq-def-S_D}. Therefore, the inner product over $\mathbb{S}^2$
$$(Y_{l'm'}^p(\widehat{x_t}), \mathcal{S}_{D_s\to D_t}^{\alpha,0}[g_{s}](x))=\sum\limits_{l,m,k}\overline{F_{l,m,(s)}^{k}(\alpha)}(Y_{l'm'}^p(\widehat{x_t}),\mathcal{S}_D^{\alpha,0}[Y_{lm}^{k}](x_t)),$$ 
in which $(Y_{l'm'}^p(\widehat{x_t}),\mathcal{S}_D^{\alpha,0}[Y_{lm}^{k}](x_t))$ is given by Eq~\eqref{eq-M-entries}.\\
When $s\neq t$, the periodic Green function 
$$\boldsymbol{G}^{\alpha,0}(\mathbf{c}_t+\rho \widehat{x_t}-\mathbf{c}_s-\rho\widehat{y_s})=\sum\limits_{n\in\mathbb{Z}}\boldsymbol{G}(x_t-y_s-(\mathbf{c}_s-\mathbf{c}_t+(n,0,0)))e^{-in\alpha}$$
So there are only two changes in the previous $M(\alpha)$:
\begin{enumerate}
\item Delete the diagonal terms $( Y_{l'm'}^p,\mathcal{S}_D Y_{lm}^q)$ and the summation should involve the case $n=0$, i.e.,
\begin{equation}\label{eq-M-ts-entries}M^{st}_{(l'm'p),(lmq)}(\alpha) = \sum_{n \in\mathbb{Z}}e^{-in\alpha}(Y_{l'm'}^p,\mathcal{S}^{st}_{D+n} Y_{lm}^q).\end{equation}
\item Replace $\mathbf{b}=(-n,0,0)$ by $\mathbf{b}_{st}=(-(c_s-c_t+n),0,0)$ in the explicit expression of $M^{st}_{(l'm'p),(lmq)}(\alpha)$. Meanwhile, the Lerch transcendent is introduced to deal with the infinite lattice sums.
\end{enumerate}
The exact expressions of $\mathbf{M}^{21}$ and $\mathbf{M}^{12}$ will be given in \ref{app1}.
 Therefore, the system of equations for solving $g_{i'}$ is given by
 $$\begin{pmatrix} \mathbf{M}^{11}(\alpha)&\mathbf{M}^{21}(\alpha)\\
 \mathbf{M}^{12}(\alpha)&\mathbf{M}^{22}(\alpha)\end{pmatrix}\begin{pmatrix}\overline{\mathbf{F}_{1}}\\\overline{\mathbf{F}_{2}}\end{pmatrix}=\begin{pmatrix}\mathbf{b}_1\\\mathbf{b}_2\end{pmatrix}.$$
 
 \begin{remark}Since $W_{0,0}=X_{0,0}=0$, they are discarded from the basis. So the corresponding components of these vectors and the rows and columns of these matrices should be deleted simultaneously.
 \end{remark}

\section{Conclusion}

This paper has established addition theorems for all three families of real vector spherical harmonics under translation, and used them to compute the full matrix representation of the quasi-periodic elastic single layer potential in closed form for one-dimensional arrays. A key feature of the approach is that the infinite lattice sums, which arise inevitably from the periodic geometry, are handled exactly via polylogarithm functions and the Lerch transcendent rather than by numerical truncation, so the resulting matrix entries carry no discretization error. The framework covers both the single-ball and dimer geometries, and reduces the operator equation to an explicit linear system that can be solved directly. Taken together, these results complete the elastic analogue of the multipole expansion method developed for the acoustic Helmholtz equation in Ammari's work, and remove what has been a missing analytical ingredient for the rigorous study of elastic subwavelength resonator chains. The explicit matrix representations obtained here are expected to serve as a practical tool in the spectral analysis of periodic elastic structures, including the investigation of band gaps, resonant frequencies, and topological phenomena in the elastic setting.

\newpage

\appendix
\section{Precise expressions of $\mathbf{M}^{21}(\alpha)$ and $\mathbf{M}^{12}(\alpha)$}
\label{app1}

\begin{definition}The {\bf Lerch transcendent} is defined by a power series in $z$ generalizing the polylogarithm:
$$\mathrm{\Phi}(z,s,D)=\sum\limits_{k=0}^\infty\dfrac{z^k}{(k+D)^s}.$$
\end{definition}

\subsection{The entries for $\mathbf{M}^{21}$}
\begin{lemma}Denote the series $\sum\limits_{n\in\mathbb{Z}}H_{l,\lambda}^{m,\mu}(\mathbf{b}_{21})e^{-in\alpha}$ by $\boxed{\mathscr{H}_{l,\lambda}^{m,\mu,(21)}(\alpha)}$. 
$$\begin{aligned}
\sum\limits_{n\in\mathbb{Z}}H_{l,\lambda}^{m,\mu}(\mathbf{b}_{21})e^{-in\alpha}=&(-1)^{\lambda+\mu}
\sqrt{\dfrac{2l+1}{2\lambda+1}}
\sqrt{{l+\lambda+\mu-m\choose \lambda+\mu}{l+\lambda+m-\mu\choose \lambda-\mu}} \sqrt{\dfrac{4\pi}{2(l+\lambda)+1}}\\&
\big(\mathrm{\Phi}(e^{-i\alpha},l+\lambda+1,2d)Y_{l+\lambda}^{m-\mu}(\dfrac{\pi}{2},\pi)+e^{i\alpha}\cdot \mathrm{\Phi}(e^{i\alpha},l+\lambda+1,1-2d)Y_{l+\lambda}^{m-\mu}(\dfrac{\pi}{2},0)\big).
\end{aligned}.$$

\end{lemma}

\begin{theorem}The expression of $\boxed{M^{21}_{(l'm'2),(lm1)}(\alpha)}$($l'\geqslant 1$, otherwise this term vanishes) is given by:
\\
If $m>0$ and $m'>0$,
$$\begin{aligned}M^{21}_{(l'm'2),(lm1)}(\alpha)= \rho^{l+3}a_{11}[r^{l'-1}\dfrac{1}{2}(\mathscr{H}_{l,l'}^{-m,-m',(21)}(\alpha)+(-1)^m\mathscr{H}_{l,l'}^{m,-m',(21)}(\alpha))l'(2l'+1)\\
+r^{l'-1}\dfrac{(-1)^{m'}}{2}(\mathscr{H}_{l,l'}^{-m,m',(21)}(\alpha)+(-1)^m\mathscr{H}_{l,l'}^{m,m',(21)}(\alpha))l'(2l'+1)]
\end{aligned}$$
If $m>0$ and $m'<0$,
$$\begin{aligned}M^{21}_{(l'm'2),(lm1)}(\alpha)= \rho^{l+3}a_{11}[r^{l'-1}\dfrac{1}{2}(\mathscr{H}_{l,l'}^{-m,m',(21)}(\alpha)+(-1)^m\mathscr{H}_{l,l'}^{m,m',(21)}(\alpha))(-i \cdot l'(2l'+1))\\
+r^{l'-1}\dfrac{(-1)^{m'}}{2}(\mathscr{H}_{l,l'}^{-m,-m',(21)}(\alpha)+(-1)^m\mathscr{H}_{l,l'}^{m,-m',(21)}(\alpha))(i\cdot l'(2l'+1))]
\end{aligned}$$
If $m>0$ and $m'=0$,
$$\begin{aligned}M^{21}_{(l'm'2),(lm1)}(\alpha)= \rho^{l+3}a_{11}r^{l'-1}\dfrac{1}{\sqrt{2}}(\mathscr{H}_{l,l'}^{-m,0,(21)}(\alpha)+(-1)^m\mathscr{H}_{l,l'}^{m,0,(21)}(\alpha))l'(2l'+1).
\end{aligned}$$

If $m<0$ and $m'>0$,
$$\begin{aligned}M^{21}_{(l'm'2),(lm1)}(\alpha)= \rho^{l+3}a_{11}[r^{l'-1}\dfrac{i}{2}(\mathscr{H}_{l,l'}^{m,-m',(21)}(\alpha)-(-1)^m \mathscr{H}_{l,l'}^{-m,-m',(21)}(\alpha))l'(2l'+1)\\
+r^{l'-1}\dfrac{(-1)^{m'}}{2}\cdot i (\mathscr{H}_{l,l'}^{m,m',(21)}(\alpha)-(-1)^m\mathscr{H}_{l,l'}^{-m,m',(21)}(\alpha))l'(2l'+1)]
\end{aligned}$$

If $m<0$ and $m'<0$,
$$\begin{aligned}M^{21}_{(l'm'2),(lm1)}(\alpha)= \rho^{l+3}a_{11}[r^{l'-1}\dfrac{i}{2}(\mathscr{H}_{l,l'}^{m,m',(21)}(\alpha)-(-1)^m \mathscr{H}_{l,l'}^{-m,m',(21)}(\alpha))(-i\cdot l'(2l'+1))\\
+r^{l'-1}\dfrac{(-1)^{m'}}{2}\cdot i (\mathscr{H}_{l,l'}^{m,-m',(21)}(\alpha)-(-1)^m\mathscr{H}_{l,l'}^{-m,-m',(21)}(\alpha))(i\cdot l'(2l'+1))]
\end{aligned}$$

If $m<0$ and $m'=0$,
$$\begin{aligned}M^{21}_{(l'm'2),(lm1)}(\alpha)= \rho^{l+3}a_{11}r^{l'-1}\dfrac{i}{\sqrt{2}}(\mathscr{H}_{l,l'}^{m,0,(21)}(\alpha)-(-1)^m\mathscr{H}_{l,l'}^{-m,0,(21)}(\alpha))l'(2l'+1).\end{aligned}$$

If $m=0$ and $m'>0$,
$$\begin{aligned}M^{21}_{(l'm'2),(lm1)}(\alpha)= \rho^{l+3}a_{11}r^{l'-1}
[\dfrac{1}{\sqrt{2}} \mathscr{H}_{l,l'}^{0,-m',(21)}(\alpha)l'(2l'+1)+\dfrac{(-1)^{m'}}{\sqrt{2}}\mathscr{H}_{l,l'}^{0,m',(21)}(\alpha)l'(2l'+1)]
\end{aligned}$$

If $m=0$ and $m'<0$,
$$\begin{aligned}M^{21}_{(l'm'2),(lm1)}(\alpha)= \rho^{l+3}a_{11}r^{l'-1}
[\dfrac{1}{\sqrt{2}} \mathscr{H}_{l,l'}^{0,m',(21)}(\alpha)(-i\cdot l'(2l'+1))+\dfrac{(-1)^{m'}}{\sqrt{2}}\mathscr{H}_{l,l'}^{0,-m',(21)}(\alpha)(i\cdot l'(2l'+1))]
\end{aligned}$$

If $m=0$ and $m'=0$,
$$\begin{aligned}M^{21}_{(l'm'2),(lm1)}(\alpha)= \rho^{l+3}a_{11}r^{l'-1}
[\mathscr{H}_{l,l'}^{0,0,(21)}(\alpha)l'(2l'+1)]
\end{aligned}$$
where $a_{ij}$ stands for the coefficient of the entries of $A^{out}_{\mathcal{S},l}(x)$.
\end{theorem}

\begin{lemma}Denote $\sum\limits_{n\in\mathbb{Z}}L_{l,j,\lambda}^{m,\mu,q,m_1}(\mathbf{b}_{21})e^{-in\alpha}$ by $\boxed{\mathscr{L}_{l,j,\lambda}^{m,\mu,q,m_1,(21)}(\alpha)}$, then
$$\begin{aligned}&\mathscr{L}_{l,j,\lambda}^{m,\mu,q,m_1,(21)}(\alpha)=
i (-1)^{\lambda+\mu+q}\text{sgn}(q-m_1)\epsilon_{q}\sqrt{\lambda(2l+1)} 
\sqrt{{l+\lambda+\mu-m\choose \lambda+\mu}{l+\lambda+m-\mu\choose \lambda-\mu}}\\
&\langle\lambda-1,\mu-m_1;1,m_1\mid \lambda,\mu\rangle
\langle \lambda-1,\mu-m_1;1,q+m_1\mid j,\mu+q \rangle \sqrt{\dfrac{4\pi}{2(l+\lambda)+1}}\\
&\big(\mathrm{\Phi}(e^{-i\alpha},l+\lambda,2d)Y_{l+\lambda}^{m-\mu}(\dfrac{\pi}{2},\pi)-e^{i\alpha}\mathrm{\Phi}(e^{i\alpha},l+\lambda,1-2d)Y_{l+\lambda}^{m-\mu}(\dfrac{\pi}{2},0)\big)
\end{aligned}$$
where $$\epsilon_q=\left\{
\begin{aligned}
\dfrac{1}{\sqrt{2}}&,\quad q=-1\\
0&,\quad q=0\\
-\dfrac{1}{\sqrt{2}}&,\quad,q=1 \end{aligned}\right.$$
\end{lemma}

\begin{theorem}
The expression of {$\boxed{M^{21}_{(l'm',2),(lm,3)}(\alpha)}$}($l'\geqslant 1,l\geqslant 1$,otherwise this term vanishes) is given by
When $m>0, m'>0$,
$$
\begin{aligned}
&M^{21}_{(l'm',2),(lm,3)}(\alpha)\\
=&\ \rho^{l+2}a_{33}\,r^{l'-1}\,\frac{\sqrt{l'(2l'+1)}}{2}\sum\limits_{q=-1}^{1}\sum\limits_{m_1=-1}^{1}\bigl[(-1)^{m'}\bigl(\mathscr{L}_{l,l',l'}^{-m,\,m'-q,\,q,\,m_1,(21)}(\alpha)+(-1)^m\mathscr{L}_{l,l',l'}^{m,\,m'-q,\,q,\,m_1,(21)}(\alpha)\bigr)\\
&\qquad\qquad\qquad\qquad+\bigl(\mathscr{L}_{l,l',l'}^{-m,\,-m'-q,\,q,\,m_1,(21)}(\alpha)+(-1)^m\mathscr{L}_{l,l',l'}^{m,\,-m'-q,\,q,\,m_1,(21)}(\alpha)\bigr)\bigr]
\end{aligned}
$$
When  $m>0, m'<0$,
$$
\begin{aligned}
&M^{21}_{(l'm',2),(lm,3)}(\alpha)\\
=&\ \rho^{l+2}a_{33}\,r^{l'-1}\,\frac{\sqrt{l'(2l'+1)}}{2}\sum\limits_{q=-1}^{1}\sum\limits_{m_1=-1}^{1}\bigl[i(-1)^{m'}\bigl(\mathscr{L}_{l,l',l'}^{-m,\,-m'-q,\,q,\,m_1,(21)}(\alpha)+(-1)^m\mathscr{L}_{l,l',l'}^{m,\,-m'-q,\,q,\,m_1,(21)}(\alpha)\bigr)\\
&\qquad\qquad\qquad\qquad-i\bigl(\mathscr{L}_{l,l',l'}^{-m,\,m'-q,\,q,\,m_1,(21)}(\alpha)+(-1)^m\mathscr{L}_{l,l',l'}^{m,\,m'-q,\,q,\,m_1,(21)}(\alpha)\bigr)\bigr]
\end{aligned}
$$

When $m>0, m'=0$,
$$
\begin{aligned}
&M^{21}_{(l',0,2),(lm,3)}(\alpha)\\
=&\ \rho^{l+2}a_{33}\,r^{l'-1}\,\frac{\sqrt{l'(2l'+1)}}{\sqrt{2}}\sum\limits_{q=-1}^{1}\sum\limits_{m_1=-1}^{1}\bigl(\mathscr{L}_{l,l',l'}^{-m,\,-q,\,q,\,m_1,(21)}(\alpha)+(-1)^m\mathscr{L}_{l,l',l'}^{m,\,-q,\,q,\,m_1,(21)}(\alpha)\bigr)
\end{aligned}
$$

When  $m<0, m'>0$,
$$
\begin{aligned}
&M^{21}_{(l'm',2),(lm,3)}(\alpha)\\
=&\ \rho^{l+2}a_{33}\,r^{l'-1}\,i\cdot\frac{\sqrt{l'(2l'+1)}}{2}\sum\limits_{q=-1}^{1}\sum\limits_{m_1=-1}^{1}\bigl[(-1)^{m'}\bigl(\mathscr{L}_{l,l',l'}^{m,\,m'-q,\,q,\,m_1,(21)}(\alpha)-(-1)^m\mathscr{L}_{l,l',l'}^{-m,\,m'-q,\,q,\,m_1,(21)}(\alpha)\bigr)\\
&\qquad\qquad\qquad\qquad+\bigl(\mathscr{L}_{l,l',l'}^{m,\,-m'-q,\,q,\,m_1,(21)}(\alpha)-(-1)^m\mathscr{L}_{l,l',l'}^{-m,\,-m'-q,\,q,\,m_1,(21)}(\alpha)\bigr)\bigr]
\end{aligned}
$$
When $m<0, m'<0$,
$$
\begin{aligned}
&M^{21}_{(l'm',2),(lm,3)}(\alpha)\\
=&\ \rho^{l+2}a_{33}\,r^{l'-1}\,\frac{\sqrt{l'(2l'+1)}}{2}\sum\limits_{q=-1}^{1}\sum\limits_{m_1=-1}^{1}\bigl[\bigl(\mathscr{L}_{l,l',l'}^{m,\,m'-q,\,q,\,m_1,(21)}(\alpha)-(-1)^m\mathscr{L}_{l,l',l'}^{-m,\,m'-q,\,q,\,m_1,(21)}(\alpha)\bigr)\\
&\qquad\qquad\qquad\qquad-(-1)^{m'}\bigl(\mathscr{L}_{l,l',l'}^{m,\,-m'-q,\,q,\,m_1,(21)}(\alpha)-(-1)^m\mathscr{L}_{l,l',l'}^{-m,\,-m'-q,\,q,\,m_1,(21)}(\alpha)\bigr)\bigr]
\end{aligned}
$$
When $m<0, m'=0$,
$$
\begin{aligned}
&M^{21}_{(l',0,2),(lm,3)}(\alpha)\\
=&\ \rho^{l+2}a_{33}\,r^{l'-1}\,\frac{i\sqrt{l'(2l'+1)}}{\sqrt{2}}\sum\limits_{q=-1}^{1}\sum\limits_{m_1=-1}^{1}\bigl(\mathscr{L}_{l,l',l'}^{m,\,-q,\,q,\,m_1,(21)}(\alpha)-(-1)^m\mathscr{L}_{l,l',l'}^{-m,\,-q,\,q,\,m_1,(21)}(\alpha)\bigr)
\end{aligned}
$$
When $m=0, m'>0$,
$$
\begin{aligned}
&M^{21}_{(l'm',2),(l,0,3)}(\alpha)\\
=&\ \rho^{l+2}a_{33}\,r^{l'-1}\,\frac{\sqrt{l'(2l'+1)}}{\sqrt{2}}\sum\limits_{q=-1}^{1}\sum\limits_{m_1=-1}^{1}\bigl[(-1)^{m'}\mathscr{L}_{l,l',l'}^{0,\,m'-q,\,q,\,m_1,(21)}(\alpha)+\mathscr{L}_{l,l',l'}^{0,\,-m'-q,\,q,\,m_1,(21)}(\alpha)\bigr]
\end{aligned}
$$

When  $m=0, m'<0$,
$$
\begin{aligned}
&M^{21}_{(l'm',2),(l,0,3)}(\alpha)\\
=&\ \rho^{l+2}a_{33}\,r^{l'-1}\,\frac{\sqrt{l'(2l'+1)}}{\sqrt{2}}\sum\limits_{q=-1}^{1}\sum\limits_{m_1=-1}^{1}\bigl[i(-1)^{m'}\mathscr{L}_{l,l',l'}^{0,\,-m'-q,\,q,\,m_1,(21)}(\alpha)-i\,\mathscr{L}_{l,l',l'}^{0,\,m'-q,\,q,\,m_1,(21)}(\alpha)\bigr]
\end{aligned}
$$

When $m=0, m'=0$,
$$
\begin{aligned}
&M^{21}_{(l',0,2),(l,0,3)}(\alpha)\\
=&\ \rho^{l+2}a_{33}\,r^{l'-1}\,\sqrt{l'(2l'+1)}\sum\limits_{q=-1}^{1}\sum\limits_{m_1=-1}^{1}\mathscr{L}_{l,l',l'}^{0,\,-q,\,q,\,m_1,(21)}(\alpha)
\end{aligned}
$$
\end{theorem}

\begin{lemma}Denote $\sum\limits_{n\in\mathbb{Z}}b^{(21)}_{-q}H_{l,\lambda}^{m,\mu}(\mathbf{b}_{21})e^{-in\alpha}$ by $\boxed{\mathscr{A}_{l,\lambda}^{m,\mu,(21)}(\alpha,q)}$, then 
$$\begin{aligned}
\sum\limits_{n\in\mathbb{Z}}b^{(21)}_{-q}H_{l,\lambda}^{m,\mu}(\mathbf{b}_{21})e^{-in\alpha}=&(-1)^{\lambda+\mu}\epsilon_q
\sqrt{\dfrac{2l+1}{2\lambda+1}}
\sqrt{{l+\lambda+\mu-m\choose \lambda+\mu}{l+\lambda+m-\mu\choose \lambda-\mu}} \sqrt{\dfrac{4\pi}{2(l+\lambda)+1}}\\&
\big(\mathrm{\Phi}(e^{-i\alpha},l+\lambda,2d)Y_{l+\lambda}^{m-\mu}(\dfrac{\pi}{2},\pi)-e^{i\alpha}\mathrm{\Phi}(e^{i\alpha},l+\lambda,1-2d)Y_{l+\lambda}^{m-\mu}(\dfrac{\pi}{2},0)\big)
\end{aligned}.$$

\end{lemma}

\begin{theorem}
The expression of {$\boxed{M^{21}_{(l'm',1),(lm,2)}(\alpha)}$}($l\geqslant 1$, otherwise this term vanishes) is given by\\
When $m>0$, $m'>0$,
$$
\begin{aligned}
&M^{21}_{(l'm',1),(lm,2)}(\alpha)\\
=&\sum\limits_{\lambda\in\{l'+1,l'+3\}}\sum\limits_{q=-1}^{1}\sum\limits_{m_1=-1}^{1}\rho^{l+1} (a_{22}-a_{12})  r^\lambda \sqrt{(l'+1)(2l'+1)}
[(-1)^{m'}(\mathscr{A}_{l,\lambda}^{-m,m'-q,(21)}(\alpha,q)\\
&+(-1)^m\mathscr{A}_{l,\lambda}^{m,m'-q,(21)}(\alpha,q))K_{l'+1,l',\lambda}^{m_1,m'-q,q}
+(\mathscr{A}_{l,\lambda}^{-m,-m'-q,(21)}(\alpha,q)+\\
&(-1)^m\mathscr{A}_{l,\lambda}^{m,-m'-q,(21)}(\alpha,q))K_{l'+1,l',\lambda}^{m_1,-m'-q,q}]\\
&-\rho^{l+1} a_{22}r^{l'+1}\dfrac{(2l+1)(l'+1)}{2}[(-1)^{m'}(\mathscr{H}_{l,l'}^{-m,m',(21)}(\alpha)\\
&+(-1)^m\mathscr{H}_{l,l'}^{m,m',(21)}(\alpha))+(\mathscr{H}_{l,l'}^{-m,-m',(21)}(\alpha)+(-1)^m\mathscr{H}_{l,l'}^{m,-m',(21)}(\alpha))]\\
&+\rho^{l+1} a_{22}r^{l'+1}(2l+1)\dfrac{\sqrt{(l'+1)(2l'+1)}}{2} \sum\limits_{q=-1}^1(-1)^q [(-1)^{m'}(\mathscr{A}_{l,l'+1}^{-m,m'-q,(21)}(\alpha,q)\\
&+(-1)^m\mathscr{A}_{l,l'+1}^{m,m'-q,(21)}(\alpha,q))\langle l'+1,m'-q;1,q\mid l',m'\rangle+(\mathscr{A}_{l,l'+1}^{-m,-m'-q,(21)}(\alpha,q)\\
&+(-1)^m\mathscr{A}_{l,l'+1}^{m,-m'-q,(21)}(\alpha,q))\langle l'+1,-m'-q;1,q\mid l',-m'\rangle]
\end{aligned}$$

When $m>0$, $m'<0$,

$$
\begin{aligned}
&M^{21}_{(l'm',1),(lm,2)}(\alpha)\\
=&\sum\limits_{\lambda\in\{l'+1,l'+3\}}\sum\limits_{q=-1}^{1}\sum\limits_{m_1=-1}^{1}\rho^{l+1} (a_{22}-a_{12})  r^\lambda \sqrt{(l'+1)(2l'+1)}
[i\cdot (-1)^{m'}(\mathscr{A}_{l,\lambda}^{-m,-m'-q,(21)}(\alpha,q)\\
&+(-1)^m\mathscr{A}_{l,\lambda}^{m,-m'-q,(21)}(\alpha,q))K_{l'+1,l',\lambda}^{m_1,-m'-q,q}
-i\cdot (\mathscr{A}_{l,\lambda}^{-m,m'-q,(21)}(\alpha,q)\\
&+(-1)^m\mathscr{A}_{l,\lambda}^{m,m'-q,(21)}(\alpha,q))K_{l'+1,l',\lambda}^{m_1,m'-q,q}]\\
&-\rho^{l+1} a_{22}r^{l'+1}\dfrac{(2l+1)(l'+1)}{2}[i\cdot (-1)^{m'}(\mathscr{H}_{l,l'}^{-m,-m',(21)}(\alpha)+(-1)^m\mathscr{H}_{l,l'}^{m,-m',(21)}(\alpha))\\
&-i\cdot (\mathscr{H}_{l,l'}^{-m,m',(21)}(\alpha)+(-1)^m\mathscr{H}_{l,l'}^{m,m',(21)}(\alpha))]\\
&+\rho^{l+1} a_{22}r^{l'+1}(2l+1)\dfrac{\sqrt{(l'+1)(2l'+1)}}{2} \sum\limits_{q=-1}^1(-1)^q [i\cdot (-1)^{m'}(\mathscr{A}_{l,l'+1}^{-m,-m'-q,(21)}(\alpha,q)\\
&+(-1)^m\mathscr{A}_{l,l'+1}^{m,-m'-q,(21)}(\alpha,q))\\
&\langle l'+1,-m'-q;1,q\mid l',-m'\rangle
-i\cdot (\mathscr{A}_{l,l'+1}^{-m,m'-q,(21)}(\alpha,q)\\
&+(-1)^m\mathscr{A}_{l,l'+1}^{m,m'-q,(21)}(\alpha,q))\langle l'+1,m'-q;1,q\mid l',m'\rangle]
\end{aligned}$$

When $m>0$, $m'=0$,
$$
\begin{aligned}
&M^{21}_{(l',0,1),(lm,2)}(\alpha)\\
=&\sum\limits_{\lambda\in\{l'+1,l'+3\}}\sum\limits_{q=-1}^{1}\sum\limits_{m_1=-1}^{1}\rho^{l+1} (a_{22}-a_{12})  r^\lambda \sqrt{2(l'+1)(2l'+1)}
[(\mathscr{A}_{l,\lambda}^{-m,-q,(21)}(\alpha,q)\\
&+(-1)^m\mathscr{A}_{l,\lambda}^{m,-q,(21)}(\alpha,q))K_{l'+1,l',\lambda}^{m_1,-q,q}]\\
&-\rho^{l+1} a_{22}r^{l'+1}\dfrac{(2l+1)(l'+1)}{\sqrt{2}}[\mathscr{H}_{l,l'}^{-m,0,(21)}(\alpha)+(-1)^m\mathscr{H}_{l,l'}^{m,0,(21)}(\alpha)]\\
&+\rho^{l+1} a_{22}r^{l'+1}(2l+1)\dfrac{\sqrt{(l'+1)(2l'+1)}}{\sqrt{2}} \sum\limits_{q=-1}^1(-1)^q (\mathscr{A}_{l,l'+1}^{-m,-q,(21)}(\alpha,q)+(-1)^m\mathscr{A}_{l,l'+1}^{m,-q,(21)}(\alpha,q))\\
&\langle l'+1,-q;1,q\mid l',0\rangle
\end{aligned}$$

When $m<0$, $m'>0$,
$$
\begin{aligned}
&M^{21}_{(l'm',1),(lm,2)}(\alpha)\\
=&\sum\limits_{\lambda\in\{l'+1,l'+3\}}\sum\limits_{q=-1}^{1}\sum\limits_{m_1=-1}^{1}\rho^{l+1} (a_{22}-a_{12})  r^\lambda i\cdot\sqrt{(l'+1)(2l'+1)}
[(-1)^{m'}(\mathscr{A}_{l,\lambda}^{m,m'-q,(21)}(\alpha,q)\\
&-(-1)^m\mathscr{A}_{l,\lambda}^{-m,m'-q,(21)}(\alpha,q))K_{l'+1,l',\lambda}^{m_1,m'-q,q}
+(\mathscr{A}_{l,\lambda}^{m,-m'-q,(21)}(\alpha,q)\\
&-(-1)^m\mathscr{A}_{l,\lambda}^{-m,-m'-q,(21)}(\alpha,q))K_{l'+1,l',\lambda}^{m_1,-m'-q,q}]\\
&-\rho^{l+1} a_{22}r^{l'+1}i\cdot \dfrac{(2l+1)(l'+1)}{2}[(-1)^{m'}(\mathscr{H}_{l,l'}^{m,m',(21)}(\alpha)\\
&-(-1)^m\mathscr{H}_{l,l'}^{-m,m',(21)}(\alpha))+(\mathscr{H}_{l,l'}^{m,-m',(21)}(\alpha)-(-1)^m\mathscr{H}_{l,l'}^{-m,-m',(21)}(\alpha))]\\
&+\rho^{l+1} a_{22}r^{l'+1}i\cdot (2l+1)\dfrac{\sqrt{(l'+1)(2l'+1)}}{2} \sum\limits_{q=-1}^1(-1)^q [(-1)^{m'}(\mathscr{A}_{l,l'+1}^{m,m'-q,(21)}(\alpha,q)\\
&-(-1)^m\mathscr{A}_{l,l'+1}^{-m,m'-q,(21)}(\alpha,q))\langle l'+1,m'-q;1,q\mid l',m'\rangle+(\mathscr{A}_{l,l'+1}^{m,-m'-q,(21)}(\alpha,q)\\
&-(-1)^m\mathscr{A}_{l,l'+1}^{-m,-m'-q,(21)}(\alpha,q))\langle l'+1,-m'-q;1,q\mid l',-m'\rangle]
\end{aligned}$$

When $m<0$, $m'<0$,

$$
\begin{aligned}
&M^{21}_{(l'm',1),(lm,2)}(\alpha)\\
=&\sum\limits_{\lambda\in\{l'+1,l'+3\}}\sum\limits_{q=-1}^{1}\sum\limits_{m_1=-1}^{1}\rho^{l+1} (a_{22}-a_{12})  r^\lambda \sqrt{(l'+1)(2l'+1)}
[-(-1)^{m'}(\mathscr{A}_{l,\lambda}^{m,-m'-q,(21)}(\alpha,q)\\
&-(-1)^m\mathscr{A}_{l,\lambda}^{-m,-m'-q,(21)}(\alpha,q))K_{l'+1,l',\lambda}^{m_1,-m'-q,q}
+ (\mathscr{A}_{l,\lambda}^{m,m'-q,(21)}(\alpha,q)\\
&-(-1)^m\mathscr{A}_{l,\lambda}^{-m,m'-q,(21)}(\alpha,q))K_{l'+1,l',\lambda}^{m_1,m'-q,q}]\\
&-\rho^{l+1} a_{22}r^{l'+1}\dfrac{(2l+1)(l'+1)}{2}
[- (-1)^{m'}(\mathscr{H}_{l,l'}^{m,-m',(21)}(\alpha)-(-1)^m\mathscr{H}_{l,l'}^{-m,-m',(21)}(\alpha))\\
&+
 (\mathscr{H}_{l,l'}^{m,m',(21)}(\alpha)-(-1)^m\mathscr{H}_{l,l'}^{-m,m',(21)}(\alpha))]\\
&+\rho^{l+1} a_{22}r^{l'+1}(2l+1)\dfrac{\sqrt{(l'+1)(2l'+1)}}{2} \sum\limits_{q=-1}^1(-1)^q [-(-1)^{m'}(\mathscr{A}_{l,l'+1}^{m,-m'-q,(21)}(\alpha,q)\\
&-(-1)^m\mathscr{A}_{l,l'+1}^{-m,-m'-q,(21)}(\alpha,q))\langle l'+1,-m'-q;1,q\mid l',-m'\rangle
+ (\mathscr{A}_{l,l'+1}^{m,m'-q,(21)}(\alpha,q)\\
&-(-1)^m\mathscr{A}_{l,l'+1}^{-m,m'-q,(21)}(\alpha,q))\langle l'+1,m'-q;1,q\mid l',m'\rangle]
\end{aligned}$$

When $m<0$, $m'=0$,
$$
\begin{aligned}
&M^{21}_{(l',0,1),(lm,2)}(\alpha)\\
=&\sum\limits_{\lambda\in\{l'+1,l'+3\}}\sum\limits_{q=-1}^{1}\sum\limits_{m_1=-1}^{1}\rho^{l+1} (a_{22}-a_{12})  r^\lambda i\cdot \sqrt{2(l'+1)(2l'+1)}
[(\mathscr{A}_{l,\lambda}^{m,-q,(21)}(\alpha,q)\\
&-(-1)^m\mathscr{A}_{l,\lambda}^{-m,-q,(21)}(\alpha,q))K_{l'+1,l',\lambda}^{m_1,-q,q}]\\
&-\rho^{l+1} a_{22}r^{l'+1}i\cdot \dfrac{(2l+1)(l'+1)}{\sqrt{2}}[(\mathscr{H}_{l,l'}^{m,0,(21)}(\alpha)-(-1)^m\mathscr{H}_{l,l'}^{-m,0,(21)}(\alpha)]\\
&+\rho^{l+1} a_{22}r^{l'+1}i\cdot (2l+1)\dfrac{\sqrt{(l'+1)(2l'+1)}}{\sqrt{2}} \sum\limits_{q=-1}^1(-1)^q [ (\mathscr{A}_{l,l'+1}^{m,-q,(21)}(\alpha,q)-(-1)^m\mathscr{A}_{l,l'+1}^{-m,-q,(21)}(\alpha,q))\\
&\langle l'+1,-q;1,q\mid l',0\rangle]
\end{aligned}$$

When $m=0$, $m'>0$,
$$
\begin{aligned}
&M^{21}_{(l'm',1),(l,0,2)}(\alpha)\\
=&\sum\limits_{\lambda\in\{l'+1,l'+3\}}\sum\limits_{q=-1}^{1}\sum\limits_{m_1=-1}^{1}\rho^{l+1} (a_{22}-a_{12})  r^\lambda \sqrt{2(l'+1)(2l'+1)}
[(-1)^{m'}\mathscr{A}_{l,\lambda}^{0,m'-q,(21)}(\alpha,q) K_{l'+1,l',\lambda}^{m_1,m'-q,q}\\
&+\mathscr{A}_{l,\lambda}^{0,-m'-q,(21)}(\alpha,q)K_{l'+1,l',\lambda}^{m_1,-m'-q,q}]\\
&-\rho^{l+1} a_{22}r^{l'+1}\dfrac{(2l+1)(l'+1)}{\sqrt{2}}[(-1)^{m'}\mathscr{H}_{l,l'}^{0,m',(21)}(\alpha)+\mathscr{H}_{l,l'}^{0,-m',(21)}(\alpha)]\\
&+\rho^{l+1} a_{22}r^{l'+1}(2l+1)\dfrac{\sqrt{(l'+1)(2l'+1)}}{\sqrt{2}} \sum\limits_{q=-1}^1(-1)^q [(-1)^{m'}\mathscr{A}_{l,l'+1}^{0,m'-q,(21)}(\alpha,q)\\
&\langle l'+1,m'-q;1,q\mid l',m'\rangle+(\mathscr{A}_{l,l'+1}^{0,-m'-q,(21)}(\alpha,q)\langle l'+1,-m'-q;1,q\mid l',-m'\rangle]
\end{aligned}$$

When $m=0$, $m'<0$,

$$
\begin{aligned}
&M^{21}_{(l'm',1),(l,0,2)}(\alpha)\\
=&\sum\limits_{\lambda\in\{l'+1,l'+3\}}\sum\limits_{q=-1}^{1}\sum\limits_{m_1=-1}^{1}\rho^{l+1} (a_{22}-a_{12})  r^\lambda \sqrt{2(l'+1)(2l'+1)}
[i\cdot (-1)^{m'}\mathscr{A}_{l,\lambda}^{0,-m'-q,(21)}(\alpha,q)K_{l'+1,l',\lambda}^{m_1,-m'-q,q}\\
&
-i\cdot \mathscr{A}_{l,\lambda}^{0,m'-q,(21)}(\alpha,q)K_{l'+1,l',\lambda}^{m_1,m'-q,q}]\\
&-\rho^{l+1} a_{22}r^{l'+1}\dfrac{(2l+1)(l'+1)}{\sqrt{2}}[i\cdot (-1)^{m'}\mathscr{H}_{l,l'}^{0,-m',(21)}(\alpha)-i\cdot \mathscr{H}_{l,l'}^{0,m',(21)}(\alpha)]\\
&+\rho^{l+1} a_{22}r^{l'+1}(2l+1)\dfrac{\sqrt{(l'+1)(2l'+1)}}{\sqrt{2}} \sum\limits_{q=-1}^1(-1)^q [i\cdot (-1)^{m'} \mathscr{A}_{l,l'+1}^{0,-m'-q,(21)}(\alpha,q)\\
&\langle l'+1,-m'-q;1,q\mid l',-m'\rangle-i\cdot \mathscr{A}_{l,l'+1}^{0,m'-q,(21)}(\alpha,q)\langle l'+1,m'-q;1,q\mid l',m'\rangle]
\end{aligned}$$

When $m=0$, $m'=0$,
$$
\begin{aligned}
&M^{21}_{(l',0,1),(l,0,2)}(\alpha)\\
=&\sum\limits_{\lambda\in\{l'+1,l'+3\}}\sum\limits_{q=-1}^{1}\sum\limits_{m_1=-1}^{1}\rho^{l+1} (a_{22}-a_{12})  r^\lambda \sqrt{4(l'+1)(2l'+1)}
[\mathscr{A}_{l,\lambda}^{0,-q,(21)}(\alpha,q)K_{l'+1,l',\lambda}^{m_1,-q,q}]\\
&-\rho^{l+1} a_{22}r^{l'+1}(2l+1)(l'+1)[\mathscr{H}_{l,l'}^{0,0,(21)}(\alpha)]\\
&+\rho^{l+1} a_{22}r^{l'+1}(2l+1)\sqrt{(l'+1)(2l'+1)} \sum\limits_{q=-1}^1(-1)^q [ \mathscr{A}_{l,l'+1}^{0,-q,(21)}(\alpha,q)\langle l'+1,-q;1,q\mid l',0\rangle]
\end{aligned}$$

\end{theorem}

\begin{theorem}
Entries for {$\boxed{M^{21}_{(l'm',3),(lm,2)}(\alpha)}$}($l'\geqslant 1,l\geqslant 1$, otherwise this term vanishes):\\
When $m>0$, $m'>0$,
$$
\begin{aligned}
&M^{21}_{(l'm',3),(lm,2)}(\alpha)\\
=&\sum\limits_{\lambda\in\{l',l'+2\mid l'\geqslant 1\}}\sum\limits_{q=-1}^{1}\sum\limits_{m_1=-1}^{1}-i\cdot \rho^{l+1} (a_{22}-a_{12})  r^\lambda \sqrt{l'(l'+1)}
[(-1)^{m'}(\mathscr{A}_{l,\lambda}^{-m,m'-q,(21)}(\alpha,q)\\
&+(-1)^m\mathscr{A}_{l,\lambda}^{m,m'-q,(21)}(\alpha,q))K_{l',l',\lambda}^{m_1,m'-q,q}
+(\mathscr{A}_{l,\lambda}^{-m,-m'-q,(21)}(\alpha,q)\\
&+(-1)^m\mathscr{A}_{l,\lambda}^{m,-m'-q,(21)}(\alpha,q))K_{l',l',\lambda}^{m_1,-m'-q,q}]\\
&- i\cdot \rho^{l+1} a_{22}r^{l'}(2l+1)\dfrac{\sqrt{l'(l'+1)}}{2} \sum\limits_{q=-1}^1(-1)^q [(-1)^{m'}(\mathscr{A}_{l,l'}^{-m,m'-q,(21)}(\alpha,q)+(-1)^m\mathscr{A}_{l,l'}^{m,m'-q,(21)}(\alpha,q))\\
&\langle l',m'-q;1,q\mid l',m'\rangle+(\mathscr{A}_{l,l'}^{-m,-m'-q,(21)}(\alpha,q)+(-1)^m\mathscr{A}_{l,l'}^{m,-m'-q,(21)}(\alpha,q))\langle l',-m'-q;1,q\mid l',-m'\rangle]
\end{aligned}$$

When $m>0$, $m'<0$,
$$
\begin{aligned}
&M^{21}_{(l'm',3),(lm,2)}(\alpha)\\
=&\sum\limits_{\lambda\in\{l',l'+2\mid l'\geqslant 1\}}\sum\limits_{q=-1}^{1}\sum\limits_{m_1=-1}^{1}\rho^{l+1} (a_{22}-a_{12})  r^\lambda \sqrt{l'(l'+1)}
[ (-1)^{m'}(\mathscr{A}_{l,\lambda}^{-m,-m'-q,(21)}(\alpha,q)\\
&+(-1)^m\mathscr{A}_{l,\lambda}^{m,-m'-q,(21)}(\alpha,q))K_{l',l',\lambda}^{m_1,-m'-q,q}
- (\mathscr{A}_{l,\lambda}^{-m,m'-q,(21)}(\alpha,q)+(-1)^m\mathscr{A}_{l,\lambda}^{m,m'-q,(21)}(\alpha,q))K_{l',l',\lambda}^{m_1,m'-q,q}]\\
&+\rho^{l+1} a_{22}r^{l'}(2l+1)\dfrac{\sqrt{l'(l'+1)}}{2} \sum\limits_{q=-1}^1(-1)^q [ (-1)^{m'}(\mathscr{A}_{l,l'}^{-m,-m'-q,(21)}(\alpha,q)+(-1)^m\mathscr{A}_{l,l'}^{m,-m'-q,(21)}(\alpha,q))\\
&\langle l',-m'-q;1,q\mid l',-m'\rangle
- (\mathscr{A}_{l,l'}^{-m,m'-q,(21)}(\alpha,q)+(-1)^m\mathscr{A}_{l,l'}^{m,m'-q,(21)}(\alpha,q))\langle l',m'-q;1,q\mid l',m'\rangle]\end{aligned}$$

When $m>0$, $m'=0$,
$$
\begin{aligned}
&M^{21}_{(l',0,3),(lm,2)}(\alpha)\\
=&\sum\limits_{\lambda\in\{l',l'+2\mid l'\geqslant 1\}}\sum\limits_{q=-1}^{1}\sum\limits_{m_1=-1}^{1} -i\cdot \rho^{l+1} (a_{22}-a_{12})  r^\lambda \sqrt{2l'(l'+1)}
[(\mathscr{A}_{l,\lambda}^{-m,-q,(21)}(\alpha,q)+\\
&(-1)^m\mathscr{A}_{l,\lambda}^{m,-q,(21)}(\alpha,q))K_{l',l',\lambda}^{m_1,-q,q}]\\
&-i\cdot \rho^{l+1} a_{22}r^{l'}(2l+1)\dfrac{\sqrt{l'(l'+1)}}{\sqrt{2}} \sum\limits_{q=-1}^1(-1)^q [ (\mathscr{A}_{l,l'}^{-m,-q,(21)}(\alpha,q)+(-1)^m\mathscr{A}_{l,l'}^{m,-q,(21)}(\alpha,q))\\
&\langle l',-q;1,q\mid l',0\rangle]
\end{aligned}$$

When $m<0$, $m'>0$,
$$
\begin{aligned}
&M^{21}_{(l'm',3),(lm,2)}(\alpha)\\
=&\sum\limits_{\lambda\in\{l',l'+2\mid l'\geqslant 1\}}\sum\limits_{q=-1}^{1}\sum\limits_{m_1=-1}^{1}\rho^{l+1} (a_{22}-a_{12})  r^\lambda \sqrt{l'(l'+1)}
[(-1)^{m'}(\mathscr{A}_{l,\lambda}^{m,m'-q,(21)}(\alpha,q)\\
&-(-1)^m\mathscr{A}_{l,\lambda}^{-m,m'-q,(21)}(\alpha,q))K_{l',l',\lambda}^{m_1,m'-q,q}
+(\mathscr{A}_{l,\lambda}^{m,-m'-q,(21)}(\alpha,q)\\
&-(-1)^m\mathscr{A}_{l,\lambda}^{-m,-m'-q,(21)}(\alpha,q))K_{l',l',\lambda}^{m_1,-m'-q,q}]\\
&+\rho^{l+1} a_{22}r^{l'} (2l+1)\dfrac{\sqrt{l'(l'+1)}}{2} \sum\limits_{q=-1}^1(-1)^q [(-1)^{m'}(\mathscr{A}_{l,l'}^{m,m'-q,(21)}(\alpha,q)-(-1)^m\mathscr{A}_{l,l'}^{-m,m'-q,(21)}(\alpha,q))\\
&\langle l',m'-q;1,q\mid l',m'\rangle+(\mathscr{A}_{l,l'}^{m,-m'-q,(21)}(\alpha,q)-(-1)^m\mathscr{A}_{l,l'}^{-m,-m'-q,(21)}(\alpha,q))\langle l',-m'-q;1,q\mid l',-m'\rangle]
\end{aligned}$$

When $m<0$, $m'<0$,
$$
\begin{aligned}
&M^{21}_{(l'm',3),(lm,2)}(\alpha)\\
=&\sum\limits_{\lambda\in\{l',l'+2\mid l'\geqslant 1\}}\sum\limits_{q=-1}^{1}\sum\limits_{m_1=-1}^{1}\rho^{l+1} (a_{22}-a_{12})  r^\lambda i\cdot \sqrt{l'(l'+1)}
[(-1)^{m'}(\mathscr{A}_{l,\lambda}^{m,-m'-q,(21)}(\alpha,q)\\
&-(-1)^m\mathscr{A}_{l,\lambda}^{-m,-m'-q,(21)}(\alpha,q))K_{l',l',\lambda}^{m_1,-m'-q,q}
- (\mathscr{A}_{l,\lambda}^{m,m'-q,(21)}(\alpha,q)\\
&-(-1)^m\mathscr{A}_{l,\lambda}^{-m,m'-q,(21)}(\alpha,q))K_{l',l',\lambda}^{m_1,m'-q,q}]\\
&-i\cdot \rho^{l+1} a_{22}r^{l'}(2l+1)\dfrac{\sqrt{l'(l'+1)}}{2} \sum\limits_{q=-1}^1(-1)^q [-(-1)^{m'}(\mathscr{A}_{l,l'}^{m,-m'-q,(21)}(\alpha,q)\\
&-(-1)^m\mathscr{A}_{l,l'}^{-m,-m'-q,(21)}(\alpha,q))
\langle l',-m'-q;1,q\mid l',-m'\rangle
+ (\mathscr{A}_{l,l'}^{m,m'-q,(21)}(\alpha,q)\\
&-(-1)^m\mathscr{A}_{l,l'}^{-m,m'-q,(21)}(\alpha,q))\langle l',m'-q;1,q\mid l',m'\rangle]
\end{aligned}$$

When $m<0$, $m'=0$,
$$
\begin{aligned}
&M^{21}_{(l',0,3),(lm,2)}(\alpha)\\
=&\sum\limits_{\lambda\in\{l',l'+2\mid l'\geqslant 1\}}\sum\limits_{q=-1}^{1}\sum\limits_{m_1=-1}^{1}  \rho^{l+1} (a_{22}-a_{12})  r^\lambda  \sqrt{2l'(l'+1)}
[(\mathscr{A}_{l,\lambda}^{m,-q,(21)}(\alpha,q)\\
&-(-1)^m\mathscr{A}_{l,\lambda}^{-m,-q,(21)}(\alpha,q))K_{l',l',\lambda}^{m_1,-q,q}]\\
&+\rho^{l+1} a_{22}r^{l'} (2l+1)\dfrac{\sqrt{l'(l'+1)}}{\sqrt{2}} \sum\limits_{q=-1}^1(-1)^q [ (\mathscr{A}_{l,l'}^{m,-q,(21)}(\alpha,q)-(-1)^m\mathscr{A}_{l,l'}^{-m,-q,(21)}(\alpha,q))\\
&\langle l',-q;1,q\mid l',0\rangle]
\end{aligned}$$

When $m=0$, $m'>0$,
$$
\begin{aligned}
&M^{21}_{(l'm',3),(l,0,2)}(\alpha)\\
=&\sum\limits_{\lambda\in\{l',l'+2\mid l'\geqslant 1\}}\sum\limits_{q=-1}^{1}\sum\limits_{m_1=-1}^{1} -i\cdot \rho^{l+1} (a_{22}-a_{12})  r^\lambda \sqrt{2l'(l'+1)}
[(-1)^{m'}\mathscr{A}_{l,\lambda}^{0,m'-q,(21)}(\alpha,q) K_{l',l',\lambda}^{m_1,m'-q,q}\\
&+\mathscr{A}_{l,\lambda}^{0,-m'-q,(21)}(\alpha,q)K_{l',l',\lambda}^{m_1,-m'-q,q}]\\
&-i\cdot \rho^{l+1} a_{22}r^{l'}(2l+1)\dfrac{\sqrt{l'(l'+1)}}{\sqrt{2}} \sum\limits_{q=-1}^1(-1)^q [(-1)^{m'}\mathscr{A}_{l,l'}^{0,m'-q,(21)}(\alpha,q)\\
&\langle l',m'-q;1,q\mid l',m'\rangle+(\mathscr{A}_{l,l'}^{0,-m'-q,(21)}(\alpha,q)\langle l',-m'-q;1,q\mid l',-m'\rangle]
\end{aligned}$$

When $m=0$, $m'<0$,

$$
\begin{aligned}
&M^{21}_{(l'm',3),(l,0,2)}(\alpha)\\
=&\sum\limits_{\lambda\in\{l',l'+2\mid l'\geqslant 1\}}\sum\limits_{q=-1}^{1}\sum\limits_{m_1=-1}^{1}\rho^{l+1} (a_{22}-a_{12})  r^\lambda \sqrt{2l'(l'+1)}
[(-1)^{m'}\mathscr{A}_{l,\lambda}^{0,-m'-q,(21)}(\alpha,q)K_{l',l',\lambda}^{m_1,-m'-q,q}\\
&
- \mathscr{A}_{l,\lambda}^{0,m'-q,(21)}(\alpha,q)K_{l',l',\lambda}^{m_1,m'-q,q}]\\
&+\rho^{l+1} a_{22}r^{l'}(2l+1)\dfrac{\sqrt{l'(l'+1)}}{\sqrt{2}} \sum\limits_{q=-1}^1(-1)^q [ (-1)^{m'} \mathscr{A}_{l,l'}^{0,-m'-q,(21)}(\alpha,q)\\
&\langle l',-m'-q;1,q\mid l',-m'\rangle- \mathscr{A}_{l,l'}^{0,m'-q,(21)}(\alpha,q)\langle l',m'-q;1,q\mid l',m'\rangle]
\end{aligned}$$

When $m=0$, $m'=0$,
$$
\begin{aligned}
&M^{21}_{(l',0,3),(l,0,2)}(\alpha)\\
=&\sum\limits_{\lambda\in\{l',l'+2\mid l'\geqslant 1\}}\sum\limits_{q=-1}^{1}\sum\limits_{m_1=-1}^{1} -i\cdot \rho^{l+1} (a_{22}-a_{12})  r^\lambda \sqrt{4l'(l'+1)}
[\mathscr{A}_{l,\lambda}^{0,-q,(21)}(\alpha,q)K_{l',l',\lambda}^{m_1,-q,q}]\\
&-i\cdot \rho^{l+1} a_{22}r^{l'}(2l+1)\sqrt{l'(l'+1)} \sum\limits_{q=-1}^1(-1)^q [ \mathscr{A}_{l,l'}^{0,-q,(21)}(\alpha,q)\langle l',-q;1,q\mid l',0\rangle]
\end{aligned}$$

\end{theorem}

\begin{theorem}
Entries for {$\boxed{M^{21}_{(l'm',3),(lm,3)}(\alpha)}$}($l'\geqslant 1, l\geqslant 1$, otherwise this term vanishes):\\

When $m>0,\ m'>0$:
$$
\begin{aligned}
&M^{21}_{(l'm',3),(lm,3)}(\alpha)\\
=&\ \rho^{l+2}a_{33}r^{l'}\Biggl\{
\frac{l'(l'+1)}{2}\Bigl[(-1)^{m'}\bigl(\mathscr{H}_{l,l'}^{-m,m',(21)}(\alpha)+(-1)^m \mathscr{H}_{l,l'}^{m,m',(21)}(\alpha)\bigr)
+\bigl(\mathscr{H}_{l,l'}^{-m,-m',(21)}(\alpha)\\
&+(-1)^m \mathscr{H}_{l,l'}^{m,-m',(21)}(\alpha)\bigr)\Bigr]\\
&-\frac{i\sqrt{l'(l'+1)}}{2}\sum_{q=-1}^{1}\sum_{m_1=-1}^{1}\Bigl[(-1)^{m'}\bigl(\mathscr{L}_{l,l',l'+1}^{-m,m'-q,q,m_1,(21)}(\alpha)+(-1)^m \mathscr{L}_{l,l',l'+1}^{m,m'-q,q,m_1,(21)}(\alpha)\bigr)\\
&\qquad\qquad\qquad\qquad\quad+\bigl(\mathscr{L}_{l,l',l'+1}^{-m,-m'-q,q,m_1,(21)}(\alpha)+(-1)^m \mathscr{L}_{l,l',l'+1}^{m,-m'-q,q,m_1,(21)}(\alpha)\bigr)\Bigr]\Biggr\}\end{aligned}
$$

When $m>0,\ m'<0$:
$$
\begin{aligned}
&M^{21}_{(l'm',3),(lm,3)}(\alpha)\\
=&\ \rho^{l+2}a_{33}r^{l'}\Biggl\{
\frac{l'(l'+1)}{2}\Bigl[i(-1)^{m'}\bigl(\mathscr{H}_{l,l'}^{-m,-m',(21)}(\alpha)+(-1)^m \mathscr{H}_{l,l'}^{m,-m',(21)}(\alpha)\bigr)
-i\bigl(\mathscr{H}_{l,l'}^{-m,m',(21)}(\alpha)\\
&+(-1)^m \mathscr{H}_{l,l'}^{m,m',(21)}(\alpha)\bigr)\Bigr]\\
&+\frac{\sqrt{l'(l'+1)}}{2}\sum_{q=-1}^{1}\sum_{m_1=-1}^{1}\Bigl[(-1)^{m'}\bigl(\mathscr{L}_{l,l',l'+1}^{-m,-m'-q,q,m_1,(21)}(\alpha)+(-1)^m \mathscr{L}_{l,l',l'+1}^{m,-m'-q,q,m_1,(21)}(\alpha)\bigr)\\
&\qquad\qquad\qquad\qquad\quad-\bigl(\mathscr{L}_{l,l',l'+1}^{-m,m'-q,q,m_1,(21)}(\alpha)+(-1)^m \mathscr{L}_{l,l',l'+1}^{m,m'-q,q,m_1,(21)}(\alpha)\bigr)\Bigr]\Biggr\}
\end{aligned}
$$

When $m>0,\ m'=0$:
$$
\begin{aligned}
&M^{21}_{(l',0,3),(lm,3)}(\alpha)\\
=&\ \rho^{l+2}a_{33}r^{l'}\Biggl\{
\frac{l'(l'+1)}{\sqrt{2}}\bigl(\mathscr{H}_{l,l'}^{-m,0,(21)}(\alpha)+(-1)^m \mathscr{H}_{l,l'}^{m,0,(21)}(\alpha)\bigr)\\
&-\frac{i\sqrt{l'(l'+1)}}{\sqrt{2}}\sum_{q=-1}^{1}\sum_{m_1=-1}^{1}\bigl(\mathscr{L}_{l,l',l'+1}^{-m,-q,q,m_1,(21)}(\alpha)+(-1)^m \mathscr{L}_{l,l',l'+1}^{m,-q,q,m_1,(21)}(\alpha)\bigr)\Biggr\}
\end{aligned}
$$

When $m<0,\ m'>0$:
$$
\begin{aligned}
&M^{21}_{(l'm',3),(lm,3)}(\alpha)\\
=&\ \rho^{l+2}a_{33}r^{l'}\Biggl\{
\frac{l'(l'+1)}{2}\Bigl[i(-1)^{m'}\bigl(\mathscr{H}_{l,l'}^{m,m',(21)}(\alpha)-(-1)^m \mathscr{H}_{l,l'}^{-m,m',(21)}(\alpha)\bigr)
+i\bigl(\mathscr{H}_{l,l'}^{m,-m',(21)}(\alpha)\\
&-(-1)^m \mathscr{H}_{l,l'}^{-m,-m',(21)}(\alpha)\bigr)\Bigr]\\
&+\frac{\sqrt{l'(l'+1)}}{2}\sum_{q=-1}^{1}\sum_{m_1=-1}^{1}\Bigl[(-1)^{m'}\bigl(\mathscr{L}_{l,l',l'+1}^{m,m'-q,q,m_1,(21)}(\alpha)-(-1)^m \mathscr{L}_{l,l',l'+1}^{-m,m'-q,q,m_1,(21)}(\alpha)\bigr)\\
&\qquad\qquad\qquad\qquad\quad+\bigl(\mathscr{L}_{l,l',l'+1}^{m,-m'-q,q,m_1,(21)}(\alpha)-(-1)^m \mathscr{L}_{l,l',l'+1}^{-m,-m'-q,q,m_1,(21)}(\alpha)\bigr)\Bigr]\Biggr\}
\end{aligned}
$$

When $m<0,\ m'<0$:
$$
\begin{aligned}
&M^{21}_{(l'm',3),(lm,3)}(\alpha)\\
=&\ \rho^{l+2}a_{33}r^{l'}\Biggl\{
\frac{l'(l'+1)}{2}\Bigl[\bigl(\mathscr{H}_{l,l'}^{m,m',(21)}(\alpha)-(-1)^m \mathscr{H}_{l,l'}^{-m,m',(21)}(\alpha)\bigr)\\
&
-(-1)^{m'}\bigl(\mathscr{H}_{l,l'}^{m,-m',(21)}(\alpha)-(-1)^m \mathscr{H}_{l,l'}^{-m,-m',(21)}(\alpha)\bigr)\Bigr]\\
&+\frac{i\sqrt{l'(l'+1)}}{2}\sum_{q=-1}^{1}\sum_{m_1=-1}^{1}\Bigl[(-1)^{m'}\bigl(\mathscr{L}_{l,l',l'+1}^{m,-m'-q,q,m_1,(21)}(\alpha)-(-1)^m \mathscr{L}_{l,l',l'+1}^{-m,-m'-q,q,m_1,(21)}(\alpha)\bigr)\\
&\qquad\qquad\qquad\qquad\quad-\bigl(\mathscr{L}_{l,l',l'+1}^{m,m'-q,q,m_1,(21)}(\alpha)-(-1)^m \mathscr{L}_{l,l',l'+1}^{-m,m'-q,q,m_1,(21)}(\alpha)\bigr)\Bigr]\Biggr\}
\end{aligned}
$$

When $m<0,\ m'=0$:
$$
\begin{aligned}
&M^{21}_{(l',0,3),(lm,3)}(\alpha)\\
=&\ \rho^{l+2}a_{33}r^{l'}\Biggl\{
\frac{il'(l'+1)}{\sqrt{2}}\bigl(\mathscr{H}_{l,l'}^{m,0,(21)}(\alpha)-(-1)^m \mathscr{H}_{l,l'}^{-m,0,(21)}(\alpha)\bigr)\\
&+\frac{\sqrt{l'(l'+1)}}{\sqrt{2}}\sum_{q=-1}^{1}\sum_{m_1=-1}^{1}\bigl(\mathscr{L}_{l,l',l'+1}^{m,-q,q,m_1,(21)}(\alpha)-(-1)^m \mathscr{L}_{l,l',l'+1}^{-m,-q,q,m_1,(21)}(\alpha)\bigr)\Biggr\}\end{aligned}
$$

When $m=0,\ m'>0$:
$$
\begin{aligned}
&M^{21}_{(l'm',3),(l,0,3)}(\alpha)\\
=&\ \rho^{l+2}a_{33}r^{l'}\Biggl\{
\frac{l'(l'+1)}{\sqrt{2}}\Bigl[(-1)^{m'}\mathscr{H}_{l,l'}^{0,m',(21)}(\alpha)+\mathscr{H}_{l,l'}^{0,-m',(21)}(\alpha)\Bigr]\\
&-\frac{i\sqrt{l'(l'+1)}}{\sqrt{2}}\sum_{q=-1}^{1}\sum_{m_1=-1}^{1}\Bigl[(-1)^{m'}\mathscr{L}_{l,l',l'+1}^{0,m'-q,q,m_1,(21)}(\alpha)+\mathscr{L}_{l,l',l'+1}^{0,-m'-q,q,m_1,(21)}(\alpha)\Bigr]\Biggr\}\end{aligned}
$$

When $m=0,\ m'<0$:
$$
\begin{aligned}
&M^{21}_{(l'm',3),(l,0,3)}(\alpha)\\
=&\ \rho^{l+2}a_{33}r^{l'}\Biggl\{
\frac{l'(l'+1)}{\sqrt{2}}\Bigl[i(-1)^{m'}\mathscr{H}_{l,l'}^{0,-m',(21)}(\alpha)-i\mathscr{H}_{l,l'}^{0,m',(21)}(\alpha)\Bigr]\\
&+\frac{\sqrt{l'(l'+1)}}{\sqrt{2}}\sum_{q=-1}^{1}\sum_{m_1=-1}^{1}\Bigl[(-1)^{m'}\mathscr{L}_{l,l',l'+1}^{0,-m'-q,q,m_1,(21)}(\alpha)-\mathscr{L}_{l,l',l'+1}^{0,m'-q,q,m_1,(21)}(\alpha)\Bigr]\Biggr\}\end{aligned}
$$

When $m=0,\ m'=0$:
$$
\begin{aligned}
&M^{21}_{(l',0,3),(l,0,3)}(\alpha)\\
=&\ \rho^{l+2}a_{33}r^{l'}\Biggl\{
l'(l'+1)\mathscr{H}_{l,l'}^{0,0,(21)}(\alpha)
-i\sqrt{l'(l'+1)}\sum_{q=-1}^{1}\sum_{m_1=-1}^{1}\mathscr{L}_{l,l',l'+1}^{0,-q,q,m_1,(21)}(\alpha)\Biggr\}
\end{aligned}
$$
\end{theorem}

\begin{lemma}Define $\mathscr{D}_{l,\lambda}^{m,\mu,(21)}(\alpha)=\sum\limits_{n\in\mathbb{Z}} (n+2d)^2H_{l,\lambda}^{m,\mu}(\mathbf{b}_{21})e^{-in\alpha}$. Then
$$\begin{aligned}
&\mathscr{D}_{l,\lambda}^{m,\mu,(21)}(\alpha)\\
=&(-1)^{\lambda+\mu}
\sqrt{\dfrac{2l+1}{2\lambda+1}}
\sqrt{{l+\lambda+\mu-m\choose \lambda+\mu}{l+\lambda+m-\mu\choose \lambda-\mu}} \sqrt{\dfrac{4\pi}{2(l+\lambda)+1}}\\
&
\big(\mathrm{\Phi}(e^{-i\alpha},l+\lambda-1,2d)Y_{l+\lambda}^{m-\mu}(\dfrac{\pi}{2},\pi)+e^{i\alpha}\mathrm{\Phi}(e^{i\alpha},l+\lambda-1,1-2d)Y_{l+\lambda}^{m-\mu}(\dfrac{\pi}{2},0)\big)
\end{aligned}$$

\end{lemma}

\begin{theorem}The entries  $\boxed{M^{21}_{(l'm',2),(lm,2)}(\alpha)}$($l'\geqslant 1,l\geqslant 1$, otherwise this term vanishes) are given by:
\\
If $m>0$ and $m'>0$,
$$\begin{aligned}&M^{21}_{(l'm',2),(lm,2)}(\alpha)\\
&=\dfrac{l'(2l'+1)}{2}[ \rho^{l+3}a_{12}r^{l'-1}+\rho^{l+1}(a_{22}-a_{12})(r^{l'+1} )+\rho^{l+1}a_{22}r^{l'+1}\dfrac{2l+1}{2l'+1}]\\
&[(-1)^{m'}(\mathscr{H}_{l,l'}^{-m,m',(21)}(\alpha)+(-1)^m\mathscr{H}_{l,l'}^{m,m',(21)}(\alpha))+(\mathscr{H}_{l,l'}^{-m,-m',(21)}(\alpha)+(-1)^m\mathscr{H}_{l,l'}^{m,-m',(21)}(\alpha))]\\
&+\dfrac{l'(2l'+1)}{2}[\rho^{l+1}(a_{22}-a_{12})r^{l'-1}][(-1)^{m'}(\mathscr{D}_{l,l'}^{-m,m',(21)}(\alpha)+(-1)^m\mathscr{D}_{l,l'}^{m,m',(21)}(\alpha))\\
&+(\mathscr{D}_{l,l'}^{-m,-m',(21)}(\alpha)+(-1)^m\mathscr{D}_{l,l'}^{m,-m',(21)}(\alpha))]\\
&+\rho^{l+1}(a_{22}-a_{12})\sqrt{l'(2l'+1)}\sum\limits_{\lambda\in\{l'+1,l'-1\}}\sum\limits_{q=-1}^{q=1}\sum\limits_{m_1=-1}^{1}r^\lambda [(-1)^{m'}(\mathscr{A}_{l,\lambda}^{-m,m'-q,(21)}(\alpha,q)\\
&+(-1)^m\mathscr{A}_{l,\lambda}^{m,m'-q,(21)}(\alpha,q))K_{l'-1,l',\lambda}^{m_1,m'-q,q}+(\mathscr{A}_{l,\lambda}^{-m,-m'-q,(21)}(\alpha,q)+(-1)^m\mathscr{A}_{l,\lambda}^{m,-m'-q,(21)}(\alpha,q))K_{l'-1,l',\lambda}^{m_1,-m'-q,q}]\\
&+\rho^{l+1}a_{22}(2l+1)\dfrac{\sqrt{l'(2l'+1)}}{2}\sum\limits_{q=-1}^1(-1)^q r^{l'-1}[(-1)^{m'}(\mathscr{A}_{l,l'-1}^{-m,m'-q,(21)}(\alpha,q)+(-1)^m\mathscr{A}_{l,l'-1}^{m,m'-q,(21)}(\alpha,q))\\
&\langle l'-1,m'-q;1,q\mid l',m' \rangle+(\mathscr{A}_{l,l'-1}^{-m,-m'-q,(21)}(\alpha,q)\\
&+(-1)^m\mathscr{A}_{l,l'-1}^{m,-m'-q,(21)}(\alpha,q))\langle l'-1,-m'-q;1,q\mid l',-m' \rangle]
\end{aligned}$$

If $m>0$ and $m'<0$,
$$
\begin{aligned}
&M^{21}_{(l'm',2),(lm,2)}(\alpha)\\
&=\dfrac{l'(2l'+1)}{2}\Bigl[ \rho^{l+3}a_{12}r^{l'-1}+\rho^{l+1}(a_{22}-a_{12})(r^{l'+1} )+\rho^{l+1}a_{22}r^{l'+1}\dfrac{2l+1}{2l'+1}\Bigr]\\
&\quad\Bigl[i\cdot(-1)^{m'}(\mathscr{H}_{l,l'}^{-m,-m',(21)}(\alpha)+(-1)^m\mathscr{H}_{l,l'}^{m,-m',(21)}(\alpha))-i\cdot(\mathscr{H}_{l,l'}^{-m,m',(21)}(\alpha)+(-1)^m\mathscr{H}_{l,l'}^{m,m',(21)}(\alpha))\Bigr]\\
&+\dfrac{l'(2l'+1)}{2}[\rho^{l+1}(a_{22}-a_{12})r^{l'-1}]\Bigl[i\cdot(-1)^{m'}(\mathscr{D}_{l,l'}^{-m,-m',(21)}(\alpha)+(-1)^m\mathscr{D}_{l,l'}^{m,-m',(21)}(\alpha))\\
&-i\cdot(\mathscr{D}_{l,l'}^{-m,m',(21)}(\alpha)+(-1)^m\mathscr{D}_{l,l'}^{m,m',(21)}(\alpha))\Bigr]
&+\rho^{l+1}(a_{22}-a_{12})\sqrt{l'(2l'+1)}\sum_{\lambda\in\{l'+1,l'-1\}}\sum_{q=-1}^{1}\sum_{m_1=-1}^{1}r^\lambda \Bigl[i\cdot(-1)^{m'}(\mathscr{A}_{l,\lambda}^{-m,-m'-q,(21)}(\alpha,q)\\
&\quad+(-1)^m\mathscr{A}_{l,\lambda}^{m,-m'-q,(21)}(\alpha,q))K_{l'-1,l',\lambda}^{m_1,-m'-q,q}-i\cdot(\mathscr{A}_{l,\lambda}^{-m,m'-q,(21)}(\alpha,q)+\\
&(-1)^m\mathscr{A}_{l,\lambda}^{m,m'-q,(21)}(\alpha,q))K_{l'-1,l',\lambda}^{m_1,m'-q,q}\Bigr]\\
&+\rho^{l+1}a_{22}(2l+1)\dfrac{\sqrt{l'(2l'+1)}}{2}\sum_{q=-1}^{1}(-1)^q r^{l'-1}\Bigl[i\cdot(-1)^{m'}(\mathscr{A}_{l,l'-1}^{-m,-m'-q,(21)}(\alpha,q)\\
&+(-1)^m\mathscr{A}_{l,l'-1}^{m,-m'-q,(21)}(\alpha,q))\langle l'-1,-m'-q;1,q\mid l',-m' \rangle-i\cdot(\mathscr{A}_{l,l'-1}^{-m,m'-q,(21)}(\alpha,q)\\
&+(-1)^m\mathscr{A}_{l,l'-1}^{m,m'-q,(21)}(\alpha,q))\langle l'-1,m'-q;1,q\mid l',m' \rangle\Bigr]
\end{aligned}
$$

If $m>0$ and $m'=0$,
$$
\begin{aligned}
&M^{21}_{(l',0,2),(lm,2)}(\alpha)\\
&=\dfrac{l'(2l'+1)}{\sqrt{2}}\Bigl[ \rho^{l+3}a_{12}r^{l'-1}+\rho^{l+1}(a_{22}-a_{12})(r^{l'+1} )+\rho^{l+1}a_{22}r^{l'+1}\dfrac{2l+1}{2l'+1}\Bigr]\\
&\quad\Bigl[(\mathscr{H}_{l,l'}^{-m,0,(21)}(\alpha)+(-1)^m\mathscr{H}_{l,l'}^{m,0,(21)}(\alpha))\Bigr]\\
&+\dfrac{l'(2l'+1)}{\sqrt{2}}[\rho^{l+1}(a_{22}-a_{12})r^{l'-1}]\Bigl[(\mathscr{D}_{l,l'}^{-m,0,(21)}(\alpha)+(-1)^m\mathscr{D}_{l,l'}^{m,0,(21)}(\alpha))\Bigr]\\
&+\rho^{l+1}(a_{22}-a_{12})\sqrt{2l'(2l'+1)}\sum_{\lambda\in\{l'+1,l'-1\}}\sum_{q=-1}^{1}\sum_{m_1=-1}^{1}r^\lambda \Bigl[(\mathscr{A}_{l,\lambda}^{-m,-q,(21)}(\alpha,q)\\
&+(-1)^m\mathscr{A}_{l,\lambda}^{m,-q,(21)}(\alpha,q))K_{l'-1,l',\lambda}^{m_1,-q,q}\Bigr]\\
&+\rho^{l+1}a_{22}(2l+1)\dfrac{\sqrt{l'(2l'+1)}}{\sqrt{2}}\sum_{q=-1}^{1}(-1)^q r^{l'-1}\Bigl[(\mathscr{A}_{l,l'-1}^{-m,-q,(21)}(\alpha,q)\\
&+(-1)^m\mathscr{A}_{l,l'-1}^{m,-q,(21)}(\alpha,q))\langle l'-1,-q;1,q\mid l',0 \rangle\Bigr].
\end{aligned}
$$

If $m<0$ and $m'>0$,
$$
\begin{aligned}
&M^{21}_{(l'm',2),(lm,2)}(\alpha)\\
&=\dfrac{l'(2l'+1)}{2}\Bigl[ \rho^{l+3}a_{12}r^{l'-1}+\rho^{l+1}(a_{22}-a_{12})(r^{l'+1} )+\rho^{l+1}a_{22}r^{l'+1}\dfrac{2l+1}{2l'+1}\Bigr]\\
&\quad\cdot i\Bigl[(-1)^{m'}(\mathscr{H}_{l,l'}^{m,m',(21)}(\alpha)-(-1)^m\mathscr{H}_{l,l'}^{-m,m',(21)}(\alpha))+(\mathscr{H}_{l,l'}^{m,-m',(21)}(\alpha)-(-1)^m\mathscr{H}_{l,l'}^{-m,-m',(21)}(\alpha))\Bigr]\\
&+\dfrac{l'(2l'+1)}{2}[\rho^{l+1}(a_{22}-a_{12})r^{l'-1}]\cdot i\Bigl[(-1)^{m'}(\mathscr{D}_{l,l'}^{m,m',(21)}(\alpha)-(-1)^m\mathscr{D}_{l,l'}^{-m,m',(21)}(\alpha))\\
&+(\mathscr{D}_{l,l'}^{m,-m',(21)}(\alpha)-(-1)^m\mathscr{D}_{l,l'}^{-m,-m',(21)}(\alpha))\Bigr]\\
&+\rho^{l+1}(a_{22}-a_{12})\sqrt{l'(2l'+1)}\sum_{\lambda\in\{l'+1,l'-1\}}\sum_{q=-1}^{1}\sum_{m_1=-1}^{1}r^\lambda \;i\Bigl[(-1)^{m'}(\mathscr{A}_{l,\lambda}^{m,m'-q,(21)}(\alpha,q)\\
&\quad-(-1)^m\mathscr{A}_{l,\lambda}^{-m,m'-q,(21)}(\alpha,q))K_{l'-1,l',\lambda}^{m_1,m'-q,q}+(\mathscr{A}_{l,\lambda}^{m,-m'-q,(21)}(\alpha,q)\\
&-(-1)^m\mathscr{A}_{l,\lambda}^{-m,-m'-q,(21)}(\alpha,q))K_{l'-1,l',\lambda}^{m_1,-m'-q,q}\Bigr]\\
&+\rho^{l+1}a_{22}(2l+1)\dfrac{\sqrt{l'(2l'+1)}}{2}\sum_{q=-1}^{1}(-1)^q r^{l'-1}\;i\Bigl[(-1)^{m'}(\mathscr{A}_{l,l'-1}^{m,m'-q,(21)}(\alpha,q)-(-1)^m\mathscr{A}_{l,l'-1}^{-m,m'-q,(21)}(\alpha,q))\\
&\quad\langle l'-1,m'-q;1,q\mid l',m' \rangle+(\mathscr{A}_{l,l'-1}^{m,-m'-q,(21)}(\alpha,q)\\
&-(-1)^m\mathscr{A}_{l,l'-1}^{-m,-m'-q,(21)}(\alpha,q))\langle l'-1,-m'-q;1,q\mid l',-m' \rangle\Bigr].
\end{aligned}
$$

If $m<0$ and $m'<0$,
$$
\begin{aligned}
&M^{21}_{(l'm',2),(lm,2)}(\alpha)\\
&=\dfrac{l'(2l'+1)}{2}\Bigl[ \rho^{l+3}a_{12}r^{l'-1}+\rho^{l+1}(a_{22}-a_{12})(r^{l'+1} )+\rho^{l+1}a_{22}r^{l'+1}\dfrac{2l+1}{2l'+1}\Bigr]\\
&\quad\Bigl[-(-1)^{m'}(\mathscr{H}_{l,l'}^{m,-m',(21)}(\alpha)-(-1)^m\mathscr{H}_{l,l'}^{-m,-m',(21)}(\alpha))+(\mathscr{H}_{l,l'}^{m,m',(21)}(\alpha)-(-1)^m\mathscr{H}_{l,l'}^{-m,m',(21)}(\alpha))\Bigr]\\
&+\dfrac{l'(2l'+1)}{2}[\rho^{l+1}(a_{22}-a_{12})r^{l'-1}]\Bigl[-(-1)^{m'}(\mathscr{D}_{l,l'}^{m,-m',(21)}(\alpha)-(-1)^m\mathscr{D}_{l,l'}^{-m,-m',(21)}(\alpha))\\
&+(\mathscr{D}_{l,l'}^{m,m',(21)}(\alpha)-(-1)^m\mathscr{D}_{l,l'}^{-m,m',(21)}(\alpha))\Bigr]\\
&+\rho^{l+1}(a_{22}-a_{12})\sqrt{l'(2l'+1)}\sum_{\lambda\in\{l'+1,l'-1\}}\sum_{q=-1}^{1}\sum_{m_1=-1}^{1}r^\lambda \Bigl[-(-1)^{m'}(\mathscr{A}_{l,\lambda}^{m,-m'-q,(21)}(\alpha,q)\\
&\quad-(-1)^m\mathscr{A}_{l,\lambda}^{-m,-m'-q,(21)}(\alpha,q))K_{l'-1,l',\lambda}^{m_1,-m'-q,q}+(\mathscr{A}_{l,\lambda}^{m,m'-q,(21)}(\alpha,q)\\
&-(-1)^m\mathscr{A}_{l,\lambda}^{-m,m'-q,(21)}(\alpha,q))K_{l'-1,l',\lambda}^{m_1,m'-q,q}\Bigr]\\
&+\rho^{l+1}a_{22}(2l+1)\dfrac{\sqrt{l'(2l'+1)}}{2}\sum_{q=-1}^{1}(-1)^q r^{l'-1}\Bigl[-(-1)^{m'}(\mathscr{A}_{l,l'-1}^{m,-m'-q,(21)}(\alpha,q)\\
&-(-1)^m\mathscr{A}_{l,l'-1}^{-m,-m'-q,(21)}(\alpha,q))\\
&\quad\langle l'-1,-m'-q;1,q\mid l',-m' \rangle+(\mathscr{A}_{l,l'-1}^{m,m'-q,(21)}(\alpha,q)\\
&-(-1)^m\mathscr{A}_{l,l'-1}^{-m,m'-q,(21)}(\alpha,q))\langle l'-1,m'-q;1,q\mid l',m' \rangle\Bigr].\end{aligned}
$$

If $m<0$ and $m'=0$,
$$
\begin{aligned}
&M^{21}_{(l',0,2),(lm,2)}(\alpha)\\
&=\dfrac{l'(2l'+1)}{\sqrt{2}}\Bigl[ \rho^{l+3}a_{12}r^{l'-1}+\rho^{l+1}(a_{22}-a_{12})(r^{l'+1} )+\rho^{l+1}a_{22}r^{l'+1}\dfrac{2l+1}{2l'+1}\Bigr]
\cdot i\Bigl[(\mathscr{H}_{l,l'}^{m,0,(21)}(\alpha)\\
&-(-1)^m\mathscr{H}_{l,l'}^{-m,0,(21)}(\alpha))\Bigr]\\
&+\dfrac{l'(2l'+1)}{\sqrt{2}}[\rho^{l+1}(a_{22}-a_{12})r^{l'-1}]\cdot i\Bigl[(\mathscr{D}_{l,l'}^{m,0,(21)}(\alpha)-(-1)^m\mathscr{D}_{l,l'}^{-m,0,(21)}(\alpha))\Bigr]\\
&+\rho^{l+1}(a_{22}-a_{12})\sqrt{2l'(2l'+1)}\sum_{\lambda\in\{l'+1,l'-1\}}\sum_{q=-1}^{1}\sum_{m_1=-1}^{1}r^\lambda \;i\Bigl[(\mathscr{A}_{l,\lambda}^{m,-q,(21)}(\alpha,q)\\
&-(-1)^m\mathscr{A}_{l,\lambda}^{-m,-q,(21)}(\alpha,q))K_{l'-1,l',\lambda}^{m_1,-q,q}\Bigr]\\
&+\rho^{l+1}a_{22}(2l+1)\dfrac{\sqrt{l'(2l'+1)}}{\sqrt{2}}\sum_{q=-1}^{1}(-1)^q r^{l'-1}\;i\Bigl[(\mathscr{A}_{l,l'-1}^{m,-q,(21)}(\alpha,q)\\
&-(-1)^m\mathscr{A}_{l,l'-1}^{-m,-q,(21)}(\alpha,q))\langle l'-1,-q;1,q\mid l',0 \rangle\Bigr].
\end{aligned}
$$

If $m=0$ and $m'>0$,
$$
\begin{aligned}
&M^{21}_{(l'm',2),(l,0,2)}(\alpha)\\
&=\dfrac{l'(2l'+1)}{\sqrt{2}}\Bigl[ \rho^{l+3}a_{12}r^{l'-1}+\rho^{l+1}(a_{22}-a_{12})(r^{l'+1} )+\rho^{l+1}a_{22}r^{l'+1}\dfrac{2l+1}{2l'+1}\Bigr]\\
&\quad\Bigl[(-1)^{m'}\mathscr{H}_{l,l'}^{0,m',(21)}(\alpha)+\mathscr{H}_{l,l'}^{0,-m',(21)}(\alpha)\Bigr]\\
&+\dfrac{l'(2l'+1)}{\sqrt{2}}[\rho^{l+1}(a_{22}-a_{12})r^{l'-1}]\Bigl[(-1)^{m'}\mathscr{D}_{l,l'}^{0,m',(21)}(\alpha)+\mathscr{D}_{l,l'}^{0,-m',(21)}(\alpha)\Bigr]\\
&+\rho^{l+1}(a_{22}-a_{12})\sqrt{2l'(2l'+1)}\sum_{\lambda\in\{l'+1,l'-1\}}\sum_{q=-1}^{1}\sum_{m_1=-1}^{1}r^\lambda \Bigl[(-1)^{m'}\mathscr{A}_{l,\lambda}^{0,m'-q,(21)}(\alpha,q)K_{l'-1,l',\lambda}^{m_1,m'-q,q}\\
&+\mathscr{A}_{l,\lambda}^{0,-m'-q,(21)}(\alpha,q)K_{l'-1,l',\lambda}^{m_1,-m'-q,q}\Bigr]\\
&+\rho^{l+1}a_{22}(2l+1)\dfrac{\sqrt{l'(2l'+1)}}{\sqrt{2}}\sum_{q=-1}^{1}(-1)^q r^{l'-1}\Bigl[(-1)^{m'}\mathscr{A}_{l,l'-1}^{0,m'-q,(21)}(\alpha,q)\langle l'-1,m'-q;1,q\mid l',m' \rangle\\
&+\mathscr{A}_{l,l'-1}^{0,-m'-q,(21)}(\alpha,q)\langle l'-1,-m'-q;1,q\mid l',-m' \rangle\Bigr].
\end{aligned}
$$

If $m=0$ and $m'<0$,
$$
\begin{aligned}
&M^{21}_{(l'm',2),(l,0,2)}(\alpha)\\
&=\dfrac{l'(2l'+1)}{\sqrt{2}}\Bigl[ \rho^{l+3}a_{12}r^{l'-1}+\rho^{l+1}(a_{22}-a_{12})(r^{l'+1} )+\rho^{l+1}a_{22}r^{l'+1}\dfrac{2l+1}{2l'+1}\Bigr]\\
&\quad i\Bigl[(-1)^{m'}\mathscr{H}_{l,l'}^{0,-m',(21)}(\alpha)-\mathscr{H}_{l,l'}^{0,m',(21)}(\alpha)\Bigr]\\
&+\dfrac{l'(2l'+1)}{\sqrt{2}}[\rho^{l+1}(a_{22}-a_{12})r^{l'-1}]\cdot i\Bigl[(-1)^{m'}\mathscr{D}_{l,l'}^{0,-m',(21)}(\alpha)-\mathscr{D}_{l,l'}^{0,m',(21)}(\alpha)\Bigr]\\
&+\rho^{l+1}(a_{22}-a_{12})\sqrt{2l'(2l'+1)}\sum_{\lambda\in\{l'+1,l'-1\}}\sum_{q=-1}^{1}\sum_{m_1=-1}^{1}r^\lambda \;i\Bigl[(-1)^{m'}\mathscr{A}_{l,\lambda}^{0,-m'-q,(21)}(\alpha,q)K_{l'-1,l',\lambda}^{m_1,-m'-q,q}\\
&-\mathscr{A}_{l,\lambda}^{0,m'-q,(21)}(\alpha,q)K_{l'-1,l',\lambda}^{m_1,m'-q,q}\Bigr]\\
&+\rho^{l+1}a_{22}(2l+1)\dfrac{\sqrt{l'(2l'+1)}}{\sqrt{2}}\sum_{q=-1}^{1}(-1)^q r^{l'-1}\;i\Bigl[(-1)^{m'}\mathscr{A}_{l,l'-1}^{0,-m'-q,(21)}(\alpha,q)\langle l'-1,-m'-q;1,q\mid l',-m' \rangle\\
&-\mathscr{A}_{l,l'-1}^{0,m'-q,(21)}(\alpha,q)\langle l'-1,m'-q;1,q\mid l',m' \rangle\Bigr].
\end{aligned}
$$

If $m=0$ and $m'=0$,
$$
\begin{aligned}
&M^{21}_{(l',0,2),(l,0,2)}(\alpha)\\
&=l'(2l'+1)\Bigl[ \rho^{l+3}a_{12}r^{l'-1}+\rho^{l+1}(a_{22}-a_{12})(r^{l'+1} )+\rho^{l+1}a_{22}r^{l'+1}\dfrac{2l+1}{2l'+1}\Bigr]\mathscr{H}_{l,l'}^{0,0,(21)}(\alpha)\\
&+{l'(2l'+1)}[\rho^{l+1}(a_{22}-a_{12})r^{l'-1}]\mathscr{D}_{l,l'}^{0,0,(21)}(\alpha)\\
&+\rho^{l+1}(a_{22}-a_{12})\cdot 2\sqrt{l'(2l'+1)}\sum_{\lambda\in\{l'+1,l'-1\}}\sum_{q=-1}^{1}\sum_{m_1=-1}^{1}r^\lambda \;\mathscr{A}_{l,\lambda}^{0,-q,(21)}(\alpha,q)K_{l'-1,l',\lambda}^{m_1,-q,q}\\
&+\rho^{l+1}a_{22}(2l+1)\sqrt{l'(2l'+1)}\sum_{q=-1}^{1}(-1)^q r^{l'-1}\;\mathscr{A}_{l,l'-1}^{0,-q,(21)}(\alpha,q)\langle l'-1,-q;1,q\mid l',0 \rangle.
\end{aligned}
$$
\end{theorem}

\begin{theorem}The remaining terms are simple:
$$M^{21}_{(l'm',3),(lm,1)}(\alpha)=M^{21}_{(l'm',1),(lm,3)}(\alpha)=0,$$
$$M^{21}_{(l'm',1),(lm1)}(\alpha)=0.$$
\end{theorem}

\subsection{The entries for $\mathbf{M}^{12}$}
\begin{lemma}Denote the series $\sum\limits_{n\in\mathbb{Z}}H_{l,\lambda}^{m,\mu}(\mathbf{b}_{12})e^{-in\alpha}$ by $\boxed{\mathscr{H}_{l,\lambda}^{m,\mu,(12)}(\alpha)}$. Then 
$$\begin{aligned}
\sum\limits_{n\in\mathbb{Z}}H_{l,\lambda}^{m,\mu}(\mathbf{b}_{12})e^{-in\alpha}=&(-1)^{\lambda+\mu}
\sqrt{\dfrac{2l+1}{2\lambda+1}}
\sqrt{{l+\lambda+\mu-m\choose \lambda+\mu}{l+\lambda+m-\mu\choose \lambda-\mu}} \sqrt{\dfrac{4\pi}{2(l+\lambda)+1}}\\&
\big(e^{-i\alpha}\cdot\mathrm{\Phi}(e^{-i\alpha},l+\lambda+1,1-2d)Y_{l+\lambda}^{m-\mu}(\dfrac{\pi}{2},\pi)+ \mathrm{\Phi}(e^{i\alpha},l+\lambda+1,2d)Y_{l+\lambda}^{m-\mu}(\dfrac{\pi}{2},0)\big)
\end{aligned}.$$

\end{lemma}

\begin{theorem}The expression of $\boxed{M^{12}_{(l'm'2),(lm1)}(\alpha)}$($l'\geqslant 1$, otherwise this term vanishes ) is given by:
\\
If $m>0$ and $m'>0$,
$$\begin{aligned}M^{12}_{(l'm'2),(lm1)}(\alpha)= \rho^{l+3}a_{11}[r^{l'-1}\dfrac{1}{2}(\mathscr{H}_{l,l'}^{-m,-m',(12)}(\alpha)+(-1)^m\mathscr{H}_{l,l'}^{m,-m',(12)}(\alpha))l'(2l'+1)\\
+r^{l'-1}\dfrac{(-1)^{m'}}{2}(\mathscr{H}_{l,l'}^{-m,m',(12)}(\alpha)+(-1)^m\mathscr{H}_{l,l'}^{m,m',(12)}(\alpha))l'(2l'+1)]
\end{aligned}$$
If $m>0$ and $m'<0$,
$$\begin{aligned}M^{12}_{(l'm'2),(lm1)}(\alpha)= \rho^{l+3}a_{11}[r^{l'-1}\dfrac{1}{2}(\mathscr{H}_{l,l'}^{-m,m',(12)}(\alpha)+(-1)^m\mathscr{H}_{l,l'}^{m,m',(12)}(\alpha))(-i \cdot l'(2l'+1))\\
+r^{l'-1}\dfrac{(-1)^{m'}}{2}(\mathscr{H}_{l,l'}^{-m,-m',(12)}(\alpha)+(-1)^m\mathscr{H}_{l,l'}^{m,-m',(12)}(\alpha))(i\cdot l'(2l'+1))]
\end{aligned}$$
If $m>0$ and $m'=0$,
$$\begin{aligned}M^{12}_{(l'm'2),(lm1)}(\alpha)= \rho^{l+3}a_{11}r^{l'-1}\dfrac{1}{\sqrt{2}}(\mathscr{H}_{l,l'}^{-m,0,(12)}(\alpha)+(-1)^m\mathscr{H}_{l,l'}^{m,0,(12)}(\alpha))l'(2l'+1).
\end{aligned}$$

If $m<0$ and $m'>0$,
$$\begin{aligned}M^{12}_{(l'm'2),(lm1)}(\alpha)= \rho^{l+3}a_{11}[r^{l'-1}\dfrac{i}{2}(\mathscr{H}_{l,l'}^{m,-m',(12)}(\alpha)-(-1)^m \mathscr{H}_{l,l'}^{-m,-m',(12)}(\alpha))l'(2l'+1)\\
+r^{l'-1}\dfrac{(-1)^{m'}}{2}\cdot i (\mathscr{H}_{l,l'}^{m,m',(12)}(\alpha)-(-1)^m\mathscr{H}_{l,l'}^{-m,m',(12)}(\alpha))l'(2l'+1)]
\end{aligned}$$

If $m<0$ and $m'<0$,
$$\begin{aligned}M^{12}_{(l'm'2),(lm1)}(\alpha)= \rho^{l+3}a_{11}[r^{l'-1}\dfrac{i}{2}(\mathscr{H}_{l,l'}^{m,m',(12)}(\alpha)-(-1)^m \mathscr{H}_{l,l'}^{-m,m',(12)}(\alpha))(-i\cdot l'(2l'+1))\\
+r^{l'-1}\dfrac{(-1)^{m'}}{2}\cdot i (\mathscr{H}_{l,l'}^{m,-m',(12)}(\alpha)-(-1)^m\mathscr{H}_{l,l'}^{-m,-m',(12)}(\alpha))(i\cdot l'(2l'+1))]
\end{aligned}$$

If $m<0$ and $m'=0$,
$$\begin{aligned}M^{12}_{(l'm'2),(lm1)}(\alpha)= \rho^{l+3}a_{11}r^{l'-1}\dfrac{i}{\sqrt{2}}(\mathscr{H}_{l,l'}^{m,0,(12)}(\alpha)-(-1)^m\mathscr{H}_{l,l'}^{-m,0,(12)}(\alpha))l'(2l'+1).\end{aligned}$$

If $m=0$ and $m'>0$,
$$\begin{aligned}M^{12}_{(l'm'2),(lm1)}(\alpha)= \rho^{l+3}a_{11}r^{l'-1}
[\dfrac{1}{\sqrt{2}} \mathscr{H}_{l,l'}^{0,-m',(12)}(\alpha)l'(2l'+1)+\dfrac{(-1)^{m'}}{\sqrt{2}}\mathscr{H}_{l,l'}^{0,m',(12)}(\alpha)l'(2l'+1)]
\end{aligned}$$

If $m=0$ and $m'<0$,
$$\begin{aligned}M^{12}_{(l'm'2),(lm1)}(\alpha)= \rho^{l+3}a_{11}r^{l'-1}
[\dfrac{1}{\sqrt{2}} \mathscr{H}_{l,l'}^{0,m',(12)}(\alpha)(-i\cdot l'(2l'+1))+\dfrac{(-1)^{m'}}{\sqrt{2}}\mathscr{H}_{l,l'}^{0,-m',(12)}(\alpha)(i\cdot l'(2l'+1))]
\end{aligned}$$

If $m=0$ and $m'=0$,
$$\begin{aligned}M^{12}_{(l'm'2),(lm1)}(\alpha)= \rho^{l+3}a_{11}r^{l'-1}
[\mathscr{H}_{l,l'}^{0,0,(12)}(\alpha)l'(2l'+1)]
\end{aligned}$$
where $a_{ij}$ stands for the coefficient of the entries of $A^{out}_{\mathcal{S},l}(x)$.
\end{theorem}

\begin{lemma}Denote $\sum\limits_{n\in\mathbb{Z}}L_{l,j,\lambda}^{m,\mu,q,m_1}(\mathbf{b}_{12})e^{-in\alpha}$ by $\boxed{\mathscr{L}_{l,j,\lambda}^{m,\mu,q,m_1,(12)}(\alpha)}$, then
$$\begin{aligned}&\mathscr{L}_{l,j,\lambda}^{m,\mu,q,m_1,(12)}(\alpha)=
i (-1)^{\lambda+\mu+q}\text{sgn}(q-m_1)\epsilon_{q}\sqrt{\lambda(2l+1)} 
\sqrt{{l+\lambda+\mu-m\choose \lambda+\mu}{l+\lambda+m-\mu\choose \lambda-\mu}}\\
&\langle\lambda-1,\mu-m_1;1,m_1\mid \lambda,\mu\rangle
\langle \lambda-1,\mu-m_1;1,q+m_1\mid j,\mu+q \rangle \sqrt{\dfrac{4\pi}{2(l+\lambda)+1}}\\
&\big(e^{-i\alpha}\cdot \mathrm{\Phi}(e^{-i\alpha},l+\lambda,1-2d)Y_{l+\lambda}^{m-\mu}(\dfrac{\pi}{2},\pi)-\mathrm{\Phi}(e^{i\alpha},l+\lambda,2d)Y_{l+\lambda}^{m-\mu}(\dfrac{\pi}{2},0)\big)
\end{aligned}$$
where $$\epsilon_q=\left\{
\begin{aligned}
\dfrac{1}{\sqrt{2}}&,\quad q=-1\\
0&,\quad q=0\\
-\dfrac{1}{\sqrt{2}}&,\quad,q=1 \end{aligned}\right.$$
\end{lemma}

\begin{theorem}
The expression of {$\boxed{M^{12}_{(l'm',2),(lm,3)}(\alpha)}$}($l'\geqslant 1, l\geqslant 1$, otherwise this term vanishes) is given by
When $m>0, m'>0$,
$$
\begin{aligned}
&M^{12}_{(l'm',2),(lm,3)}(\alpha)\\
=&\ \rho^{l+2}a_{33}\,r^{l'-1}\,\frac{\sqrt{l'(2l'+1)}}{2}\sum\limits_{q=-1}^{1}\sum\limits_{m_1=-1}^{1}\bigl[(-1)^{m'}\bigl(\mathscr{L}_{l,l',l'}^{-m,\,m'-q,\,q,\,m_1,(12)}(\alpha)+(-1)^m\mathscr{L}_{l,l',l'}^{m,\,m'-q,\,q,\,m_1,(12)}(\alpha)\bigr)\\
&\qquad\qquad\qquad\qquad+\bigl(\mathscr{L}_{l,l',l'}^{-m,\,-m'-q,\,q,\,m_1,(12)}(\alpha)+(-1)^m\mathscr{L}_{l,l',l'}^{m,\,-m'-q,\,q,\,m_1,(12)}(\alpha)\bigr)\bigr]
\end{aligned}
$$
When  $m>0, m'<0$,
$$
\begin{aligned}
&M^{12}_{(l'm',2),(lm,3)}(\alpha)\\
=&\ \rho^{l+2}a_{33}\,r^{l'-1}\,\frac{\sqrt{l'(2l'+1)}}{2}\sum\limits_{q=-1}^{1}\sum\limits_{m_1=-1}^{1}\bigl[i(-1)^{m'}\bigl(\mathscr{L}_{l,l',l'}^{-m,\,-m'-q,\,q,\,m_1,(12)}(\alpha)+(-1)^m\mathscr{L}_{l,l',l'}^{m,\,-m'-q,\,q,\,m_1,(12)}(\alpha)\bigr)\\
&\qquad\qquad\qquad\qquad-i\bigl(\mathscr{L}_{l,l',l'}^{-m,\,m'-q,\,q,\,m_1,(12)}(\alpha)+(-1)^m\mathscr{L}_{l,l',l'}^{m,\,m'-q,\,q,\,m_1,(12)}(\alpha)\bigr)\bigr]
\end{aligned}
$$

When $m>0, m'=0$,
$$
\begin{aligned}
&M^{12}_{(l',0,2),(lm,3)}(\alpha)\\
=&\ \rho^{l+2}a_{33}\,r^{l'-1}\,\frac{\sqrt{l'(2l'+1)}}{\sqrt{2}}\sum\limits_{q=-1}^{1}\sum\limits_{m_1=-1}^{1}\bigl(\mathscr{L}_{l,l',l'}^{-m,\,-q,\,q,\,m_1,(12)}(\alpha)+(-1)^m\mathscr{L}_{l,l',l'}^{m,\,-q,\,q,\,m_1,(12)}(\alpha)\bigr)
\end{aligned}
$$

When  $m<0, m'>0$,
$$
\begin{aligned}
&M^{12}_{(l'm',2),(lm,3)}(\alpha)\\
=&\ \rho^{l+2}a_{33}\,r^{l'-1}\,i\cdot\frac{\sqrt{l'(2l'+1)}}{2}\sum\limits_{q=-1}^{1}\sum\limits_{m_1=-1}^{1}\bigl[(-1)^{m'}\bigl(\mathscr{L}_{l,l',l'}^{m,\,m'-q,\,q,\,m_1,(12)}(\alpha)-(-1)^m\mathscr{L}_{l,l',l'}^{-m,\,m'-q,\,q,\,m_1,(12)}(\alpha)\bigr)\\
&\qquad\qquad\qquad\qquad+\bigl(\mathscr{L}_{l,l',l'}^{m,\,-m'-q,\,q,\,m_1,(12)}(\alpha)-(-1)^m\mathscr{L}_{l,l',l'}^{-m,\,-m'-q,\,q,\,m_1,(12)}(\alpha)\bigr)\bigr]
\end{aligned}
$$
When $m<0, m'<0$,
$$
\begin{aligned}
&M^{12}_{(l'm',2),(lm,3)}(\alpha)\\
=&\ \rho^{l+2}a_{33}\,r^{l'-1}\,\frac{\sqrt{l'(2l'+1)}}{2}\sum\limits_{q=-1}^{1}\sum\limits_{m_1=-1}^{1}\bigl[\bigl(\mathscr{L}_{l,l',l'}^{m,\,m'-q,\,q,\,m_1,(12)}(\alpha)-(-1)^m\mathscr{L}_{l,l',l'}^{-m,\,m'-q,\,q,\,m_1,(12)}(\alpha)\bigr)\\
&\qquad\qquad\qquad\qquad-(-1)^{m'}\bigl(\mathscr{L}_{l,l',l'}^{m,\,-m'-q,\,q,\,m_1,(12)}(\alpha)-(-1)^m\mathscr{L}_{l,l',l'}^{-m,\,-m'-q,\,q,\,m_1,(12)}(\alpha)\bigr)\bigr]
\end{aligned}
$$
When $m<0, m'=0$,
$$
\begin{aligned}
&M^{12}_{(l',0,2),(lm,3)}(\alpha)\\
=&\ \rho^{l+2}a_{33}\,r^{l'-1}\,\frac{i\sqrt{l'(2l'+1)}}{\sqrt{2}}\sum\limits_{q=-1}^{1}\sum\limits_{m_1=-1}^{1}\bigl(\mathscr{L}_{l,l',l'}^{m,\,-q,\,q,\,m_1,(12)}(\alpha)-(-1)^m\mathscr{L}_{l,l',l'}^{-m,\,-q,\,q,\,m_1,(12)}(\alpha)\bigr)
\end{aligned}
$$
When $m=0, m'>0$,
$$
\begin{aligned}
&M^{12}_{(l'm',2),(l,0,3)}(\alpha)\\
=&\ \rho^{l+2}a_{33}\,r^{l'-1}\,\frac{\sqrt{l'(2l'+1)}}{\sqrt{2}}\sum\limits_{q=-1}^{1}\sum\limits_{m_1=-1}^{1}\bigl[(-1)^{m'}\mathscr{L}_{l,l',l'}^{0,\,m'-q,\,q,\,m_1,(12)}(\alpha)+\mathscr{L}_{l,l',l'}^{0,\,-m'-q,\,q,\,m_1,(12)}(\alpha)\bigr]
\end{aligned}
$$

When  $m=0, m'<0$,
$$
\begin{aligned}
&M^{12}_{(l'm',2),(l,0,3)}(\alpha)\\
=&\ \rho^{l+2}a_{33}\,r^{l'-1}\,\frac{\sqrt{l'(2l'+1)}}{\sqrt{2}}\sum\limits_{q=-1}^{1}\sum\limits_{m_1=-1}^{1}\bigl[i(-1)^{m'}\mathscr{L}_{l,l',l'}^{0,\,-m'-q,\,q,\,m_1,(12)}(\alpha)-i\,\mathscr{L}_{l,l',l'}^{0,\,m'-q,\,q,\,m_1,(12)}(\alpha)\bigr]
\end{aligned}
$$

When $m=0, m'=0$,
$$
\begin{aligned}
&M^{12}_{(l',0,2),(l,0,3)}(\alpha)\\
=&\ \rho^{l+2}a_{33}\,r^{l'-1}\,\sqrt{l'(2l'+1)}\sum\limits_{q=-1}^{1}\sum\limits_{m_1=-1}^{1}\mathscr{L}_{l,l',l'}^{0,\,-q,\,q,\,m_1,(12)}(\alpha)
\end{aligned}
$$
\end{theorem}

\begin{lemma}Denote $\sum\limits_{n\in\mathbb{Z}}b^{(12)}_{-q}H_{l,\lambda}^{m,\mu}(\mathbf{b}_{12})e^{-in\alpha}$ by $\boxed{\mathscr{A}_{l,\lambda}^{m,\mu,(12)}(\alpha,q)}$, then 
$$\begin{aligned}
\sum\limits_{n\in\mathbb{Z}}b^{(12)}_{-q}H_{l,\lambda}^{m,\mu}(\mathbf{b}_{12})e^{-in\alpha}
=&(-1)^{\lambda+\mu}\epsilon_q
\sqrt{\dfrac{2l+1}{2\lambda+1}}
\sqrt{{l+\lambda+\mu-m\choose \lambda+\mu}{l+\lambda+m-\mu\choose \lambda-\mu}} \sqrt{\dfrac{4\pi}{2(l+\lambda)+1}}\\&
\big(e^{-i\alpha}\cdot \mathrm{\Phi}(e^{-i\alpha},l+\lambda,1-2d)Y_{l+\lambda}^{m-\mu}(\dfrac{\pi}{2},\pi)-\mathrm{\Phi}(e^{i\alpha},l+\lambda,2d)Y_{l+\lambda}^{m-\mu}(\dfrac{\pi}{2},0)\big)
\end{aligned}.$$

\end{lemma}

\begin{theorem}
The expression of {$\boxed{M^{12}_{(l'm',1),(lm,2)}(\alpha)}$}($l\geqslant 1$, otherwise this term vanishes) is given by\\
When $m>0$, $m'>0$,
$$
\begin{aligned}
&M^{12}_{(l'm',1),(lm,2)}(\alpha)\\
=&\sum\limits_{\lambda\in\{l'+1,l'+3\}}\sum\limits_{q=-1}^{1}\sum\limits_{m_1=-1}^{1}\rho^{l+1} (a_{22}-a_{12})  r^\lambda \sqrt{(l'+1)(2l'+1)}
[(-1)^{m'}(\mathscr{A}_{l,\lambda}^{-m,m'-q,(12)}(\alpha,q)\\
&+(-1)^m\mathscr{A}_{l,\lambda}^{m,m'-q,(12)}(\alpha,q))K_{l'+1,l',\lambda}^{m_1,m'-q,q}\\
&
+(\mathscr{A}_{l,\lambda}^{-m,-m'-q,(12)}(\alpha,q)+(-1)^m\mathscr{A}_{l,\lambda}^{m,-m'-q,(12)}(\alpha,q))K_{l'+1,l',\lambda}^{m_1,-m'-q,q}]\\
&-\rho^{l+1} a_{22}r^{l'+1}\dfrac{(2l+1)(l'+1)}{2}[(-1)^{m'}(\mathscr{H}_{l,l'}^{-m,m',(12)}(\alpha)\\
&+(-1)^m\mathscr{H}_{l,l'}^{m,m',(12)}(\alpha))+(\mathscr{H}_{l,l'}^{-m,-m',(12)}(\alpha)+(-1)^m\mathscr{H}_{l,l'}^{m,-m',(12)}(\alpha))]\\
&+\rho^{l+1} a_{22}r^{l'+1}(2l+1)\dfrac{\sqrt{(l'+1)(2l'+1)}}{2} \sum\limits_{q=-1}^1(-1)^q [(-1)^{m'}(\mathscr{A}_{l,l'+1}^{-m,m'-q,(12)}(\alpha,q)\\
&+(-1)^m\mathscr{A}_{l,l'+1}^{m,m'-q,(12)}(\alpha,q))\langle l'+1,m'-q;1,q\mid l',m'\rangle+(\mathscr{A}_{l,l'+1}^{-m,-m'-q,(12)}(\alpha,q)\\
&+(-1)^m\mathscr{A}_{l,l'+1}^{m,-m'-q,(12)}(\alpha,q))\langle l'+1,-m'-q;1,q\mid l',-m'\rangle]
\end{aligned}$$

When $m>0$, $m'<0$,

$$
\begin{aligned}
&M^{12}_{(l'm',1),(lm,2)}(\alpha)\\
=&\sum\limits_{\lambda\in\{l'+1,l'+3\}}\sum\limits_{q=-1}^{1}\sum\limits_{m_1=-1}^{1}\rho^{l+1} (a_{22}-a_{12})  r^\lambda \sqrt{(l'+1)(2l'+1)}
[i\cdot (-1)^{m'}(\mathscr{A}_{l,\lambda}^{-m,-m'-q,(12)}(\alpha,q)\\
&+(-1)^m\mathscr{A}_{l,\lambda}^{m,-m'-q,(12)}(\alpha,q))K_{l'+1,l',\lambda}^{m_1,-m'-q,q}
-i\cdot (\mathscr{A}_{l,\lambda}^{-m,m'-q,(12)}(\alpha,q)\\
&+(-1)^m\mathscr{A}_{l,\lambda}^{m,m'-q,(12)}(\alpha,q))K_{l'+1,l',\lambda}^{m_1,m'-q,q}]\\
&-\rho^{l+1} a_{22}r^{l'+1}\dfrac{(2l+1)(l'+1)}{2}[i\cdot (-1)^{m'}(\mathscr{H}_{l,l'}^{-m,-m',(12)}(\alpha)\\
&+(-1)^m\mathscr{H}_{l,l'}^{m,-m',(12)}(\alpha))\\
&-i\cdot (\mathscr{H}_{l,l'}^{-m,m',(12)}(\alpha)+(-1)^m\mathscr{H}_{l,l'}^{m,m',(12)}(\alpha))]\\
&+\rho^{l+1} a_{22}r^{l'+1}(2l+1)\dfrac{\sqrt{(l'+1)(2l'+1)}}{2} \sum\limits_{q=-1}^1(-1)^q [i\cdot (-1)^{m'}(\mathscr{A}_{l,l'+1}^{-m,-m'-q,(12)}(\alpha,q)\\
&+(-1)^m\mathscr{A}_{l,l'+1}^{m,-m'-q,(12)}(\alpha,q))\\
&\langle l'+1,-m'-q;1,q\mid l',-m'\rangle
-i\cdot (\mathscr{A}_{l,l'+1}^{-m,m'-q,(12)}(\alpha,q)\\
&+(-1)^m\mathscr{A}_{l,l'+1}^{m,m'-q,(12)}(\alpha,q))\langle l'+1,m'-q;1,q\mid l',m'\rangle]
\end{aligned}$$

When $m>0$, $m'=0$,
$$
\begin{aligned}
&M^{12}_{(l',0,1),(lm,2)}(\alpha)\\
=&\sum\limits_{\lambda\in\{l'+1,l'+3\}}\sum\limits_{q=-1}^{1}\sum\limits_{m_1=-1}^{1}\rho^{l+1} (a_{22}-a_{12})  r^\lambda \sqrt{2(l'+1)(2l'+1)}
[(\mathscr{A}_{l,\lambda}^{-m,-q,(12)}(\alpha,q)\\
&+(-1)^m\mathscr{A}_{l,\lambda}^{m,-q,(12)}(\alpha,q))K_{l'+1,l',\lambda}^{m_1,-q,q}]\\
&-\rho^{l+1} a_{22}r^{l'+1}\dfrac{(2l+1)(l'+1)}{\sqrt{2}}[\mathscr{H}_{l,l'}^{-m,0,(12)}(\alpha)+(-1)^m\mathscr{H}_{l,l'}^{m,0,(12)}(\alpha)]\\
&+\rho^{l+1} a_{22}r^{l'+1}(2l+1)\dfrac{\sqrt{(l'+1)(2l'+1)}}{\sqrt{2}} \sum\limits_{q=-1}^1(-1)^q (\mathscr{A}_{l,l'+1}^{-m,-q,(12)}(\alpha,q)+(-1)^m\mathscr{A}_{l,l'+1}^{m,-q,(12)}(\alpha,q))\\
&\langle l'+1,-q;1,q\mid l',0\rangle
\end{aligned}$$

When $m<0$, $m'>0$,
$$
\begin{aligned}
&M^{12}_{(l'm',1),(lm,2)}(\alpha)\\
=&\sum\limits_{\lambda\in\{l'+1,l'+3\}}\sum\limits_{q=-1}^{1}\sum\limits_{m_1=-1}^{1}\rho^{l+1} (a_{22}-a_{12})  r^\lambda i\cdot\sqrt{(l'+1)(2l'+1)}
[(-1)^{m'}(\mathscr{A}_{l,\lambda}^{m,m'-q,(12)}(\alpha,q)\\
&-(-1)^m\mathscr{A}_{l,\lambda}^{-m,m'-q,(12)}(\alpha,q))K_{l'+1,l',\lambda}^{m_1,m'-q,q}
+(\mathscr{A}_{l,\lambda}^{m,-m'-q,(12)}(\alpha,q)\\
&-(-1)^m\mathscr{A}_{l,\lambda}^{-m,-m'-q,(12)}(\alpha,q))K_{l'+1,l',\lambda}^{m_1,-m'-q,q}]\\
&-\rho^{l+1} a_{22}r^{l'+1}i\cdot \dfrac{(2l+1)(l'+1)}{2}[(-1)^{m'}(\mathscr{H}_{l,l'}^{m,m',(12)}(\alpha)-(-1)^m\mathscr{H}_{l,l'}^{-m,m',(12)}(\alpha))+(\mathscr{H}_{l,l'}^{m,-m',(12)}(\alpha)\\
&-(-1)^m\mathscr{H}_{l,l'}^{-m,-m',(12)}(\alpha))]\\
&+\rho^{l+1} a_{22}r^{l'+1}i\cdot (2l+1)\dfrac{\sqrt{(l'+1)(2l'+1)}}{2} \sum\limits_{q=-1}^1(-1)^q [(-1)^{m'}(\mathscr{A}_{l,l'+1}^{m,m'-q,(12)}(\alpha,q)\\
&-(-1)^m\mathscr{A}_{l,l'+1}^{-m,m'-q,(12)}(\alpha,q))\langle l'+1,m'-q;1,q\mid l',m'\rangle+(\mathscr{A}_{l,l'+1}^{m,-m'-q,(12)}(\alpha,q)\\
&-(-1)^m\mathscr{A}_{l,l'+1}^{-m,-m'-q,(12)}(\alpha,q))\langle l'+1,-m'-q;1,q\mid l',-m'\rangle]
\end{aligned}$$

When $m<0$, $m'<0$,

$$
\begin{aligned}
&M^{12}_{(l'm',1),(lm,2)}(\alpha)\\
=&\sum\limits_{\lambda\in\{l'+1,l'+3\}}\sum\limits_{q=-1}^{1}\sum\limits_{m_1=-1}^{1}\rho^{l+1} (a_{22}-a_{12})  r^\lambda \sqrt{(l'+1)(2l'+1)}
[-(-1)^{m'}(\mathscr{A}_{l,\lambda}^{m,-m'-q,(12)}(\alpha,q)\\
&-(-1)^m\mathscr{A}_{l,\lambda}^{-m,-m'-q,(12)}(\alpha,q))K_{l'+1,l',\lambda}^{m_1,-m'-q,q}
+ (\mathscr{A}_{l,\lambda}^{m,m'-q,(12)}(\alpha,q)\\
&-(-1)^m\mathscr{A}_{l,\lambda}^{-m,m'-q,(12)}(\alpha,q))K_{l'+1,l',\lambda}^{m_1,m'-q,q}]\\
&-\rho^{l+1} a_{22}r^{l'+1}\dfrac{(2l+1)(l'+1)}{2}
[- (-1)^{m'}(\mathscr{H}_{l,l'}^{m,-m',(12)}(\alpha)-(-1)^m\mathscr{H}_{l,l'}^{-m,-m',(12)}(\alpha))\\
&+
 (\mathscr{H}_{l,l'}^{m,m',(12)}(\alpha)-(-1)^m\mathscr{H}_{l,l'}^{-m,m',(12)}(\alpha))]\\
&+\rho^{l+1} a_{22}r^{l'+1}(2l+1)\dfrac{\sqrt{(l'+1)(2l'+1)}}{2} \sum\limits_{q=-1}^1(-1)^q [-(-1)^{m'}(\mathscr{A}_{l,l'+1}^{m,-m'-q,(12)}(\alpha,q)\\
&-(-1)^m\mathscr{A}_{l,l'+1}^{-m,-m'-q,(12)}(\alpha,q))\\
&\langle l'+1,-m'-q;1,q\mid l',-m'\rangle
+ (\mathscr{A}_{l,l'+1}^{m,m'-q,(12)}(\alpha,q)\\
&-(-1)^m\mathscr{A}_{l,l'+1}^{-m,m'-q,(12)}(\alpha,q))\langle l'+1,m'-q;1,q\mid l',m'\rangle]
\end{aligned}$$

When $m<0$, $m'=0$,
$$
\begin{aligned}
&M^{12}_{(l',0,1),(lm,2)}(\alpha)\\
=&\sum\limits_{\lambda\in\{l'+1,l'+3\}}\sum\limits_{q=-1}^{1}\sum\limits_{m_1=-1}^{1}\rho^{l+1} (a_{22}-a_{12})  r^\lambda i\cdot \sqrt{2(l'+1)(2l'+1)}
[(\mathscr{A}_{l,\lambda}^{m,-q,(12)}(\alpha,q)\\
&-(-1)^m\mathscr{A}_{l,\lambda}^{-m,-q,(12)}(\alpha,q))K_{l'+1,l',\lambda}^{m_1,-q,q}]\\
&-\rho^{l+1} a_{22}r^{l'+1}i\cdot \dfrac{(2l+1)(l'+1)}{\sqrt{2}}[(\mathscr{H}_{l,l'}^{m,0,(12)}(\alpha)-(-1)^m\mathscr{H}_{l,l'}^{-m,0,(12)}(\alpha)]\\
&+\rho^{l+1} a_{22}r^{l'+1}i\cdot (2l+1)\dfrac{\sqrt{(l'+1)(2l'+1)}}{\sqrt{2}} \sum\limits_{q=-1}^1(-1)^q [ (\mathscr{A}_{l,l'+1}^{m,-q,(12)}(\alpha,q)-(-1)^m\mathscr{A}_{l,l'+1}^{-m,-q,(12)}(\alpha,q))\\
&\langle l'+1,-q;1,q\mid l',0\rangle]
\end{aligned}$$

When $m=0$, $m'>0$,
$$
\begin{aligned}
&M^{12}_{(l'm',1),(l,0,2)}(\alpha)\\
=&\sum\limits_{\lambda\in\{l'+1,l'+3\}}\sum\limits_{q=-1}^{1}\sum\limits_{m_1=-1}^{1}\rho^{l+1} (a_{22}-a_{12})  r^\lambda \sqrt{2(l'+1)(2l'+1)}
[(-1)^{m'}\mathscr{A}_{l,\lambda}^{0,m'-q,(12)}(\alpha,q) K_{l'+1,l',\lambda}^{m_1,m'-q,q}\\
&+\mathscr{A}_{l,\lambda}^{0,-m'-q,(12)}(\alpha,q)K_{l'+1,l',\lambda}^{m_1,-m'-q,q}]\\
&-\rho^{l+1} a_{22}r^{l'+1}\dfrac{(2l+1)(l'+1)}{\sqrt{2}}[(-1)^{m'}\mathscr{H}_{l,l'}^{0,m',(12)}(\alpha)+\mathscr{H}_{l,l'}^{0,-m',(12)}(\alpha)]\\
&+\rho^{l+1} a_{22}r^{l'+1}(2l+1)\dfrac{\sqrt{(l'+1)(2l'+1)}}{\sqrt{2}} \sum\limits_{q=-1}^1(-1)^q [(-1)^{m'}\mathscr{A}_{l,l'+1}^{0,m'-q,(12)}(\alpha,q)\\
&\langle l'+1,m'-q;1,q\mid l',m'\rangle+(\mathscr{A}_{l,l'+1}^{0,-m'-q,(12)}(\alpha,q)\langle l'+1,-m'-q;1,q\mid l',-m'\rangle]
\end{aligned}$$

When $m=0$, $m'<0$,

$$
\begin{aligned}
&M^{12}_{(l'm',1),(l,0,2)}(\alpha)\\
=&\sum\limits_{\lambda\in\{l'+1,l'+3\}}\sum\limits_{q=-1}^{1}\sum\limits_{m_1=-1}^{1}\rho^{l+1} (a_{22}-a_{12})  r^\lambda \sqrt{2(l'+1)(2l'+1)}
[i\cdot (-1)^{m'}\mathscr{A}_{l,\lambda}^{0,-m'-q,(12)}(\alpha,q)K_{l'+1,l',\lambda}^{m_1,-m'-q,q}\\
&
-i\cdot \mathscr{A}_{l,\lambda}^{0,m'-q,(12)}(\alpha,q)K_{l'+1,l',\lambda}^{m_1,m'-q,q}]\\
&-\rho^{l+1} a_{22}r^{l'+1}\dfrac{(2l+1)(l'+1)}{\sqrt{2}}[i\cdot (-1)^{m'}\mathscr{H}_{l,l'}^{0,-m',(12)}(\alpha)-i\cdot \mathscr{H}_{l,l'}^{0,m',(12)}(\alpha)]\\
&+\rho^{l+1} a_{22}r^{l'+1}(2l+1)\dfrac{\sqrt{(l'+1)(2l'+1)}}{\sqrt{2}} \sum\limits_{q=-1}^1(-1)^q [i\cdot (-1)^{m'} \mathscr{A}_{l,l'+1}^{0,-m'-q,(12)}(\alpha,q)\\
&\langle l'+1,-m'-q;1,q\mid l',-m'\rangle-i\cdot \mathscr{A}_{l,l'+1}^{0,m'-q,(12)}(\alpha,q)\langle l'+1,m'-q;1,q\mid l',m'\rangle]
\end{aligned}$$

When $m=0$, $m'=0$,
$$
\begin{aligned}
&M^{12}_{(l',0,1),(l,0,2)}(\alpha)\\
=&\sum\limits_{\lambda\in\{l'+1,l'+3\}}\sum\limits_{q=-1}^{1}\sum\limits_{m_1=-1}^{1}\rho^{l+1} (a_{22}-a_{12})  r^\lambda \sqrt{4(l'+1)(2l'+1)}
[\mathscr{A}_{l,\lambda}^{0,-q,(12)}(\alpha,q)K_{l'+1,l',\lambda}^{m_1,-q,q}]\\
&-\rho^{l+1} a_{22}r^{l'+1}(2l+1)(l'+1)[\mathscr{H}_{l,l'}^{0,0,(12)}(\alpha)]\\
&+\rho^{l+1} a_{22}r^{l'+1}(2l+1)\sqrt{(l'+1)(2l'+1)} \sum\limits_{q=-1}^1(-1)^q [ \mathscr{A}_{l,l'+1}^{0,-q,(12)}(\alpha,q)\langle l'+1,-q;1,q\mid l',0\rangle]
\end{aligned}$$

\end{theorem}

\begin{theorem}
Entries for {$\boxed{M^{12}_{(l'm',3),(lm,2)}(\alpha)}$}($l'\geqslant 1,l\geqslant 1$, otherwise this term vanishes):\\
When $m>0$, $m'>0$,
$$
\begin{aligned}
&M^{12}_{(l'm',3),(lm,2)}(\alpha)\\
=&\sum\limits_{\lambda\in\{l',l'+2\mid l'\geqslant 1\}}\sum\limits_{q=-1}^{1}\sum\limits_{m_1=-1}^{1}-i\cdot \rho^{l+1} (a_{22}-a_{12})  r^\lambda \sqrt{l'(l'+1)}
[(-1)^{m'}(\mathscr{A}_{l,\lambda}^{-m,m'-q,(12)}(\alpha,q)\\
&+(-1)^m\mathscr{A}_{l,\lambda}^{m,m'-q,(12)}(\alpha,q))K_{l',l',\lambda}^{m_1,m'-q,q}
+(\mathscr{A}_{l,\lambda}^{-m,-m'-q,(12)}(\alpha,q)\\
&+(-1)^m\mathscr{A}_{l,\lambda}^{m,-m'-q,(12)}(\alpha,q))K_{l',l',\lambda}^{m_1,-m'-q,q}]\\
&- i\cdot \rho^{l+1} a_{22}r^{l'}(2l+1)\dfrac{\sqrt{l'(l'+1)}}{2} \sum\limits_{q=-1}^1(-1)^q [(-1)^{m'}(\mathscr{A}_{l,l'}^{-m,m'-q,(12)}(\alpha,q)+(-1)^m\mathscr{A}_{l,l'}^{m,m'-q,(12)}(\alpha,q))\\
&\langle l',m'-q;1,q\mid l',m'\rangle+(\mathscr{A}_{l,l'}^{-m,-m'-q,(12)}(\alpha,q)+(-1)^m\mathscr{A}_{l,l'}^{m,-m'-q,(12)}(\alpha,q))\langle l',-m'-q;1,q\mid l',-m'\rangle]
\end{aligned}$$

When $m>0$, $m'<0$,
$$
\begin{aligned}
&M^{12}_{(l'm',3),(lm,2)}(\alpha)\\
=&\sum\limits_{\lambda\in\{l',l'+2\mid l'\geqslant 1\}}\sum\limits_{q=-1}^{1}\sum\limits_{m_1=-1}^{1}\rho^{l+1} (a_{22}-a_{12})  r^\lambda \sqrt{l'(l'+1)}
[ (-1)^{m'}(\mathscr{A}_{l,\lambda}^{-m,-m'-q,(12)}(\alpha,q)\\
&+(-1)^m\mathscr{A}_{l,\lambda}^{m,-m'-q,(12)}(\alpha,q))K_{l',l',\lambda}^{m_1,-m'-q,q}
- (\mathscr{A}_{l,\lambda}^{-m,m'-q,(12)}(\alpha,q)+(-1)^m\mathscr{A}_{l,\lambda}^{m,m'-q,(12)}(\alpha,q))K_{l',l',\lambda}^{m_1,m'-q,q}]\\
&+\rho^{l+1} a_{22}r^{l'}(2l+1)\dfrac{\sqrt{l'(l'+1)}}{2} \sum\limits_{q=-1}^1(-1)^q [ (-1)^{m'}(\mathscr{A}_{l,l'}^{-m,-m'-q,(12)}(\alpha,q)+(-1)^m\mathscr{A}_{l,l'}^{m,-m'-q,(12)}(\alpha,q))\\
&\langle l',-m'-q;1,q\mid l',-m'\rangle
- (\mathscr{A}_{l,l'}^{-m,m'-q,(12)}(\alpha,q)+(-1)^m\mathscr{A}_{l,l'}^{m,m'-q,(12)}(\alpha,q))\langle l',m'-q;1,q\mid l',m'\rangle]\end{aligned}$$

When $m>0$, $m'=0$,
$$
\begin{aligned}
&M^{12}_{(l',0,3),(lm,2)}(\alpha)\\
=&\sum\limits_{\lambda\in\{l',l'+2\mid l'\geqslant 1\}}\sum\limits_{q=-1}^{1}\sum\limits_{m_1=-1}^{1} -i\cdot \rho^{l+1} (a_{22}-a_{12})  r^\lambda \sqrt{2l'(l'+1)}
[(\mathscr{A}_{l,\lambda}^{-m,-q,(12)}(\alpha,q)\\
&+(-1)^m\mathscr{A}_{l,\lambda}^{m,-q,(12)}(\alpha,q))K_{l',l',\lambda}^{m_1,-q,q}]\\
&-i\cdot \rho^{l+1} a_{22}r^{l'}(2l+1)\dfrac{\sqrt{l'(l'+1)}}{\sqrt{2}} \sum\limits_{q=-1}^1(-1)^q [ (\mathscr{A}_{l,l'}^{-m,-q,(12)}(\alpha,q)+(-1)^m\mathscr{A}_{l,l'}^{m,-q,(12)}(\alpha,q))\\
&\langle l',-q;1,q\mid l',0\rangle]
\end{aligned}$$

When $m<0$, $m'>0$,
$$
\begin{aligned}
&M^{12}_{(l'm',3),(lm,2)}(\alpha)\\
=&\sum\limits_{\lambda\in\{l',l'+2\mid l'\geqslant 1\}}\sum\limits_{q=-1}^{1}\sum\limits_{m_1=-1}^{1}\rho^{l+1} (a_{22}-a_{12})  r^\lambda \sqrt{l'(l'+1)}
[(-1)^{m'}(\mathscr{A}_{l,\lambda}^{m,m'-q,(12)}(\alpha,q)\\
&-(-1)^m\mathscr{A}_{l,\lambda}^{-m,m'-q,(12)}(\alpha,q))K_{l',l',\lambda}^{m_1,m'-q,q}
+(\mathscr{A}_{l,\lambda}^{m,-m'-q,(12)}(\alpha,q)\\
&-(-1)^m\mathscr{A}_{l,\lambda}^{-m,-m'-q,(12)}(\alpha,q))K_{l',l',\lambda}^{m_1,-m'-q,q}]\\
&+\rho^{l+1} a_{22}r^{l'} (2l+1)\dfrac{\sqrt{l'(l'+1)}}{2} \sum\limits_{q=-1}^1(-1)^q [(-1)^{m'}(\mathscr{A}_{l,l'}^{m,m'-q,(12)}(\alpha,q)-(-1)^m\mathscr{A}_{l,l'}^{-m,m'-q,(12)}(\alpha,q))\\
&\langle l',m'-q;1,q\mid l',m'\rangle+(\mathscr{A}_{l,l'}^{m,-m'-q,(12)}(\alpha,q)-(-1)^m\mathscr{A}_{l,l'}^{-m,-m'-q,(12)}(\alpha,q))\langle l',-m'-q;1,q\mid l',-m'\rangle]
\end{aligned}$$

When $m<0$, $m'<0$,
$$
\begin{aligned}
&M^{12}_{(l'm',3),(lm,2)}(\alpha)\\
=&\sum\limits_{\lambda\in\{l',l'+2\mid l'\geqslant 1\}}\sum\limits_{q=-1}^{1}\sum\limits_{m_1=-1}^{1}\rho^{l+1} (a_{22}-a_{12})  r^\lambda i\cdot \sqrt{l'(l'+1)}
[(-1)^{m'}(\mathscr{A}_{l,\lambda}^{m,-m'-q,(12)}(\alpha,q)\\
&-(-1)^m\mathscr{A}_{l,\lambda}^{-m,-m'-q,(12)}(\alpha,q))K_{l',l',\lambda}^{m_1,-m'-q,q}
- (\mathscr{A}_{l,\lambda}^{m,m'-q,(12)}(\alpha,q)\\
&-(-1)^m\mathscr{A}_{l,\lambda}^{-m,m'-q,(12)}(\alpha,q))K_{l',l',\lambda}^{m_1,m'-q,q}]\\
&-i\cdot \rho^{l+1} a_{22}r^{l'}(2l+1)\dfrac{\sqrt{l'(l'+1)}}{2} \sum\limits_{q=-1}^1(-1)^q [-(-1)^{m'}(\mathscr{A}_{l,l'}^{m,-m'-q,(12)}(\alpha,q)\\
&-(-1)^m\mathscr{A}_{l,l'}^{-m,-m'-q,(12)}(\alpha,q))\\
&\langle l',-m'-q;1,q\mid l',-m'\rangle
+ (\mathscr{A}_{l,l'}^{m,m'-q,(12)}(\alpha,q)-(-1)^m\mathscr{A}_{l,l'}^{-m,m'-q,(12)}(\alpha,q))\langle l',m'-q;1,q\mid l',m'\rangle]
\end{aligned}$$

When $m<0$, $m'=0$,
$$
\begin{aligned}
&M^{12}_{(l',0,3),(lm,2)}(\alpha)\\
=&\sum\limits_{\lambda\in\{l',l'+2\mid l'\geqslant 1\}}\sum\limits_{q=-1}^{1}\sum\limits_{m_1=-1}^{1}  \rho^{l+1} (a_{22}-a_{12})  r^\lambda  \sqrt{2l'(l'+1)}
[(\mathscr{A}_{l,\lambda}^{m,-q,(12)}(\alpha,q)\\
&-(-1)^m\mathscr{A}_{l,\lambda}^{-m,-q,(12)}(\alpha,q))K_{l',l',\lambda}^{m_1,-q,q}]\\
&+\rho^{l+1} a_{22}r^{l'} (2l+1)\dfrac{\sqrt{l'(l'+1)}}{\sqrt{2}} \sum\limits_{q=-1}^1(-1)^q [ (\mathscr{A}_{l,l'}^{m,-q,(12)}(\alpha,q)-(-1)^m\mathscr{A}_{l,l'}^{-m,-q,(12)}(\alpha,q))\\
&\langle l',-q;1,q\mid l',0\rangle]
\end{aligned}$$

When $m=0$, $m'>0$,
$$
\begin{aligned}
&M^{12}_{(l'm',3),(l,0,2)}(\alpha)\\
=&\sum\limits_{\lambda\in\{l',l'+2\mid l'\geqslant 1\}}\sum\limits_{q=-1}^{1}\sum\limits_{m_1=-1}^{1} -i\cdot \rho^{l+1} (a_{22}-a_{12})  r^\lambda \sqrt{2l'(l'+1)}
[(-1)^{m'}\mathscr{A}_{l,\lambda}^{0,m'-q,(12)}(\alpha,q) K_{l',l',\lambda}^{m_1,m'-q,q}\\
&+\mathscr{A}_{l,\lambda}^{0,-m'-q,(12)}(\alpha,q)K_{l',l',\lambda}^{m_1,-m'-q,q}]\\
&-i\cdot \rho^{l+1} a_{22}r^{l'}(2l+1)\dfrac{\sqrt{l'(l'+1)}}{\sqrt{2}} \sum\limits_{q=-1}^1(-1)^q [(-1)^{m'}\mathscr{A}_{l,l'}^{0,m'-q,(12)}(\alpha,q)\\
&\langle l',m'-q;1,q\mid l',m'\rangle+(\mathscr{A}_{l,l'}^{0,-m'-q,(12)}(\alpha,q)\langle l',-m'-q;1,q\mid l',-m'\rangle]
\end{aligned}$$

When $m=0$, $m'<0$,

$$
\begin{aligned}
&M^{12}_{(l'm',3),(l,0,2)}(\alpha)\\
=&\sum\limits_{\lambda\in\{l',l'+2\mid l'\geqslant 1\}}\sum\limits_{q=-1}^{1}\sum\limits_{m_1=-1}^{1}\rho^{l+1} (a_{22}-a_{12})  r^\lambda \sqrt{2l'(l'+1)}
[(-1)^{m'}\mathscr{A}_{l,\lambda}^{0,-m'-q,(12)}(\alpha,q)K_{l',l',\lambda}^{m_1,-m'-q,q}\\
&
- \mathscr{A}_{l,\lambda}^{0,m'-q,(12)}(\alpha,q)K_{l',l',\lambda}^{m_1,m'-q,q}]\\
&+\rho^{l+1} a_{22}r^{l'}(2l+1)\dfrac{\sqrt{l'(l'+1)}}{\sqrt{2}} \sum\limits_{q=-1}^1(-1)^q [ (-1)^{m'} \mathscr{A}_{l,l'}^{0,-m'-q,(12)}(\alpha,q)\\
&\langle l',-m'-q;1,q\mid l',-m'\rangle- \mathscr{A}_{l,l'}^{0,m'-q,(12)}(\alpha,q)\langle l',m'-q;1,q\mid l',m'\rangle]
\end{aligned}$$

When $m=0$, $m'=0$,
$$
\begin{aligned}
&M^{12}_{(l',0,3),(l,0,2)}(\alpha)\\
=&\sum\limits_{\lambda\in\{l',l'+2\mid l'\geqslant 1\}}\sum\limits_{q=-1}^{1}\sum\limits_{m_1=-1}^{1} -i\cdot \rho^{l+1} (a_{22}-a_{12})  r^\lambda \sqrt{4l'(l'+1)}
[\mathscr{A}_{l,\lambda}^{0,-q,(12)}(\alpha,q)K_{l',l',\lambda}^{m_1,-q,q}]\\
&-i\cdot \rho^{l+1} a_{22}r^{l'}(2l+1)\sqrt{l'(l'+1)} \sum\limits_{q=-1}^1(-1)^q [ \mathscr{A}_{l,l'}^{0,-q,(12)}(\alpha,q)\langle l',-q;1,q\mid l',0\rangle]
\end{aligned}$$

\end{theorem}

\begin{theorem}
Entries for {$\boxed{M^{12}_{(l'm',3),(lm,3)}(\alpha)}$}($l'\geqslant 1,l\geqslant 1$, otherwise this term vanishes):\\

When $m>0,\ m'>0$:
$$
\begin{aligned}
&M^{12}_{(l'm',3),(lm,3)}(\alpha)\\
=&\ \rho^{l+2}a_{33}r^{l'}\Biggl\{
\frac{l'(l'+1)}{2}\Bigl[(-1)^{m'}\bigl(\mathscr{H}_{l,l'}^{-m,m',(12)}(\alpha)+(-1)^m \mathscr{H}_{l,l'}^{m,m',(12)}(\alpha)\bigr)\\
&
+\bigl(\mathscr{H}_{l,l'}^{-m,-m',(12)}(\alpha)+(-1)^m \mathscr{H}_{l,l'}^{m,-m',(12)}(\alpha)\bigr)\Bigr]\\
&-\frac{i\sqrt{l'(l'+1)}}{2}\sum_{q=-1}^{1}\sum_{m_1=-1}^{1}\Bigl[(-1)^{m'}\bigl(\mathscr{L}_{l,l',l'+1}^{-m,m'-q,q,m_1,(12)}(\alpha)+(-1)^m \mathscr{L}_{l,l',l'+1}^{m,m'-q,q,m_1,(12)}(\alpha)\bigr)\\
&\qquad\qquad\qquad\qquad\quad+\bigl(\mathscr{L}_{l,l',l'+1}^{-m,-m'-q,q,m_1,(12)}(\alpha)+(-1)^m \mathscr{L}_{l,l',l'+1}^{m,-m'-q,q,m_1,(12)}(\alpha)\bigr)\Bigr]\Biggr\}\end{aligned}
$$

When $m>0,\ m'<0$:
$$
\begin{aligned}
&M^{12}_{(l'm',3),(lm,3)}(\alpha)\\
=&\ \rho^{l+2}a_{33}r^{l'}\Biggl\{
\frac{l'(l'+1)}{2}\Bigl[i(-1)^{m'}\bigl(\mathscr{H}_{l,l'}^{-m,-m',(12)}(\alpha)+(-1)^m \mathscr{H}_{l,l'}^{m,-m',(12)}(\alpha)\bigr)\\
&
-i\bigl(\mathscr{H}_{l,l'}^{-m,m',(12)}(\alpha)+(-1)^m \mathscr{H}_{l,l'}^{m,m',(12)}(\alpha)\bigr)\Bigr]\\
&+\frac{\sqrt{l'(l'+1)}}{2}\sum_{q=-1}^{1}\sum_{m_1=-1}^{1}\Bigl[(-1)^{m'}\bigl(\mathscr{L}_{l,l',l'+1}^{-m,-m'-q,q,m_1,(12)}(\alpha)+(-1)^m \mathscr{L}_{l,l',l'+1}^{m,-m'-q,q,m_1,(12)}(\alpha)\bigr)\\
&\qquad\qquad\qquad\qquad\quad-\bigl(\mathscr{L}_{l,l',l'+1}^{-m,m'-q,q,m_1,(12)}(\alpha)+(-1)^m \mathscr{L}_{l,l',l'+1}^{m,m'-q,q,m_1,(12)}(\alpha)\bigr)\Bigr]\Biggr\}
\end{aligned}
$$

When $m>0,\ m'=0$:
$$
\begin{aligned}
&M^{12}_{(l',0,3),(lm,3)}(\alpha)\\
=&\ \rho^{l+2}a_{33}r^{l'}\Biggl\{
\frac{l'(l'+1)}{\sqrt{2}}\bigl(\mathscr{H}_{l,l'}^{-m,0,(12)}(\alpha)+(-1)^m \mathscr{H}_{l,l'}^{m,0,(12)}(\alpha)\bigr)\\
&-\frac{i\sqrt{l'(l'+1)}}{\sqrt{2}}\sum_{q=-1}^{1}\sum_{m_1=-1}^{1}\bigl(\mathscr{L}_{l,l',l'+1}^{-m,-q,q,m_1,(12)}(\alpha)+(-1)^m \mathscr{L}_{l,l',l'+1}^{m,-q,q,m_1,(12)}(\alpha)\bigr)\Biggr\}
\end{aligned}
$$

When $m<0,\ m'>0$:
$$
\begin{aligned}
&M^{12}_{(l'm',3),(lm,3)}(\alpha)\\
=&\ \rho^{l+2}a_{33}r^{l'}\Biggl\{
\frac{l'(l'+1)}{2}\Bigl[i(-1)^{m'}\bigl(\mathscr{H}_{l,l'}^{m,m',(12)}(\alpha)-(-1)^m \mathscr{H}_{l,l'}^{-m,m',(12)}(\alpha)\bigr)
+i\bigl(\mathscr{H}_{l,l'}^{m,-m',(12)}(\alpha)\\
&-(-1)^m \mathscr{H}_{l,l'}^{-m,-m',(12)}(\alpha)\bigr)\Bigr]\\
&+\frac{\sqrt{l'(l'+1)}}{2}\sum_{q=-1}^{1}\sum_{m_1=-1}^{1}\Bigl[(-1)^{m'}\bigl(\mathscr{L}_{l,l',l'+1}^{m,m'-q,q,m_1,(12)}(\alpha)-(-1)^m \mathscr{L}_{l,l',l'+1}^{-m,m'-q,q,m_1,(12)}(\alpha)\bigr)\\
&\qquad\qquad\qquad\qquad\quad+\bigl(\mathscr{L}_{l,l',l'+1}^{m,-m'-q,q,m_1,(12)}(\alpha)-(-1)^m \mathscr{L}_{l,l',l'+1}^{-m,-m'-q,q,m_1,(12)}(\alpha)\bigr)\Bigr]\Biggr\}
\end{aligned}
$$

When $m<0,\ m'<0$:
$$
\begin{aligned}
&M^{12}_{(l'm',3),(lm,3)}(\alpha)\\
=&\ \rho^{l+2}a_{33}r^{l'}\Biggl\{
\frac{l'(l'+1)}{2}\Bigl[\bigl(\mathscr{H}_{l,l'}^{m,m',(12)}(\alpha)-(-1)^m \mathscr{H}_{l,l'}^{-m,m',(12)}(\alpha)\bigr)
-(-1)^{m'}\bigl(\mathscr{H}_{l,l'}^{m,-m',(12)}(\alpha)\\
&-(-1)^m \mathscr{H}_{l,l'}^{-m,-m',(12)}(\alpha)\bigr)\Bigr]\\
&+\frac{i\sqrt{l'(l'+1)}}{2}\sum_{q=-1}^{1}\sum_{m_1=-1}^{1}\Bigl[(-1)^{m'}\bigl(\mathscr{L}_{l,l',l'+1}^{m,-m'-q,q,m_1,(12)}(\alpha)-(-1)^m \mathscr{L}_{l,l',l'+1}^{-m,-m'-q,q,m_1,(12)}(\alpha)\bigr)\\
&\qquad\qquad\qquad\qquad\quad-\bigl(\mathscr{L}_{l,l',l'+1}^{m,m'-q,q,m_1,(12)}(\alpha)-(-1)^m \mathscr{L}_{l,l',l'+1}^{-m,m'-q,q,m_1,(12)}(\alpha)\bigr)\Bigr]\Biggr\}
\end{aligned}
$$

When $m<0,\ m'=0$:
$$
\begin{aligned}
&M^{12}_{(l',0,3),(lm,3)}(\alpha)\\
=&\ \rho^{l+2}a_{33}r^{l'}\Biggl\{
\frac{il'(l'+1)}{\sqrt{2}}\bigl(\mathscr{H}_{l,l'}^{m,0,(12)}(\alpha)-(-1)^m \mathscr{H}_{l,l'}^{-m,0,(12)}(\alpha)\bigr)\\
&+\frac{\sqrt{l'(l'+1)}}{\sqrt{2}}\sum_{q=-1}^{1}\sum_{m_1=-1}^{1}\bigl(\mathscr{L}_{l,l',l'+1}^{m,-q,q,m_1,(12)}(\alpha)-(-1)^m \mathscr{L}_{l,l',l'+1}^{-m,-q,q,m_1,(12)}(\alpha)\bigr)\Biggr\}\end{aligned}
$$

When $m=0,\ m'>0$:
$$
\begin{aligned}
&M^{12}_{(l'm',3),(l,0,3)}(\alpha)\\
=&\ \rho^{l+2}a_{33}r^{l'}\Biggl\{
\frac{l'(l'+1)}{\sqrt{2}}\Bigl[(-1)^{m'}\mathscr{H}_{l,l'}^{0,m',(12)}(\alpha)+\mathscr{H}_{l,l'}^{0,-m',(12)}(\alpha)\Bigr]\\
&-\frac{i\sqrt{l'(l'+1)}}{\sqrt{2}}\sum_{q=-1}^{1}\sum_{m_1=-1}^{1}\Bigl[(-1)^{m'}\mathscr{L}_{l,l',l'+1}^{0,m'-q,q,m_1,(12)}(\alpha)+\mathscr{L}_{l,l',l'+1}^{0,-m'-q,q,m_1,(12)}(\alpha)\Bigr]\Biggr\}\end{aligned}
$$

When $m=0,\ m'<0$:
$$
\begin{aligned}
&M^{12}_{(l'm',3),(l,0,3)}(\alpha)\\
=&\ \rho^{l+2}a_{33}r^{l'}\Biggl\{
\frac{l'(l'+1)}{\sqrt{2}}\Bigl[i(-1)^{m'}\mathscr{H}_{l,l'}^{0,-m',(12)}(\alpha)-i\mathscr{H}_{l,l'}^{0,m',(12)}(\alpha)\Bigr]\\
&+\frac{\sqrt{l'(l'+1)}}{\sqrt{2}}\sum_{q=-1}^{1}\sum_{m_1=-1}^{1}\Bigl[(-1)^{m'}\mathscr{L}_{l,l',l'+1}^{0,-m'-q,q,m_1,(12)}(\alpha)-\mathscr{L}_{l,l',l'+1}^{0,m'-q,q,m_1,(12)}(\alpha)\Bigr]\Biggr\}\end{aligned}
$$

When $m=0,\ m'=0$:
$$
\begin{aligned}
&M^{12}_{(l',0,3),(l,0,3)}(\alpha)\\
=&\ \rho^{l+2}a_{33}r^{l'}\Biggl\{
l'(l'+1)\mathscr{H}_{l,l'}^{0,0,(12)}(\alpha)
-i\sqrt{l'(l'+1)}\sum_{q=-1}^{1}\sum_{m_1=-1}^{1}\mathscr{L}_{l,l',l'+1}^{0,-q,q,m_1,(12)}(\alpha)\Biggr\}
\end{aligned}
$$
\end{theorem}

\begin{lemma}Define $\mathscr{D}_{l,\lambda}^{m,\mu,(12)}(\alpha)=\sum\limits_{n\in\mathbb{Z}} (n+2d)^2H_{l,\lambda}^{m,\mu}(\mathbf{b}_{12})e^{-in\alpha}$. Then
$$\begin{aligned}
&\mathscr{D}_{l,\lambda}^{m,\mu,(12)}(\alpha)\\
=&(-1)^{\lambda+\mu}
\sqrt{\dfrac{2l+1}{2\lambda+1}}
\sqrt{{l+\lambda+\mu-m\choose \lambda+\mu}{l+\lambda+m-\mu\choose \lambda-\mu}} \sqrt{\dfrac{4\pi}{2(l+\lambda)+1}}\\
&
\big(e^{-i\alpha}\cdot \mathrm{\Phi}(e^{-i\alpha},l+\lambda-1,1-2d)Y_{l+\lambda}^{m-\mu}(\dfrac{\pi}{2},\pi)+\mathrm{\Phi}(e^{i\alpha},l+\lambda-1,2d)Y_{l+\lambda}^{m-\mu}(\dfrac{\pi}{2},0)\big)
\end{aligned}$$

\end{lemma}

\begin{theorem}The entries for  $\boxed{M^{12}_{(l'm',2),(lm,2)}(\alpha)}$($l'\geqslant 1,l\geqslant 1$, otherwise this term vanishes) are given by:
\\
If $m>0$ and $m'>0$,
$$\begin{aligned}&M^{12}_{(l'm',2),(lm,2)}(\alpha)\\
&=\dfrac{l'(2l'+1)}{2}[ \rho^{l+3}a_{12}r^{l'-1}+\rho^{l+1}(a_{22}-a_{12})(r^{l'+1} )+\rho^{l+1}a_{22}r^{l'+1}\dfrac{2l+1}{2l'+1}]\\
&[(-1)^{m'}(\mathscr{H}_{l,l'}^{-m,m',(12)}(\alpha)+(-1)^m\mathscr{H}_{l,l'}^{m,m',(12)}(\alpha))+(\mathscr{H}_{l,l'}^{-m,-m',(12)}(\alpha)+(-1)^m\mathscr{H}_{l,l'}^{m,-m',(12)}(\alpha))]\\
&+\dfrac{l'(2l'+1)}{2}[\rho^{l+1}(a_{22}-a_{12})r^{l'-1}][(-1)^{m'}(\mathscr{D}_{l,l'}^{-m,m',(12)}(\alpha)+(-1)^m\mathscr{D}_{l,l'}^{m,m',(12)}(\alpha))+(\mathscr{D}_{l,l'}^{-m,-m',(12)}(\alpha)\\
&+(-1)^m\mathscr{D}_{l,l'}^{m,-m',(12)}(\alpha))]\\
&+\rho^{l+1}(a_{22}-a_{12})\sqrt{l'(2l'+1)}\sum\limits_{\lambda\in\{l'+1,l'-1\}}\sum\limits_{q=-1}^{q=1}\sum\limits_{m_1=-1}^{1}r^\lambda [(-1)^{m'}(\mathscr{A}_{l,\lambda}^{-m,m'-q,(12)}(\alpha,q)\\
&+(-1)^m\mathscr{A}_{l,\lambda}^{m,m'-q,(12)}(\alpha,q))K_{l'-1,l',\lambda}^{m_1,m'-q,q}+(\mathscr{A}_{l,\lambda}^{-m,-m'-q,(12)}(\alpha,q)+(-1)^m\mathscr{A}_{l,\lambda}^{m,-m'-q,(12)}(\alpha,q))K_{l'-1,l',\lambda}^{m_1,-m'-q,q}]\\
&+\rho^{l+1}a_{22}(2l+1)\dfrac{\sqrt{l'(2l'+1)}}{2}\sum\limits_{q=-1}^1(-1)^q r^{l'-1}[(-1)^{m'}(\mathscr{A}_{l,l'-1}^{-m,m'-q,(12)}(\alpha,q)+(-1)^m\mathscr{A}_{l,l'-1}^{m,m'-q,(12)}(\alpha,q))\\
&\langle l'-1,m'-q;1,q\mid l',m' \rangle+(\mathscr{A}_{l,l'-1}^{-m,-m'-q,(12)}(\alpha,q)\\
&+(-1)^m\mathscr{A}_{l,l'-1}^{m,-m'-q,(12)}(\alpha,q))\langle l'-1,-m'-q;1,q\mid l',-m' \rangle]
\end{aligned}$$

If $m>0$ and $m'<0$,
$$
\begin{aligned}
&M^{12}_{(l'm',2),(lm,2)}(\alpha)\\
&=\dfrac{l'(2l'+1)}{2}\Bigl[ \rho^{l+3}a_{12}r^{l'-1}+\rho^{l+1}(a_{22}-a_{12})(r^{l'+1} )+\rho^{l+1}a_{22}r^{l'+1}\dfrac{2l+1}{2l'+1}\Bigr]\\
&\quad\Bigl[i\cdot(-1)^{m'}(\mathscr{H}_{l,l'}^{-m,-m',(12)}(\alpha)+(-1)^m\mathscr{H}_{l,l'}^{m,-m',(12)}(\alpha))-i\cdot(\mathscr{H}_{l,l'}^{-m,m',(12)}(\alpha)+(-1)^m\mathscr{H}_{l,l'}^{m,m',(12)}(\alpha))\Bigr]\\
&+\dfrac{l'(2l'+1)}{2}[\rho^{l+1}(a_{22}-a_{12})r^{l'-1}]\Bigl[i\cdot(-1)^{m'}(\mathscr{D}_{l,l'}^{-m,-m',(12)}(\alpha)+(-1)^m\mathscr{D}_{l,l'}^{m,-m',(12)}(\alpha))\\
&-i\cdot(\mathscr{D}_{l,l'}^{-m,m',(12)}(\alpha)+(-1)^m\mathscr{D}_{l,l'}^{m,m',(12)}(\alpha))\Bigr]\\
&+\rho^{l+1}(a_{22}-a_{12})\sqrt{l'(2l'+1)}\sum_{\lambda\in\{l'+1,l'-1\}}\sum_{q=-1}^{1}\sum_{m_1=-1}^{1}r^\lambda \Bigl[i\cdot(-1)^{m'}(\mathscr{A}_{l,\lambda}^{-m,-m'-q,(12)}(\alpha,q)\\
&+(-1)^m\mathscr{A}_{l,\lambda}^{m,-m'-q,(12)}(\alpha,q))K_{l'-1,l',\lambda}^{m_1,-m'-q,q}-i\cdot(\mathscr{A}_{l,\lambda}^{-m,m'-q,(12)}(\alpha,q)\\
&+(-1)^m\mathscr{A}_{l,\lambda}^{m,m'-q,(12)}(\alpha,q))K_{l'-1,l',\lambda}^{m_1,m'-q,q}\Bigr]\\
&+\rho^{l+1}a_{22}(2l+1)\dfrac{\sqrt{l'(2l'+1)}}{2}\sum_{q=-1}^{1}(-1)^q r^{l'-1}\Bigl[i\cdot(-1)^{m'}(\mathscr{A}_{l,l'-1}^{-m,-m'-q,(12)}(\alpha,q)\\
&+(-1)^m\mathscr{A}_{l,l'-1}^{m,-m'-q,(12)}(\alpha,q))\langle l'-1,-m'-q;1,q\mid l',-m' \rangle-i\cdot(\mathscr{A}_{l,l'-1}^{-m,m'-q,(12)}(\alpha,q)\\
&+(-1)^m\mathscr{A}_{l,l'-1}^{m,m'-q,(12)}(\alpha,q))\langle l'-1,m'-q;1,q\mid l',m' \rangle\Bigr]
\end{aligned}
$$

If $m>0$ and $m'=0$,
$$
\begin{aligned}
&M^{12}_{(l',0,2),(lm,2)}(\alpha)\\
&=\dfrac{l'(2l'+1)}{\sqrt{2}}\Bigl[ \rho^{l+3}a_{12}r^{l'-1}+\rho^{l+1}(a_{22}-a_{12})(r^{l'+1} )+\rho^{l+1}a_{22}r^{l'+1}\dfrac{2l+1}{2l'+1}\Bigr]\\
&\quad\Bigl[(\mathscr{H}_{l,l'}^{-m,0,(12)}(\alpha)+(-1)^m\mathscr{H}_{l,l'}^{m,0,(12)}(\alpha))\Bigr]\\
&+\dfrac{l'(2l'+1)}{\sqrt{2}}[\rho^{l+1}(a_{22}-a_{12})r^{l'-1}]\Bigl[(\mathscr{D}_{l,l'}^{-m,0,(12)}(\alpha)+(-1)^m\mathscr{D}_{l,l'}^{m,0,(12)}(\alpha))\Bigr]\\
&+\rho^{l+1}(a_{22}-a_{12})\sqrt{2l'(2l'+1)}\sum_{\lambda\in\{l'+1,l'-1\}}\sum_{q=-1}^{1}\sum_{m_1=-1}^{1}r^\lambda \Bigl[(\mathscr{A}_{l,\lambda}^{-m,-q,(12)}(\alpha,q)\\
&+(-1)^m\mathscr{A}_{l,\lambda}^{m,-q,(12)}(\alpha,q))K_{l'-1,l',\lambda}^{m_1,-q,q}\Bigr]\\
&+\rho^{l+1}a_{22}(2l+1)\dfrac{\sqrt{l'(2l'+1)}}{\sqrt{2}}\sum_{q=-1}^{1}(-1)^q r^{l'-1}\Bigl[(\mathscr{A}_{l,l'-1}^{-m,-q,(12)}(\alpha,q)\\
&+(-1)^m\mathscr{A}_{l,l'-1}^{m,-q,(12)}(\alpha,q))\langle l'-1,-q;1,q\mid l',0 \rangle\Bigr].
\end{aligned}
$$

If $m<0$ and $m'>0$,
$$
\begin{aligned}
&M^{12}_{(l'm',2),(lm,2)}(\alpha)\\
&=\dfrac{l'(2l'+1)}{2}\Bigl[ \rho^{l+3}a_{12}r^{l'-1}+\rho^{l+1}(a_{22}-a_{12})(r^{l'+1} )+\rho^{l+1}a_{22}r^{l'+1}\dfrac{2l+1}{2l'+1}\Bigr]\\
&\quad\cdot i\Bigl[(-1)^{m'}(\mathscr{H}_{l,l'}^{m,m',(12)}(\alpha)-(-1)^m\mathscr{H}_{l,l'}^{-m,m',(12)}(\alpha))+(\mathscr{H}_{l,l'}^{m,-m',(12)}(\alpha)-(-1)^m\mathscr{H}_{l,l'}^{-m,-m',(12)}(\alpha))\Bigr]\\
&+\dfrac{l'(2l'+1)}{2}[\rho^{l+1}(a_{22}-a_{12})r^{l'-1}]\cdot i\Bigl[(-1)^{m'}(\mathscr{D}_{l,l'}^{m,m',(12)}(\alpha)-(-1)^m\mathscr{D}_{l,l'}^{-m,m',(12)}(\alpha))\\
&+(\mathscr{D}_{l,l'}^{m,-m',(12)}(\alpha)-(-1)^m\mathscr{D}_{l,l'}^{-m,-m',(12)}(\alpha))\Bigr]\\
&+\rho^{l+1}(a_{22}-a_{12})\sqrt{l'(2l'+1)}\sum_{\lambda\in\{l'+1,l'-1\}}\sum_{q=-1}^{1}\sum_{m_1=-1}^{1}r^\lambda \;i\Bigl[(-1)^{m'}(\mathscr{A}_{l,\lambda}^{m,m'-q,(12)}(\alpha,q)\\
&\quad-(-1)^m\mathscr{A}_{l,\lambda}^{-m,m'-q,(12)}(\alpha,q))K_{l'-1,l',\lambda}^{m_1,m'-q,q}+(\mathscr{A}_{l,\lambda}^{m,-m'-q,(12)}(\alpha,q)\\
&-(-1)^m\mathscr{A}_{l,\lambda}^{-m,-m'-q,(12)}(\alpha,q))K_{l'-1,l',\lambda}^{m_1,-m'-q,q}\Bigr]\\
&+\rho^{l+1}a_{22}(2l+1)\dfrac{\sqrt{l'(2l'+1)}}{2}\sum_{q=-1}^{1}(-1)^q r^{l'-1}\;i\Bigl[(-1)^{m'}(\mathscr{A}_{l,l'-1}^{m,m'-q,(12)}(\alpha,q)-(-1)^m\mathscr{A}_{l,l'-1}^{-m,m'-q,(12)}(\alpha,q))\\
&\quad\langle l'-1,m'-q;1,q\mid l',m' \rangle+(\mathscr{A}_{l,l'-1}^{m,-m'-q,(12)}(\alpha,q)\\
&-(-1)^m\mathscr{A}_{l,l'-1}^{-m,-m'-q,(12)}(\alpha,q))\langle l'-1,-m'-q;1,q\mid l',-m' \rangle\Bigr].
\end{aligned}
$$

If $m<0$ and $m'<0$,
$$
\begin{aligned}
&M^{12}_{(l'm',2),(lm,2)}(\alpha)\\
&=\dfrac{l'(2l'+1)}{2}\Bigl[ \rho^{l+3}a_{12}r^{l'-1}+\rho^{l+1}(a_{22}-a_{12})(r^{l'+1} )+\rho^{l+1}a_{22}r^{l'+1}\dfrac{2l+1}{2l'+1}\Bigr]\\
&\quad\Bigl[-(-1)^{m'}(\mathscr{H}_{l,l'}^{m,-m',(12)}(\alpha)-(-1)^m\mathscr{H}_{l,l'}^{-m,-m',(12)}(\alpha))+(\mathscr{H}_{l,l'}^{m,m',(12)}(\alpha)-(-1)^m\mathscr{H}_{l,l'}^{-m,m',(12)}(\alpha))\Bigr]\\
&+\dfrac{l'(2l'+1)}{2}[\rho^{l+1}(a_{22}-a_{12})r^{l'-1}]\Bigl[-(-1)^{m'}(\mathscr{D}_{l,l'}^{m,-m',(12)}(\alpha)\\
&-(-1)^m\mathscr{D}_{l,l'}^{-m,-m',(12)}(\alpha))+(\mathscr{D}_{l,l'}^{m,m',(12)}(\alpha)-(-1)^m\mathscr{D}_{l,l'}^{-m,m',(12)}(\alpha))\Bigr]\\
&+\rho^{l+1}(a_{22}-a_{12})\sqrt{l'(2l'+1)}\sum_{\lambda\in\{l'+1,l'-1\}}\sum_{q=-1}^{1}\sum_{m_1=-1}^{1}r^\lambda \Bigl[-(-1)^{m'}(\mathscr{A}_{l,\lambda}^{m,-m'-q,(12)}(\alpha,q)\\
&\quad-(-1)^m\mathscr{A}_{l,\lambda}^{-m,-m'-q,(12)}(\alpha,q))K_{l'-1,l',\lambda}^{m_1,-m'-q,q}+(\mathscr{A}_{l,\lambda}^{m,m'-q,(12)}(\alpha,q)\\
&-(-1)^m\mathscr{A}_{l,\lambda}^{-m,m'-q,(12)}(\alpha,q))K_{l'-1,l',\lambda}^{m_1,m'-q,q}\Bigr]\\
&+\rho^{l+1}a_{22}(2l+1)\dfrac{\sqrt{l'(2l'+1)}}{2}\sum_{q=-1}^{1}(-1)^q r^{l'-1}\Bigl[-(-1)^{m'}(\mathscr{A}_{l,l'-1}^{m,-m'-q,(12)}(\alpha,q)\\
&-(-1)^m\mathscr{A}_{l,l'-1}^{-m,-m'-q,(12)}(\alpha,q))\\
&\quad\langle l'-1,-m'-q;1,q\mid l',-m' \rangle+(\mathscr{A}_{l,l'-1}^{m,m'-q,(12)}(\alpha,q)\\
&-(-1)^m\mathscr{A}_{l,l'-1}^{-m,m'-q,(12)}(\alpha,q))\langle l'-1,m'-q;1,q\mid l',m' \rangle\Bigr].\end{aligned}
$$

If $m<0$ and $m'=0$,
$$
\begin{aligned}
&M^{12}_{(l',0,2),(lm,2)}(\alpha)\\
&=\dfrac{l'(2l'+1)}{\sqrt{2}}\Bigl[ \rho^{l+3}a_{12}r^{l'-1}+\rho^{l+1}(a_{22}-a_{12})(r^{l'+1} )+\rho^{l+1}a_{22}r^{l'+1}\dfrac{2l+1}{2l'+1}\Bigr]
\cdot i\Bigl[(\mathscr{H}_{l,l'}^{m,0,(12)}(\alpha)\\
&-(-1)^m\mathscr{H}_{l,l'}^{-m,0,(12)}(\alpha))\Bigr]\\
&+\dfrac{l'(2l'+1)}{\sqrt{2}}[\rho^{l+1}(a_{22}-a_{12})r^{l'-1}]\cdot i\Bigl[(\mathscr{D}_{l,l'}^{m,0,(12)}(\alpha)-(-1)^m\mathscr{D}_{l,l'}^{-m,0,(12)}(\alpha))\Bigr]\\
&+\rho^{l+1}(a_{22}-a_{12})\sqrt{2l'(2l'+1)}\sum_{\lambda\in\{l'+1,l'-1\}}\sum_{q=-1}^{1}\sum_{m_1=-1}^{1}r^\lambda \;i\Bigl[(\mathscr{A}_{l,\lambda}^{m,-q,(12)}(\alpha,q)\\
&-(-1)^m\mathscr{A}_{l,\lambda}^{-m,-q,(12)}(\alpha,q))K_{l'-1,l',\lambda}^{m_1,-q,q}\Bigr]\\
&+\rho^{l+1}a_{22}(2l+1)\dfrac{\sqrt{l'(2l'+1)}}{\sqrt{2}}\sum_{q=-1}^{1}(-1)^q r^{l'-1}\;i\Bigl[(\mathscr{A}_{l,l'-1}^{m,-q,(12)}(\alpha,q)\\
&-(-1)^m\mathscr{A}_{l,l'-1}^{-m,-q,(12)}(\alpha,q))\langle l'-1,-q;1,q\mid l',0 \rangle\Bigr].
\end{aligned}
$$

If $m=0$ and $m'>0$,
$$
\begin{aligned}
&M^{12}_{(l'm',2),(l,0,2)}(\alpha)\\
&=\dfrac{l'(2l'+1)}{\sqrt{2}}\Bigl[ \rho^{l+3}a_{12}r^{l'-1}+\rho^{l+1}(a_{22}-a_{12})(r^{l'+1} )+\rho^{l+1}a_{22}r^{l'+1}\dfrac{2l+1}{2l'+1}\Bigr]\\
&\quad\Bigl[(-1)^{m'}\mathscr{H}_{l,l'}^{0,m',(12)}(\alpha)+\mathscr{H}_{l,l'}^{0,-m',(12)}(\alpha)\Bigr]\\
&+\dfrac{l'(2l'+1)}{\sqrt{2}}[\rho^{l+1}(a_{22}-a_{12})r^{l'-1}]\Bigl[(-1)^{m'}\mathscr{D}_{l,l'}^{0,m',(12)}(\alpha)+\mathscr{D}_{l,l'}^{0,-m',(12)}(\alpha)\Bigr]\\
&+\rho^{l+1}(a_{22}-a_{12})\sqrt{2l'(2l'+1)}\sum_{\lambda\in\{l'+1,l'-1\}}\sum_{q=-1}^{1}\sum_{m_1=-1}^{1}r^\lambda \Bigl[(-1)^{m'}\mathscr{A}_{l,\lambda}^{0,m'-q,(12)}(\alpha,q)K_{l'-1,l',\lambda}^{m_1,m'-q,q}\\
&+\mathscr{A}_{l,\lambda}^{0,-m'-q,(12)}(\alpha,q)K_{l'-1,l',\lambda}^{m_1,-m'-q,q}\Bigr]\\
&+\rho^{l+1}a_{22}(2l+1)\dfrac{\sqrt{l'(2l'+1)}}{\sqrt{2}}\sum_{q=-1}^{1}(-1)^q r^{l'-1}\Bigl[(-1)^{m'}\mathscr{A}_{l,l'-1}^{0,m'-q,(12)}(\alpha,q)\langle l'-1,m'-q;1,q\mid l',m' \rangle\\
&+\mathscr{A}_{l,l'-1}^{0,-m'-q,(12)}(\alpha,q)\langle l'-1,-m'-q;1,q\mid l',-m' \rangle\Bigr].
\end{aligned}
$$

If $m=0$ and $m'<0$,
$$
\begin{aligned}
&M^{12}_{(l'm',2),(l,0,2)}(\alpha)\\
&=\dfrac{l'(2l'+1)}{\sqrt{2}}\Bigl[ \rho^{l+3}a_{12}r^{l'-1}+\rho^{l+1}(a_{22}-a_{12})(r^{l'+1} )+\rho^{l+1}a_{22}r^{l'+1}\dfrac{2l+1}{2l'+1}\Bigr]\\
&\quad i\Bigl[(-1)^{m'}\mathscr{H}_{l,l'}^{0,-m',(12)}(\alpha)-\mathscr{H}_{l,l'}^{0,m',(12)}(\alpha)\Bigr]\\
&+\dfrac{l'(2l'+1)}{\sqrt{2}}[\rho^{l+1}(a_{22}-a_{12})r^{l'-1}]\cdot i\Bigl[(-1)^{m'}\mathscr{D}_{l,l'}^{0,-m',(12)}(\alpha)-\mathscr{D}_{l,l'}^{0,m',(12)}(\alpha)\Bigr]\\
&+\rho^{l+1}(a_{22}-a_{12})\sqrt{2l'(2l'+1)}\sum_{\lambda\in\{l'+1,l'-1\}}\sum_{q=-1}^{1}\sum_{m_1=-1}^{1}r^\lambda \;i\Bigl[(-1)^{m'}\mathscr{A}_{l,\lambda}^{0,-m'-q,(12)}(\alpha,q)K_{l'-1,l',\lambda}^{m_1,-m'-q,q}\\
&-\mathscr{A}_{l,\lambda}^{0,m'-q,(12)}(\alpha,q)K_{l'-1,l',\lambda}^{m_1,m'-q,q}\Bigr]\\
&+\rho^{l+1}a_{22}(2l+1)\dfrac{\sqrt{l'(2l'+1)}}{\sqrt{2}}\sum_{q=-1}^{1}(-1)^q r^{l'-1}\;i\Bigl[(-1)^{m'}\mathscr{A}_{l,l'-1}^{0,-m'-q,(12)}(\alpha,q)\langle l'-1,-m'-q;1,q\mid l',-m' \rangle\\
&-\mathscr{A}_{l,l'-1}^{0,m'-q,(12)}(\alpha,q)\langle l'-1,m'-q;1,q\mid l',m' \rangle\Bigr].
\end{aligned}
$$

If $m=0$ and $m'=0$,
$$
\begin{aligned}
&M^{12}_{(l',0,2),(l,0,2)}(\alpha)\\
&=l'(2l'+1)\Bigl[ \rho^{l+3}a_{12}r^{l'-1}+\rho^{l+1}(a_{22}-a_{12})(r^{l'+1} )+\rho^{l+1}a_{22}r^{l'+1}\dfrac{2l+1}{2l'+1}\Bigr]\mathscr{H}_{l,l'}^{0,0,(12)}(\alpha)\\
&+{l'(2l'+1)}[\rho^{l+1}(a_{22}-a_{12})r^{l'-1}]\mathscr{D}_{l,l'}^{0,0,(12)}(\alpha)\\
&+\rho^{l+1}(a_{22}-a_{12})\cdot 2\sqrt{l'(2l'+1)}\sum_{\lambda\in\{l'+1,l'-1\}}\sum_{q=-1}^{1}\sum_{m_1=-1}^{1}r^\lambda \;\mathscr{A}_{l,\lambda}^{0,-q,(12)}(\alpha,q)K_{l'-1,l',\lambda}^{m_1,-q,q}\\
&+\rho^{l+1}a_{22}(2l+1)\sqrt{l'(2l'+1)}\sum_{q=-1}^{1}(-1)^q r^{l'-1}\;\mathscr{A}_{l,l'-1}^{0,-q,(12)}(\alpha,q)\langle l'-1,-q;1,q\mid l',0 \rangle.
\end{aligned}
$$
\end{theorem}

\begin{theorem}The remaining terms are simple:
$$M^{12}_{(l'm',3),(lm,1)}(\alpha)=M^{12}_{(l'm',1),(lm,3)}(\alpha)=0,$$
$$M^{12}_{(l'm',1),(lm,1)}(\alpha)=0.$$
\end{theorem}

\section{The first few vector spherical harmonics}
\label{app2}
First, consider the table of vector spherical harmonics up to the second order as listed below:

\[l = 0: \quad Y_{0,0} = \frac{1}{2} \sqrt{\frac{1}{\pi}}, \quad l = 1: \quad Y_{1,-1} = \sqrt{\frac{3}{4\pi}} y, \quad Y_{1,0} = \sqrt{\frac{3}{4\pi}} z, \quad Y_{1,1} = \sqrt{\frac{3}{4\pi}} x,\]

\[l = 2: \quad Y_{2,-2} = \frac{1}{2} \sqrt{\frac{15}{\pi}} xy, \quad Y_{2,-1} = \frac{1}{2} \sqrt{\frac{15}{\pi}} yz, \quad Y_{2,0} = \frac{1}{4} \sqrt{\frac{5}{\pi}} (-x^2 - y^2 + 2z^2),\]

\[Y_{2,1} = \frac{1}{2} \sqrt{\frac{15}{\pi}} xz, \quad Y_{2,2} = \frac{1}{4} \sqrt{\frac{15}{\pi}} (x^2 - y^2).\]

This gives first the obvious result that

\[l = 0: \quad V_{0,0} = -\frac{1}{2} \sqrt{\frac{1}{\pi}} (x, y, z)^\top \quad \text{and} \quad W_{00} = X_{00} = 0.\]

The spherical harmonics \( V_{l m} \) up to order 2 are then given as follows

\[l = 0: \quad V_{0,0} = -\frac{1}{2} \sqrt{\frac{1}{\pi}} (x, y, z)^\top,\]

\[l = 1: \quad V_{1,-1} = \sqrt{\frac{3}{4\pi}} ((0, 1, 0) - 2y(x, y, z))^\top,\]

\[V_{1,0} = \sqrt{\frac{3}{4\pi}} ((0, 0, 1) - 2z(x, y, z))^\top,\]

\[V_{1,1} = \sqrt{\frac{3}{4\pi}} ((1, 0, 0) - 2x(x, y, z))^\top,\]

\[l = 2: \quad V_{2,-2} = \frac{1}{2} \sqrt{\frac{15}{\pi}} ((y, x, 0) - 5xy(x, y, z))^\top,\]

\[V_{2,-1} = \frac{1}{2} \sqrt{\frac{15}{\pi}} ((0, z, y) - 5yz(x, y, z))^\top,\]

\[V_{2,0} = \frac{1}{2} \sqrt{\frac{5}{\pi}} ((-x, -y, 2z) - \frac{5}{2}(-x^2 - y^2 + 2z^2)(x, y, z))^\top,\]

\[V_{2,1} = \frac{1}{2} \sqrt{\frac{15}{\pi}} ((z, 0, x) - 5xz(x, y, z))^\top,\]

\[V_{2,2} = \frac{1}{2} \sqrt{\frac{15}{\pi}} ((x, -y, 0) - \frac{5}{2}(x^2 - y^2)(x, y, z))^\top.\]

The spherical harmonics \( W_{l m} \) up to order 2 are given by

\[l = 0: \quad W_{0,0} = 0,\]

\[l = 1: \quad W_{1,-1} = \sqrt{\frac{3}{4\pi}} (0, 1, 0)^\top, \quad W_{1,0} = \sqrt{\frac{3}{4\pi}} (0, 0, 1)^\top, \quad W_{1,1} = \sqrt{\frac{3}{4\pi}} (1, 0, 0)^\top,\]

\[l = 2: \quad W_{2,-2} = \frac{1}{2} \sqrt{\frac{15}{\pi}} (y, x, 0)^\top, \quad W_{2,-1} = \frac{1}{2} \sqrt{\frac{15}{\pi}} (0, z, y)^\top, \quad W_{2,0} = \frac{1}{2} \sqrt{\frac{5}{\pi}} (-x, -y, 2z)^\top,\]

\[W_{2,1} = \frac{1}{2} \sqrt{\frac{15}{\pi}} (z, 0, x)^\top, \quad W_{2,2} = \frac{1}{2} \sqrt{\frac{15}{\pi}} (x, -y, 0)^\top.\]

And finally, the spherical harmonics \( X_{l m} \) up to order 2 are given by

\[l = 0: \quad X_{0,0} = 0,\]

\[l = 1: \quad X_{1,-1} = \sqrt{\frac{3}{4\pi}} (-z, 0, x)^\top, \quad X_{1,0} = \sqrt{\frac{3}{4\pi}} (y, -x, 0)^\top,\]

\[X_{1,1} = \sqrt{\frac{3}{4\pi}} (0, z, -y)^\top,\]

\[l = 2: \quad X_{2,-2} = \frac{1}{2} \sqrt{\frac{15}{\pi}} (-xz, yz, x^2 - y^2)^\top, \quad X_{2,-1} = \frac{1}{2} \sqrt{\frac{15}{\pi}} (y^2 - z^2, -xy, xz)^\top,\]

\[X_{2,0} = \frac{1}{2} \sqrt{\frac{5}{\pi}} (3yz, -3xz, 0)^\top, \quad X_{2,1} = \frac{1}{2} \sqrt{\frac{15}{\pi}} (xy, z^2 - x^2, -yz)^\top,\]

\[X_{2,2} = \frac{1}{2} \sqrt{\frac{15}{\pi}} (yz, xz, -2xy)^\top.\]


\bibliographystyle{plainnat}
\bibliography{MP}
\end{document}